\documentclass{amsart}
\usepackage{amscd,amssymb,amsxtra}
\usepackage[mathscr]{eucal}
\usepackage[all]{xy}
\usepackage{bbm}
\usepackage{graphicx}

\usepackage[usenames, dvipsnames]{color}
\usepackage{accents}
\usepackage{pdfsync}

\definecolor{lightgray}{gray}{.5}
\definecolor{lightred}{rgb}{.6,.6,1}
\definecolor{purple}{rgb}{.7,0,1}
\definecolor{softyellow}{rgb}{.8745,.8588,.7647}
\definecolor{softblue}{rgb}{.6862,.8078,.8745}

\usepackage{comment}
\usepackage{datetime}

\usepackage{hyperref}
\hypersetup{%
  bookmarksnumbered=true,%
  bookmarks=true,%
  colorlinks=true,%
allcolors=lightgray,%
}

\usepackage{mathtools}
\usepackage{ushort}
\excludecomment{prem}

\title[Kervaire invariant one]{On the non-existence of elements of
 Kervaire invariant one}

\author{M.~A.~Hill}
\address{Department of Mathematics \\ University of Virginia
\\Charlottesville, VA 22904}
\email{michael.a.hill@math.uva.edu}

\author{M.~J.~Hopkins}
\address{Department of Mathematics \\ Harvard University
\\Cambridge, MA 02138}
\email{mjh@math.harvard.edu}

\author{D.~C.~Ravenel}
\address{Department of Mathematics \\ Rochester University
\\Rochester, NY}
\email{doug@math.rochester.edu}

\newtheorem{thm}[equation]{Theorem}
\newtheorem{cor}[equation]{Corollary}
\newtheorem{lem}[equation]{Lemma}
\newtheorem{prop}[equation]{Proposition}
\newtheorem*{thm*}{Theorem}
\newtheorem*{cor*}{Corollary}
\newtheorem*{lem*}{Lemma}
\newtheorem*{prop*}{Proposition}
\newtheorem{cond}[equation]{Condition}

\ifx\undefined\theoremstyle
	\newtheorem{definition}[equation]{Definition}
	\newenvironment{defin}{\begin{definition}\rm}{\end{definition}}
	\newtheorem{conj}[equation]{Conjecture}
	\newtheorem{rem}[equation]{Remark}
	\newtheorem*{rem*}{Remark}
	\newtheorem{rems}[equation]{Remarks}

\else
	\theoremstyle{definition}

	\newtheorem{defin}[equation]{Definition}

	\theoremstyle{remark}

	\newtheorem{rem}[equation]{Remark}
	\newtheorem*{rem*}{Remark}

\ifx\undefined\prem	\newtheorem{prem}{\bf Working Remark} \fi

\fi

\newtheorem{eg}[equation]{Example}

\ifx\undefined\pf
\newenvironment{pf}{\bigskip{\em Proof:\/}}{\qed\medskip}
\newenvironment{pf*}[1]{\bigskip{\em #1:\/}}{\qed\medskip}
\fi

\ifx\undefined\numberwithin
\def\numberwithin#1#2{\makeatletter\@ifundefined{c@#1}{\@nocnterrr}{%
  \@ifundefined{c@#2}{\@nocnterr}{%
  \@addtoreset{#1}{#2}%
  \toks@\expandafter\expandafter\expandafter{\csname the#1\endcsname}%
  \expandafter\xdef\csname the#1\endcsname
    {\expandafter\noexpand\csname the#2\endcsname
     .\the\toks@}}}\makeatother}\fi

\ifx\undefined\qed\newcommand{\qed}{\hfil\rule{4pt}{6pt}\bigskip}\fi
\ifx\undefined\\operatorname\newcommand{\operatorname}[1]{\mathop{\mbox{\rm #1}}}\fi

\newcommand{\ext}{\operatorname{Ext}}

\newcommand{\Q}{{\mathbb Q}}
\newcommand{\Z}{{\mathbb Z}}
\newcommand{\R}{{\mathbb R}}

\newcommand{\zerowidth}[1]{\hbox to 0pt{\hss$\displaystyle #1$\hss}}
\newcommand{\LL}{[\mkern-2mu[}
\newcommand{\RR}{]\mkern-2mu]}

\newcommand{\F}{{\mathbb F}}

\newcommand{\C}{\mathbb C}

\ifx\undefined\eqref\newcommand{\eqref}[1]{\rm (\ref{#1})}\fi

\newcounter{thmItem}

\newcommand{\thmItemref}[1]
 	 {{\rm \ref{#1})}}

\newcommand{\thmListItem}[1]{\rm \romannumeral#1)}

\newenvironment{thmList}{\begin{list}%
{\rm \roman{thmItem})}{\usecounter{thmItem}
\setlength{\labelwidth}{2em}
\setlength{\itemindent}{2em}
\setlength{\leftmargin}{0pt}
\setlength{\listparindent}{0pt}
\setlength{\parsep}{0pt}
\setlength{\partopsep}{0pt}
\setlength{\itemsep}{\medskipamount}
\setlength{\topsep}{\medskipamount}
}}{\end{list}}

\newcounter{textItem}

\newcounter{condItem}

\ifx\undefined\slot
\newcommand{\slot}{\,-\,}
\fi

\newcounter{probi}

\newif\ifsols

\newcounter{problemi}

\newcounter{multListi}

\newcounter{nmultListi}

\DeclareMathOperator{\sym}{Sym}

\DeclareMathOperator{\Endo}{End}

\newcommand{\cp}{\mathbf{CP}}
\newcommand{\rp}{\mathbf{RP}}

\newcommand{\extmu}[3]{\ext^{#1,#2}_{MU_{\ast}MU}(MU_{\ast},#3)}

\newcommand{\spectra}{\EuScript S}

\newcommand{\swg}[1]{\EuScript {SW}^{#1}}
\newcommand{\mur}{MU_{\R}}
\newcommand{\MUR}{\mathcal{MU}_{\R}}
\newcommand{\kr}{K_{\R}}
\newcommand{\ko}{KO}
\DeclareMathOperator{\spin}{Spin}

\newcommand{\mutn}[1]{MU^{(\!(#1)\!)}}
\newcommand{\ssmutn}[1]{MU#1}
\DeclareMathOperator{\map}{Map}
\DeclareMathOperator{\fil}{fil}

\newcommand{\rmod}[1]{\mathcal M_{#1}}

\DeclareMathOperator{\slicesymbol}{\EuScript{S}\EuScript{l}}

\newcommand{\gslice}[1]{\slicesymbol^{G}_{>{#1}}}
\newcommand{\geslice}[1]{\slicesymbol^{G}_{\ge{#1}}}
\newcommand{\lslice}[2]{\slicesymbol^{#1}_{<{#2}}}
\newcommand{\leslice}[2]{\slicesymbol^{#1}_{\le{#2}}}

\newcommand{\slc}[2]{P_{#1}^{#1}#2}

\DeclareMathOperator{\ho}{ho}
\DeclareMathOperator{\spaces}{\mathcal T}

\newcommand{\smashove}[1]{\underset{#1}{\wedge}}

\newcommand{\norm}{N}

\newcommand{\tg}[1]{\tilde{E}#1}

\newcommand{\piu}{\pi^{u}}
\newcommand{\pim}{\underline{\pi}}

\newcommand{\hu}{H^{{u}}}
\newcommand{\gfp}[1]{\Phi^{#1}}
\newcommand{\gfpm}[1]{\Phi_{M}^{#1}}
\newcommand{\phig}{\gfp{G}}
\newcommand{\phizt}{\gfp{\zt}}
\newcommand{\phiztm}{\gfpm{\zt}}
\newcommand{\phigm}{\gfpm{G}}
\newcommand{\phihm}{\gfpm{H}}
\newcommand{\phih}{\gfp{H}}

\newcommand{\epfam}{E\pfamily_{+}}
\newcommand{\tepfam}{\tilde E\pfamily}

\newcommand{\einfty}{E_{\infty}}
\newcommand{\ainfty}{A_{\infty}}

\newcommand{\Zm}{\underline{\Z}}
\newcommand{\Am}{\underline{A}}
\newcommand{\Mm}{\ushort{M}}

\DeclareMathOperator{\ind}{ind}

\newcommand{\einftycat}{\mathbf{Comm}}
\newcommand{\ainftycat}{\mathbf{Alg}}
\newcommand{\geinftycat}[1]{\einftycat_{#1}}
\newcommand{\gainftycat}[1]{\ainftycat_{#1}}
\newcommand{\ugeinftycat}[1]{\einftycat^{#1}}
\newcommand{\ugainftycat}[1]{\ainftycat^{#1}}

\newcommand{\mm}{m}
\newcommand{\mbm}{\gbar{m}}

\newcommand{\zlt}{\Z_{(2)}}
\newcommand{\zltm}{\underline{\Z}_{(2)}}

\newcommand{\gbar}{\bar}
\newcommand{\uc}{h}

\newcommand{\tensj}[1]{\if{displaystyle} #1\underset{H}{\otimes} \else
#1\otimes_{H}\fi}

\DeclareMathAlphabet{\scr}{U}{rsfs}{m}{n}
\DeclareMathOperator{\thom}{Thom}

\newcommand{\icat}{\scr{J}}
\newcommand{\igcat}[1]{\icat_{#1}}
\newcommand{\uigcat}[1]{\icat^{#1}}
\newcommand{\gspaces}[1]{\underline{\spaces\mkern-4mu}\mkern4mu_{\!#1}}
\newcommand{\ugspaces}[1]{\spaces^{#1}}

\newcommand{\gspectra}[1]{\underaccent{\bar}{\spectra}_{#1}} 
\newcommand{\ugspectra}[1]{\spectra^{#1}}
\newcommand{\flugspectra}[1]{\spectra^{#1}_{\text{fl}}}

\newcommand{\spectrar}{\spectra_{\R}}

\newcommand{\hocolim}{\ho\!\varinjlim}
\newcommand{\holim}{\ho\!\varprojlim}

\newcommand{\cat}[1]{{\mathcal #1}}

\newcommand{\normrbar}{\gbar{\mathfrak d}}
\newcommand{\normr}{\Delta}

\newcommand{\rr}{\gbar{r}}

\newcommand{\zt}{C_{2}}
\newcommand{\ze}{C_{8}}
\newcommand{\zf}{C_{4}}
\newcommand{\ztn}{C_{2^{n}}}
\newcommand{\ztnm}[1]{C_{2^{#1}}}
\newcommand{\slicecell}{\widehat{S}}

\newcommand{\D}{D} 
\newcommand{\refine}{\widehat{W}}
\newcommand{\magic}{\Omega}
\newcommand{\premagic}{\Omega_{\mathbb{O}}}
\newcommand{\ie}{i.e.,}

\newcommand{\Ext}{\mbox{\rm Ext}}

\newcommand{\ckp}{c_{k}}
\newcommand{\ck}{c_{k}}

\newcommand{\gcat}[1]{\mathfrak{Cat}_{#1}}
\newcommand{\vcat}[1]{\mathfrak{Cat}_{\cat #1}}
\newcommand{\categories}{\mathfrak{Cat}}

\newcommand{\isosets}{\mathfrak{Sets}_{\text{iso}}}

\DeclareMathOperator{\comm}{\mathbf{comm}}
\DeclareMathOperator{\ass}{\mathbf{ass}}
\newcommand{\commcat}[1]{\comm\cat #1}
\newcommand{\asscat}[1]{\ass\cat #1}
\newcommand{\one}{\mathbf{1}}
\newcommand{\tensorp}[1]{#1^{\otimes}_{\ast}}
\DeclareMathOperator{\ev}{ev}
\DeclareMathOperator{\ob}{ob}
\newcommand{\nat}{\mathbb N_{0}}

\newcommand{\fixg}[1]{#1^{G}}
\newcommand{\tensove}[1]{\underset{#1}{\otimes}}

\DeclareMathOperator{\RO}{\mathit {RO}}

\newcommand{\icatr}{\icat_{\R}}
\newcommand{\icatc}{\icat_{\C}}

\DeclareMathOperator{\cyl}{cyl}
\newcommand{\rspectra}{\spectra_{\R}}
\newcommand{\cspectra}{\spectra_{\C}}

\newcommand{\eggroup}[2]{E_{#1}#2}

\newcommand{\egsigma}[1]{\eggroup{#1}{\Sigma}}
\newcommand{\pfamily}{\mathcal P}
\newcommand{\tgp}{\tg{\pfamily}}

\newcommand{\bpstarbp}{BP_{\ast}BP}
\newcommand{\bpstar}{BP_{\ast}}

\newcommand{\mustarmu}{MU_{\ast}MU}

\newcommand{\mustar}{MU_{\ast}}

\newcommand{\xyhookar}{\ar@{^{(}->}@<-.3ex>}

\newcommand{\mfg}{\mathcal{M}_{FG}}
\newcommand{\mfgh}{\mathcal{M}_{FG}^{h}}
\newcommand{\funiv}{F_{\text{univ}}}
\newcommand{\theClass}{\upsilon}
\newcommand{\sta}{\mathcal A}
\newcommand{\ltgroup}{F_{f}}
\DeclareMathOperator{\tot}{Tot}

\newcommand{\lder}{\mathbf{L}}
\newcommand{\rder}{\mathbf{R}}
\newcommand{\lwedge}{\overset{\lder}{\wedge}}

\newcommand{\weq}{\mathcal W}
\newcommand{\sets}{\mathbf{Sets}}
\newcommand{\pscofib}[1]{\spectra^{#1}_{\mathbf{cof}}}

\newcommand{\pnalg}[1]{P_{\text{alg}}^{#1}}

\newcommand{\gswdes}[1]{\mathbf{Sp^{G}_{#1}}}

\newcommand{\augspectra}[1]{\spectra_{1}^{#1}}
\newcommand{\aspectra}[1]{\spectra_{1}^{#1}}
\newcommand{\gencof}{\mathcal A_{\text{cof}}}
\newcommand{\genacyclic}{\mathcal B_{\text{acyclic}}}

\newcommand{\pist}{\pi^{\text{st}}}
\DeclareMathOperator{\id}{Id}

\begin{document}

\begin{abstract}
We show that the Kervaire invariant one elements
$\theta_{j}\in\pi_{2^{j+1}-2}S^{0}$ exist only for $j\le 6$.  By
Browder's Theorem, this means that smooth framed manifolds of Kervaire
invariant one exist only in dimensions $2$, $6$, $14$, $30$, $62$, and
possibly $126$.  Except for dimension $126$ this resolves a
longstanding problem in algebraic topology.
\end{abstract}

\dedicatory{Dedicated to Mark Mahowald}
\thanks{M.~A.~Hill was partially supported by NSF grants DMS-0905160 , DMS-1307896 and
the Sloan foundation}
\thanks{M.~J.~Hopkins was partially supported the  NSF grant
DMS-0906194}
\thanks{D.~C.~Ravenel was partially supported by the NSF grants DMS-1307896 and  DMS-0901560}
\thanks{All three authors received support from the DARPA grants
HR0011-10-1-0054-DOD35CAP and FA9550-07-1-0555}
\maketitle

\tableofcontents

\section{Introduction}
\numberwithin{equation}{section}

The existence of smooth framed manifolds of Kervaire invariant one is
one of the oldest unresolved issues in differential and algebraic
topology.  The question originated in the work of Pontryagin in the
1930's.  It took a definitive form in the paper~\cite{kervaire60} of
Kervaire in which he constructed a combinatorial $10$-manifold with no
smooth structure, and in the work of
Kervaire-Milnor~\cite{kervaire63:_group} on $h$-cobordism classes of
manifolds homeomorphic to a sphere.  The question was connected to
homotopy theory by Browder in his fundamental
paper~\cite{browder69:_kervair} where he showed that smooth framed
manifolds of Kervaire invariant one exist only in dimensions of the
form $(2^{j+1}-2)$, and that a manifold exists in that dimension if
and only if the class
\[
h_{j}^{2}\in\ext_{\sta}^{2,2^{j+1}}(\Z/2,\Z/2)
\]
in the $E_{2}$-term of the classical Adams spectral represents an element
\[
\theta_{j}\in \pi_{2^{j+1}-2}S^{0}
\]
in the stable homotopy groups of spheres.  The classes $h_{j}^{2}$ for
$j\le 3$ represent the squares of the Hopf maps.  The element
$\theta_{4}\in\pi_{30}S^{0}$ had been observed in existing
computations~\cite{may64:_lie_hopf,MR0214072,MR0304463}, and was
constructed explicitly as a framed manifold by Jones~\cite{MR508888}.
The element $\theta_{5}\in\pi_{60}S^{0}$ was constructed by
Barratt-Mahowald, and Barratt-Jones-Mahowald (see~\cite{BJM} and the
discussion therein).

The purpose of this paper is to prove the following theorem
\begin{thm}
\label{thm:20} For $j\ge 7$ the class $h_{j}^{2}\in
\ext_{\sta}^{2,2^{j+1}}(\Z/2,\Z/2)$ does not represent an element of the
stable homotopy groups of spheres.  In other words, the Kervaire
invariant elements $\theta_{j}$ do not exist for $j\ge 7$.
\end{thm}

Smooth framed manifolds of Kervaire invariant one therefore exist only
in dimensions $2$, $6$, $14$, $30$, $62$, and possibly $126$.  At the
time of writing, our methods still leave open the existence of
$\theta_{6}$.

Many open issues in algebraic and differential topology depend on
knowing whether or not the Kervaire invariant one elements
$\theta_{j}$ exist for $j\ge 6$.  The following results represent some
of the issues now settled by Theorem~\ref{thm:20}.  In the statements,
the phrase ``exceptional dimensions'' refers to the dimensions $2$,
$6$, $14$, $30$, $62$, and $126$.  In all cases the situation in the
dimension $126$ is unresolved.  By Browder's
work~\cite{browder69:_kervair} the results listed below were known
when the dimension in question was not $2$ less than a power of $2$.
Modulo Browder's result~\cite{browder69:_kervair} the reduction of the
statements to Theorem~\ref{thm:20} can be found in the references
cited.

\begin{thm}[\cite{kervaire63:_group,MR802786}] 
\label{thm:69} Except in the six
exceptional dimensions, every stably framed smooth manifold is framed
cobordant to a homotopy sphere. \qed
\end{thm}

In the first five of the exceptional dimensions it is known that not
every stably framed manifold is framed cobordant to a homotopy
sphere.  The situation is unresolved in dimension $126$.

\begin{thm}[\cite{kervaire63:_group}]
\label{thm:110} Let $M^{m}$ be the manifold with boundary constructed
by plumbing together two copies of the unit tangent bundle to
$S^{2k+1}$ (so $m=4k+2$), and set $\Sigma^{m-1}=\partial M^{m}$.
Unless $m$ is one of the six exceptional dimensions, the space
$M^{m}/\Sigma^{m-1}$ is a triangulable manifold which does not admit
any smooth structure, and the manifold $\Sigma^{m-1}$ (the Kervaire
sphere) is homoemorphic but not diffeomorphic to $S^{m-1}$.\qed
\end{thm}

In the first five of the exceptional cases, the Kervaire sphere is
known to be diffeomorphic to the ordinary sphere, and the Kervaire
manifold can be smoothed.

\begin{thm}[\cite{kervaire63:_group,MR802786}] 
Let $\Theta_{n}$ be the Kervaire-Milnor group of $h$-cobordism classes
of homotopy $n$-spheres.  Unless $(4k+2)$ is one of the six
exceptional dimensions, 
\[
\Theta_{4k+2} \approx \pi_{4k+2}S^{0} 
\]
and 
\[
\left|\Theta_{4k+1}\right| = a_{k}\left|\pi_{4k+1}S^{0}\right|,
\]
where $a_{k}$ is $1$ if $k$ is even, and $2$ if $k$ is odd.
\end{thm}

\begin{thm}[\cite{BJM3}]
Unless $n$ is $1$, or one of the six exceptional dimensions, the
Whitehead square $[\iota_{n+1},\iota_{n+1}]\in \pi_{2n+1}S^{n+1}$ is not
divisible by $2$.  \qed
\end{thm}

\subsection{Outline of the argument} Our proof builds on the
strategy used by the third author in~\cite{Rav:Arf} and on the
homotopy theoretic refinement developed by the second author and
Haynes Miller (see~\cite{MR1642902}).  

We construct a multiplicative cohomology theory $\magic$ and establish
the following results:
\begin{thm} [The Detection Theorem]\label{thm:25}
If $\theta_{j}\in\pi_{2^{j+1}-2}S^{0}$ is an element of Kervaire
invariant $1$, and $j>2$, then the ``Hurewicz'' image of $\theta_{j}$ in
$\Omega^{2-2^{j+1}}(\text{pt})$ is non-zero.
\end{thm}

\begin{thm}[The Periodicity Theorem]
\label{thm:39}
The cohomology theory $\magic$ is $256$-fold periodic: 
For all $X$, 
\[
\magic^{\ast}(X) \approx \magic^{\ast+256}(X).
\]
\end{thm}

\begin{thm}[The Gap Theorem]
\label{thm:42}
The groups $\magic^{i}(\text{pt})$ are zero for $0<i<4$.
\end{thm}

These three results easily imply Theorem~\ref{thm:20}.  The
Periodicity Theorem and the Gap Theorem imply that the groups
$\magic^{i}(\text{pt})$ are zero for $i\equiv 2\mod 256$.  By the Detection
Theorem, if $\theta_{j}$ exists it has a non-zero Hurewicz image in
$\magic^{2-2^{j+1}}(\text{pt})$.  But this latter group is zero if
$j\ge 7$.

\subsection{The cohomology theory $\magic$}
\label{sec:cohom-theory-omega}

Write $C_{n}$ for the cyclic group of order $n$.  Our cohomology
theory $\magic$ is part of a pair $(\magic,\premagic)$ analogous to
the orthogonal and unitary $K$-theory spectra $KO$ and $KU$.  The role
of complex conjugation on $KU$ is played by an action of $\ze$ on
$\premagic$, and $\magic$ arises as its fixed points.  It is better to
think of $\premagic$ as generalizing Atiyah's $\zt$-equivariant
$K_{\R}$-theory~\cite{atiyah66:_k}, and in fact $\premagic$ is
constructed from the corresponding real bordism spectrum, as we now
describe.

Let $\mur$ be the $\zt$-equivariant {\em real bordism} spectrum of
Landweber~\cite{MR0222890} and Fujii~\cite{MR0420597}. Roughly
speaking one can think of $\mur$ as describing the cobordism theory of
{\em real manifolds}, which are stably almost complex manifolds
equipped with a conjugate linear action of $\zt$, such as the space of
complex points of a smooth variety defined over $\R$.  A real manifold
of real dimension $2n$ determines a homotopy class of maps
\[
S^{n\rho_{2}}\to \mur
\]
where $n\rho_{2}$ is the direct sum of $n$ copies of the real regular
representation of $\zt$, and $S^{n\rho_{2}}$ is its one point
compactification.  

Write
\[
\mutn{\ze}=\mur\wedge\mur\wedge\mur\wedge\mur
\]
for the $\ze$-equivariant spectrum gotten by smashing $4$ copies of
$\mur$ together and letting $\ze$ act by
\[
(a,b,c,d)\mapsto (\bar{d},a,b,c).
\]
Very roughly speaking, $\mutn{\ze}$ can be thought of as the cobordism theory
of stably almost manifolds equipped with a $\ze$-action, with the
property that the restriction of the action to $\zt\subset \ze$
determines a real structure.   If $M$ is a real manifold
then $M\times M\times M\times M$ with the $\ze$-action 
\[
(a,b,c,d) \mapsto (\bar d, a,b,c)
\]
is an example.   A suitable $\ze$-manifold $M$ of real dimension $8n$
determines a homotopy class of maps 
\[
S^{n\rho_{8}}\to \mutn{\ze},
\]
where $n\rho_{8}$ is the direct sum of $n$ copies of the real regular
representation of $\ze$, and $S^{n\rho_{8}}$ is its one point
compactification.  

To define $\magic$ we invert an 
equivariant analogue 
\[
\D:S^{\ell\rho_{8}}\to \mutn{\ze}
\]
of the Bott periodicity class and form the $\ze$-equivariant spectrum
$\premagic=\D^{-1}\mutn{\ze}$ (in fact $\ell$ works out to be $19$).
The cohomology theory $\magic$ is defined to be the homotopy fixed
point spectrum of the $\ze$-action on $\premagic$.

There is some flexibility in the choice of $\D$, but it needs to be
chosen in order that the Periodicity Theorem holds, and in order that
the map from the fixed point spectrum of $\premagic$ to the homotopy
fixed point spectrum be a weak equivalence.  It also needs to be
chosen in such a way that the Detection Theorem is preserved (see
Remark~\ref{rem:59}).  That such an $\D$ can be chosen with these
properties is a relatively easy fact, albeit mildly technical.  It is
specified in Corollary~\ref{thm:140}.  It can be described in the form
$M\times M\times M\times M$ for a suitable real manifold $M$, though
we do not do so.

\subsection{The Detection Theorem}
\label{sec:detection-theorem}

Since the non-equivariant spectrum $\premagic$ underlying $\premagic$
is complex orientable, the inclusion of the unit $S^{0}\to \magic$
induces a map
\[
\xymatrix{
*+{\extmu{s}{t}{MU_{\ast}}} \ar[d] \ar@{}[r]|-{\implies} &
\pi_{t-s}{S^{0}} \ar[d]\\
*+{H^{s}(\ze;\pi_{t}{\premagic})} \ar@{}[r]|-{\implies} &
\pi_{t-s}\magic
}
\]
from the Adams-Novikov spectral sequence to the $\ze$ homotopy fixed
point spectral sequence for $\pi_{\ast}\magic$.
In~\S\ref{sec:group-cohom-adams} we give an ad hoc construction of
this spectral sequence, conveniently adapted to describing the map of
$E_{2}$-terms.  It gives the horizontal arrow in the diagram of
spectral sequences below.
\[
\xymatrix{
*++\txt<8pc>{Adams-Novikov spectral sequence} \ar@{->}[d]\ar@{->}[r] & *\txt<8pc>{$\ze$ homotopy fixed point spectral sequence}\\
*++\txt<8pc>{Classical Adams spectral sequence}  & \\
}
\]
The Detection Theorem is proved by investigating this diagram, and
follows from a purely algebraic result.   

\begin{thm}[Algebraic Detection Theorem]
If 
\[
x\in \Ext_{MU_{*} (MU)}^{2, 2^{j+1}}\left(MU_{*}, MU_{*}\right)
\]
is any element mapping to $h_{j}^{2}$ in the $E_{2}$-term of the
classical Adams spectral sequence, and $j>2$, then the image of $x$ in
$H^{2} (\ze; \pi_{2^{j+1}}\premagic)$ is nonzero.
\end{thm}

The restriction $j>2$ is not actually necessary, but the other values
of $j$ require separate arguments.  Since we do not need them, we have
chosen to leave them to the interested reader.

To deduce the Detection Theorem from the Algebraic Detection Theorem
suppose that $\theta_{j}:S^{2^{j+1}-2}\to S^{0}$ is a map represented
by $h_{j}^{2}$ in the classical Adams spectral sequence.  Then
$\theta_{j}$ has Adams filtration $0$, $1$ or $2$ in the Adams-Novikov
spectral sequence, since the Adams filtration can only increase under
a map.  Since both
\[
\ext^{0,2^{j+1}-2}_{MU_{\ast}MU}\left(MU_{*}, MU_{*}\right)\quad\text{
and }\quad\ext^{{1,2^{j+1}-1}}_{MU_{\ast}MU}\left(MU_{*},
MU_{*}\right) 
\]
are zero, the class $\theta_{j}$ must be represented in Adams
filtration $2$ by some element $x$ which is a permanent cycle.  By the
Algebraic Detection Theorem, the element $x$ has a non-trivial image
$b_{j}\in H^{2}(\ze;\pi_{2^{j+1}}\premagic)$, representing the
image of $\theta_{j}$ in $\pi_{2^{j+1}-2}\magic$.  If this 
image is zero then the class $b_{j}$ must be in the image of the 
differential
\[
d_{2}:H^{0}(\ze;\pi_{2^{j+1}-1}\premagic) \to 
H^{2}(\ze;\pi_{2^{j+1}}\premagic).
\]
But $\pi_{\text{odd}}\premagic=0$, so this cannot happen.

The proof of the Algebraic Detection Theorem is given in
\S\ref{sec-detect}.  The method of proof is similar to that used
in~\cite{Rav:Arf}, where an analogous result is established at primes
greater than $3$.

\subsection{The slice filtration and the Gap Theorem}

While the Detection Theorem and the Periodicity Theorem involve the
homotopy fixed point spectral sequence for $\magic$, the Gap Theorem
results from studying $\premagic$ as an honest equivariant spectrum.
What permits the mixing of the two approaches is the following result,
which is part of Theorem~\ref{thm:52}.

\begin{thm}[Homotopy Fixed Point Theorem] The map from the fixed point
spectrum of $\premagic$ to the homotopy fixed point spectrum of
$\premagic$ is a weak equivalence.
\end{thm}

In particular, for all $n$, the map
\[
\pi_{n}^{\ze}\premagic\to \pi_{n}\premagic^{h\ze}=\pi_{n}\magic
\]
is an isomorphism, in which the symbol $\pi_{n}^{\ze}\premagic$
denotes the group of equivariant homotopy classes of maps from $S^{n}$
(with the trivial action) to $\premagic$.

We study the equivariant homotopy type of $\premagic$ using an
analogue of the Postnikov tower.  We call this tower the {\em slice
tower}.  Versions of it have appeared in work of Dan
Dugger~\cite{MR2240234}, Hopkins-Morel (unpublished),
Voevodsky\cite{MR2101286,MR1977582,MR1890744}, and
Hu-Kriz~\cite{MR1808224}.   

The slice tower is defined for any finite group $G$.  For a subgroup
$K\subset G$, let $\rho_{K}$ denote its regular representation and
write
\[
\slicecell(m,K)=G_{+}\smashove{K}S^{m\rho_{K}}\qquad m\in\Z.
\]

\begin{defin}
\label{def:18} The set of {\em slice cells} (for $G$) is
\[
\{\slicecell(m,K), \Sigma^{-1}\slicecell(m,K)\mid m\in \Z,
K\subset G \}.
\]
\end{defin}

\begin{defin}
\label{def:19}
A slice cell $\slicecell$ is {\em free} if it is of the form $G_{+}\wedge
S^{m}$ for some $m$.  An {\em isotropic} slice cell is one which is not free.
\end{defin}

We define the {\em dimension} of a slice cell $\slicecell$ by
\begin{align*}
\dim \slicecell(m,K) &=m|K| \\
\dim\Sigma^{-1}\slicecell(m,K) &= m |K|-1.
\end{align*}
Finally the {\em slice section} $P^{n}X$ is constructed by attaching
cones on slice cells $\slicecell$ with $\dim \slicecell>n$ to kill all
maps $\slicecell\to X$ with $\dim \slicecell>n$.  There is a natural
map
\[
P^{n}X\to P^{n-1}X
\]
The {\em $n$-slice of $X$} is defined to be its homotopy fiber
$\slc{n}{X}$.

In this way a tower $\{P^{n}X\}$, $n\in \Z$ is associated to each
equivariant spectrum $X$.  The homotopy colimit $\hocolim_{n}P^{n}X$ is
contractible, and $\holim_{n}P^{n}X$ is just $X$.  The {\em slice
spectral sequence} for $X$ is the spectral sequence of the slice
tower, relating $\pi_{\ast}\slc{n}{X}$ to $\pi_{\ast}X$.

The key technical result of the whole paper is the following.

\begin{thm}[The Slice Theorem] 
\label{thm:62} 
The $\ze$-spectrum
$\slc{n}{\mutn{\ze}}$ is contractible if $n$ is odd.  If $n$ is even
then $\slc{n}{\mutn{\ze}}$ is weakly equivalent to $H\Zm\wedge W$, where
$H\Zm$ is the Eilenberg-Mac~Lane spectrum associated to the constant
Mackey functor $\Zm$, and $W$ is a wedge of isotropic slice cells of
dimension $n$.
\end{thm}

The Slice Theorem actually holds more generally for the spectra
$\mutn{C_{2^{k}}}$ formed like $\mutn{\ze}$, using the smash product
of $2^{k-1}$ copies of $\mur$.  The more general statement is
Theorem~\ref{thm:89}

The Gap Theorem depends on the following result.

\begin{lem}[The Cell Lemma] Let $G=\ztn$ for some $n\ne 0$.  If
$\slicecell$ is an isotropic slice cell of even dimension, then the
groups $\pi_{k}^{G}H\Zm\wedge \slicecell$ are zero for $-4<k<0$.
\end{lem}

This is an easy explicit computation, and reduces to the fact that the
orbit space $S^{m\rho_{G}}/G$ is simply connected, being the
suspension of a connected space. 

Since the restriction of $\rho_{G}$ to a subgroup $K\subset G$ is
isomorphic to $(|G/K|)\,\rho_{K}$ there is an equivalence
\[
S^{m\rho_{G}}\wedge (G_{+}\smashove{K} S^{n\rho_{K}})\approx
G_{+}\smashove{K} S^{(n+m|G/K|)\rho_{K}}.
\]
It follows that if $\slicecell$ is a slice cell of dimension $d$, then
for any $m$, $S^{m\rho_{G}}\wedge \slicecell$ is a slice cell of
dimension $d+m|G|$.  Moreover, if $\slicecell$ is isotropic, then so
is $S^{m\rho_{G}}\wedge \slicecell$.  The Cell Lemma and the Slice
Theorem then imply that for any $m$, the group
\[
\pi^{\ze}_{i}S^{m\rho_{\ze}}\wedge\mutn{\ze}
\]
is zero for $-4<i<0$.   Since
\[
\pi^{\ze}_{i}\premagic = \varinjlim \pi_{i}S^{-m\ell\rho_{\ze}} \mutn{\ze}
\]
this implies that
\[
\pi^{\ze}_{i}\premagic=\pi_{i}\magic = 0
\]
for $-4<i<0$, which is the Gap Theorem.

The Periodicity Theorem is proved with a small amount of computation
in the $RO(\ze)$-graded slice spectral sequence for $\premagic$.  It
makes use of the fact that $\premagic$ is an equivariant commutative
ring spectrum.  Using the nilpotence machinery of~\cite{DHS,HS}
instead of explicit computation, it can be shown that the groups
$\pi_{\ast}\magic$ are periodic with {\em some} period which a power
of $2$.  This would be enough to show that only finitely many of the
$\theta_{j}$ can exist.  Some computation is necessary to get the
actual period stated in the Periodicity Theorem.

All of the results are fairly easy consequences of the Slice Theorem,
which in turn reduces to a single computational fact: that the
quotient of $\mutn{\ze}$ by the analogue of the ``Lazard ring'' is the
Eilenberg-Mac~Lane spectrum $H\Zm$ associated to the constant Mackey
functor $\Zm$.  We call this the {\em Reduction Theorem} and its
generalization to $\ztn$ appears as Theorem~\ref{thm:21}.  It is
proved for $G=\zt$ in Hu-Kriz~\cite{MR1808224}, and the analogue in
motivic homotopy theory is the main result of the (unpublished) work
of the second author and Morel mentioned earlier, where it is used to
identify the Voevodsky slices of $MGL$.  It would be very interesting
to find a proof of Theorem~\ref{thm:21} along the lines of Quillen's
argument in~\cite{MR0290382}.

During the long period between revisions of this paper, Haynes
Miller's Bourbaki talk on this material has appeared~\cite{MR3050712}.
We refer the reader there for a incisive overview.

\subsection{Summary of the contents}

We now turn to a more detailed summary of the contents of this paper.
In \S\ref{sec:some-equiv-stable} we recall the basics of equivariant
stable homotopy theory, establish many conventions and explain some
simple computations.  One of our main new constructions, introduced in
\S\ref{sec:change-group-2} is the multiplicative {\em norm functor}.
We merely state our main results about the norm, deferring the details
of the proofs to the appendices.  Another useful technique, the {\em
method of twisted monoid rings}, is described
in~\S\ref{sec:meth-poly-algebr}.  It is used in constructing
convenient filtrations of rings, and in forming the quotient of an
equivariant commutative ring spectrum by a regular sequence, in the
situation in which the group is acting non-trivially on the sequence.

Section~\ref{sec:ths-slice-filtration} introduces the slice
filtration, and establishes many of its basic properties, including
the strong convergence of the slice spectral sequence
(Theorem~\ref{thm:34}), and an important result on the distribution of
groups in the $E_{2}$-term (Corollary~\ref{thm:36}).  The notions of
{\em pure spectra}, {\em isotropic spectra}, and {\em spectra with
cellular slices} are introduced in \S\ref{sec:even-spectra}.  In these
terms, the Slice Theorem states that $\mutn{\ztn}$ is both pure and
isotropic.  Most of the material of these first sections makes no
restriction on the group $G$.  

From \S\ref{sec:case-mu} forward we restrict attention to the case in
which $G$ is cyclic of order a power of $2$, and we localize all
spectra at the prime $2$.  The spectra $\mutn{G}$ are introduced and
some of the basic properties are established.  The groundwork is laid
for the proof of the Slice Theorem.  The Reduction Theorem
(Theorem~\ref{thm:21}) is stated in \S\ref{sec:slice-theorem-i}.  The
Reduction Theorem is the backbone of the Slice Theorem, and is the
only part that is not ``formal'' in the sense that it depends on the
outcome of certain computations.

The Slice Theorem is also proved in \S\ref{sec:slice-theorem-i},
assuming that the Reduction Theorem holds.  The proof of the Reduction
Theorem is in \S\ref{sec:reduction-theorem}.  The Gap Theorem in
proved in \S\ref{sec:gap-theorem}, the Periodicity theorem in
\S\ref{sec:periodicity-theorem}.  The Homotopy Fixed Point Theorem is
proved in \S\ref{sec:fixed-points-homot}, and the Detection Theorem in
\S\ref{sec-detect}.

The paper concludes with two appendices devoted to
foundations of equivariant stable homotopy theory.  Two factors
contribute to the length of this material.  One is simply the wish to
make this paper as self-contained as possible and to collect material
central to our investigation in one place.  The other reason is that our
methods rely on multiplicative aspects of equivariant stable
homotopy theory that do not appear in the existing literature.
Establishing the basic properties of these structures involve details
of the foundations and cannot be done at the level of user interface.
Because of this, a relatively complete account of equivariant
orthogonal spectra is required.

\subsection{Acknowledgments} First and foremost the authors would like
to thank Ben Mann and the support of DARPA through the grant number
FA9550-07-1-0555.  It was the urging of Ben and the opportunity
created by this funding that brought the authors together in
collaboration in the first place.  Though the results described in
this paper were an unexpected outcome of our program, it's safe to say
they would not have come into being without Ben's prodding.  As it
became clear that the techniques of equivariant homotopy theory were
relevant to our project we drew heavily on the paper~\cite{MR1808224}
of Po Hu and Igor Kriz.  We'd like to acknowledge a debt of influence
to that paper, and to thank the authors for writing it.  We were also
helped by the thesis of Dan Dugger (which appears as
\cite{MR2240234}).  The second author would like to thank Dan Dugger,
Marc Levine, Jacob Lurie, and Fabien Morel for several useful
conversations.  Early drafts of this manuscript were read by Mark
Hovey, Tyler Lawson, and Peter Landweber, and the authors would like
to express their gratitude for their many detailed comments.  We also
owe thanks to Haynes Miller for a very thoughtful and careful reading
of our earlier drafts, and for his helpful suggestions for
terminology.  Thanks are due to Stefan Schwede for sharing with us his
construction of $\mur$, to Mike Mandell for diligently manning the
hotline for questions about the foundations of equivariant orthogonal
spectra, to Andrew Blumberg for his many valuable comments on the
second revision, and to Anna Marie Bohmann and Emily Riehl for
valuable comments on our description of ``working fiberwise.''

Finally, and most importantly, the authors would like to thank Mark
Mahowald for a lifetime of mathematical ideas and inspiration, and for
many helpful discussions in the early stages of this project.

\section{Equivariant stable homotopy theory}
\label{sec:some-equiv-stable}

We will work in the category of equivariant orthogonal
spectra~\cite{MR1806878,MR1922205}.  In this section we survey some of
the main properties of the theory and establish some notation.  The
definitions, proofs, constructions, and other details are explained in
Appendices~\ref{sec:equiv-orth-spectra}
and~\ref{sec:homot-theory-equiv}.  The reader is also referred to the
books of tom Dieck~\cite{MR889050,TtD:TGRT}, and the
survey of Greenlees and May~\cite{MR1361893} for an overview of
equivariant stable homotopy theory, and for further references.

We set up the basics of equivariant stable homotopy theory in the
framework of {\em homotopical category} in the sense
of~\cite{dwyer04:_homot}.  A homotopical category is a pair $(\cat C,
\cat W)$ consisting of a category $\cat C$ and a collection $\cat W$
of morphisms in $\cat C$ called {\em weak equivalences} containing all
identity maps, and satisfying the ``two out of six property'' that in
the situation
\[
\bullet \xrightarrow{u} \bullet \xrightarrow{v} \bullet
\xrightarrow{w} \bullet
\]
if $v u$ and $w v$ are in $\mathcal W$ then so are $u$, $v$, $w$, and
$v w u$.  Any class $\cat W$ defined as the collection of morphisms
$u$ taken to isomorphisms by some fixed functor automatically
satisfies this property.  This holds in particular when $\cat W$
consists of the weak equivalences in a model category structure.  In
this situation we will say that the model structure {\em refines} the
homotopical category structure, and that the homotopical category is
{\em completed} to a model category structure.  

Associated to a homotopical category $(\cat C, \cat W)$ is the homotopy category
$\ho\cat C$ and the functor $\cat C\to\ho\cat C$, 
characterized uniquely up to unique isomorphism by the following
universal property: for every category $\cat D$, and every functor
$F:\cat C\to \cat D$ taking the stable weak equivalences as
isomorphisms, there is a unique functor $\ho \cat C\to D$
making the diagram
\[
\xymatrix{
\cat C \ar[dr]_{F}\ar[r] & \ho\cat C \ar[d] \\
& \cat D 
}
\]
commute.    See~\S\ref{sec:categories-with-weak-2} for more on the theory of
homotopical categories, for a description of the issues that arise
when doing homotopy theory in a homotopical category, the techniques
for dealing with them, and for an explanation of the notion of left
($\lder$) and right ($\rder$) derived functors appearing in the
discussion below.

\subsection{$G$-spaces} 
\label{sec:g-spaces} 

We begin with unstable equivariant homotopy theory.  Let $G$ be a
finite group, and $\ugspaces{G}$ the topological category of pointed
compactly generated, weak Hausdorff left $G$-spaces and spaces of
equivariant maps.  The category $\ugspaces{G}$ is a closed symmetric
monoidal category under the smash product operation.  The tensor unit
is the $0$-sphere $S^{0}$ equipped with the trivial $G$-action.

We call a category enriched over $\ugspaces{G}$ a {\em $G$-equivariant
topological category}.  Since it is closed monoidal, $\ugspaces{G}$
may be regarded as enriched over itself.  We denote the enriched
category by $\gspaces{G}$.  Thus $\gspaces{G}$ is the $G$-equivariant
topological category of $G$-spaces and $G$-spaces of continuous, not
necessarily equivariant maps, on which $G$ acts by conjugation.  There
is an isomorphism
\[
\ugspaces{G}(X,Y)=\gspaces{G}(X,Y)^{G}.
\]
See \S\ref{sec:enriched-categories} and~\S\ref{sec:equivariant-spaces}
for further background and discussion.

The homotopy set (group, for $n>0$) $\pi^{H}_{n}(X)$ of a pointed
$G$-space is defined for $H\subset G$ and $n\ge 0$ to be the set of
$H$-equivariant homotopy classes of pointed maps
\[
S^{n}\to X.
\]
This is the same as the ordinary homotopy set (group) $\pi_{n}(X^{H})$
of the space of $H$ fixed-points in $X$.

A map $f:X\to Y$ in $\ugspaces{G}$ is a {\em weak equivalence} if for
all $H\subset G$ the map $X^{H}\to Y^{H}$ of $H$-fixed point spaces is
an ordinary weak equivalence.  Equivalently, $f:X\to Y$ is a weak
equivalence if for all $H\subset G$ and all choices of base point
$x_{0}\in X^{H}$ the induced map $\pi^{H}_{n}(X,x_{0})\to
\pi^{H}_{n}(Y,f(x_{0}))$ is an isomorphism.  Equipped with the weak
equivalences, the category underlying $\ugspaces{G}$ becomes a
homotopical category.  It can be completed to a topological model
category in which a {\em fibration} is a map $X\to Y$ which for every
$H\subset G$ is a Serre fibration on fixed points $X^{H}\to Y^{H}$.
The smash product of $G$-spaces makes $\ugspaces{G}$ into a symmetric
monoidal category in the sense of Schwede-Shipley~\cite[Definition
3.1]{MR1997322}, and $\gspaces{G}$ into an {\em enriched model
category}.

Every pointed $G$-space is weakly equivalent to a $G$-CW complex
constructed inductively from the basepoint by attaching equivariant
cells of the form $G/H\times D^{n}$ along maps from $G/H \times
S^{n-1}$.

We will write both
\[
\ho\ugspaces{G}(X,Y)\text{ and } [X,Y]^{G}
\]
for the set of maps from $X$ to $Y$ in the homotopy category of
$\ugspaces{G}$.  When $X$ is cofibrant and $Y$ is fibrant this can be
calculated as the set of homotopy classes of maps from $X$ to $Y$ in
$\ugspaces{G}$
\[
[X,Y]^{G} = \pi_{0}\ugspaces{G}(X,Y) = \pi_{0}^{G}\gspaces{G}(X,Y).
\]

We will make frequent use of {\em finite dimensional real orthogonal
representations of $G$}.   To keep the terminology simple these will
be referred to as {\em representations} of $G$.   

An important role is played by the equivariant spheres $S^{V}$ arising
as the one point compactification of representations $V$ of $G$.  When
$V$ is the trivial representation of dimension $n$, $S^{V}$ is just
the $n$-sphere $S^{n}$ with the trivial $G$-action.  We combine these
two notations and write
\[
S^{V+n}=S^{V\oplus \R^{n}}.
\]

Associated to $S^{V}$ is the equivariant homotopy set
\[
\pi_{V}^{G}X = [S^{V},X]^{G}
\]
defined to be the set of homotopy classes of $G$-equivariant maps
from $S^{V}$ to $X$.  The set $\pi_{V}^{G}X$ is a group if $\dim
V>0$ and an abelian group if $\dim V^{G}>1$, where $V^{G}$ is the
space of $G$-invariant vectors in $V$.   

Also associated to the sphere $S^{V}$ one has the equivariant suspension
$\Sigma^{V}X = S^{V}\wedge X$ and the equivariant loop space
$\Omega^{V}X= \gspaces{G}(S^{V},X)$.

Now suppose that $V_{1}$ and $V_{2}$ are two orthogonal
representations of $G$ and that for each irreducible representation
$U$ of $G$ occurring in $V_{1}$ one has
\begin{equation}
\label{eq:161}
\dim \hom^{G}(U,V_{2})\ge \dim \hom^{G}(U, V_{1}).
\end{equation}
Then one may choose an equivariant linear isometric embedding $t:V_{1}\to
V_{2}$ and form
\begin{equation}
\label{eq:154}
\pi_{V_{2}-t(V_{1})}^{G}(X),
\end{equation}
in which $V_{2}-t(V_{1})$ denotes the orthogonal complement of the
image of $V_{1}$ in $V_{2}$.  The groups~\eqref{eq:154} form a local
system over the Stiefel manifold $O(V_{1},V_{2})^{G}$ of equivariant linear isometric 
embeddings.   If instead of~\eqref{eq:161} the one has
$V_{2}>V_{1}$ in the sense of Definition~\ref{def:28} below, then the
Stiefel manifold $O(V_{1},V_{2})$ is simply connected and one may
define 
\[
\pi_{V_{2}-V_{2}}^{G}(X)
\]
to be the group of global sections of this local system.  For any
$t\in O(V_{1},V_{2})^{G}$ the restriction map gives a canonical
isomorphism $\pi_{V_{2}-V_{2}}^{G}(X)\to \pi_{V_{2}-t(V_{1})}^{G}(X)$.

\begin{defin}
\label{def:28} Let $V_{1}$ and $V_{2}$ be two non-zero $G$-representations.  We
write $V_{1}<V_{2}$ if for every irreducible $G$-representation $U$,
\[
\dim \hom^{G}(U,V_{1}) < \dim \hom^{G}(U, V_{2})-1.
\]
\end{defin}

This relation makes the set of $G$-representations into a (large) partially
ordered set.   

We will shortly (\S\ref{sec:homot-theory-ugsp}) be interested in the
special case in which $V_{1}$ is a trivial representation of dimension
$k$.  As above we will write
\[
\pi_{V_{2}-k}^{G}(X)
\]
for this group.   In this way, for any $n\in\Z$ there is a
well-defined group
\[
\pi_{V+n}^{G}(X)
\]
provided $\dim V^{G}\ge -n+2$.

\subsection{Equivariant stable homotopy theory}

There is a choice to be made when stabilizing equivariant homotopy
theory.  If one only seeks that fibration sequences and cofibration
sequences become weakly equivalent, then one stabilizes in the usual
way, using suspensions by spheres with trivial $G$-action.  But if one
wants to have Spanier-Whitehead duals of finite $G$-CW complexes, one
needs to stabilize with respect to the spheres $S^{V}$ where $V$ is a
finite dimensional representation of $G$.   

We will do equivariant stable homotopy theory in the category of
equivariant orthogonal spectra, equipped with the stable weak
equivalences.  In order for this to be considered viable some
properties must be established that guarantee computations made with
equivariant orthogonal spectra ultimately reduce to computations in
$\ho\ugspaces{G}$ in the expected manner.   We therefore begin by
discussing the equivariant Spanier-Whitehead category, and formulate
six properties an equivariant stable homotopy should satisfy in order
that it faithfully extend the Spanier-Whitehead category.   These properties
aren't enough for all of our purposes, so after establishing them for
equivariant orthogonal spectra we turn the more refined structures
(indexed products, coproducts, and smash products) that we require.

\subsubsection{Spanier-Whitehead stabilization}
\label{sec:g-spectra}

The {\em $G$-equivariant Spanier-Whitehead} category $\swg{G}$ is the
category whose objects are finite pointed $G$-CW complexes and with
maps
\[
\{X,Y \}^{G} = \varinjlim_{V}[S^{V}\wedge X,S^{V}\wedge Y]^{G},
\]
in which the colimit is taken over the partially ordered set of
$G$-representations.  For an informative discussion of this category
the reader is referred to~\cite{MR764596}.

There is a direct analogue~\cite{MR764596,MR0343260} of
Spanier-Whitehead duality in $\swg{G}$, in which a finite based $G$-CW
complex embedded in a representation sphere $S^{V}$ is ``$V$-dual'' to
the unreduced suspension of its complement.

\begin{eg}
\label{eg:27} Suppose that $X$ is a finite pointed $G$-set $B$.  If
there is an equivariant embedding $B\subset S^{V}$ (for instance when
$V$ is the $G$-representation with basis $B$) the $V$-dual of $B$
works out to be $S^{V}\wedge B$.
\end{eg}

If one wants finite $G$-CW complexes to have actual duals, in the
sense of objects in a symmetric monoidal category, then one must
enlarge the category $\swg{G}$ by formally adding, for each finite
$G$-CW complex $Y$ and each finite dimensional representation $V$ of
$G$, an object $S^{-V}\wedge Y$ defined by
\begin{equation}
\label{eq:58}
\{X,S^{-V}\wedge Y \}^{G} = 
\{S^{V}\wedge X,Y \}^{G}.
\end{equation}
Since $\{S^{V}\wedge (\slot),Y \}^{G}$ is a functor on $\swg{G}$, this
amounts to simply working in an enlargement of the Yoneda embedding of
$\swg{G}$.  One checks that for any $Z$, the map $Z\to S^{-V}\wedge
S^{V}\wedge Z$ corresponding to the identity map of $S^{V}\wedge Z$
under~\eqref{eq:58} is an isomorphism, and that symmetric monoidal
structure given by the smash product extends to this enlarged
category.   If $X$ and $Y$ are $V$-duals in $\swg{G}$, then $X$ and
$S^{-V}\wedge Y$ are duals in the enlarged equivariant
Spanier-Whitehead category. 

\begin{eg}
\label{eg:28}
From Example~\ref{eg:27}, $B$ is self-dual in the enlarged equivariant
Spanier-Whitehead category.
\end{eg}

As in the non-equivariant case, the equivariant Spanier-Whitehead
category still suffers the defect that it is also not quite set up for
doing stable homotopy theory.  What one wants is a complete closed
symmetric monoidal category $\ugspectra{G}$ of {\em $G$-equivariant
spectra}, equipped with the structure of a homotopical category (or
even a Quillen model category), and related to $\ugspaces{G}$ by a
pair of adjoint (suspension spectrum and zero space) functors
\[
\Sigma^{\infty}:\ugspaces{G} \leftrightarrows \ugspectra{G}: \Omega^{\infty}.
\]
In order to know that computations made in this category reduce in the
expected manner to those in classical stable homotopy theory, one
would like this data to satisfy

\begin{itemize}
\item[$\gswdes{1}$]
The functors $\Sigma^{\infty}$ and
$\Omega^{\infty}$ induce adjoint functors 
\[
\lder\Sigma^{\infty}:\ho\ugspaces{G} \leftrightarrows \ho \ugspectra{G}: \rder\Omega^{\infty}.
\]
on the homotopy categories. 

\item[$\gswdes{2}$] The symmetric monoidal structure on
$\ugspectra{G}$ induces a closed symmetric monoidal structure on the
homotopy category $\ho\ugspectra{G}$ and the functor
$\lder\Sigma^{\infty}$ is symmetric monoidal.

\item[$\gswdes{3}$]   The functor
$\lder\Sigma^{\infty}$ extends to a fully faithful, symmetric monoidal
embedding of $\swg{G}$ into $\ho\ugspectra{G}$.

\item [$\gswdes{4}$] The objects $S^{V}$ are invertible in
$\ho\ugspectra{G}$ under the smash product, so in particular the above
embedding of $\swg{G}$ extends to an embedding of the extended
Spanier-Whitehead category.

\item [$\gswdes{5}$] Arbitrary coproducts (denoted $\vee$) exist in
$\ho\ugspectra{G}$ and can be computed by the formation of wedges.  If
$\{X_{\alpha} \}$ is a collection of objects of $\ugspectra{G}$ and
$K$ is a finite $G$-CW complex, then the map
\[
\bigoplus_{\alpha}\ho\ugspectra{G}(K,X_{\alpha}) \to 
\ho\ugspectra{G}(K,\bigvee_{\alpha}X_{\alpha}) 
\]
is an isomorphism.

\item[$\gswdes{6}$]
Up to weak equivalence every
object $X$ is presentable in $\ugspectra{G}$ as a homotopy colimit
\[
\cdots\to S^{-V_{n}}\wedge X_{V_{n}}\to S^{-V_{n+1}}\wedge
X_{V_{n+1}}\to\cdots,
\]
in which $\{V_{n}\}$ is a fixed increasing sequence of representations
eventually containing every finite dimensional representation of $G$,
and each $X_{V_{n}}$ is weakly equivalent to an object of the form
$\Sigma^{\infty}K_{V_{n}}$, with $K_{V_{n}}$ a $G$-CW complex.
\end{itemize}

These properties aren't meant to constitute a characterization of
$\ugspectra{G}$, though they nearly do.  The first five insist that
$\ugspectra{G}$ not be too small, and the last that it not be too big.
Combined they show that, any computation one wishes to make in
$\ho\ugspectra{G}$ can, in principle, be reduced to a computation in
$\swg{G}$.

In all of the common models, and in particular in equivariant
orthogonal spectra, the presentation $\gswdes{6}$ is functorial.  We
call this the {\em canonical homotopy presentation}.  It is described
in detail in \S\ref{sec:canon-pres}.  For many purposes one can ignore
most of the technical details of equivariant spectra, and just think
in terms of the canonical homotopy presentation.

Finally, unless the emphasis is on foundations, we will drop the
$\lder$ and $\rder$ and implicitly assume that all of the functors
have been derived, unless otherwise specified.

\subsubsection{Equivariant orthogonal spectra}
\label{sec:equiv-orth-spectra-1}

An {\em orthogonal $G$-spectrum} consists of a collection of pointed
$G$-spaces $X_{V}$ indexed by the finite dimensional orthogonal
representations $V$ of $G$, an action of the orthogonal group $O(V)$
(of non-equivariant maps) on $X_{V}$, and for each (not necessarily
$G$-equivariant) orthogonal inclusion $t:V\subset W$ a map
$S^{W-t(V)}\wedge X_{V}\to X_{W}$, in which $W-t(V)$ denotes the
orthogonal complement of the image of $V$ in $W$.  These maps are
required to be compatible with the actions of $G$ and $O(V)$.  Maps of
equivariant orthogonal spectra are defined in the evident manner.  For
a more careful and detailed description
see~\ref{sec:orthogonal-spectra}.

Depending on the context, we will refer to orthogonal $G$-spectra as
``equivariant orthogonal spectra,'' ``orthogonal spectra,''
``$G$-spectra,'' and sometimes just as ``spectra''.

As with $G$-spaces, there are two useful ways of making the collection
of $G$-spectra into a category.  There is the topological category
$\ugspectra{G}$ just described, and the $G$-equivariant topological
category $\gspectra{G}$ of equivariant orthogonal spectra and
$G$-spaces of non-equivariant maps.  Thus for equivariant orthogonal
spectra $X$ and $Y$ there is an identification
\[
\ugspectra{G}(X,Y) = \gspectra{G}(X,Y)^{G}.
\]
We will use the abbreviated notation $\spectra$ to denote
$\ugspectra{G}$ when $G$ is the trivial group.

If $V$ and $W$ are two orthogonal representations of $G$ the same
dimension, and $O(V,W)$ is the $G$-space of (not necessarily
equivariant) orthogonal maps, then
\[
O(V,W)_{+}\underset{O(V)}{\wedge}X_{V}\to X_{W}
\]
is a $G$-equivariant homeomorphism.  In particular an orthogonal
$G$-spectrum $X$ is determined by the $X_{V}$ with $V$ a trivial
$G$-representation.  This implies that the category $\ugspectra{G}$ is
equivalent to the category of objects in $\spectra$ equipped with a
$G$-action (Proposition~\ref{thm:143}).

Both $\ugspectra{G}$ and $\gspectra{G}$ are tensored and cotensored
over $G$-spaces:
\begin{align*}
(X\wedge K)_{V} &= X_{V}\wedge K \\
\left(X^{K} \right)_{V} &= (X_{V})^{K}.
\end{align*}
Both categories are complete and cocomplete.

\begin{defin}
\label{def:6} The {\em suspension} and {\em $0$-space} functors are
defined by
\begin{align*}
(\Sigma^{\infty}K)_{V} &=S^{V}\wedge K \\
\Omega^{\infty}X &= X_{\{0\}}
\end{align*}
where $\{0 \}$ is the zero vector space.
\end{defin}

The suspension spectrum functor is left adjoint to the $0$-space
functor.  One has $\Sigma^{\infty}K = S^{0}\wedge K$ and more
generally $\Sigma^{\infty}(K\wedge L) = \left(\Sigma^\infty
K\right)\wedge L$.  The functors $\Sigma^{\infty}$ and
$\Omega^{\infty}$ may be regarded as topological functors between
$\ugspaces{G}$ and $\ugspectra{G}$ or as $\ugspaces{G}$-enriched
functors relating $\gspaces{G}$ and $\gspectra{G}$.

For each $G$-representation $V$ there is a $G$-spectrum $S^{-V}$
characterized by the existence of a functorial equivariant isomorphism
\begin{equation}
\label{eq:167}
\gspectra{G}(S^{-V},X) \approx X_{V}
\end{equation}
(see \S\ref{sec:orthogonal-spectra}).  By the enriched Yoneda Lemma,
every equivariant orthogonal $G$-spectrum $X$ is functorially
expressed as a reflexive coequalizer
\begin{equation}
\bigvee_{V,W} S^{-W}\wedge \gspectra{G}(S^{-W},S^{-V})\wedge
X_{V}\rightrightarrows
\bigvee_{V} S^{-V}\wedge X_{V} \to X.
\end{equation}
We call this the {\em tautological presentation} of $X$.

The category $\ugspectra{G}$ is a closed symmetric monoidal category
under the {\em smash product operation.}  The tensor unit is the sphere
spectrum $S^{0}$.    
There are canonical identifications 
\[
S^{-V}\wedge S^{-W}\approx S^{-V\oplus W},
\]
and in fact the association 
\[
V\mapsto S^{-V}
\]
is a symmetric monoidal functor from the category of finite
dimensional representations of $G$ (and isomorphisms) to
$\ugspectra{G}$.  Because of the tautological presentation, this
actually determines the smash product functor (see
\S\ref{sec:smash-prod-orth}).

Regarding the adjoint functors
\[
\Sigma^{\infty}:\ugspaces{G} \leftrightarrows \ugspectra{G}:
\Omega^{\infty},
\]
the left adjoint $\Sigma^{\infty}$
is symmetric monoidal.  We will usually drop the $\Sigma^{\infty}$ and
either not distinguish in notation between the suspension spectrum of
a $G$-space and the $G$-space itself, or use $S^{0}\wedge K$.

\subsubsection{Change of group and indexed monoidal products}
\label{sec:change-group-2} 

The fact that the category $\ugspectra{G}$ is equivalent to the
category of objects in $\spectra$ equipped with a $G$-action has an
important and useful consequence.  It means that if a construction
involving spectra happens to produce something with a $G$-action, it
defines a functor with values in $G$-spectra.  For example, if
$H\subset G$ is a subgroup, there is a restriction functor
$i_{H}^{\ast}:\ugspectra{G}\to \ugspectra{H}$ given by simply
restricting the action to $H$.  This functor has both a left and a
right adjoint.  The left adjoint is given by
\[
X \mapsto  G_{+}\underset{H}{\wedge}X
\]
and may be written as a ``wedge''
\[
\bigvee_{i\in G/H}X_{i}
\]
where $X_{i} = (H_{i})_{+}\underset{H}{\wedge}{X}$ with $H_{i}\subset
G$ the coset indexed by $i$.  Similarly, the right adjoint is given by
the $H$-fixed points of the internal function spectrum from $G$ to
$X$, and may be written as a kind of product 
\[
\prod_{i\in H\backslash G}X^{i} \approx
\prod_{i\in G/H}X_{i}.
\]
where $X^{i}= \hom^{H}(H^{i},X)$ and $H^{i}$ is the left $H$-coset
with index $i$.  The identification of the two expressions is made
using the map $g\mapsto g^{-1}$.  There is also an analogous
construction involving the smash product
\[
\norm_{H}^{G}X=\bigwedge_{i\in G/H}X_{i}
\]
These are special cases of a more general construction.

Suppose that $G$ is a finite group, and $J$ is a finite set on which
$G$ acts.  Write $\cat{B}_{J}G$ for the category with object set $J$,
in which a map from $j$ to $j'$ is an element $g\in G$ with $g\cdot
j=j'$.  We abbreviate this to $\cat{B}G$ in case $J=\text{pt}$.  Given
a functor
\[
X:\cat{B}_{J}G\to \spectra
\]
define the {\em indexed wedge}, {\em indexed product} and {\em indexed
smash product} of $X$ to be
\[
\bigvee_{j\in J} X_{j}, \quad\quad
\prod_{j\in J} X_{j}, \quad\text{and}\quad
\bigwedge_{j\in J}X_{j}
\]
respectively.  The group $G$ acts naturally on the indexed wedge and
indexed smash product and so they define functors from the category of
$\cat{B}_{J}G$-diagrams of spectra to $\ugspectra{G}$.  For more
details, see \S\ref{sec:index-mono-prod-1}.

Suppose that $H$ is a subgroup of $G$ and $J=G/H$.  In this case the
inclusion $\cat{B}_{\{e \}}H\to \cat{B}_{J}G$ of the full subcategory
containing the identity coset is an equivalence.  The restriction
functor and its left Kan extension therefore give an equivalence of
the category of $\cat{B}_{J}G$-diagrams of spectra with
$\ugspectra{H}$.  Under this equivalence, the indexed wedge works out
to be the functor
\[
G_{+}\smashove{H}(\slot).
\]
The indexed smash product is the {\em norm functor} 
\[
\norm_{H}^{G}:\ugspectra{H}\to \ugspectra{G},
\]
sending an $H$-spectrum $X$ to the $G$-spectrum
\[
\bigwedge_{j\in G/H}X_{j}.
\]

\begin{rem}
\label{rem:56} When the context is clear, we will sometimes abbreviate
the $\norm_{H}^{G}$ simply to $\norm$ in order to avoid clustering of
symbols.
\end{rem}

The norm distributes over wedges in much the same way as the iterated
smash product.  A precise statement of the general ``distributive
law'' appears in \S\ref{sec:distributive-laws}.

The functor $\norm_{H}^{G}$ is symmetric monoidal, commutes with
sifted colimits, and so filtered colimits and reflexive coequalizers
(Proposition~\ref{thm:150}).  The fact that $V\mapsto S^{-V}$ is
symmetric monoidal implies that
\begin{equation}
\label{eq:103}
\norm_{H}^{G}S^{-V} = S^{-\ind_{H}^{G}V},
\end{equation}
where $\ind_{H}^{G}V$ is the induced representation.   From the
definition, one also concludes that for a pointed $G$-space $T$,
\[
\norm_{H}^{G}\left(S^{-V}\wedge T\right) = S^{-\ind_{H}^{G}V}\wedge \norm_{H}^{G}T,
\]
where $\norm_{H}^{G}T$ is the analogous norm functor on spaces.

Th norm first appeared in group cohomology (Evens~\cite{MR0153725}),
and is often referred to as the ``Evens transfer'' or the ``norm
transfer.''  The analogue in stable homotopy theory originates in
Greenlees-May~\cite{MR1491447}.

\subsubsection{Stable weak equivalences} 
\label{sec:homot-theory-ugsp}

The inequality of Definition~\ref{def:28} gives the collection of
finite dimensional orthogonal $G$-representations the structure of a
(large) partially ordered set.  When $V_{1}$ is the trivial
representation of dimension $k$ the condition $V_{2}>V_{1}$ means that
\begin{equation}
\label{eq:15}
\dim V_{2}^{G}>k+1,
\end{equation}
and we will use instead the abbreviation $V_{2}>k$.
Using~\eqref{eq:15} we extend this to all $k\in\Z$.

Suppose we are given $X\in\ugspectra{G}$, $K\in\ugspaces{G}$, and two
representations $V_{1}<V_{2}$.  Choose an equivariant isometric
embedding $t:V_{1}\to V_{2}$ and let $W$ be the orthogonal complement
of $t(V_{1})$ in $V_{2}$.  Define
\begin{equation}
\label{eq:4}
[S^{V_{1}}\wedge K,X_{V_{1}}]^{G} \to 
[S^{V_{2}}\wedge K,X_{V_{2}}]^{G}
\end{equation}
by using the identification $S^{W}\wedge S^{V_{1}}\approx S^{V_{2}}$
and the structure map $S^{W}\wedge X_{V_{1}}\to X_{V_{2}}$ to form the composite 
\[
[S^{V_{1}}\wedge K,X_{V_{1}}]^{G}\to
[S^{W}\wedge S^{V_{1}}\wedge K,S^{W}\wedge X_{V_{1}}]^{G}\to
[S^{V_{2}}\wedge K,X_{V_{2}}]^{G}.
\]
This map depends only on the path component of $t$ in
$O(V_{1},V_{2})^{G}$, so the condition $V_{1}<V_{2}$ implies
that~\eqref{eq:4} is independent of the choice of $t$.

\begin{defin}
Let $X$ be a $G$-spectrum and $k\in\Z$.  For $H\subset G$ the {\em
$H$-equivariant $k^{\text{th}}$ stable homotopy group of $X$} is the
group
\begin{equation}
\label{eq:23}
\pi_{k}^{H}X = \varinjlim_{V> -k}\pi^{H}_{V+k}X_{V},
\end{equation}
in which the colimit is taken over the partial ordered set of
orthogonal $G$-representations $V$ satisfying $V>-k$.
\end{defin}

The poset of $G$-representations is a class, not a set, so one must
check that the colimit~\eqref{eq:23} actually exists.

\begin{defin}
\label{def:8}
An increasing sequence $V_{n}\subset V_{n+1}\subset\cdots$ of
finite dimensional representations of $G$ is {\em exhausting} if any
finite dimensional representation $V$ of $G$ admits an equivariant
embedding in some $V_{n}$.
\end{defin}

Any exhausting sequence $\cdots\subset V_{n}\subset
V_{n+1}\subset\cdots$ is final in the poset of $G$-representations, so
the map
\[
\varinjlim_{n}\pi^{H}_{V_{n}+k}X_{V_{n}} \to 
\varinjlim_{V> -k}\pi^{H}_{V+k}X_{V}
\]
is an isomorphism.   This gives the existence of the
colimit~\eqref{eq:23}, and shows that $\pi_{k}^{H}X$ can be computed as 
\[
\pi_{k}^{H}X = \varinjlim_{n}\pi^{H}_{V_{n}+k}X_{V_{n}} 
\]
in which $\cdots\subset V_{n}\subset V_{n+1}\subset\cdots$ is any
choice of exhausting sequence.

\begin{defin}
A {\em stable weak equivalence} (or just {\em weak
equivalence}, for short) is a map $X\to Y$ in $\ugspectra{G}$
inducing an isomorphism of stable homotopy groups $\pi^{H}_{k}$ for
all $k\in \Z$ and $H\subset G$.
\end{defin}

Equipped with the stable weak equivalences, the category
$\ugspectra{G}$ becomes a {\em homotopical category} in the sense
of~\cite{dwyer04:_homot}, and so both the homotopy category
$\ho\ugspectra{G}$ and the functor $\ugspectra{G}\to\ho\ugspectra{G}$
are defined.  As with $G$-spaces, we will often employ the notation
\[
[X,Y]^{G}
\]
for $\ho\ugspectra{G}(X,Y)$.  See~\S\ref{sec:categories-with-weak-2}
for more on the theory of homotopical categories, and for an
explanation of the notion of left ($\lder$) and right ($\rder$)
derived functors appearing in the discussion below.

\subsubsection{Properties $\gswdes{1}$-$\gswdes{6}$}
\label{sec:prop-gswd-gswd} 

We now describe how properties $\gswdes{1}$-$\gswdes{6}$ are verified,
deferring most of the technical details to
Appendix~\ref{sec:homot-theory-equiv}.  The first five properties
assert things only about the homotopy category and, save the fact that
the symmetric monoidal structure is closed, they can be established
using only the language of homotopical categories.

For $\gswdes{1}$, one checks directly from the definition that the
functor $\Sigma^{\infty}$ preserves weak
equivalences between $G$-spaces with non-degenerate base points, so
that $\lder\Sigma^{\infty}X$ can be computed as $\Sigma^{\infty}X$ if
$X$ has a non-degenerate base point, or as $\Sigma^{\infty}\tilde X$
in general, where $\tilde X$ is formed from $X$ by adding a whisker at
the base point.  The right derived functor $\rder\Omega^{\infty}$ is
given by choosing any exhausting sequence and forming
\[
\rder\Omega^{\infty}X = \hocolim \Omega^{V_{n}}X_{V_{n}}
\]
(See Proposition~\ref{thm:258}).  Verifying that
$\lder\Sigma^{\infty}$ and $\rder\Omega^{\infty}$ are adjoint functors
makes use of formula~\eqref{eq:6} below.

Regarding the symmetric monoidal structure ($\gswdes{2}$), the smash
product is not known to preserve weak equivalences between all objects
but it does so on the full subcategory of
$\ugspectra{G}\times\ugspectra{G}$ consisting of pairs $(X,Y)$ for
which one of $X$ or $Y$ is {\em cellular} in the sense that it
constructed inductively, starting with $\ast$ and attaching cells of
the form $G_{+}\underset{H}{\wedge}S^{-V}\wedge D^{m}_{+}$, with $V$ a
representation of $H$.  Every $G$-spectrum receives a functorial weak
equivalence from a cellular object, so this is enough to induce a
symmetric monoidal smash product on $\ho\ugspectra{G}$.  See
\S\ref{sec:weak-equiv-smash}.  The fact that the symmetric monoidal
structure is closed is best understood in the context of model
categories.  See \S\ref{sec:furth-homot-prop}, and especially
Corollary~\ref{thm:113}.

For $\gswdes{3}$, there is a useful formula for maps
in $\ho\ugspectra{G}$ in good cases.  Choose an exhausting sequence
$\{V_{n} \}$.  For $K$ a finite $G$-CW complex, $\ell\in\Z$, and any
$Y\in\ugspectra{G}$ the definition of stable weak equivalence and some
elementary facts about homotopical categories lead to the formula
(Proposition~\ref{thm-copied:244})
\begin{equation}
\label{eq:6}
\ho\ugspectra{G}(S^{\ell}\wedge K,Y) = \varinjlim_{n} [S^{V_{n}+\ell}\wedge
K, Y_{V_{n}}]^{G},
\end{equation}
Using this one easily checks that functor $K\mapsto S^{0}\wedge K$
extends to a symmetric monoidal functor $\swg{G}\to\ho\ugspectra{G}$.
A little more work gives the generalization
(Proposition~\ref{thm:161})
\begin{equation}
\label{eq:180}
\ho\ugspectra{G}(S^{-V}\wedge K,Y) = \varinjlim_{n} [S^{V_{n}}\wedge
K, Y_{V\oplus V_{n}}]^{G},
\end{equation}
in which $V$ a representation of $G$, 

For any representations $V$, $W$ of $G$, and any $X\in\ugspectra{G}$,
the map
\[
S^{-V} \wedge S^{V}\wedge X\to X
\]
is a weak equivalence (Proposition~\ref{thm:159}).   This ultimately
implies that $S^{V}$ is invertible in $\ho\ugspectra{G}$
(Corollary~\ref{thm:160}).  This establishes $\gswdes{4}$.

The fact that the formation of arbitrary wedges preserves weak
equivalences gives the first part of Property~$\gswdes{5}$
(Corollary~\ref{thm:90}).  The second part follows from~\eqref{eq:6}.

The canonical homotopy presentation of
$\gswdes{6}$ is constructed by choosing an exhausting sequence
$V=\{V_{n} \}$ and letting $X_{n}$ be an equivariant CW approximation
to $X_{V_{n}}$.  Since it involves more than just the homotopy
category, the construction is easier to describe with a model category
structure in place.  For the details see \S\ref{sec:canon-pres}.  

Indexed monoidal products have convenient homotopy properties in
$\ugspectra{G}$.  The formation of indexed wedges is {\em homotopical}
(Definition~\ref{def-copied:56}), in the sense that it preserves weak
equivalences.  This means that it need not be derived.  The same is
true of the formation of {\em finite} indexed products.  The map
\[
\bigvee_{\alpha} X_{\alpha} \to 
\prod_{\alpha} X_{\alpha} 
\]
from a finite indexed wedge to a finite indexed product is always a
weak equivalence.  This means in particular that for $H\subset G$, the
map from the left to the right adjoint of the restriction functor
\[
\ugspectra{G}\to \ugspectra{H}
\]
is always a weak equivalence.  Thus for $X\in\ugspectra{G}$ and
$Y\in\ugspectra{H}$ there are isomorphisms
\begin{equation}
\label{eq:21}
[X,G_{+}\underset{H}{\wedge}Y ]^{G} \approx
\big[X,\prod_{i\in G/H}Y_{i} \big]^{G} \approx
[i_{H}^{\ast} X,Y ]^{H}.
\end{equation}
The composite is the {\em Wirthm\"uller isomorphism}.  Because of it,
the right adjoint to the restriction functor tends not to appear
explicitly when discussing the homotopy category.

Up to weak equivalence, indexed smash products can be computed using
cellular approximations.  Combining this with the properties of the
norm listed in \S\ref{sec:change-group-2} leads to a useful
description of $\norm_{H}^{G}X$ in terms of the canonical homotopy
presentation
\[
\norm_{H}^{G}X = \hocolim_{V_{n}} S^{-\ind_{H}^{G}V_{n}}\wedge
\norm_{H}^{G}X_{V_{n}}.
\]

Finally, note that the formula~\eqref{eq:6} also implies that for any
$k\in \Z$
\[
\pi_{k}^{H}X \approx \ho \ugspectra{H}(S^{k},X) \approx \ho\ugspectra{G}(G_{+}\underset{H}{\wedge}S^{k},X)
\]
where for $k>0$ $S^{-k}$ is defined to be $S^{-\R^{k}}$ with $\R^{k}$
the trivial representation.

\subsubsection{The model structure}
\label{sec:model-structure}

Not all of the functors one wishes to consider have convenient
homotopy theoretic properties.

\begin{eg}
\label{eg:14}
For a $G$-spectrum $X$ let
\[
\sym^{n}X = X^{\wedge n}/\Sigma_{n}
\]
be the orbit spectrum of the $n$-fold iterated smash product by the action of
the symmetric group.  The map 
\[
S^{-1}\wedge S^{1}\to S^{0}
\]
is a weak equivalence.   However, the induced map 
\[
\sym^{n}(S^{-1}\wedge S^{1})\to \sym^{n}S^{0}
\]
is not.  The right side is $S^{0}$ since it is the tensor unit, while
the left side works out (Proposition~\ref{thm:200}) to be weakly
equivalent to the suspension spectrum of classifying space for
$G$-equivariant principal $\Sigma_{n}$-bundles, pointed by a disjoint
basepoint.
\end{eg}

In order to go further it is useful to refine the homotopical category
structure on $\ugspectra{G}$ to a model category.  Let $\gencof$ be
the set of maps
\begin{equation}
\gencof = \{G_{+}\smashove{H} S^{-V}\wedge S^{n-1}_{+}\to
G_{+}\smashove{H} S^{-V}\wedge D^{n}_{+} \}
\end{equation}
in which $H\subset G$ is a subgroup of $G$ and $V$ is a representation
of $H$ {\em containing a non-zero invariant vector}.     The set
$\gencof$ is the set of generating cofibrations in the {\em
positive complete} model structure on $\ugspectra{G}$.     The weak
equivalences are the stable weak equivalences and the fibrations are
the maps having the right lifting property against the acyclic
cofibrations.   See~\S\ref{sec:strong-posit-stable} for more details.

It works out that the symmetric power construction is homotopical on
the class of cofibrant objects in the positive complete model
structure (Proposition~\ref{thm:198}).

\begin{rem}
\label{rem:30} The condition that $V$ contains a non-zero invariant
vector is the {\em positivity} condition.  It is due to Jeff Smith,
and arose first in the theory of symmetric spectra (see the paragraph
following Corollary~0.6 in~\cite{MR1806878}).  The choice is dictated
by two requirements.  One is that symmetric power construction sends
weak equivalences between cofibrant objects to weak equivalences.
This is the key point in showing that the forgetful functor from
commutative algebras in $\ugspectra{G}$ to $\ugspectra{G}$ creates a
model category structure on commutative algebras in $\ugspectra{G}$
(Proposition~\ref{thm:114}).  The other is that the geometric fixed
point functor (\S\ref{sec:geom-fixed-points}) preserves (acyclic)
cofibrations.  The first requirement could be met by replacing
``positive'' with $\dim V>0$.  The second requires $\dim V^{G}>0$,
once one is using a positive model structure on $\spectra$.
\end{rem}

\subsubsection{Virtual representation spheres and \texorpdfstring{$\RO(G)$}{RO(G)}-graded cohomology}
\label{sec:smash-product-stable}

Using the spectra $S^{-V_{1}}$ and the spaces $S^{V_{0}}$ one can
associate a stable ``sphere'' to each virtual representation $V$ of
$G$.  To do so, one first represents $V$ as difference
$[V_{0}]-[V_{1}]$ of representations, and then sets
\[
S^{V} = S^{-V_{1}}\wedge S^{V_{0}}.  
\]
If $(V_{0}, V_{1})$ and $(W_{0}, W_{1})$ are two pairs of
orthogonal representations representing the same virtual
representation
\[
V=[V_{0}]-[V_{1}] = [W_{0}]-[W_{1}]\in \RO(G),
\]
then there is a representation $U$, and an equivariant orthogonal
isomorphism
\[
W_{1}\oplus V_{0}\oplus U\approx
V_{1}\oplus W_{0}\oplus U. 
\]
A choice of such data gives weak equivalences 
\begin{multline*}
S^{-W_{1}}\wedge S^{W_{0}} \leftarrow
S^{-W_{1}\oplus V_{0}\oplus U}\wedge S^{W_{0}\oplus V_{0}\oplus U}  \\
\approx
S^{-V_{1}\oplus W_{0}\oplus U}\wedge S^{W_{0}\oplus V_{0}\oplus U} \rightarrow
S^{-V_{1}}\wedge S^{V_{0}}
\end{multline*}
Thus, up to weak equivalence 
\[
S^{V}=S^{-V_{1}}\wedge S^{V_{0}} 
\]
depends only on $V$.  However, the weak equivalence between the
spheres arising from different choices depends on data not specified
in the notation.  This leads to some subtleties in grading equivariant
stable homotopy groups over the real representation ring $\RO(G)$.
See~\cite[\S6]{MR764596} and~\cite[Chapter XIII]{MR1413302}.  The
virtual representation spheres arising in this paper always occur as
explicit differences of actual representations.

In the positive complete model structure, the spectrum
$S^{-V_{1}}\wedge S^{V_{0}}$ will be cofibrant if and only if the
dimension of the fixed point space $V_{1}^{G}$ is positive.

\begin{defin}
\label{def:39} Suppose that $V$ is a virtual representation of $G$.  A
{\em positive representative} of $V$ consists of a pair of
representations $(V_{0}, V_{1})$ with $\dim V_{1}^{G}>0$ and for which
\[
V=[V_{0}]-[V_{1}] \in \RO(G).
\]
\end{defin}

Associated to every subgroup $H\subset G$ and every representation
$V\in \RO(H)$ is the group
\[
\pi_{V}^{H}(X) = [S^{V},X]^{H}.
\]
An equivariant cohomology theory is associated to every equivariant
orthogonal spectrum $E$, by
\begin{align*}
E^{k}(X) &= [X,S^{k}\wedge E]^{G} \\
E_{k}(X) &= [S^{k},E\wedge X]^{G} = \pi^{G}_{k}(E\wedge X). 
\end{align*}
There is also an $RO(G)$-graded version, defined by 
\begin{align*}
E^{V}(X) &= [X,S^{V}\wedge E]^{G} \\
E_{V}(X) &= [S^{V},E\wedge X]^{G} = \pi^{G}_{V}(E\wedge X).
\end{align*}

\subsection{Multiplicative properties}

\subsubsection{Commutative and associative
algebras} \label{sec:equiv-einfty-rings}

\begin{defin}
\label{def:15} An {\em equivariant commutative algebra} (or just {\em
commutative algebra}) is a unital commutative monoid in
$\ugspectra{G}$ with respect to the smash product operation.  An {\em
equivariant associative algebra} ({\em associative algebra}) is a
unital associative monoid with respect to the smash product.
\end{defin}

There is a weaker ``up to homotopy notion'' that sometimes comes up.

\begin{defin}
\label{def:16} An {\em equivariant homotopy associative algebra} (or
just {\em homotopy associative algebra}) is an associative algebra in
$\ho\ugspectra{G}$.  An {\em equivariant homotopy commutative algebra}
(or just {\em homotopy commutative algebra}) is a commutative algebra
in $\ho\ugspectra{G}$.
\end{defin}

The category of commutative algebras in $\ugspectra{G}$ and spaces of
equivariant multiplicative maps will be denoted $\ugeinftycat{G}$.
The $\ugspaces{G}$-enriched category of $G$-equivariant commutative
algebras and $G$-spaces of non-equivariant multiplicative maps will be
denoted $\geinftycat{G}$.  The categories $\ugeinftycat{G}$ and
$\geinftycat{G}$ are tensored and cotensored over $\spaces$ and
$\ugspaces{G}$ respectively.  The tensor product of an equivariant
commutative algebra $R$ and a $G$-space $T$ will be denoted
\[
R\otimes T
\]
to distinguish it from the smash product.  By definition
\[
\geinftycat{G}(R\otimes T, E) = \gspaces{G}(T,\geinftycat{G}(R,E)).
\]

We make $\ugeinftycat{G}$ into a homotopical category by defining a
map to be a weak equivalence if the underlying map of orthogonal
$G$-spectra is.  The {\em free commutative algebra functor}
\[
X\mapsto \sym(X)=\bigvee_{n\ge 0}\sym^{n}X,
\]
is left adjoint to the forgetful functor.  It takes weak equivalences
between cofibrant spectra to weak equivalences
(Proposition~\ref{thm:198}).  This is the key point in showing that
the forgetful functor
\[
\ugeinftycat{G}\to \ugspectra{G}
\]
creates a ($\spaces$-enriched) model category from the positive
complete model structure on $\ugspectra{G}$
(Proposition~\ref{thm:114}), and that
\[
\geinftycat{G}\to \gspectra{G}
\]
creates a $\gspaces{G}$-enriched model structure.    For $H\subset G$
the forgetful functor $\ugeinftycat{G}\to \ugeinftycat{H}$ and its
left adjoint form a Quillen morphism.   A similar set of results
applies to associative algebras.

Modules over an equivariant commutative ring are defined in the
evident way using the smash product.  The category of left modules
over $R$ and equivariant maps will be denoted $\rmod{R}$.  A map of
$R$-modules is defined to be a {\em weak equivalence} if the
underlying map of spectra is a weak equivalence.  The adjoint ``free
module'' and ``forgetful'' functors
\[
X\mapsto R\wedge X: \ugspectra{G} \leftrightarrows \rmod{R}: M\mapsto M
\]
create a model category structure on $\rmod{R}$.  It becomes an
enriched symmetric monoidal model category under the operation
\[
M\smashove{R}N
\]
where $M$ is regarded as a right $R$-module via
\[
M\wedge R\xrightarrow{\text{flip}} R\wedge M\to M,
\]
and $M\smashove{R}N$ is defined by the coequalizer diagram 
\[
M\wedge R\wedge N\rightrightarrows M\wedge N\to M\smashove{R}N.
\]

There are also the related notions of $\einfty$ and $\ainfty$
algebras.  It can be shown that the categories of $\einfty$ and
commutative algebras are Quillen equivalent, as are those of $\ainfty$
and associative algebras.

\subsubsection{Commutative algebras and indexed monoidal products}
\label{sec:norm-induction}

Because it is symmetric monoidal, the functor $\norm$ take commutative
algebras to commutative algebras, and so induces a functor
\[
\norm=\norm_{H}^{G}:\ugeinftycat{H}\to \ugeinftycat{G}.
\]
The following result is proved in the Appendices, as
Corollaries~\ref{thm:184} and~\ref{thm:17}.

\begin{prop}
\label{thm:119}
The functor
\[
\norm:\ugeinftycat{H}\to \ugeinftycat{G}.
\]
is left adjoint to the restriction functor $i^{\ast}$.  Together they
form a Quillen morphism of model categories. \qed
\end{prop}

\begin{cor}
\label{thm:122}
There is a natural
isomorphism
\[
\norm_{H}^{G}(i_{H}^{\ast} R) \to R\otimes (G/H),
\]
under which the counit of the adjunction is identified with the map
\[
R\otimes (G/H)\to R\otimes(\text{pt})
\]
given by the unique $G$-map $G/H\to \text{pt}$.  
\end{cor}

\begin{pf}
Since both $R\otimes (G/H)$ and the left adjoint to restriction
corepresent the same functor, this follows from
Proposition~\ref{thm:119}
\end{pf}

A useful consequence Corollary~\ref{thm:122} is that the group
$N(H)/H$ of $G$-automor\-phisms of $G/H$ acts naturally on $\norm_{H}^{G}
(i_{H}^{\ast}R)$.  The result below is used in the main computational
assertion of Proposition~\ref{thm:124}.

\begin{cor}
\label{thm:123}
For $\gamma\in N(H)/H$ the following diagram commutes:
\[
\xymatrix@C=1ex{
\norm_{H}^{G}(i_{H}^{\ast}R)\ar[dr] \ar[rr]^{\gamma} & & \norm_{H}^{G}(i_{H}^{\ast}R)\ar[dl] \\
&R&
}
\]
\end{cor}

\begin{pf}
Immediate from Corollary~\ref{thm:122}.
\end{pf}

At this point a serious technical issue arises.  The spectra
underlying commutative rings are almost never cofibrant.  This means
that there is no guarantee that the norm of a commutative ring has the
correct homotopy type.  The fact that it does is one of the main
results of Appendix~B.  The following is a consequence of
Proposition~\ref{thm:16}.

\begin{prop}
\label{thm:56}
Suppose that $R$ is a cofibrant commutative $H$-algebra, and $\tilde
R\to R$ is a cofibrant approximation of the underlying $H$-spectrum.
If $\tilde Z\to Z$ is a weak equivalence of $G$-spectra then
\[
\norm_{H}^{G}(\tilde R)\wedge \tilde Z\to \norm_{H}^{G}(R) \wedge Z
\]
is a weak equivalence.  \qed
\end{prop}

We refer to the property exhibited in Proposition~\ref{thm:56} by
saying that cofibrant commutative rings are {\em very flat}.

\subsubsection{Other uses of the norm}
\label{sec:other-uses-norm}

There are several important constructions derived from the norm
functor which also go by the name of ``the norm.''

Suppose that $R$ is a $G$-equivariant commutative ring spectrum, and
$X$ is an $H$-spectrum for a subgroup $H\subset G$.  Write
\[
R^{0}_{H}(X) = [X,i_{H}^{\ast}R]^{H}.
\]
There is a norm map
\[
\norm_{H}^{G}:R^{0}_{H}(X)  \to R^{0}_{G}(\norm_{H}^{G} X)
\]
defined by sending an $H$-equivariant map $X\to R$ to the composite
\[
\norm_{H}^{G} X\to \norm_{H}^{G} (i_{H}^{\ast}R)\to R,
\]
in which the second map is the counit of the restriction-norm
adjunction.  This is the {\em norm map on equivariant spectrum
cohomology}, and is the form in which the norm is described in
Greenlees-May~\cite{MR1491447}.   For an explicit comparison
with~\cite{MR1491447} see~\cite{MR3162261}.

When $V$ is a representation of $H$ and $X=S^{V}$ the above gives a
map
\[
\norm = \norm_{H}^{G}: \pi^{H}_{V} R \to \pi^{G}_{\ind V} R
\]
in which $\ind V$ is the induced representation which we call the {\em
norm map on the $RO(G)$-graded homotopy groups of commutative rings}

Now suppose that $X$ is a pointed $G$-space.  There is a norm map
\[
\norm_{H}^{G}:R^{0}_{H}(X)\to R^{0}_{G}(X)
\]
sending
\[
x\in R^{0}_{H}(X) = [S^{0}\wedge X, i_{H}^{\ast}R]^{H}
\]
to the composite
\[
S^{0}\wedge X \to S^{0}\wedge \norm(X) \approx \norm(S^{0}\wedge X)\to
\norm(i_{H}^{\ast}R)\to R,
\]
in which the equivariant map of pointed $G$-spaces
\[
X\to \norm_{H}^{G}(X)
\]
is the ``diagonal'' 
\[
X\to \prod_{j\in G/H}X_{j}\to \bigwedge_{j\in G/H} X_{j}
\]
whose $j^{\text{th}}$ component is the inverse to the isomorphism
\[
X_{j}= (H_{j})_{+}\smashove{H} X\to X
\]
given by the action map.   That this is actually equivariant is probably
most easily seen by making the identification
\[
X_{j} \approx \hom_{H}(H_{j}^{-1},X)
\]
in which $H_{j}^{-1}$ denotes the {\em left} $H$-coset consisting of
the inverses of the elements of $H_{j}$, and then writing
\[
\prod_{j\in G/H} X_{j} \approx \hom_{H}(G,X).
\]
Under this identification, the ``diagonal'' map  is the map
\[
X\to \hom_{H}(G,X)
\]
adjoint to the action map 
\[
G\underset{H}{\times}X\to X,
\]
which is clearly equivariant.   

One can combine these construction to define the {\em norm
on $RO(G)$-graded cohomology} of a $G$-space $X$
\[
\norm_{H}^{G}: R^{V}_{H}(X)\to R^{\ind V}_{G}(X)
\]
sending
\[
S^{0}\wedge X \xrightarrow{a}{} S^{V}\wedge i_{H}^{\ast}R
\]
to the composite
\[
S^{0}\wedge X\xrightarrow{}{} S^{0}\wedge \norm X \xrightarrow{\norm a}{} S^{\ind V}\wedge
\norm i_{H}^{\ast}R \to S^{\ind V}\wedge R.
\]

\subsection{The method of twisted monoid rings}
\label{sec:meth-poly-algebr}

In this section we describe the method of twisted monoid rings.   The
basic constructions are categorical, and in
\S\S\ref{sec:polynomial-algebras}-\ref{sec:weakly-comm-algebr} we do
not make any homotopy theoretic considerations.   In
\S\ref{sec:quotient-rings} we take up the homotopy theoretic aspects
of our constructions.

\subsubsection{Twisted monoid rings}
\label{sec:polynomial-algebras}

We start with a subgroup $H$ of $G$, and a positive representative
$(V_{0},V_{1})$ of a virtual representation $V$ of $H$.  Let
\[
S^{0}[S^{V}] = \bigvee_{k\ge 0} (S^{V})^{\wedge k}
\]
be the free $H$-equivariant associative algebra generated by $S^{V}=
S^{-V_{1}}\wedge S^{V_{0}}$, and
\[
\gbar{x} \in \pi^{H}_{V}S^{0}[S^{V}]
\]
the homotopy class of the generating inclusion.  When $|\gbar{x}|=0$,
the spectrum $S^{0}[S^{V}]$ is the monoid ring of the free monoid on
one generator, and is in fact commutative.  For general $\gbar{x}$ it
is the Thom spectrum of an associative monoid map from the free monoid
on one generator to the classifying space for $H$-equivariant real
vector bundles, hence a {\em twisted monoid ring}.  It
is not typically a commutative algebra, though the
$RO(H)$-equivariant homotopy groups make it appear so, since
$\pi^{H}_{\star}S^{0}[S^{V}]$ is a free module over
$\pi^{H}_{\star}S^{0}$ with basis $\{1,\gbar x, \gbar x^{2},\dots \}$.
It will be convenient to use the notation
\[
S^{0}[\gbar{x}] = S^{0}[S^{V}].
\]

Using the norm functor we can form the $G$-equivariant
twisted monoid ring
\[
\norm_{H}^{G}
(S^{0}[S^{V}]) = S^{0}[G_{+}\smashove{H}S^{V}].
\]
This spectrum can also be described as a Thom spectrum over the free
commutative monoid generated by $G/H$.  Things will look cleaner, and
better resemble the (polynomial) algebras we are modeling if we use
the alternate notation
\[
S^{0}[G\cdot S^{V}] \text{ and } S^{0}[G\cdot\gbar{x}].
\]
Though the symbol $H$ is omitted in this notation, it is still
referenced.  The representation $V$ is representation of $H$, and
$\gbar{x}$ is an $H$-equivariant map with domain $S^{V}$.

By smashing examples like these together we can make associative
algebras that are twisted forms of free commutative
monoid algebras over $S^{0}$, in which the group $G$ is allowed to act
on the monoid.  More explicitly, suppose we are given a sequence
(possibly infinite) of subgroups $H_{i}\subset G$ and for each $i$ a
positive representative $((V_{i})_{0}, (V_{i})_{1})$ of a virtual
representation $V_{i}$ of $H_{i}$.  For each $i$ form
\[
S^{0}[G \cdot \gbar{x}_{i}]
\]
as described above, smash the first $m$ together to make
\[
S^{0}[G\cdot\gbar{x}_{1},\dots,G\cdot\gbar{x}_{m}],
\]
and then pass to the  colimit to construct the $G$-equivariant
associative algebra
\[
T=S^{0}[G\cdot\gbar{x}_{1},G\cdot\gbar{x}_{2},\dots].
\]

The twisted monoid ring $T$ can also be described as
a Thom spectrum over the free commutative monoid generated by the
$G$-set
\[
J=\coprod_{i=1}^{\infty} G/H_{i}.
\]
By construction, it is an indexed smash product of an indexed wedge
\begin{equation}
\label{eq:8}
T = \bigwedge_{j\in J}\bigvee_{n=0}^{\infty} S^{n V(j)}
\end{equation}
where for $j= g H_{i}$, $V(j)$ is the virtual representation of $g
H_{i}g^{-1}$ with positive representative
\[
V(j) = \big(V(j)_{0},
V(j)_{1}\big) = 
(g H_{i})\underset{H_{i}}{\times}\big((V_{i})_{0}, (V_{i})_{1}\big)
\]

All of this can be done relative to an equivariant commutative algebra
$R$ by defining
\[
R[G\cdot\gbar{x}_{1},G\cdot\gbar{x}_{2},\dots]
\]
to be
\[
R\wedge S^{0}[G\cdot\gbar{x}_{1},G\cdot\gbar{x}_{2},\dots].
\]

Because they can fail to be commutative, these twisted
monoid rings do not have all of the algebraic properties one might
hope for.  But it is possible to naturally construct all of the
equivariant {\em monomial ideals}.  Here is how.

Applied to~\eqref{eq:8} the distributive law of
\S\ref{sec:distributive-laws} gives an isomorphism of $T$ with the
indexed wedge
\[
T=\bigvee_{f\in \nat^{J}} S^{V_{f}}
\]
in which $f$ is running through the set of 
functions 
\[
J\to \nat = \{0,1,2,\dots \}
\]
taking non-zero values on only finitely many elements ({\em finitely
supported} functions).  The group $G$ acts on the set $\nat^{J}$
through its action on $J$, and $V_{f}$ is the virtual representation
\[
V_{f}=\sum_{j\in J}f(j)\cdot V(j)
\]
of the stabilizer $H_{f}$ of $f$, with the evident positive representation
\[
\left(\bigoplus_{j\in J}V(j)_{0}^{f(j)},\bigoplus_{j\in J}V(j)_{1}^{f(j)}\right).
\]

The $G$-set $\nat^{J}$ is a commutative monoid under addition of
functions, and the ring structure on $T$ is the indexed sum of the
obvious isomorphisms
\[
S^{V_{f}}\wedge S^{V_{g}} \approx S^{V_{f}\oplus V_{g}}\approx
S^{V_{f+g}}.
\]

Recall (for example from~\cite{MR0132791}) that an ideal in a
commutative monoid $L$ is a subset $I\subset L$ with the property that
$L+I\subset I$.   Given a $G$-stable ideal $I\subset \nat^{J}$ form the
$G$-spectrum
\[
T_{I} = \bigvee_{f\in I} S^{V_{f}}.
\]
The formula for the multiplication in $T$ implies that $T_{I}$ is an
equivariant sub bimodule of $T$, and that the association $I\mapsto
T_{I}$ is an inclusion preserving function from the set of ideals in
$\nat^{J}$ to the set of sub-bimodules of $T$.  For a more general and
systematic discussion see \S\ref{sec:monomial-ideals}.  

\begin{eg}
\label{eg:11} The monomial ideal corresponding to the set $I$ of all
non-zero elements of $\nat^{J}$ is the augmentation ideal (up to
homotopy it is the fiber of the map $T\to S^{0}$).  It is convenient
to denote this $T$ bimodule as $(G\cdot\gbar{x}_{1},
G\cdot\gbar{x}_{2},\dots )$.  More generally, for an integer $n>0$ the
set $nI = I+\cdots+I$ of $n$-fold sums of elements of $I$ is a monoid
ideal.  It corresponds to the monomial ideal given by the
$n^{\text{th}}$ power of the augmentation ideal.
\end{eg}

\begin{eg}
\label{eg:9}
Let $\dim:\nat^{J}\to
\nat$ be the function given by
\[
\dim f = \dim V_{f}=\sum_{j\in J} f(j)\dim V_{j}.
\]
If for all $i$, $\dim V_{i}>0$, then the set $\{f\mid \dim f\ge d\}$
is a monoid ideal in $\nat^{J}$ and corresponds to the monomial ideal
$M\subset T$ consisting the wedge of spheres of dimension greater
than $d$.  The quotient bimodule $M_{d}/M_{d-1}$ can be identified
with the indexed coproduct
\[
\bigvee_{\dim f=d} S^{V_{f}}
\]
on which $T$ is acting through the augmentation $T\to S^{0}$.  These
monomial ideals play an important role in the proof of the Slice
Theorem in \S\ref{sec:slice-theorem-i}.
\end{eg}

\subsubsection{The method of twisted monoid rings}
\label{sec:weakly-comm-algebr}

\begin{defin}
\label{def:32}
Suppose that
\[
f_{i}:B_{i}\to R,\quad i=1,\dots m
\]
are algebra maps from associative algebra $B_{i}$ to a commutative
algebra $R$.  The {\em smash product} of the $f_{i}$ is the algebra
map
\[
\bigwedge^{m} f_{i} : \bigwedge^{m} B_{i} \to \bigwedge^{m}R \to R,
\]
in which the right most map is the iterated multiplication.  If $B$ is
an $H$-equivariant associative algebra, and $f:B\to i_{H}^{\ast}R$ is
an algebra map, we define the {\em norm of $f$} to be the $G$-equivariant algebra
map
\[
\norm_{H}^{G}B\to R
\]
given by
\[
\norm_{H}^{G} B\to \norm_{H}^{G}(i_{H}^{\ast}R)\to R,
\]
in which the rightmost map is the counit of the adjunction described
in Proposition~\ref{thm:119}.
\end{defin}

These constructions make it easy to map a twisted monoid ring to a
commutative algebra.  Suppose that $R$ is a fibrant $G$-equivariant
commutative algebra, and we're given a sequence
\[
\gbar{x}_{i}\in \pi^{H_{i}}_{V_{i}} R,\qquad i = 1,2,\dots.
\]
A choice of positive representative $((V_{0})_{i},(V_{1})_{i})$ of
$V_{i}$ and a map
\[
S^{V_{i}}\to R
\]
representing $\gbar{x}_{i}$ determines an associative algebra map
\[
S^{0}[\gbar{x}_{i}] \to R.
\]
Applying the norm gives a $G$-equivariant associative algebra map
\[
S^{0}[G\cdot\gbar{x}_{i}] \to R.
\]
By smashing these together we can make a sequence of equivariant
algebra maps
\[
S^{0}[G\cdot\gbar{x}_{1},\dots,G\cdot\gbar{x}_{m}] \to R.
\]
Passing to the colimit gives an equivariant algebra map
\begin{equation}\label{eq:56}
S^{0}[G\cdot\gbar{x}_{1},G\cdot\gbar{x}_{2},\dots] \to R
\end{equation}
representing the sequence $\gbar{x}_{i}$.  We will refer to this
process by saying that the map~\eqref{eq:56} is constructed by the
{\em method of twisted monoid rings}.  The whole construction can also
be made relative to a commutative algebra $S$, leading to an
$S$-algebra map
\begin{equation}\label{eq:37}
S[G\cdot\gbar{x}_{1},G\cdot\gbar{x}_{2},\dots] \to R
\end{equation}
when $R$ is a commutative $S$-algebra.

\subsubsection{Quotient modules}
\label{sec:quotient-rings}

One important construction in ordinary stable homotopy theory is the
formation of the quotient of a module $M$ over a commutative algebra
$R$ by the ideal generated by a sequence $\{x_{1},x_{2},\dots
\}\subset \pi_{\ast}R$.  This is done by inductively forming the
cofibration sequence of $R$-modules
\begin{equation}
\label{eq:57}
\Sigma^{|x_{n}|}M/(x_{1},\dots,x_{n-1}) \to
M/(x_{1},\dots,x_{n-1}) \to
M/(x_{1},\dots,x_{n})
\end{equation}
and passing to the homotopy colimit in the end.  There is an evident
equivalence
\[
M/(x_{1},\dots) \approx M\smashove{R} R/(x_{1},\dots)
\]
in case $M$ is a cofibrant $R$-module.  The situation is slightly
trickier in equivariant stable homotopy theory, where the group $G$
might act on the elements $x_{i}$, and prevent the inductive approach
described above.  The method of twisted monoid rings
(\S\ref{sec:polynomial-algebras}) can be used to get around this
difficulty.

Suppose that $R$ is a fibrant equivariant commutative algebra, and
that
\[
\gbar x_{i}\in \pi_{V_{i}}^{H_{i}}(R) \qquad i=1,2,\dots
\]
is a sequence of equivariant homotopy classes.  Using the method of
twisted monoid rings, construct an associative $R$-algebra map
\begin{equation}\label{eq:38}
T=R[G\cdot\gbar{x}_{1},G\cdot\gbar{x}_{2},\dots] \to R.
\end{equation}
Using this map, we may regard an equivariant $R$-module $M$ as a
$T$-module.   In addition to~\eqref{eq:38} we will make use of the
augmentation $\epsilon:T\to R$ sending the $\gbar{x}_{i}$ to zero.

\begin{defin}
The {\em quotient module $M/(G\cdot\gbar x_{1},\dots)$} is the
$R$-module 
\[
M\underset{T}{\overset\lder{\wedge}}R
\]
in which $T$ acts on $M$ through the map~\eqref{eq:38} and on $R$
through the augmentation.
\end{defin}

The symbol $\overset{\lder}{\wedge}$ denotes derived smash product.
By Proposition~\ref{thm:48} it can be computed
by taking a cofibrant approximation in either variable.  

Let us check that this construction reduces to the usual one when the
group acting is the trivial group and $M$ is a cofibrant $R$-module.
For ease of notation, write
\begin{align*}
T &=R[x_{1},\dots] \\
T_{n} &=R[x_{1},\dots, x_{n}].
\end{align*}
Using the isomorphism
\[
R[x_{1},\dots] \approx R[x_{1},\dots, x_{n}]\smashove{R} R[x_{n+1},\dots]
\]
one can construct an associative algebra map
\[
T \to  R[x_{n+1},\dots]
\]
by smashing the augmentation
\[
R[x_{1},\dots, x_{n}] \to R
\]
sending each $x_{i}$ to $0$, with the identity map of
$R[x_{n+1},\dots]$.  By construction, the evident map of $T$-algebras
\[
\varinjlim R[x_{n+1}, \dots] \to R
\]
is an isomorphism, and hence so is the map 
\[
\varinjlim M\underset{T}{\wedge}R[x_{n+1}, \dots] \to
M\underset{T}{\wedge}R.
\]

In fact this isomorphism is also a derived
equivalence.   To see this construct a sequence
\[
\to\cdots N_{n+1}\to N_{n+2}\cdots
\]
of cofibrations of cofibrant left $T$-module approximations to 
\[
\to\cdots R[x_{n+1},\dots] \to R[x_{n+2},\dots] \to \cdots.
\]
We have 
\[
\pi_{\ast}\varinjlim N_{k} \approx \varinjlim \pi_{\ast}N_{k} \approx
\varinjlim (\pi_{\ast}R)[x_{k},\dots] \approx R
\]
from which one concludes that the map 
\[
\varinjlim N_{k} \to \varinjlim R[x_{k},\dots]
\]
is a cofibrant approximation.  It follows that
\[
M/(x_{1},\dots) \approx 
\hocolim M/(x_{1},\dots, x_{n}).
\]

To compare $M/(x_{1},\dots, x_{n-1})$ with $M/(x_{1},\dots, x_{n})$
let $T_{n}\to R[x_{n}]$ be associative algebra map constructed from
the isomorphism
\[
T_{n}\approx T_{n-1}\underset{R}{\wedge}R[x_{n}].
\]
by smashing the augmentation of $T_{n-1}$ with the identity map of
$R[x_{n}]$.   We have
\[
M/(x_{1},\dots,x_{n-1}) \sim 
M\underset{T_{n-1}}{\wedge}R \approx
M\underset{T_{n}}{\wedge}T_{n}\underset{T_{n-1}}{\wedge}R \approx
M\underset{T_{n}}{\wedge}R[x_{n}].
\]
By Proposition~\ref{thm:48} $M\underset{T_{n}}{\wedge}R[x_{n}]$ is a
cofibrant $R[x_{n}]$-module.  The cofibration sequence~\eqref{eq:57}
is now constructed by applying the functor
\begin{equation}
\label{eq:80}
M/(x_{1},\dots, x_{n-1})\smashove{R[x_{n}]} (\slot).
\end{equation}
to the pushout diagram of $R[x_{n}]$ bimodules 
\begin{equation}
\label{eq:40}
\xymatrix{
(x_{n}) \ar[r]\ar[d]  & R[x_{n}]  \ar[d] \\
\ast \ar[r] & R
}
\end{equation}
and appealing to Corollary~\ref{thm:102}.

A similar discussion applies to the equivariant situation, giving 
\[
M/(G\cdot \gbar{x}_{1},\dots) \approx \varinjlim M/(G\cdot
\gbar{x}_{1},\dots, G\cdot \gbar{x}_{n}),
\]
a relation 
\[
M/(G\cdot \gbar{x}_{1},\dots, G\cdot \gbar{x}_{n}) \approx M/(G\cdot
\gbar{x}_{1},\dots, G\cdot \gbar{x}_{n-1})
\smashove{R[G\cdot\gbar{x}_{n}]}R,
\]
and a cofibration sequence 
\[
(G\cdot \gbar x_{n})\cdot M/(G\cdot \gbar{x}_{1},\dots, G\cdot
\gbar{x}_{n-1})\to M/(G\cdot \gbar{x}_{1},\dots, G\cdot
\gbar{x}_{n-1}) \to M/(G\cdot \gbar{x}_{1},\dots, G\cdot
\gbar{x}_{n}),
\]
derived by applying the functor 
\[
M/(G\cdot \gbar{x}_{1},\dots, G\cdot \gbar{x}_{n-1})
\smashove{R[G\cdot\gbar{x}_{n}]}(\slot)
\]
to
\[
(G\cdot \gbar x_{n})\to R[G\cdot\gbar{x}_{n}]\to R.
\]

One can also easily deduce the equivalences
\begin{equation}
\label{eq:64} 
R/(G\cdot\gbar{x}_{1},\dots,G\cdot\gbar{x}_{n}) \approx
R/(G\cdot\gbar{x}_{1})\smashove{R}\cdots\smashove{R}R/(G\cdot\gbar{x}_{1})
\end{equation}
and
\begin{equation}
\label{eq:65} R/(G\cdot\gbar{x}_{1},\dots) \approx \varinjlim
R/(G\cdot\gbar{x}_{1})\smashove{R}\cdots\smashove{R}
R/(G\cdot\gbar{x}_{n}).
\end{equation}
These expressions play an important role in the proof
Lemma~\ref{thm:181}, which is a key step in the proof the Reduction
Theorem.

\subsection{Fixed points, isotropy separation and geometric fixed points}

\subsubsection{Fixed point spectra}
\label{sec:fixed-point-spectra}

The {\em fixed point spectrum} of a $G$-spectrum $X$ is defined to be
the spectrum of $G$ fixed points in the underlying, non-equivariant
spectrum $i_{0}^{\ast}X$.    In other words it is given by
\[
X\mapsto (i_{0}^{\ast}X)^{G}.
\]
The notation $i_{0}^{\ast}X^{G}$ can get clumsy and we will usually
abbreviate it to $X^{G}$.  

The functor of fixed points has a left adjoint which sends
$S^{-V}\wedge X_{V}\in \spectra$ to $S^{-V}\wedge X_{V}\in
\ugspectra{G}$, where in the latter expression $V$ is regarded as a
representation of $G$ with trivial $G$-action and $X_{V}$ is regarded
as a space with trivial $G$-action.  It can be computed for general
$X$ in terms of the tautological presentation
\[
\bigvee_{V,W} \icat (V,W)_{+}\wedge S^{-W}\wedge
X_{V}\rightrightarrows
\bigvee_{V} S^{-V}\wedge X_{V} \to X
\]
for the
trivial group (see~\eqref{eq:32}), once one observes that 
\[
\icat (V,W) = 
\igcat G(V,W)
\] 
when $V$ and $W$ have trivial $G$-action.  

Under the equivalence between $\ugspectra{G}$ and the category of
objects in $\spectra$ equipped with a $G$-action, the fixed point
spectrum functor is formed by passing to objectwise fixed points, and
its left adjoint is given by regarding a non-equivariant spectrum as a
$G$-object with trivial $G$-action.

The fixed point functor and its left adjoint form a Quillen morphism
in the positive complete model structures.  Neither the fixed point
functor nor its left adjoint is homotopical and so both need to be
derived.  As one can easily check from the definition, if $X$ is
fibrant (or more generally has the property that for some exhausting
sequence $\{V_{n} \}$, the map $X_{V_{n}}\to
\Omega^{V_{n+1}-V_{n}}X_{V_{n+1}}$ is a weak equivalence) there is an
isomorphism
\[
\pi_{\ast}(X^{G}) \approx \pi^{G}_{\ast}X.
\]

The (derived) fixed point functor on spectra doesn't always have the
properties one might be led to expect by analogy with spaces.  For
example even though the composition of the fixed point functor with
its left adjoint is the identity, the composition of the derived
functors is not.  The derived fixed point functor does not generally
commute with smash products, or with the formation of suspension
spectra.

\subsubsection{Isotropy separation and geometric fixed points}
\label{sec:isotropy-separation}

A standard approach to getting at the equivariant homotopy type of a
$G$-spectrum $X$ is to nest $X$ between two pieces, one an aggregate
of information about the spectra $i_{H}^{\ast}X$ for all proper
subgroups $H\subset G$, and the other a localization of $X$ at a
``purely $G$'' part.  This is the {\em isotropy separation sequence}
of $X$.

More formally, let $\pfamily$ denote the family of proper subgroups of
$G$, and $E\pfamily$ the ``classifying space'' for $\pfamily$,
characterized up to equivariant weak equivalence by the property that
the space of fixed points $E\pfamily^{G}$ is empty, while for any
proper $H\subset G$, $E\pfamily^{H}$ is weakly contractible.  For
convenience we will assume that $E\pfamily$ has been chosen to be a
$G$-CW complex.  Such a model can be constructed as the join of
infinitely many copies of $G/H$ with $H$ ranging through the proper
subgroups of $G$.  It can also be constructed as the unit sphere in
the sum of infinitely many copies of the reduced regular
representation of $G$.  Any $G$-CW complex $E\pfamily$ admits an
equivariant cell decomposition into cells of the form $(G/H)_{+}\wedge
D^{n}_{+}$ with $H$ a proper subgroup of $G$.

Let $\tilde E\pfamily$ be the mapping cone of $E\pfamily\to
\text{pt}$, with the cone point taken as base point.
The $G$-CW complex $\tilde E\pfamily$ is characterized up to equivariant
homotopy equivalence by the property
\[
\big(\tilde E\pfamily\big)^{H} \sim 
\begin{cases}
S^{0}& H=G \\
\ast & H\ne G.
\end{cases}
\]
The important {\em isotropy separation sequence} is constructed by
smashing a $G$ spectrum $X$ with the defining cofibration
sequence for $\tilde E\pfamily$
\begin{equation}
\label{eq:3}
E\pfamily_{+}\wedge X \to X \to \tilde E\pfamily \wedge X.
\end{equation}
The term on the left can be described in terms of the action of proper
subgroups $H\subset G$ on $X$.  The homotopy type of the term on the
right is determined by its right derived fixed point spectrum
\[
\phig(X) = \big((\tilde E\pfamily\wedge X)_{f}\big)^{G},
\]
in which the subscript $f$ indicates a functorial fibrant replacement.
The functor $\phig(X)$ is the {\em geometric fixed point functor} and
has many remarkable properties.

\begin{prop}
\label{thm:202}
\begin{thmList}
\item The functor $\phig$ sends weak equivalences to weak equivalences.
\item The functor $\phig$ commutes with filtered homotopy colimits.
\item For a $G$-space $A$ and a representation $V$ of $G$ there is a
weak equivalence $\phig(S^{-V}\wedge A) \approx S^{-V^{G}}\wedge A^{G}$
where $V^{G}\subset V$ is the subspace of $G$-invariant vectors.
\item For $G$-spectra $X$ and $Y$ the spectra
\[
\phig(X\wedge Y) \quad\text{and}\quad \phig(X)\wedge \phig(Y).
\]
are related by a natural chain of weak equivalences.
\end{thmList}
\end{prop}

\begin{rem}
\label{rem:2}
Note that in terms of the canonical homotopy presentation
\[
X\approx \hocolim_{V} S^{-V}\wedge X_{V}
\]
properties~\thmListItem{1}-\thmListItem{3} of
Proposition~\ref{thm:202} imply that
\begin{equation}\label{eq:35}
\phig X \approx \hocolim_{V} S^{-V^{G}} \wedge X_{V}^{G}.
\end{equation}
\end{rem}

\begin{pf*}{Sketch of proof}
The first assertion follows from the fact that smashing with $\tilde
E\pfamily$ is homotopical (\S\ref{sec:weak-equiv-smash}), so need not
be derived, and that the fixed point functor is homotopical on the
full subcategory of fibrant objects.  The second is straightforward.
Part~\thmListItem{3} is Corollary~\ref{thm:205}.  By the remark
above, the canonical
homotopy presentation reduces part~\thmListItem{4} to the case
$X=S^{-V}\wedge A$, $Y=S^{-W}\wedge B$, with $A$ and $B$ $G$-CW
complexes.  One easily checks the assertion in this case using
part~\thmListItem{3}.
\end{pf*}

\begin{rem}
\label{rem:42} 
When $G=\ztn$, the space $E\pfamily$ is the space
$E\zt$ with $G$ acting through the epimorphism $G\to\zt$.   Taking
$S^{\infty}$ with the antipodal action as a model of $E\zt$, this
leads to an identification
\[
\tepfam \approx \lim_{n\to\infty} S^{n\sigma},
\]
in which $S^{n\sigma}$ denotes the one point compactification of the
direct sum of $n$ copies of the real sign representation of $G$.
\end{rem}

\begin{rem}
\label{rem:8} The isotropy separation sequence often leads to the
situation of needing to show that a map $X\to Y$ of cofibrant
$G$-spectra induces a weak equivalence
\[
\tilde E\pfamily\wedge X\to \tilde E\pfamily\wedge Y.
\]
Since for every proper $H\subset G$, $\pi^{H}_{\ast}\tilde E\pfamily
\wedge X=\pi^{H}_{\ast}\tilde E\pfamily\wedge Y=0$, this is
equivalent to showing that the map of geometric fixed point spectra
$\phig X\to \phig Y$ is a weak equivalence.
\end{rem}

\begin{rem}
\label{rem:20} Since for every proper $H\subset G$,
$\pi^{H}_{\ast}\tilde E\pfamily\wedge X=0$, it is also true that
\[
[T,\tilde E\pfamily\wedge X]^{G}_{\ast} = 0
\]
for every $G$-CW spectrum $T$ built entirely from $G$-cells of the
form $G_{+}\smashove{H}D^{n}$ with $H$ a proper subgroup of $G$.
Similarly, if $T$ is gotten from $T_{0}$ by attaching $G$-cells
induced from proper subgroups, then the restriction map
\[
[T,\tilde E\pfamily\wedge X]^{G}_{\ast} \to
[T_{0},\tilde E\pfamily\wedge X]^{G}_{\ast}
\]
is an isomorphism.  This holds, for example, if $T$ is the suspension
spectrum of a $G$-CW complex, and $T_{0}\subset T$ is the subcomplex
of $G$-fixed points.
\end{rem}

\begin{rem}
\label{rem:6} 
For a subgroup $H\subset G$ and a $G$-spectrum $X$, it
will be convenient to use the abbreviation
\[
\Phi^{H}X
\]
for the more correct $\Phi^{H}i_{H}^{\ast}X$.  This situation comes up
in our proof of the ``homotopy fixed point'' property of
Theorem~\ref{thm:52}, where the more compound notation becomes a
little unwieldy.
\end{rem}

We end this section with a simple result whose proof illustrates a
typical use of the geometric fixed point spectra.   

\begin{prop}
\label{thm:229}
Suppose that $X$ is a $G$-spectrum with the property that for all
$H\subset G$, the geometric fixed point spectrum $\phih X$ is
contractible.   Then $X$ is contractible as a $G$-spectrum.
\end{prop}

\begin{pf}
By induction on $|G|$ we may assume that for proper $H\subset G$, the
spectrum $i_{H}^{\ast}X$ is contractible.  Since both
$G_{+}\underset{H}{\wedge}(\slot)$ and the formation of mapping cones
are homotopical it follows that $T\wedge X$ is contractible for any
$G$-CW complex built entirely from cells of the form
$G_{+}\smashove{H}D^{n}$ with $H\subset G$ proper.  This applies in
particular to $T=\epfam$.  The isotropy separation sequence then shows
that
\[
X\to \tepfam\wedge X
\]
is a weak equivalence.  But Remark~\ref{rem:8} and our assumption that
$\phig X$ is contractible imply that $\tepfam\wedge X$ is contractible.
\end{pf}

\subsubsection{Monoidal geometric fixed points}
\label{sec:mono-geom-fixed} 

For some purposes it is useful to have a version of the geometric
fixed point functor which is lax symmetric monoidal.  For example,
such a functor automatically takes (commutative) algebras to
(commutative) algebras.  The geometric fixed point functor defined
in~\cite[\S{V.4}]{MR1922205} has this property.  We denote it $\phigm$
and refer to it as the {\em monoidal geometric fixed point functor} in
order to distinguish it from $\phig$.  The following is a compendium
of results from~\cite[\S{V.4}]{MR1922205}.  The construction and
proofs are described in \S\ref{sec:geom-fixed-points}.

\begin{prop}
\label{thm:203}
The monoidal geometric fixed point functor has the following
properties:
\begin{thmList}
\item It preserves acyclic cofibrations.
\item It is lax symmetric monoidal.
\item If $X$ and $Y$ are cofibrant, the map 
\[
\phigm(X)\wedge \phigm(Y)\to \phigm(X\wedge Y)
\]
is an isomorphism.
\item It commutes with cobase change along a closed inclusion.
\item It commutes with directed colimits.
\end{thmList}
\end{prop}

Property \thmListItem{3} implies that $\phigm$ is weakly symmetric monoidal in
the sense of the definition below.

\begin{defin}[Schwede-Shipley~\cite{MR1997322}]
\label{def:40} A functor $F:\cat C\to \cat D$ between (symmetric)
monoidal model categories is {\em weakly (symmetric) monoidal} if it
is lax (symmetric) monoidal, and the map
\[
F(X)\wedge F(Y)\to F(X\wedge Y)
\]
is a weak equivalence when $X$ and $Y$ are cofibrant.
\end{defin}

The next result is~\cite[Proposition V.4.17]{MR1922205}, and is
discussed in more detail as Proposition~\ref{thm:14}.

\begin{prop}
\label{thm:220} 
The left derived functor of $\phigm$ is $\phig$.  More
specifically, there are natural transformations
\[
\phig(X)\to \tilde{\Phi}^{G}_{M}(X)\xleftarrow{\sim} \phigm(X)
\]
in which the rightmost arrow is a always weak equivalence and the
leftmost arrow is a weak equivalence when $X$ is cofibrant.  \qed
\end{prop}

Because $\phig$ is lax monoidal, it determines functors
\begin{align*}
\phigm&:\ugainftycat{G}\to\ainftycat \\
\phigm&:\ugeinftycat{G}\to\einftycat,
\end{align*}
and for each associative algebra $R$ a functor 
\[
\phigm: \rmod{R}\to\rmod{\phigm R}.
\]
In addition, if $R$ is an associative algebra, $M$ a right $R$-module
and $N$ a left $R$-module there is a natural map
\begin{equation}
\label{eq:118}
\phigm(M\smashove{R}N) \to \phigm M\smashove{\phigm R} \phigm N.
\end{equation}
The argument for~\cite[Proposition V.4.7]{MR1922205} shows
that~\eqref{eq:118} is a weak equivalence (in fact an isomorphism) if
$M$ and $N$ are cofibrant and $R$ is ``cellular.'' See
Proposition~\ref{thm:174}.  (Recently, Blumberg and
Mandell~\cite[Appendix~A]{blumberg_cylcotomic} have shown that one
need only require one of $M$ or $N$ to be cofibrant in order to
guarantee that this map is an isomorphism.)

While these properties if $\phigm$ are very convenient, they must be
used with caution.  The value $\phigm(X)$ is only guaranteed to have
the ``correct'' homotopy type on cofibrant objects.  The spectrum
underlying a commutative algebra is rarely known to be cofibrant,
making the monoidal geometric fixed point functor difficult to use in
that context.  The situation is a little better with associative
algebras.  The weak equivalence~\eqref{eq:118} leads to an expression
for the geometric fixed point spectrum of a quotient module which we
will use in~\S\ref{sec:proof-theor-refthm:21}.  In order to do so, we
need criteria guaranteeing that the monoidal geometric fixed point
functor realizes the correct homotopy type.  Such criteria are
described in \S\ref{sec:geom-fixed-points-1}.

\subsubsection{Geometric fixed points and the norm}

The geometric fixed point construction interacts well with the norm.
Suppose $H\subset G$ is a subgroup, and that $X$ is an $H$-spectrum.
The following result is proved as Proposition~\ref{thm:101}.  Our
original version merely concluded that the transformation in question
was a weak equivalence on cofibrant objects.  Andrew Blumberg and Mike
Mandell pointed out that it is in fact an isomorphism on cofibrant
objects, and at their request we have included a proof of the stronger
statement.

\begin{prop}
\label{thm:120}
There is a natural map
\[
\phihm X \to \phigm (\norm_{H}^{G} X)
\]
which is an isomorphism, hence a weak equivalence, on cofibrant objects.
\end{prop}

Because of Proposition~\ref{thm:220} and the fact that the norm
preserves cofibrant objects (Proposition~\ref{thm:195}), the above
result gives a natural zig-zag of weak equivalences relating
$\phih(X)$ and $\phig(\norm_{H}^{G}X)$ when $X$ is cofibrant.  In fact
there is a natural zig-zag of maps 
\[
\phih X \leftrightarrow \phig (\norm_{H}^{G} X)
\]
which is a weak equivalence not only for cofibrant $X$, but for
suspension spectra of cofibrant $G$-spaces and for the spectra
underlying cofibrant commutative rings.  The actual statement is
somewhat technical, and is one of the main results of Appendix~B.  The
condition is described in the statement of Proposition~\ref{thm:97}.
See also Remarks~\ref{rem:33} and~\ref{rem:10}.

\begin{cor}
\label{thm:30} For the spectra satisfying the condition of
Proposition~\ref{thm:97}, the composite functor
\[
\phig\circ\norm_{H}^{G}:\ugspectra{H}\to \spectra
\]
preserves, up to weak equivalence, wedges, directed colimits along
closed inclusions and cofiber sequences.
\end{cor}

\begin{pf}
The properties obviously hold for $\phih$.
\end{pf}

There is another useful result describing the interaction of the
geometric fixed point functor with the norm map in $RO(G)$-graded
cohomology described in \S\ref{sec:other-uses-norm}.  Suppose that $R$
is a $G$-equivariant commutative algebra, $X$ is a $G$-space, and
$V\in\RO(H)$ a virtual real representation of a subgroup $H\subset G$.
In this situation one can compose the norm
\[
\norm: R^{V}_{H}(X) \to R^{\ind V}_{G}(X)
\]
with the geometric fixed point map
\[
\phig:R^{\ind V}_{G}(X) \to (\phig R)^{V^{H}}(X^{G}),
\]
where $V^{H}\subset V$ is the subspace of $H$-fixed vectors, and
$X^{G}$ is the space of $G$-fixed points in $X$.

\begin{prop}
\label{thm:125}
The composite
\[
\phig\circ \norm: R^{V}_{H}(X) \to (\phig R)^{V^{H}}(X^{G})
\]
is a ring homomorphism.
\end{prop}

\begin{pf}
Multiplicativity is a consequence of the fact that both the norm and
the geometric fixed point functors are weakly monoidal.   Additivity
follows from the fact that the composition $\phig\circ\norm$ preserves
wedges (Corollary~\ref{thm:30}).
\end{pf}

\section{Mackey functors, homology and homotopy}
\label{sec:mackey-functors}

\subsection{Mackey functors}
\label{sec:mackey-functors-1}

In equivariant homotopy theory, the role of ``abelian group'' is
played by the notion of a {\em Mackey functor}
(Dress~\cite{MR0384917}).  The following is the summary of Dress'
definition as it appears in Greenlees-May~\cite{MR1361893}.

\begin{defin}[Dress~\cite{MR0384917}]
\label{def:17}
A {\em Mackey functor} consists of a pair $\Mm=(\Mm_{\ast},\Mm^{\ast})$ of
functors from the category of finite $G$-sets to the category of
abelian groups.  The two functors have the same object function
(denoted $\Mm$) and take disjoint unions to direct sums.  The functor
$\Mm_{\ast}$ is covariant, while $\Mm^{\ast}$ is contravariant, and
together they take a pullback diagram of finite $G$-sets
\[
\xymatrix{
S\ar[r]^{\delta}\ar[d]_{\gamma} & A \ar[d]^{\alpha} \\
T\ar[r]_{\beta} & B
}
\]
to a commutative square
\[
\xymatrix{
\Mm(P) \ar[r]^{\delta_{\ast}}  & \Mm(X)\\
\Mm(Y)\ar[u]^{\gamma^{\ast}}\ar[r]_{\beta_{\ast}} & \Mm(Z) \ar[u]_{\alpha^{\ast}}
}
\]
where $\alpha^{\ast}=\Mm^{\ast}(\alpha)$,
$\beta_{\ast}=\Mm_{\ast}(\beta)$, etc.
\end{defin}

The contravariant maps $\Mm^{\ast}(\alpha)$ are called the {\em
restriction} maps, and the covariant maps $\Mm_{\ast}(\beta)$, the {\em
transfer} maps.

A Mackey functor can also be defined as a contravariant additive
functor from the full subcategory of $\ugspectra{G}$ consisting of the
suspension spectra $\Sigma^{\infty}B_{+}$ of finite $G$-sets $B$.  It
is a theorem of tom Dieck that these definitions are equivalent.
See~\cite[\S5]{MR1361893}.

The equivariant homotopy groups of a $G$-spectrum $X$ are naturally
part of the Mackey functor $\pim_{n}X$ defined by
\begin{align*}
(\pim_{n}X)^{\ast}(B) &= [S^{n}\wedge B_{+},X]^{G} \\
(\pim_{n}X)_{\ast}(B) &= [S^{n},X\wedge B_{+}]^{G}.
\end{align*}
The identification of the two object functors
\[
[S^{n}\wedge B_{+},X]^{G} \approx [S^{n},X\wedge B_{+}]^{G}
\]
comes from the self-duality of finite $G$-sets (Example~\ref{eg:28}).
For $B=G/H$ one has
\[
\pim_{n}X(B) = \pi_{n}^{H}X.
\]

The Mackey functor $\pim_{0}S^{0}$ is the {\em Burnside ring} Mackey
functor $\Am$.  It is the free Mackey functor on one generator.  For a
$G$-set $B$, the value $\Am(B)$ is the group completion of the monoid
of isomorphism classes of finite $G$-sets $T\to B$ over $B$ under
disjoint union.  The restriction maps are given by pullback, and the
transfer maps by composition.  The group $\Am(G/H)$ works out to be
isomorphic to the abelian group underlying the Burnside ring of finite
$H$-sets.  

Just as every abelian group can occur as a stable homotopy group,
every Mackey functor $\Mm$ can occur as an equivariant stable homotopy
group.  In fact associated to each Mackey functor $\Mm$ is an
equivariant Eilenberg-Mac~Lane spectrum $H\Mm$, characterized by the
property
\[
\pim_{n}H\Mm=\begin{cases}
\Mm &\quad n=0 \\
0 &\quad n\ne 0.
\end{cases}
\]
See~\cite[\S5]{MR1361893} or~\cite{MR598689}. 

The homology and cohomology groups of a $G$-spectrum $X$ with
coefficients in $\Mm$ are defined by
\begin{align*}
H_{k}^{G}(X;\Mm) &= \pi^{G}_{k}H\Mm\wedge X \\
H_{G}^{k}(X;\Mm) &= [X,\Sigma^{k}H\Mm]^{G}.
\end{align*}
For a pointed $G$-space $Y$ one defines
\begin{align*}
H^{G}_{n}(Y;\Mm) &= H^{G}_{n}(\Sigma^{\infty}Y;\Mm) \\
H_{G}^{n}(Y;\Mm) &= H_{G}^{n}(\Sigma^{\infty} Y;\Mm).
\end{align*}
We emphasize that the equivariant (co)homology groups of pointed
$G$-spaces $Y$ we consider will always be {\em reduced} (co)homology
groups.

We will have need to consider the ordinary,
non-equivariant homology and cohomology groups of the spectrum
$i_{0}^{\ast}X$ underlying a $G$-spectrum $X$.  It will be convenient
to employ the shorthand notation
\begin{align*}
H_{n}^{u}(X;\Z) &= H_{n}(i_{0}^{\ast}X;\Z) \\
H^{n}_{u}(X;\Z) &= H^{n}(i_{0}^{\ast}X;\Z)
\end{align*}
for these groups.

\subsection{Constant and permutation Mackey functors}
\label{sec:const-perm-mack}

The {\em constant Mackey functor $\Zm$} is the functor represented
on the category of finite $G$-sets by the abelian group $\Z$ with
trivial $G$-action.  The value of $\Zm$ on a finite $G$-set $B$ is the
group of functions
\[
\Zm(B)=\hom^{G}(B,\Z)=\hom(B/G,\Z).
\]
The restriction maps are given by precomposition, and the transfer
maps by summing over the fibers.  For $K\subset H\subset G$, the
transfer map associated by $\Zm$ to
\[
G/K\to G/H
\]
is the map $\Z\to\Z$ given by multiplication by the index of $K$ in
$H$.   

\begin{defin}
\label{def:23} Suppose that $S$ is a $G$-set, and write $\Z\{S \}$ for
the free abelian group generated by $S$, and $\Z^{S}$ for the ring of
functions $S$ to $\Z$.  The {\em permutation Mackey functor}
$\Zm\{S\}$ is the Mackey functor
\[
\Zm\{S\}(B)=\hom^{G}(B,\Z\{S\}),
\]
whose restriction maps are given by precomposition and transfer maps by
summing over the fibers.   
\end{defin}

The permutation Mackey functor $\Z\{S \}$ is naturally isomorphic to
the Mackey functor $\pim_{0}H\Zm\wedge S_{+}$.   To see this note that 
restricting to underlying non-equivariant spectra gives a map
\[
\pim_{0}H\Zm\wedge S_{+}(B) = [B_{+},H\Zm\wedge S_{+}]^{G}
\to [B_{+},H\Z\wedge S_{+}],
\]
whose image lies in the $G$-invariant part.   Since
\[
[B_{+},H\Zm\wedge S_{+}]  = \hom(B,\Z\{S \})
\]
this gives a natural transformation 
\[
\pim_{0} H\Zm\wedge S_{+}\to \Zm\{S \}.
\]
Since both sides take filtered colimits in $S$ to filtered colimits,
to check that it is an isomorphism, it suffices to do so when $S$ is
finite.  In that case we can use the self duality of finite
$G$-sets to compute
\[
[B_{+},H\Zm\wedge S_{+}]^{G} \approx [B_{+}\wedge S_{+},H\Zm]^{G},
\]
and then observe that by definition of the constant Mackey functor
$\Zm$, the forgetful map
\[
[B_{+}\wedge S_{+},H\Zm]^{G} \to 
[B_{+}\wedge S_{+},H\Z]
\]
is an isomorphism with the $G$-invariant part of the target.  The
claim then follows from the compatibility of equivariant
Spanier-Whitehead duality with the restriction functor to
non-equivariant spectra.  

The properties of permutation Mackey functors listed in the Lemma
below follow immediately from the definition.  They are used in
\S\ref{sec:even-spectra} to establish some of our basic tools for
investigating the slice tower.  To formulate part~\thmItemref{part:3}, note that
every $G$-set $B$ receives a functorial map from a free $G$-set,
namely $G\times B$, and the group of equivariant automorphisms of
$G\times B$ over $B$ is canonically isomorphic to $G$.  For instance,
one can give $G\times B$ the product of the left action on $G$ and the
trivial action on $B$, and take the map $G\times B\to B$ to be the
original action mapping.  With this choice the automorphisms $G\times
B$ over $B$ are of the form $(g,b)\mapsto (g\,x, x^{-1}b)$ with $x\in
G$.

\begin{lem}
\label{thm:109}
\begin{thmList}
Let $\Mm$ be a permutation Mackey functor and $B$ finite $G$-set.
\item\label{part:4} 
If $B'\to B$ is a surjective map of finite
$G$-sets, then 
\[
\Mm(B)\to \Mm(B') \rightrightarrows \Mm(B'\underset{B}{\times}B')
\]
is an equalizer.
\item\label{part:3} Restriction along the action map
$G\times B\to B$ gives an isomorphism
\[
\Mm(B)\to \Mm(G\times B)^{G}.
\]
\item\label{part:1} The restriction mapping $\Mm(G/H)\to \Mm(G)$ gives an
isomorphism
\[
\Mm(G/H) \to \Mm(G)^{H}
\]
of $\Mm(G/H)$ with the $H$-invariant part of $\Mm(G)$.   
\item \label{part:2}  A map $\Mm\to \Mm'$ of
permutation Mackey functors is an isomorphism if and only if $\Mm(G)\to
\Mm'(G)$ is an isomorphism. 
\end{thmList}\qed
\end{lem}

\subsection{Equivariant cellular chains and cochains}
\label{sec:mack-funct-cohom}

The Mackey functor homology and cohomology groups of a $G$-CW spectrum
$X$ can be computed from a chain complex analogous to the complex of
cellular chains (see, for example~\cite[\S5]{MR1361893}).  Write
$X^{(n)}$ for the $n$-skeleton of $X$ so that
\[
X^{(n)}/X^{(n-1)}\approx {X_{n}}_{+}\wedge S^{n}
\]
with ${X_{n}}$ a discrete $G$-set.  Set
\begin{align*}
C_{n}^{\text{cell}}(X;\Mm) &= \pi^{G}_{n}H\Mm\wedge X^{(n)}/X^{(n-1)} = \pi^{G}_{0}H\Mm\wedge {X_{n}}_{+} \\
C^{n}_{\text{cell}}(X;\Mm) &= [X^{(n)}/X^{(n-1)},\Sigma^{n}H\Mm]^{G} =
[\Sigma^{\infty}{X_{n}}_{+},H\Mm]^{G}.
\end{align*}
The map
\[
X^{(n)}/X^{(n-1)}\to \Sigma X^{(n-1)}/X^{(n-2)}
\]
defines boundary and coboundary maps
\begin{align*}
C_{n}^{\text{cell}}(X;\Mm) &\to C_{n-1}^{\text{cell}}(X;\Mm) \\
C^{n-1}_{\text{cell}}(X;\Mm) &\to C^{n}_{\text{cell}}(X;\Mm).
\end{align*}
The equivariant homology and cohomology groups of $X$ with
coefficients in $\Mm$ are the homology and cohomology groups of these
complexes.  By writing the $G$-set ${X_{n}}$ as a coproduct of finite
$G$-sets $X_{n}^{\alpha}$ one can express $C_{n}^{\text{cell}}(X;\Mm)$
and $C^{n}_{\text{cell}}(X;\Mm)$ in terms of the values of the Mackey
functor $\Mm$ on the $X_{n}^{\alpha}$.

\begin{eg}
\label{eg:25}
Write $\rho_{G}$ for the (real) regular representation of $G$ and
$\rho_{G}-1$ for the reduced regular representation.   
The groups
\[
H^{\ast}(S^{\rho_{G}-1};\Mm)
\]
play an important role in equivariant stable homotopy theory.  To
describe them we need an equivariant cell decomposition of
$S^{\rho_{G}-1}$.  Since $S^{\rho_{G}-1}$ is the mapping cone of the
map
\[
S(\rho_{G}-1)\to \text{pt}
\]
from the unit sphere in $(\rho_{G}-1)$ it suffices to construct an
equivariant cell decomposition of $S(\rho_{G}-1)$.  Write $g=|G|$.
Think of $\R^{G}$ as the vector space with basis the elements of $G$.
The boundary of the standard simplex in this space is equivariantly
homeomorphic to $S(\rho_{G}-1)$.  The simplicial decomposition of this
simplex is not an equivariant cell decomposition, but the barycentric
subdivision is.  Thus $S(\rho_{G}-1)$ is homeomorphic to the geometric
realization of the poset of non-empty proper subsets of $G$.  This
leads to the complex
\begin{equation}
\label{eq:108}
\Mm(G/G) \to \Mm(S_{0}) \to \Mm(S_{1})\to\dots \to \Mm(S_{g-1})
\end{equation}
in which $S_{k}$ is the $G$-set of flags $F_{0}\subset\dots\subset
F_{k}\subset G$ of proper inclusions of subsets of $G$, with $G$
acting by translation.  The coboundary map is the alternating sum of
the restriction maps derived by omitting one of the sets in a flag.
\end{eg}

\begin{cor}
\label{thm:188}
For any Mackey functor $\Mm$, the group
\[
\pi_{\rho_{G}-1}^{G}H\Mm = H^{0}_{G}(S^{\rho_{G}-1};\Mm)
\]
is given by
\[
\bigcap_{H\subsetneq G} \ker \big(\Mm(G/G)\to \Mm(G/H)\big).
\]
\end{cor}

\begin{pf}
Using the complex~\eqref{eq:108} it suffices to show
that the orbit types occurring in $S_{0}$ are precisely the transitive
$G$-sets of the form $G/H$ with $H$ a proper subgroup of $G$.  The set
$S_{0}$ is the set of non-empty proper subsets $S\subset G$.  Any
proper subgroup $H$ of $G$ occurs as the stabilizer of itself,
regarded as a subset of $G$.  Since the subsets are proper, the
group $G$ does not occur as a stabilizer.
\end{pf}

\begin{eg}
\label{eg:17} 
Let $X$ be the sphere $S^{d-1}$ with the action of
$\ztn$ given by the antipodal map, and pointed by adding a disjoint
base point.  The usual cell decomposition into hemispheres is
equivariant for the action of $\ztn$, and for this cell structure one
has $X^{(j)}/X^{(j-1)}=(\ztn/\ztnm{n-1})_{+}\wedge S^{j}$.  The
complex of cellular chains $C_{\ast}^{\text{cell}}(X;\Mm)$ works out
to be the complex of length $d$
\[
\Mm(\ztn) \xrightarrow{\phantom{1+\gamma}}{}
\cdots
\xrightarrow{1+\gamma}{}
\Mm(\ztn)\xrightarrow{1-\gamma}{} \Mm(\ztn)
\]
in which $\gamma\in \ztn$ is the generator.   
\end{eg}

\begin{eg}
\label{eg:23} Let $G=\ztn$ and $\sigma$ the sign representation of
$G$.  Suspending the cell decomposition of Example~\ref{eg:17} gives
an equivariant cell decomposition of $S^{d\sigma}$ whose $k$-skeleton
is $S^{k\sigma}$ and whose set of $k$-cells is ${\zt}_{+}\times
D^{k}$, in which $G$ acts on $\zt$ through the unique surjective map
$G\to \zt$.  The complex of cellular chains
$C_{\ast}^{\text{cell}}(S^{d\sigma};\Mm)$ works out to be the complex of
length $(d+1)$
\[
\Mm(\zt) \xrightarrow{\phantom{1+\gamma}}{}
\cdots
\xrightarrow{1-\gamma}{}
\Mm(\zt)\xrightarrow{\phantom{1+\gamma}}{}\Mm(\text{pt})
\]
in which $\gamma\in G$ is the preferred generator.   
\end{eg}

If $\Mm$ is the constant Mackey functor $\Zm$, then
$C_{n}^{\text{cell}}(X;\Mm)$ is the permutation Mackey functor
$\Zm\{{X_{n}}\}$, and associates to a finite $G$-set $B$ the group of
equivariant functions
\[
B\to \Z\{X_{n} \}=C_{n}^{\text{cell}}X.
\]
In this way the entire Mackey functor chain complex
$C_{\ast}^{\text{cell}}(X;\Zm)$ is encoded in the ordinary cellular
chain complex $C_{\ast}^{\text{cell}}(X)$ for $i_{0}^{\ast}X$,
equipped with the action of $G$.  The equivariant homology group
$H_{\ast}^{G}(X;\Zm)$ are just the homology groups of the complex
\[
\hom_{G}(G/G, C_{\ast}^{\text{cell}}(X)) = C_{\ast}^{\text{cell}}(X)^{G}
\]
of $G$-invariant cellular chains.   Similarly the equivariant
cohomology groups $H^{\ast}_{G}(X;\Z)$ are given by the cohomology
groups of the complex
\[
C^{\ast}_{\text{cell}}(X)^{G}
\]
of equivariant cochains.  The equivariant homology and cohomology
groups depend only on the equivariant chain homotopy type of these
complexes of permutation Mackey functors.

\begin{eg}
\label{eg:26} If $X$ is a $G$-{\em space} admitting the structure of a
$G$-CW complex, then the cohomology groups $H^{\ast}_{G}(X;\Zm)$ are
isomorphic to the cohomology groups
\[
H^{\ast}(X/G;\Z)
\]
of the orbit space.   Indeed the equivariant cell decomposition of $X$
induces a cell decomposition of $X/G$ and one has an isomorphism
\[
C^{\ast}_{\text{cell}}(X)^{G} \approx 
C^{\ast}_{\text{cell}}(X/G).
\]
\end{eg}

\begin{eg}
\label{eg:5}
Suppose that $V$ is a representation of $G$ of dimension $d$, and
consider the equivariant cellular chain complex
\[
C_{d}^{\text{cell}}(S^{V};\Zm) \to C_{d-1}^{\text{cell}}(S^{V};\Zm)\to \dots
\to C_{0}^{\text{cell}}(S^{V};\Zm),
\]
associated to an equivariant cell decomposition of $S^{V}$.  The
underlying homology groups are those of the sphere $S^{V}$.  In
particular, the kernel of
\[
C_{d}^{\text{cell}}(S^{V};\Zm) \to C_{d-1}^{\text{cell}}(S^{V};\Zm)
\]
is isomorphic, as a $G$-module, to $H^{u}_{d}(S^{V};\Z)$.  If $V$ is
orientable then the $G$-action is trivial, and one finds that the
restriction map
\[
H_{d}^{G}(S^{V};\Zm)\to H_{d}^{u}(S^{V};\Z)
\]
is an isomorphism.  A choice of orientation gives an equivariant
isomorphism
\[
H^{u}_{d}(S^{V};\Z)\approx \Z.
\]
Thus when $V$ is oriented there is a unique isomorphism
\[
H_{d}^{G}(S^{V};\Zm)\approx \Z
\]
extending the non-equivariant isomorphism given by the orientation.
\end{eg}

\subsection{Homology and geometric fixed points}
\label{sec:homol-geom-fixed}

In addition to the Mackey functor
homotopy groups $\pim_{\ast}X$ there are the $RO(G)$ graded homotopy
groups $\pi^{G}_{\star}X$ defined by
\[
\pi^{G}_{V}X = [S^{V},X]^{G}\qquad V\in RO(G).
\]
Here $RO(G)$ is the Grothendieck group of real representations of $G$.
The use of $\star$ for the wildcard symbol in $\pi^{G}_{\star}$ is
taken from Hu-Kriz~\cite{MR1808224}.  The $RO(G)$-graded homotopy
groups are also part of a Mackey functor $\pim_{\star}(X)$ defined by
\[
\pim_{V}X(B) = [S^{V}\wedge B_{+},X]^{G}.
\]
As with $\Z$-graded homotopy groups, we'll use the abbreviation
\[
\pi^{H}_{V}X = (\pim_{V}X)(G/H).
\]
In this section we will make use of $RO(G)$-graded homotopy groups to
describe the geometric fixed point spectrum $\phig H\Zm$ when $G=\ztn$
(Proposition~\ref{thm:94} below).

There are a few distinguished elements of $RO(G)$-graded homotopy
groups we will need.  Let $V$ be a representation of $G$ and $S^{0}\to
S^{V}$ the one point compactification of the inclusion $\{0 \}\subset
V$.
\begin{defin}
\label{def:43}
The element 
\[
a_{V}\in \pi^{G}_{-V}S^{0}
\]
is the element corresponding under the suspension isomorphism
$\pi^{G}_{-V}S^{0}\approx \pi^{G}_{0}S^{V}$ to the map
$S^{0}\hookrightarrow S^{V}$ described above.
\end{defin}

The element $a_{V}$ is the Euler class of $V$ in $RO(G)$-graded
equivariant stable cohomotopy.  If $V$ contains a trivial
representation then $a_{V}=0$.  For two representations $V$ and $W$
one has
\[
a_{V\oplus W}=a_{V}a_{W}\in \pi^{G}_{-V-W}S^{0}.
\]

When $V$ is oriented, Example~\ref{eg:5} provides a preferred
generator of $H^{G}_{d}(S^{V};\Z)$.  We give the corresponding
$RO(G)$-graded homotopy class name.

\begin{defin}
\label{def:54} Let $V$ be an oriented representation of $G$ of
dimension $d$.  The element
\[
u_{V}\in\pi^{G}_{d-V}H\Zm.
\]
is the element corresponding to the preferred generator of $\pi_{d}
H\Zm\wedge S^{V}= H^{G}_{d}(S^{V};\Zm)$ given by Example~\ref{eg:5}.   
\end{defin}

If $V$ is trivial then $u_{V}=1$.  If $V$ and $W$ are two oriented
representations of $G$, and $V\oplus W$ is given the direct sum
orientation, then
\[
u_{V\oplus W}=u_{V}u_{W}.
\]
Among other things this implies that the class $u_{V}$ is {\em stable}
in $V$ in the sense that $u_{V+1}=u_{V}$.

For any $V$, the representation $V\oplus V$ has a canonical
orientation, giving 
\[
u_{V\oplus V}\in\pi^{G}_{2d-2V} H\Zm.
\]
When $V$ is oriented this class can be identified, up to sign, with
$u_{V}^{2}$.

The elements $a_{V}$ and $u_{V}$ behave well with respect to the norm.
The following result is a simple consequence of the
fact~\eqref{eq:103} that $\norm S^{V}=S^{\ind V}$.
\begin{lem}
\label{thm:138} Suppose that $V$ is a $d$-dimensional representation
of a subgroup $H\subset G$.  Then
\begin{align*}
\norm a_{V} &= a_{\ind V} \\
u_{\ind d}\cdot \norm u_{V} &= u_{\ind V},
\end{align*}
where $\ind V=\ind_{H}^{G}V$ is the induced representation and $d$ is
the trivial representation.   \qed
\end{lem}

\begin{rem}
\label{rem:61} As is standard in algebra, we will adopt the convention
that the operation of mutiplication by an element of a ring on a
module is denoted by the element of the ring.  We will also use it in
closely related contexts.  For example, for a $G$-spectrum $X$ we will
refer to the to the maps
\begin{align*}
a_{V}\wedge 1_{X} &: S^{-V}\wedge X\to X \\
u_{V}\wedge 1_{X} &: S^{d-V}\wedge X\to H\Zm\wedge X
\end{align*}
as {\em multiplication} by $a_{V}$ and $u_{V}$ respectively, and, when
no confusion is likely, denote them simply by $a_{V}$ and $u_{V}$.
Note that $X$ might be a virtual representation sphere.  This means
that we will not usually distinguish in notation between these maps
and their suspensions.  Similarly, if $R$ is any equivariant algebra,
and $x\in \pi^{G}_{V} S^{0}$ then the product of $x$ with
$1\in\pi^{G}_{0}R$ will be denoted $x\in\pi^{G}_{V}R$.  In accordance
with this, at various places in this paper the symbol $a_{V}$ might
refer to a map $S^{-V}\to S^{0}$, or its suspension $S^{0}\to S^{V}$
or the Hurewicz image $S^{0}\to H\Zm\wedge S^{V}$ or equivalently an
element of $\pi^{G}_{0}H\Zm\wedge S^{V}$.
\end{rem}

\begin{eg}
\label{eg:15} Let $S^{\infty V}$ be the colimit of the spaces $S^{n
V}$ under the standard inclusions.  Each of these inclusions is
``multiplication by $a_{V}$.''  Smashing with a $G$-spectrum $X$ we
find that $S^{\infty V}\wedge X$ is the colimit of the sequence
\[
X\xrightarrow{a_{V}}{} S^{V}\wedge X\xrightarrow{a_{V}} S^{V\oplus V}\wedge X
\cdots \xrightarrow{a_{V}}{} S^{ n V}\wedge X \xrightarrow{a_{V}}{}\cdots.
\]
Using the suspension isomorphism to replace
$\pi^{G}_{\star}S^{ n V}\wedge X$ with
$\pi^{G}_{\star-\, n V}X$ the sequence of the $RO(G)$-graded
groups becomes
\[
\pi^{G}_{\star}X\xrightarrow{a_{V}}{} 
\pi^{G}_{\star-V}X\cdots\xrightarrow{a_{V}}{}
\pi^{G}_{\star-\, n V}X\cdots
\]
from which one gets an isomorphism
\[
\pi^{G}_{\star}S^{\infty V}\wedge X \approx
a_{V}^{-1}\pi^{G}_{\star} X. 
\]
Under this isomorphism the effect in $RO(G)$-graded homotopy groups
induced by the inclusion
\[
S^{ n V}\wedge X \to S^{\infty V}\wedge X
\]
sends $x\in \pi^{G}_{\star} X \approx \pi^{G}_{\star+V}S^{ n
V}\wedge X$ to $a_{V}^{-n}x\in a_{V}^{-1}\pi^{G}_{\star}X$.
\end{eg}

\begin{eg}
\label{eg:20} Specializing Example~\ref{eg:23}, let $G=\ztn$ and
$\sigma$ the sign representation.  Consider the equivariant homology
of $S^{d\sigma}$ with coefficients in the constant Mackey functor
$\Z$.  The complex of cellular chains works out to be
(Example~\ref{eg:23}) the complex of length $(d+1)$
\[
\Z\xrightarrow{\phantom{2}}{}\cdots
\Z\xrightarrow{2}{}\Z\xrightarrow{0}{} \Z\xrightarrow{2}{}\Z,
\]
Our conventions provide nomenclature for the homology classes.  When
$d$ is odd the group $H_{d}(S^{d\sigma};\Zm)$ is zero.  When $d$ is
even, the representation $d\sigma$ acquires a canonical orientation,
the group $H_{d}^{\zt}(S^{d\sigma};\Zm)$ is canonically isomorphic to
the integers, and the preferred generator is the class $u_{d\sigma}$
(Remark~\ref{rem:61}).  For every even $0\le k < d$ the group
$H^{G}_{k}(S^{d\sigma};\Zm)$ is cyclic of order $2$ generated by the
image of $u_{k\sigma}\in H^{G}_{k}(S^{k\sigma};\Zm)$ under the map
induced by the inclusion $S^{k\sigma}\subset S^{d\sigma}$.  As
explained in Remark~\ref{rem:61} this induced map is multiplication by
$a_{(d-k)\sigma}$, and so this generator corresponds to the element
\[
a_{(d-k)\sigma}\cdot u_{k\sigma}\in \pi^{G}_{k - d\sigma}(H\Zm)
\]
under the suspension isomorphism 
\[
\pi^{G}_{k - d\sigma}(H\Zm)  \approx
\pi^{G}_{k}(H\Zm\wedge S^{d\sigma}) = H^{G}_{k}(S^{d\sigma};\Zm).
\]
\end{eg}

\begin{eg}
\label{eg:24}
Passing to the limit as $d\to\infty$ and using the last part of
Example~\ref{eg:15} we find that $a_{(d-k)\sigma}\cdot u_{k\sigma}$ is
sent to 
\[
a_{d\sigma}^{-1}\cdot a_{(d-k)\sigma}\cdot u_{k\sigma} =
a_{k\sigma}^{-1}u_{k\sigma}\in \pi_{k}S^{\infty \sigma}.
\]
Writing $b=a_{2\sigma}^{-1}u_{2\sigma}$ we find that the homogeneous
component 
\[
\pi^{\zt}_{2n}H\Zm \wedge S^{\infty \sigma}\subset \pi^{\zt}_{\star}H\Zm
\wedge S^{\infty \sigma}= a_{2\sigma}^{-1}\pi^{\zt}H\Zm
\]
is cyclic of order $2$, generated by
$b^{n}$.
\end{eg}

We now explicitly describe the geometric fixed point spectrum of
$H\Zm$ when $G=\ztn$.  The computation plays an important role in the
proof of the Reduction Theorem.

\begin{prop}
\label{thm:94}
Let $G=\ztn$.  For any $G$-spectrum $X$, the $RO(G)$-graded
homotopy groups of $\tgp\wedge X$ are given by
\[
\pi^{G}_{\star}(\tgp\wedge  X) = a_{\sigma}^{-1} \pi^{G}_{\star} (X).
\]
The homotopy groups of the commutative algebra $\phig H\Zm$ are given
by
\[
\pi_{\ast}(\phig H\Zm)=\Z/2[b],
\]
where $b=u_{2\sigma}a_{\sigma}^{-2}\in\pi_{2}(\phig H\Zm)
=\pi^{G}_{2}(\tgp \wedge H\Zm)\subset a_{\sigma}^{-1} \pi^{G}_{\star}
H\Zm$.
\end{prop}

\begin{pf}
As mentioned in Remark~\ref{rem:42}, the space $\tgp$ can be identified
with
\[
\lim_{n\to\infty} S^{n\sigma}.
\]
The first assertion therefore follows from Example~\ref{eg:15}.  The
second assertion follows from Example~\ref{eg:24} and the fact that
the map $a_{\sigma}^{-1}\pi^{G}_{\star}X\to \pi^{G}_{\star}\tgp\wedge
X$ is a ring homomorphism when $X$ is an equivariant algebra.
\end{pf}

\subsection{A gap in homology}
\label{sec:gap-homology}

We conclude \S\ref{sec:mackey-functors} with some further observations
about $S^{\rho_{G}-1}$.  Proposition~\ref{thm:108} below constitutes
the computational part of the Gap Theorem, and contains the Cell Lemma
as a special case.

\begin{eg}
\label{eg:8} Suppose that $G$ is not the trivial group.  In
\S\ref{sec:even-spectra} we will encounter the group
\[
\pi^{G}_{\rho_{G}-2}H\Zm \approx H^{1}_{G}(S^{\rho_{G}-1};\Zm),
\]
which, by Example~\ref{eg:26}, is isomorphic to 
\[
H^{1}(S^{\rho_{G}-1}/G;\Z).
\]
The $G$-space $S^{\rho_{G}-1}$ is the unreduced suspension of the unit
sphere $S(\rho_{G}-1)$, and so the orbit space is also a suspension.
If $|G|>2$ then $S(\rho_{G}-1)$ is connected, hence so is the orbit
space. If $G=\zt$, then $S(\rho_{G}-1)\approx G$ and the orbit space
is still connected.  In all cases then, the unreduced suspension
$S^{\rho_{G}-1}/G$ is simply connected.  Thus
\[
\pi^{G}_{\rho_{G}-2}H\Zm \approx H^{2}_{G}(S^{\rho_{G}};\Zm) \approx
H^{1}_{G}(S^{\rho_{G}-1};\Zm) = 0.
\]
In fact, the same argument shows that for $n>0$  the orbit space
$S^{n(\rho_{G}-1)}/G$ is simply connected, and hence
\[
H^{0}_{G}(S^{n(\rho_{G}-1)};\Zm) = 
H^{1}_{G}(S^{n(\rho_{G}-1)};\Zm) = 0
\]
or, equivalently
\[
\pi^{G}_{n(\rho_{G}-1)}H\Zm   = \pi^{G}_{n(\rho_{G}-1)-1}H\Zm = 0.
\]
\end{eg}

Building on this, we have

\begin{prop}
\label{thm:108}
Let $G$ be any non-trivial finite group and $n\ge 0$ an integer.
Except in case $G=C_{3}$, $i=3$, $n=1$ the groups
\[
H^{i}_{G}(S^{n\rho_{G}};\Zm)
\]
are zero for $0<i<4$.  In the exceptional case one has
\[
H^{3}_{G}(S^{\rho_{C_{3}}};\Zm) = \Z. 
\]
\end{prop}

\begin{pf}
Since 
\[
H^{i}_{G}(S^{n\rho_{G}};\Zm) \approx 
H^{i-n}_{G}(S^{n(\rho_{G}-1)};\Zm),
\]
connectivity and Example~\ref{eg:8} show that
$H^{i}_{G}(S^{n\rho_{G}};\Zm) = 0$ for $i\le n+1$.   This takes care
of the cases in which $n+1\ge 3$, leaving only $n=1$, and in that case
only the group
\[
H^{2}_{G}(S^{\rho_{G}-1};\Zm)
\]
which is isomorphic to 
\[
H^{2}(S^{\rho_{G}-1}/G;\Z).
\]
Since the orbit space $S^{\rho_{G}-1}/G$ is simply connected, the
universal coefficient theorem gives an inclusion
\[
H^{2}(S^{\rho_{G}-1}/G;\Z)\to
H^{2}(S^{\rho_{G}-1}/G;\Q).
\]
It therefore suffices to show that
\[
H^{2}(S^{\rho_{G}-1}/G;\Q)=0.
\]
But since $G$ is finite, this group is just the $G$-invariant part of 
\[
H^{2}(S^{\rho_{G}-1};\Q)
\]
which is zero since $G$ does not have order $3$.  When $G$ does have
order $3$ the group is $\Q$.  The claim follows since the homology
groups are finitely generated.
\end{pf}

\section{The slice filtration}
\label{sec:ths-slice-filtration}

The slice filtration is an equivariant analogue of the Postnikov
tower, to which it reduces in the case of the trivial group.  In this
section we introduce the slice filtration and establish some of its
basic properties.  We work for the most part with a general finite
group $G$, though our application to the Kervaire invariant problem
involves only the case $G=\ztn$.  While the situation for general $G$
exhibits many remarkable properties, the reader should regard as
exploratory the apparatus of definitions at this level of generality.

From now until the end of \S\ref{sec-detect} our focus will be on
homotopy theory.  Though it will not appear in the notation, all
spectra should be replaced by cofibrant or fibrant approximations
where appropriate.  

\subsection{Slice cells}
\label{sec:slice-cells}

\subsubsection{Slice cells and their dimension}
\label{sec:slice-cells-and-dimensions}

For a subgroup $K\subset G$ let $\rho_{K}$ denote its regular
representation, and write
\[
\slicecell(m,K)=G_{+}\smashove{K}S^{m\rho_{K}}\qquad m\in\Z.
\]

\begin{defin}
\label{def:20}
The set of  {\em slice cells} is
\[
\{\slicecell(m,K), \Sigma^{-1}\slicecell(m,K)\mid m\in \Z,
K\subset G \}.
\]
\end{defin}

This brings two notions of ``cell'' into the story: the slice cells
and the cells of the form $G/H_{+}\wedge D^{m}$, used to manufacture
$G$-CW spectra.  We'll refer to the traditional equivariant cells as
``$G$-cells'' in order to distinguish them from the slice cells.

\begin{defin}
\label{def:44}
A slice cell is {\em regular} if it is of the form $\slicecell(m,K)$. 
\end{defin}

\begin{defin}
\label{def:21}
A slice cell is {\em induced} if it is of the form
\[
G_{+}\smashove{H}\slicecell,
\]
where $\slicecell$ is a slice cell for $H$ and $H\subset G$ is a
proper subgroup.  It is {\em free} if $H$ is the trivial group.  A
slice cell is {\em isotropic} if it is not free.
\end{defin}

Since
\begin{align*}
[G_{+}\smashove{H}S,X]^{G} &\approx [S,i_{H}^{\ast}X]^{H}\quad\text{ and }  \\
[X,G_{+}\smashove{H}S]^{G} &\approx [i_{H}^{\ast}X,S]^{H},
\end{align*}
induction on $|G|$ usually reduces claims about cells to the case of
those which are not induced.  The slice cells which are not induced
are those of the form $S^{m\rho_{G}}$ and $S^{m\rho_{G}-1}$.

\begin{defin}
The {\em dimension} of a slice cell is defined by
\begin{align*}
\dim \slicecell(m,K) &=m|K| \\
\dim\Sigma^{-1}\slicecell(m,K) &= m |K|-1.
\end{align*}
\end{defin}
In other words the dimension of a slice cell is that of its underlying
spheres.

\begin{rem}
Not every suspension of a slice cell is a slice cell.  Typically, the
spectrum $\Sigma^{-2}\slicecell(m,K)$ will {\em not} be a slice
cell, and will {\em not} exhibit the properties of a slice cell of dimension
$\dim \slicecell(m,K)-2$.
\end{rem}

The following is immediate from the definition.

\begin{prop}\label{thm:32}
Let $H\subset G$ be a subgroup.  If $\slicecell$ is a $G$-slice cell of
dimension $d$, then $i_{H}^{\ast}\slicecell$ is a wedge of $H$-slice cells of
dimension $d$.  If $\slicecell$ is an $H$-slice cell of dimension $d$ then
$G_{+}\smashove{H}\slicecell$ is a $G$-slice cell of dimension $d$.  \qed
\end{prop}

The regular slice cells behave well under the norm.

\begin{prop}
\label{thm:77}
Let $H\subset G$ be a subgroup.  If $\refine$ is a wedge of regular
$H$-slice cells, then $\norm_{H}^{G}\refine$ is a wedge of regular
$G$-slice cells.
\end{prop}

\begin{pf}
The wedges of regular $H$-slice cells are exactly the indexed wedges
(in the sense of \S\ref{sec:norm-induction}) of spectra of the form
$S^{m\rho_{K}}$ for $K\subset H$, and $m\in\Z$.  Since regular
representations induce to regular representations, the
identity~\eqref{eq:103} and the distribution formula
(Proposition~\ref{thm:9}) show that the norm of such an indexed wedge
is an indexed wedge of $S^{m\rho_{K}}$ with $K\subset G$ and
$m\in\Z$.  The claim follows.
\end{pf}

\subsubsection{Slice positive and slice null spectra}
\label{sec:slice-positive-slice}

Underlying the theory of the Postnikov tower is the notion of
``connectivity'' and the class of $(n-1)$-connected spectra.  In this
section we describe the slice analogues of these ideas.  There is a
simple relationship between ``connectivity'' and
``slice-positivity'' which we will describe in detail in
\S\ref{sec:conv-slice-tower}.

\begin{defin}
\label{def:41} A $G$-spectrum $Y$ is {\em slice $n$-null}, written
\[
Y < n\quad\text{or}\quad Y \le n-1
\]
if for every slice cell $\slicecell$ with $\dim \slicecell\ge n$ the
$G$-space
\[
\gspectra{G}(\slicecell,Y)
\]
is equivariantly contractible.  A $G$-spectrum $X$ is {\em slice $n$-positive}, written
\[
X> n\quad\text{or}\quad X \ge n+1
\]
if 
\[
\gspectra{G}(X,Y)
\]
is equivariantly contractible for every $Y$ with $Y\le n$.
\end{defin}

We will use the terms {\em slice-positive} and {\em slice-null}
instead of ``slice $0$-positive'' and ``slice $0$-null.''  The full
subcategory of $\ugspectra{G}$ consisting of $X$ with $X> n$ will be
denoted $\gslice{n}$ or $\geslice{n+1}$.  Similarly, the full
subcategory of $\ugspectra{G}$ consisting of $X$ with $X < n$ will be
denoted $\lslice{G}{n}$ or $\leslice{G}{n-1}$.

\begin{rem}
\label{rem:23}
The category $\gslice{n}$ is the smallest full subcategory of
$\ugspectra{G}$ containing the slice cells $\slicecell$ with
$\dim\slicecell> n$ and possessing the following properties:
\begin{itemize}
\item [i)] If $X$ is weakly equivalent to an object of
$\gslice{n}$, then $X$ is in $\gslice{n}$.
\item [ii)] Arbitrary wedges of objects of $\gslice{n}$ are in
$\gslice{n}$.
\item [iii)]  If $X\to Y\to Z$ is a cofibration sequence and $X$ and
$Y$ are in $\gslice{n}$ then so is $Z$.
\item [iv)] If $X\to Y\to Z$ is a cofibration sequence and $X$ and
$Z$ are in $\gslice{n}$ then so is $Y$.
\end{itemize}
More briefly, these properties are that $\gslice{n}$ is closed under
weak equivalences, homotopy colimits (properties \thmListItem{2} and
\thmListItem{3}), and extensions.    
\end{rem}

\begin{rem}
\label{rem:76} 
The fiber of a map of slice $n$-positive spectra is not
assumed to be slice $n$-positive, and need not be.  For example, the
fiber of $\ast \to S^{\rho_{G}}$ is $S^{\rho_{G}-1}$ which is not
slice $(|G|-1)$-positive, even though both $\ast$ and $S^{\rho_{G}}$
are.
\end{rem}

For $n=0,-1$, the notions of slice $n$-null and slice $n$-positive are
familiar.

\begin{prop}
\label{thm:65}
For a $G$-spectrum $X$ the following hold
\begin{thmList}
\item $X\ge 0 \iff X$ is $(-1)$-connected, i.e. $\pim_{k} X= 0$ for $k< 0$;
\item $X\ge -1 \iff X$ is $(-2)$-connected, i.e. $\pim_{k} X= 0$ for $k<-1$;
\item $X < 0 \iff X$ is $0$-coconnected, i.e.  $\pim_{k} X= 0$ for $k\ge 0$;
\item $X < -1 \iff X$ is $(-1)$-coconnected, i.e.  $\pim_{k} X= 0$ for $k\ge -1$;
\end{thmList}
\end{prop}

\begin{pf}
These are all straightforward consequences of the fact that $S^{0}$ is
a slice cell of dimension $0$, and $S^{-1}$ is a slice cell of
dimension $(-1)$.
\end{pf}

\begin{rem}
It is not the case that if $Y>0$ then $\pi_{0}Y=0$.  In
Proposition~\ref{thm:35} we will see that the fiber $F$ of $S^{0}\to
H\Zm$ has the property that $F>0$.  On the other hand $\pim_{0}F$ is the
augmentation ideal of the Burnside ring.   Proposition~\ref{thm:64} below
gives a characterization of slice-positive spectra.
\end{rem}

The classes of slice $n$-null and slice $n$-positive spectra are
preserved under change of group.

\begin{prop}
\label{thm:67} 
Suppose $H\subset G$, that $X$ is a $G$-spectrum and $Y$ is an $H$-spectrum.
The following implications hold
\begin{align*}
X>n &\implies i_{H}^{\ast}X>n \\
X<n &\implies i_{H}^{\ast}X<n \\
Y>n &\implies G_{+}\smashove{H}Y>n \\
Y<n &\implies G_{+}\smashove{H}Y<n.
\end{align*}
\end{prop}

\begin{pf}
The second and third implications are straightforward consequences of
Proposition~\ref{thm:32}.   The fourth implication follows from the
Wirthm\"uller isomorphism and Proposition~\ref{thm:32}, and the first
implication is an immediate consequence of the fourth.
\end{pf}

We end this section with a mild simplification of the condition that a
spectrum be slice $n$-null.

\begin{lem}
\label{thm:2}
For a $G$-spectrum $X$, the following are equivalent
\begin{thmList}
\item $X < n$;
\item For all slice cells $\slicecell$ with $\dim \slicecell \ge n$,
$[\slicecell,X]^{G}=0$.
\end{thmList}
\end{lem}

\begin{pf}
The first condition trivially implies the second.  We prove that the
second implies the first by induction on $|G|$.  By the induction
hypothesis we may assume that the $G$-space
\[
\gspectra{G}(\slicecell,X)
\]
is contractible for all induced slice cells $\slicecell$ with $\dim
\slicecell\ge n$, and that for all slice cells $\slicecell$ with
$\dim\slicecell\ge n$, and all proper $H\subset G$, the space
\[
\gspectra{G}(\slicecell,X)^{H}
\]
is contractible.   We therefore also know that the $G$-space 
\[
\gspectra{G}(T\wedge\slicecell, X)
\]
is contractible, for all slice cells $\slicecell$ with
$\dim\slicecell\ge n$ and all $(-1)$-connected $G$-CW spectra $T$
built entirely from induced $G$-cells.   We must show that for each
$t \ge 0$, the groups 
\begin{gather*}
[S^{t}\wedge S^{m\rho_{G}-1},X]^{G} \qquad m|G|-1\ge n \\
[S^{t}\wedge S^{m\rho_{G}},X]^{G} \qquad m|G|\ge n \\
\end{gather*}
are zero.  They are zero by assumption when $t=0$.  For $t>0$, the
first case is a special case of the second, since $S^{1}\wedge
S^{m\rho_{G}-1}$ is a slice cell of dimension $m|G|$.  Let $T$ be the
homotopy fiber of the map
\[
S^{t}\subset S^{t\rho_{G}}, 
\]
and consider the exact sequence
\[
[S^{t\rho_{G}}\wedge S^{m\rho_{G}},X]^{G}
\to
[S^{t}\wedge S^{m\rho_{G}},X]^{G}
\to
[T\wedge S^{m\rho_{G}},X]^{G}.
\]
The leftmost group is zero since $S^{t\rho_{G}}\wedge S^{m\rho_{G}}$
is a slice cell of dimension $(t+m)|G|\ge n$.  The rightmost group is
zero by the induction hypothesis since $T$ is $(-1)$-connected and
built entirely from induced $G$-cells.  It follows from exactness that
the middle group is zero.
\end{pf}

\subsection{The slice tower}\label{sec:slice-tower}

Let $P^{n} X$ be the Bousfield localization, or Dror Farjoun
nullification (\cite{MR1320991,hirschhorn03:_model}) of $X$ with
respect to the class $\gslice{n}$ , and $P_{n+1}X$ the homotopy fiber
of $X\to P^{n}X$.  Thus, by definition, there is a functorial
fibration sequence
\[
P_{n+1}X\to X \to P^{n}X.
\]

The functor $P^{n}X$ can be constructed as the colimit of a sequence
of functors
\[
W_{0}X\to W_{1}X\to\cdots.
\]
The $W_{i}X$ are defined inductively starting with $W_{0}X=X$, and
taking $W_{k}X$ to be the cofiber of
\[
\bigvee_{I} \Sigma^{t}\slicecell \to W_{k-1}X,
\]
in which the indexing set $I$ is the set of maps
$\Sigma^{t}\slicecell\to W_{k-1}X$ with $\slicecell>n$ a slice cell
and $t\ge 0$.  By Lemma~\ref{thm:2} the functors $P^{n}$ can also be
formed using the analogous construction using only slice cells
themselves, and not their suspensions.  

\begin{prop}
\label{thm:35} A spectrum $X$ is slice $n$-positive if an only if it
admits (up to weak equivalence) a filtration
\[
X_{0} \subset X_{1}\subset \cdots
\]
whose associated graded spectrum $\bigvee X_{k}/X_{k-1}$ is a wedge of
slice cells of dimension greater than $n$.  For any spectrum $X$,
$P_{n+1}X$ is slice $n$-positive.
\end{prop}

\begin{pf}
This follows easily from the construction of $P^{n}X$ described above.
\end{pf}

The map $P_{n+1}X\to X$ is characterized up to a contractible space of
choices by the properties
\begin{itemize}
\item [i)] for all $X$, $P_{n+1}X\in\gslice{n}$;
\item [ii)] for all $A\in \gslice{n}$ and all $X$, the map
$\gspectra{G}(A,P_{n+1}X)\to \gspectra{G}(A,X)$ is a weak equivalence
of $G$-spaces.
\end{itemize}
In other words, $P_{n+1}X\to X$ is the ``universal map'' from an
object of $\gslice{n}$ to $X$.  Similarly $X\to P^{n}X$ is the
universal map from $X$ to a slice $(n+1)$-null $G$-spectrum
$Z$.  More specifically
\begin{itemize}
\item [iii)] the spectrum $P^{n}X$ is slice $(n+1)$-null;
\item [iv)] for any slice $(n+1)$-null $Z$, the map
\[
\gspectra{G}(P^{n}X,Z)\to
\gspectra{G}(X,Z)
\]
is a weak equivalence.
\end{itemize}
These conditions lead to a useful recognition principle.
\begin{lem}
\label{thm:7}
Suppose $X$ is a $G$-spectrum and that 
\[
\tilde{P}_{n+1}\to X\to \tilde P^{n}
\]
is a fibration sequence with the property that $\tilde P^{n}\le n$ and
$\tilde P_{n+1}>n$.  Then the canonical maps $\tilde{P}_{n+1}\to
P_{n+1}X$ and $P^{n}X\to \tilde P^{n}$ are weak equivalences.
\end{lem}

\begin{pf}
We show that the map $X\to\tilde P^{n}$ satisfies the universal
property of $P^{n}X$.  Suppose that $Z\le n$, and consider the
fibration sequence of $G$-spaces
\[
\gspectra{G}(\tilde P^{n},Z)\to\gspectra{G}(X,Z)\to\gspectra{G}(\tilde
P_{n+1},Z)
\]
The rightmost space is contractible since $\tilde P_{n+1}>n$, so the
map $\gspectra{G}(\tilde P^{n},Z)\to\gspectra{G}(X,Z)$ is a weak equivalence.
\end{pf}

The following consequence of Lemma~\ref{thm:7} is used in the
proof of the Reduction Theorem.

\begin{cor}
\label{thm:197}
Suppose that $X\to Y\to Z$ is a cofibration sequence, and that the mapping
cone of $P^{n}X\to P^{n}Y$ is slice $(n+1)$-null.   Then both
\[
P^{n}X\to P^{n}Y\to P^{n}Z
\]
and
\[
P_{n+1}X\to P_{n+1}Y\to P_{n+1}Z
\]
are cofibration sequences.
\end{cor}

\begin{pf}
Consider the diagram 
\[
\xymatrix@R=1.5em{
P_{n+1}X  \ar[r]\ar[d]  & P_{n+1}Y  \ar[r]\ar[d]   & \tilde P_{n+1}Z \ar[d] \\
X  \ar[r]\ar[d]  & Y  \ar[r]\ar[d]   & Z  \ar[d]\\
P^{n}X  \ar[r]    &    P^{n}Y   \ar[r]  & \tilde P^{n}Z
}
\]
in which the rows and columns are cofibration sequences.  By
construction, $\tilde P_{n+1}Z$ is slice $n$-positive
(Remark~\ref{rem:23}).  If $\tilde P^{n}Z \le n$ then the right column
satisfies the condition of~\ref{thm:7}, and the result follows.
\end{pf}

Since $\gslice{n}\subset\gslice{n-1}$, there is a natural
transformation
\[
P^{n}X\to P^{n-1}X.
\]

\begin{defin}
\label{def:24} The {\em slice tower} of $X$ is the tower $\{P^{n}X
\}_{n\in\Z}$.  The spectrum $P^{n}X$ is the {\em $n^{\text{th}}$ slice
section} of $X$.
\end{defin}

When considering more than one group, we will write
$P^{n}X=P^{n}_{G}X$ and $P_{n}X=P_{n}^{G}X$.

Let $\slc{n}{X}$ be the fiber of the map
\[
P^{n}X\to P^{n-1}X.
\]
\begin{defin}
\label{def:1} The {\em $n$-slice} of a spectrum $X$ is $\slc{n}X$.  A
spectrum $X$ is an {\em $n$-slice} if $X=\slc{n}{X}$.
\end{defin}

The spectrum $P_{n+1}X$ is analogous to the $n$-connected cover of
$X$, and for two values of $n$ they coincide.  The following is
a straightforward consequence of Proposition~\ref{thm:65}.

\begin{prop}
\label{thm:33} For any spectrum $X$, $P_{0}X$ is the $(-1)$-connected
cover of $X$ and $P_{-1}X$ is the $(-2)$-connected cover of $X$.
The $(-1)$-slice of $X$ is given by
\[
\slc{-1}{X} = \Sigma^{-1}H\pim_{-1}X.
\]
\qed
\end{prop}

The formation of slice sections and therefore of the slices themselves
behave well with respect to change of group.

\begin{prop}
\label{thm:26} 
The functor $P^{n}$ commutes with both restriction to a
subgroup and left induction.  More precisely, given $H\subset G$ there
are natural weak equivalences
\[
i_{H}^{\ast}(P^{n}_{G}X) \to P^{n}_{H}(i_{H}^{\ast}X)
\]
and
\[
G_{+}\smashove{H} (P^{n}_{H}X)\to P^{n}_{G}(G_{+}\smashove{H}X).
\]
\end{prop}

\begin{pf}
This is an easy consequence of Lemma~\ref{thm:7} and Proposition~\ref{thm:67}.
\end{pf}

\begin{rem}
\label{rem:14}
When $G$ is the trivial group the slice cells are just ordinary cells
and the slice tower becomes the Postnikov tower.
It therefore follows from Proposition~\ref{thm:26} that the tower of
non-equivariant spectra underlying the slice tower is the Postnikov tower.
\end{rem}

\subsection{Multiplicative properties of the slice tower}
\label{sec:mult-prop-slice}

The slice filtration does not quite have the multiplicative properties
one might expect.  In this section we collect a few results describing
how things work.  One important result is Corollary~\ref{thm:177}
asserting that the slice sections of a $(-1)$-connected commutative or
associative algebra are $(-1)$-connected commutative or associative
algebras.  We'll show in \S\ref{sec:furth-mult-prop} the slice
filtration is multiplicative for the special class of ``pure''
spectra, defined in \S\ref{sec:even-spectra}.

\begin{lem}
\label{thm:71} Smashing with $S^{m\rho_{G}}$ gives a bijection of the
set of slice cells $\slicecell$ with $\dim \slicecell=k$ and those with
$\dim \slicecell=k+m|G|$.
\end{lem}

\begin{pf}
Since the restriction of $\rho_{G}$ to $K\subset G$ is $|G/K|\rho_{K}$ there is
an identity 
\[
S^{\rho_{G}}\wedge (G_{+}\smashove{K}S^{m\rho_{K}}) \approx 
G_{+}\smashove{K}(S^{\rho_{G}}\wedge S^{m\rho_{K}}) \approx 
G_{+}\smashove{K} S^{(|G/K|+m)\rho_{K}}.
\]
The result follows easily from this.
\end{pf}

\begin{cor}
\label{thm:29}
Smashing with $S^{m\rho_{G}}$ gives an equivalence
\[
\geslice{n} \to \geslice{n+m|G|}.
\]
\qed
\end{cor}

\begin{cor}
\label{thm:70}
The natural maps
\begin{align*}
S^{m\rho_{G}}\wedge P_{k+1} X &\to P_{k+m|G|+1}\left(S^{m\rho_{G}}\wedge
X\right) \\
S^{m\rho_{G}}\wedge P^{k} X &\to P^{k+m|G|}\left(S^{m\rho_{G}}\wedge
X\right) \\
\end{align*}
are weak equivalences.  \qed
\end{cor}

\begin{prop}
\label{thm:6}
If $X\ge n$,  $Y\ge m$, and $n$ is divisible by $|G|$ then $X\wedge Y\ge n+m$.
\end{prop}

\begin{pf}
By smashing $X$ with $S^{(-n/|G|)\rho_{G}}$ and using
Corollary~\ref{thm:70} we may assume $n=0$.  Suppose that $Z<m$.
Since $Y\ge m$, the zero space of function spectrum $Z^{Y}$ is
contractible, and so $Z^{Y}$ is $0$-coconnected.   Since $X$ is
$(-1)$-connected (Proposition~\ref{thm:65})
\[
\gspectra{G}(X\wedge Y,Z) \approx 
\gspectra{G}(X, Z^{Y})
\]
is contractible and so $X\wedge Y\ge m$.
\end{pf}

\begin{defin}
\label{def:25}
A map $X\to Y$ is a {\em $P^{n}$-equivalence} if
$P^{n}X\to P^{n}Y$ is an equivalence.   Equivalently, $X\to Y$ is a
$P^{n}$-equivalence if for every $Z<n$, the map 
\[
\gspectra{G}(Y,Z) \to \gspectra{G}(X,Z)
\]
is a weak equivalence.
\end{defin}

\begin{lem}
\label{thm:104} If the homotopy fiber $F$ of $f:X\to Y$ is in $\gslice{n}$,
then $f$ is a a $P^{n}$ equivalence.
\end{lem}

\begin{pf}
Immediate from the fibration sequence 
\[
\gspectra{G}(Y,Z) \to \gspectra{G}(X,Z) \to \gspectra{G}(F,Z).
\]

\end{pf}

\begin{rem}
\label{rem:9} The converse of the above result is not true.  For
instance, $\ast\to S^{0}$ is a $P^{-1}$-equivalence, but the fiber
$S^{-1}$ is not in $\gslice{-1}$.
\end{rem}

\begin{lem}
\label{thm:76}
\begin{thmList}
\item If $Y\to Z$ is a $P^{n}$-equivalence and $X\ge 0$, then $X\wedge
Y\to X\wedge Z$ is a $P^{n}$-equivalence;
\item For $X_{1},\dots,X_{k}\in\geslice{0}$, the map
\[
X_{1}\wedge
\dots \wedge X_{k}\to P^{n}X_{1}\wedge \dots\wedge P^{n}X_{k}
\]
is a $P^{n}$-equivalence.
\end{thmList}
\end{lem}

\begin{pf}
Since $P_{n+1}X$ and $P_{n+1}Y$ are both slice $n$-positive the
vertical map in the square below are $P^{n}$-equivalences by
Lemmas~\ref{thm:104} and~\ref{thm:6}
\[
\xymatrix{
X\wedge Y  \ar[r]\ar[d]  & X\wedge Z  \ar[d] \\
X\wedge P^{n}Y  \ar[r]        & X\wedge P^{n}Z\mathrlap{\ .} 
}
\]
The bottom row is a weak equivalence by assumption.  It follows that
the top row is a $P^{n}$-equivalence.  The second assertion is proved
by induction on $k$, the case $k=1$ being trivial.  For the induction
step consider
\[
\xymatrix{
X_{1}\wedge \dots\wedge X_{k-1}\wedge X_{k}  \ar[r]  & P^{n}X_{1}\wedge
\dots\wedge P^{n}X_{k-1}\wedge X_{k}  \ar[d] \\
& P^{n}X_{1}\wedge \dots\wedge P^{n}X_{k-1}\wedge P^{n}X_{k}.
}
\]
The first map is a $P^{n}$-equivalence by the induction hypothesis and
part \thmListItem{1}.    The second map is a $P^{n}$-equivalence by part \thmListItem{1}.
\end{pf}

\begin{rem}
\label{rem:29}
Lemma~\ref{thm:76} can be described as asserting that the functor
\[
P^{n}:
\left\{\text{$(-1)$-connected spectra}  \right\}
\to
\left\{\gslice{n}\text{-null spectra}  \right\}
\]
is weakly monoidal.
\end{rem}

\begin{cor}
\label{thm:177} Let $R$ be a $(-1)$-connected $G$-spectrum.  If $R$ is
a homotopy commutative or homotopy associative algebra, then so is
$P^{n}R$ for all $n$. \qed
\end{cor}

The following additional results are proved in
\S\ref{sec:slice-tower-symm}.  The first two are
Propositions~\ref{thm:98},~\ref{thm:91}, and the third is easily
deduced from Proposition~\ref{thm:79}.

\begin{prop}
\label{thm:100}
Suppose that $n\ge 0$ is an integer.  If $A$ is a slice
$(n-1)$-positive $H$-spectrum then $\norm_{H}^{G}A$ is a slice
$(n-1)$-positive $G$-spectrum.  \qed
\end{prop}

\begin{prop}
\label{thm:24}
Suppose that $n\ge 0$ is an integer.  If $A$ is a slice
$(n-1)$-positive $G$-spectrum then for every $m>0$, the symmetric
smash power $\sym^{m}A$ is slice $(n-1)$-positive. \qed
\end{prop}

\begin{prop}
\label{thm:49}
Suppose that $n\ge 0$ is an integer.  If $R$ is a
$(-1)$-connected equivariant commutative ring, then the slice section
$P^{n}R$ can be given a the structure of an equivariant commutative
ring in such a way that $R\to P^{n}R$ is a commutative ring
homomorphism.  Moreover this commutative ring structure is unique.
\qed
\end{prop}

\subsection{The slice spectral sequence}
\label{sec:conv-slice-tower}

The {\em slice spectral sequence} is the homotopy spectral of the
slice tower.  The main point of this section is to establish strong
convergence of the slice spectral sequence, and to show that for any
$X$ the $E_{2}$-term is distributed in the gray region of
Figure~\ref{fig:1}.   We begin with some results relating the slice
sections to Postnikov sections.

\subsubsection{Connectivity and the slice filtration}
\label{sec:constr-slice-cells}

Our convergence result for the slice spectral sequence depends on
knowing how slice cells are constructed from $G$-cells.  We will say
that a space or spectrum $X$ {\em decomposes} into the elements of a
collection of spectra $\{T_{\alpha} \}$ if $X$ is weakly equivalent to
a spectrum $\tilde X$ admitting an increasing filtration
\[
X_{0}\subset X_{1}\subset \cdots
\]
with the property that $X_{n}/X_{n-1}$ is weakly equivalent to a wedge
of $T_{\alpha}$.  

\begin{rem}
\label{rem:43} A $G$-spectrum $X$ decomposes into a collection of
spectra $\{G/H_{+}\wedge S^{m} \}$, with $H$ and $m$ ranging through
some indexing list, if and only if $X$ is weakly equivalent a $G$-CW
spectrum with $G$-cells of the form $G/H_{+}\wedge D^{m}$, with $H$
and $m$ ranging through the same list.
\end{rem}

\begin{rem}
\label{rem:57}
To say that $X$ decomposes into the elements of a collection of
compact objects $\{T_{\alpha} \}$ means that $X$ is in the smallest
subcategory of $\gspectra{G}$ closed under weak equivalences,
arbitrary wedges, and the formation of mapping cones and extensions
(\ie\ the properties listed in Remark~\ref{rem:23}).
\end{rem}

\begin{lem}
\label{thm:37} Let $\slicecell$ be a slice cell.  If $\dim
\slicecell=n\ge 0$, then $\slicecell$ decomposes into the spectra
$G/H_{+}\wedge S^{k}$ with $\lfloor n/|G| \rfloor \le k\le n$.  If
$\dim \slicecell=n < 0$ then $\slicecell$ decomposes into
$G/H_{+}\wedge S^{k}$ with $n\le k \le \lfloor n/|G| \rfloor$.
\end{lem}

\begin{pf}
The cell structure of $S^{\rho_{G}-1}$ described in
Example~\ref{eg:25} has $G$-cells ranging in dimension from $0$ to
$|G|-1$, and suspends to a cell decomposition of $S^{\rho_{G}}$ with
$G$-cells whose dimension ranges from $1$ to $|G|$.  The cases
$\slicecell=S^{m\rho_{G}}$ and $\slicecell=S^{m\rho_{G}-1}$ with $m\ge
0$ are handled by smashing these together and passing to suspension
spectra.  For $m<0$, Spanier-Whitehead duality gives an equivariant
cell decomposition of $S^{m\rho_{G}}$ into cells whose dimensions
range from $m|G|$ to $m$ and of $\Sigma^{-1}S^{m\rho_{G}}$ into cells
whose dimensions range from $n=m|G|-1$ to $m-1=\lfloor
n/|G|\rfloor$. Finally, the case in which $\slicecell$ is induced from
a subgroup $K\subset G$ is proved by left inducing its $K$-equivariant
cell decomposition.
\end{pf}

\begin{cor}
\label{thm:38} Let $Y\in\geslice{n}$.  If $n\ge 0$, then $Y$ can be
decomposed into the spectra $G/H_{+}\wedge S^{m}$ with $m\ge \lfloor
n/|G|\rfloor$.  If $n\le 0$ then $Y$ can be decomposed into
$G/H_{+}\wedge S^{m}$ with $m\ge n$.
\end{cor}

\begin{pf}
The class of $G$-spectra $Y$ which can be decomposed into
$G/H_{+}\wedge S^{m}$ with $m\ge \lfloor n/|G|\rfloor$ is closed under
weak equivalences, homotopy colimits, and extensions.  By
Lemma~\ref{thm:37} it contains the slice cells $\slicecell$ with $\dim
\slicecell\ge n$.  It therefore contains all $Y\in\geslice{n}$ by
Remark~\ref{rem:23}.  A similar argument handles the case $n<0$.
\end{pf}

\begin{prop}\label{thm:40}
\begin{thmList}  Write $g=|G|$.
\item\label{item:1} If $n\ge 0$ and $k>n$ then $(G/H)_{+}\wedge S^{k}
>n$.
\item\label{item:2} If $m\le -1$ and $k\ge m$ then $(G/H)_{+}\wedge
S^{k} \ge (m+1)g-1$.
\item If $Y\ge n$ with $n\ge 0$, then $\pim_{i}Y=0 \text{ for }i <
\lfloor n/g\rfloor$.
\item If $Y \ge n$ with $n\le 0$, then $\pim_{i}Y=0$ for $i<n$.
\end{thmList}
\end{prop}

\begin{pf}
We start with the first assertion.  We will prove the claim by
induction on $|G|$, the case of the trivial group being obvious.
Using Proposition~\ref{thm:67} we may assume by induction that
$(G/H)_{+}\wedge S^{k}> n$ when $n\ge 0$ and $H\subset G$ is a proper
subgroup.  This implies that if $T$ is an equivariant CW-spectrum
built from $G$-cells of the form $(G/H)_{+}\wedge S^{k}$ with $k>n$
and $H\subset G$ a proper subgroup, then $T>n$.  The homotopy fiber of
the natural inclusion
\[
S^{k}\to S^{k\rho_{G}}
\]
can be identified with the suspension spectrum of
$S(k\rho_{G}-k)_{+}\wedge S^{k}$, and so is such a $T$.  Since
$S^{k\rho_{G}}\ge k|G|\ge k g>n$ the fibration sequence
\[
T\to S^{k}\to S^{k\rho_{G}}
\] 
exhibits $S^{k}$ as an extension of two slice $n$-positive spectra,
making it slice $n$-positive.  The second assertion is trivial for
$k\ge 0$ since in that case $G/H_{+}\wedge S^{k}\ge 0$ and
$(m+1)g-1\le -1$.  The case $k\le -1$ is handled by writing
\[
(G/H)_{+}\wedge S^{k} =
\Sigma^{-1}(G/H)_{+}\wedge S^{(k+1)\rho_{G}}\wedge S^{-(k+1)(\rho_{G}-1)}.
\]
Since $-(k+1)\ge 0$, the spectrum $S^{-(k+1)(\rho_{G}-1)}$ is a suspension spectrum
and so
\[
(G/H)_{+}\wedge S^{k} \ge (k+1)g-1 \ge (m+1)g-1.
\]
The third and fourth assertions are immediate from
Corollary~\ref{thm:38}.
\end{pf}

\begin{rem}
We've stated part~\thmItemref{item:2} of Proposition~\ref{thm:40} in
the form in which it is most clearly proved.  When it comes up, it is
needed as the implication that for $n<0$,
\[
k \ge \lfloor (n+1)/g \rfloor \implies G/H_{+}\wedge S^{k} > n.
\]
To relate these, write $m=\lfloor (n+1)/g \rfloor$, so that
\[
m+1 > (n+1)/g
\]
and by part~\thmItemref{item:2} of Proposition~\ref{thm:40}
\[
G/H_{+}\wedge S^{k} \ge (m+1)g-1 > n.
\]
\end{rem}

\subsubsection{The spectral sequence}
\label{sec:spectral-sequence}

\begin{figure}
\includegraphics[]{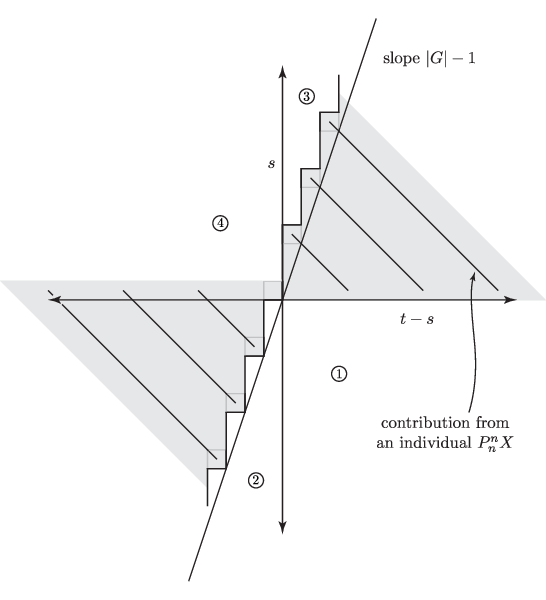}
\caption{The slice spectral sequence}
\label{fig:1}
\end{figure}

The slice spectral sequence is the spectral sequence associated to the
tower of fibration $\{P^{n}X \}$, and it takes the form
\[
E_{2}^{s,t}=\pi_{t-s}^{G}\slc{t}{X}\implies \pi_{t-s}^{G}X.
\]
It can be regarded as a spectral sequence of Mackey functors, or of
individual homotopy groups.  We have chosen our indexing so that the
display of the spectral sequence is in accord with the classical Adams
spectral sequence: the $E_{r}^{s,t}$-term is placed in the plane in
position $(t-s,s)$.  The situation is depicted in Figure~\ref{fig:1}.
The differential $d_{r}$ maps $E_{r}^{s,t}$ to $E_{r}^{s+r, t+r-1}$,
or in terms the display in the plane, the group in position $(t-s,s)$
to the group in position $(t-s-1,s+r)$.

The following is an immediate consequence of Proposition~\ref{thm:40}.
As there, we write $g=|G|$.

\begin{thm}
\label{thm:34} Let $X$ be a $G$-spectrum.  The
Mackey functor homotopy groups of $P^{n}X$ satisfy
\[
\pim_{k}P^{n}X = 0\text{ for }
\begin{cases}
k>n &\quad \text{ if }n\ge 0 \\
k \ge \lfloor (n+1)/g\rfloor &\quad \text{ if }n < 0
\end{cases}
\]
and the map $X\to P^{n}X$ induces an isomorphism
\[
\pim_{k}X\xrightarrow{\approx}{} \pim_{k}P^{n}X\text{ for }
\begin{cases}
k<\lfloor (n+1)/g \rfloor &\quad \text{ if }n\ge 0 \\
k \le n &\quad \text{ if }n < 0\mathrlap{\ .}
\end{cases}
\]
Thus for any $X$,
\[
\varinjlim_{n} P^{n}X
\]
is contractible, the map
\[
X\to \varprojlim_{n} P^{n}X
\]
is a weak equivalence, and for each $k$, the map
\[
\{\pim_{k}(X) \}\to \{\pim_{k}P^{n}X \}
\]
from the constant tower to the slice tower of Mackey functors
is a pro-isomorphism. \qed
\end{thm}

\begin{cor}
\label{thm:36}
If $M$ is an $n$-slice then
\[
\pim_{k}M=0
\]
if $n\ge 0$ and $k$ lies outside of the region $\lfloor n/g \rfloor\le
k\le n$, or if $n<0$ and $k$ lies outside of the region $n\le k <
\lfloor (n+1)/g\rfloor$. \qed
\end{cor}

Theorem~\ref{thm:34} gives the strong convergence of the slice
spectral sequence, while Corollary~\ref{thm:36} shows that the
$E_{2}$-term vanishes outside of a restricted range of dimensions.
The situation is depicted in Figure~\ref{fig:1}.  The homotopy groups
of individual slices lie along lines of slope $-1$, and the groups
contributing to $\pim_{\ast}P^{n}X$ lie to the left of a line of slope
$-1$ intersecting the $(t-s)$-axis at $(t-s)=n$.  All of the groups
outside the gray region are zero.  The vanishing in the regions
labeled $1$-$4$ correspond to the four parts of
Proposition~\ref{thm:40}.

Proposition~\ref{thm:40} also gives a relationship between the Postnikov
tower and the slice tower.

\begin{cor}
\label{thm:23}
If $X$ is an $(n-1)$-connected $G$-spectrum with $n\ge 0$ then $X\ge n$.
\end{cor}

\begin{pf}
The assumption on $X$ means it is weakly equivalent to a $G$-CW
spectrum having $G$-cells $G/H_{+}\wedge S^{m}$ only in dimensions $m\ge
n$.  By part~\thmItemref{item:1} of Proposition~\ref{thm:40} these
are in $\geslice{n}$.
\end{pf}

We end this section with an application.  The next result says that if
a tower looks like the slice tower, then it is the slice tower.

\begin{prop}
\label{thm:41}
Suppose that $X\to \{\tilde P^{n}\}$ is a map from $X$ to a tower of
fibrations with the properties
\begin{thmList}
\item\label{limit} the map $X\to \varprojlim\tilde P^{n}$ is a weak equivalence;
\item\label{colimit} the spectrum $\varinjlim_{n}\tilde P^{n}$ is contractible;
\item\label{fiber} for all $n$, the fiber of the map $\tilde P^{n}\to\tilde
P^{n-1}$ is an $n$-slice.
\end{thmList}
Then $\tilde P^{n}$ is the slice tower of $X$.   
\end{prop}

\begin{pf}
We first show that $\tilde P^{n}$ is slice $(n+1)$-null.  We
will use the criteria of Lemma~\ref{thm:2}.  Suppose that $\slicecell$
is a slice cell with $\dim\slicecell >n$.  By
condition~\thmItemref{fiber}, the maps
\[
[\slicecell,\tilde P^{n} ]^{G} \to 
[\slicecell,\tilde P^{n-1} ]^{G} \to 
[\slicecell,\tilde P^{n-2} ]^{G} \to\cdots
\]
are all monomorphisms.  Since $\slicecell$ is finite, the map 
\[
\varinjlim_{k\le n}[\slicecell,\tilde P^{k} ]^{G} \to 
[\slicecell,\varinjlim_{k\le n}\tilde P^{k} ]^{G}
\]
is an isomorphism.  It then follows from
assumption~\thmItemref{colimit} that $[\slicecell, \tilde
P^{n}]^{G}=0$.  This shows that $\tilde P^{n}$ is slice $(n+1)$-null.
Now let $\tilde P_{n+1}$ be the homotopy fiber of the map $X\to \tilde
P^{n}$.  By Lemma~\ref{thm:7}, the result will follow if we can show
$\tilde P_{n+1}>n$.  By assumption~\thmItemref{fiber}, for any
$N>n+1$, the spectrum
\[
\tilde P_{n+1}\cup C\tilde P_{N}
\]
admits a finite filtration whose layers are $m$-slices, with $m\ge
n+1$.   It follows that 
\[
\tilde P_{n+1}\cup C\tilde P_{N}>n.
\]
In view of the cofibration sequence
\[
\tilde P_{N}\to \tilde P_{n+1}\to \tilde P_{n+1}\cup C\tilde P_{N},
\]
to show that $\tilde P_{n+1}>n$ it suffices to show that $\tilde
P_{N}>n$ for {\em some} $N>n$.

Let $Z$ be any slice $(n+1)$-null spectrum.  We
need to show that the $G$-space
\[
\gspectra{G}(\tilde P_{N},Z)
\]
is contractible.  We do this by studying the Mackey functor homotopy
groups of the spectra involved, and appealing to an argument using the
usual equivariant notion of connectivity.  By Theorem~\ref{thm:34},
there is an integer $m$ with the property that for $k>m$,
\[
\pim_{k}Z= 0.
\]
By Corollary~\ref{thm:36} and assumption~\thmItemref{fiber}, for
$N\gg0$ and any $N'>N$, 
\[
\pim_{k}\tilde P_{N}\cup C\tilde P_{N'} = 0, \quad k\le m,
\]
so
\[
\pim_{k}\tilde P_{N'}\to \pim_{k}\tilde P_{N}
\]
is an isomorphism for $k\le m$.  Since $\holim_{N'}\tilde P_{N'}$ is
contractible this implies that for $N\gg0$
\[
\pim_{k}\tilde P_{N}=0, \quad k\le m.
\]
Thus for $N\gg0$,  $\tilde P_{N}$ is $m$-connected in the usual sense
and so 
\[
\gspectra{G}(\tilde P_{N},Z)
\]
is contractible.
\end{pf}

\subsection{The $RO(G)$-graded slice spectral sequence}
\label{sec:rog-graded-slice}
Applying $RO(G)$-graded homotopy groups to the slice
tower leads to an $RO(G)$-graded slice spectral sequence
\[
E_{2}^{s,V}= \pi^{G}_{V-s}\slc{\dim V}{X} \implies \pi^{G}_{V-s}X.
\]
The grading convention is chosen so that it restricts to the one of
\S\ref{sec:spectral-sequence} when $V$ is a trivial virtual
representation.  The $r^{\text{th}}$ differential is a map
\[
d_{r}:E_{2}^{s,V} \to E_{2}^{s+r,V+(r-1)}.
\]
The $RO(G)$-graded slice spectral sequence is a sum of spectral
sequences, one for each element of $RO(G)/\Z$.  We will call the
spectral sequence corresponding to the coset $V+\Z\in\RO(G)/\Z$ the
{\em slice spectral sequence for $\pi^{G}_{V+\ast}X$.}  This spectral
sequence can be displayed on the $(x,y)$-plane, and we will do so
following Adams conventions, with the term $E_{2}^{s,V+t}$ displayed
at a position with $x$-coordinate $(V+t-s)$ and $y$-coordinate $s$.
For an example, see Figures~\ref{fig:4}, \ref{fig:5} and~\ref{fig:6}
in \S\ref{sec:periodicity-theorem}.  

\subsection{Special slices} 
\label{sec:special-slices}

In this section we investigate special slices of spectra, and
introduce the notion of a {\em spectrum with cellular slices}, and of
a {\em pure} $G$-spectrum.  Our main result
(Proposition~\ref{thm:189}) asserts that a map $X\to Y$ of $G$-spectra
with cellular slices is a weak equivalence if and only if the
underlying map of non-equivariant spectra is.  This result plays an
important role in the proof of the Reduction Theorem in
\S\ref{sec:reduction-theorem}.  We also include material useful for
investigating the slices of more general spectra.

\subsubsection{Slice positive spectra, $0$-slices and $(-1)$-slices}
\label{sec:slice-posit-spectra}

In this section we will describe methods for determining the slices of
spectra, and introduce a convenient class of equivariant spectra.  Our
first results make use of the isotropy separation sequence
(\S\ref{sec:isotropy-separation}) obtained by smashing with the
cofibration sequence of pointed $G$-spaces
\[
E\pfamily_{+}\to S^{0} \to \tilde E\pfamily.
\]
The space $E\pfamily_{+}$ is an equivariant CW-complex built from
$G$-cells of the form $(G/H)_{+}\wedge S^{n}$ with $H\subset G$ a
proper subgroup.  It follows that if $W$ is a pointed $G$-space whose
$H$-fixed points are contractible for all proper $H\subset G$, then
$\gspaces{G}(E\pfamily_{+}, W)$ is contractible.

\begin{lem}
\label{thm:74}
Fix an integer $d$.  If $X$ is a $G$-spectrum with the property that
$i_{H}^{\ast}X>d$ for all proper $H\subset G$, then
$E\pfamily_{+}\wedge X>d$.
\end{lem}

\begin{pf}
Suppose that $Z\le d$.  Then 
\[
\gspectra{G}(E\pfamily_{+}\wedge X,Z) \approx
\gspaces{G}(E\pfamily_{+},\gspectra{G}(X,Z)).
\]
By the assumption on $X$, the $G$-space $\gspectra{G}(X,Z)$ has
contractible $H$ fixed points for all proper $H\subset G$.    The
Lemma now follows from the remark preceding its statement.
\end{pf}

\begin{lem}
\label{thm:75} Write $g=|G|$.  The suspension spectrum of $\tilde
E\pfamily$ is in $\geslice{g-1}$.
\end{lem}

\begin{pf}
The map $\tilde E\pfamily\wedge S^{0}\to \tilde E\pfamily\wedge
S^{\rho_{G}-1}$ is a weak equivalence (Proposition~\ref{thm:202} and
Remark~\ref{rem:8}).  The suspension spectrum of $\tilde E\pfamily$ is
in $\geslice{0}$, since it is $(-1)$-connected
(Proposition~\ref{thm:65}).  So $\tilde E\pfamily\wedge
S^{\rho_{G}-1}\ge g-1$ by Proposition~\ref{thm:6}.
\end{pf}

\begin{prop}
\label{thm:64} 
A $G$-spectrum $X$ is slice positive if and only if it
is $(-1)$ connected and $\piu_{0}X=0$ (\ie the non-equivariant
spectrum $i_{0}^{\ast}X$ underlying $X$ is $0$-connected).
\end{prop}

\begin{pf}
The only if assertion follows from the fact that the slice cells of
positive dimension are $(-1)$-connected and have $0$-connected
underlying spectra.  The ``if'' assertion is proved by induction on
$|G|$, the case of the trivial group being trivial.  For the induction
step we may assume $X$ is $(-1)$-connected and has the property that
$i_{H}^{\ast}X>0$ for all proper $H\subset G$.  Consider the isotropy
separation sequence for $X$
\[
E\pfamily_{+}\wedge X\to X\to \tilde E\pfamily \wedge X.
\]
The leftmost term is slice-positive by Lemma~\ref{thm:74}, and the
rightmost term is by Propositions~\ref{thm:65} and~\ref{thm:6}, and
Lemma~\ref{thm:75}.  It follows that $X$ is slice-positive.
\end{pf}

\begin{eg}
\label{eg:7} Suppose that $f:S\to S'$ is a surjective map of $G$-sets.
Proposition~\ref{thm:64} implies that the suspension spectrum of the
mapping cone of $f$ is slice positive.  This implies that if $H\Mm$ is
an Eilenberg-MacLane spectrum which is a zero slice then for every
surjective $S\to S'$ the map $\Mm(S')\to \Mm(S)$ is a monomorphism.  The
proposition below shows that this is also a sufficient condition.
\end{eg}

\begin{prop}
\label{thm:185} 
\begin{thmList}
\item A spectrum $X$ is a $(-1)$-slice if and only if it is of the
form $X=\Sigma^{-1}H\Mm$, with $\Mm$ an arbitrary Mackey functor.
\item\label{restriction-mono} A spectrum $X$ is a $0$-slice if and
only if it is of the form $H\Mm$ with $\Mm$ a Mackey functor all of whose
restriction maps are monomorphisms.
\end{thmList}
\end{prop}

\begin{rem}
\label{rem:44} The condition on $\Mm$ in \thmItemref{restriction-mono}
is that if $S\to S'$ is a surjective map of finite $G$-sets then
$\Mm(S')\to \Mm(S)$ is a monomorphism.  Let $G$ act on $G\times S$ and
$G\times S'$ through its left action on $G$.  Then $G\times S\to
G\times S'$ has a section, so $\Mm(G\times S')\to \Mm(G\times S)$ is
always a monomorphism.  Using this one easily checks that condition is
also equivalent to requiring that for every finite $G$-set $S'$, the
map $\Mm(S')\to \Mm(G\times S')$, induced by the action mapping $G\times
S'\to S'$, is a monomorphism.
\end{rem}

\begin{pf}
The first assertion is immediate from Proposition~\ref{thm:33}, which,
combined with part~\thmItemref{item:1} of Proposition~\ref{thm:40},
also shows that a $0$-slice is an Eilenberg-MacLane spectrum.
Example~\ref{eg:7} gives the ``only if'' part of the second assertion.
For the ``if'' part, suppose that $\Mm$ is a Mackey functor all of
whose restrictions maps are monomorphisms, and consider the sequence
\[
P_{1}H\Mm\to H\Mm\to P^{0}H\Mm.
\]
Since $P_{1}H\Mm\ge 0$ it is $(-1)$-connected, and so $P_{1}H\Mm$ is an
Eilenberg-MacLane spectrum.  For convenience, write 
\begin{align*}
\Mm' &=\pim_{0}P_{1}H\Mm \\
\Mm'' &=\pim_{0}P^{0}H\Mm
\end{align*}
so that there is a short exact sequence
\[
\Mm'\rightarrowtail \Mm\twoheadrightarrow \Mm''.
\]
Suppose that $S$ is any finite $G$-set and consider the following
diagram
\[
\xymatrix{
\Mm'(S)  \ar[r]\ar[d]  & \Mm(S)  \ar[d] \ar[r]  &   \ar[d]\Mm''(S)\\
\Mm'(G\times S)  \ar[r]        & \Mm(G\times S)\ar[r] & \Mm''(G\times S)
}
\]
in which the rows are short exact, and the vertical maps are induced
by the action mapping, as in Remark~\ref{rem:44}.  The bottom right
arrow is an isomorphism since $i_{0}^{\ast}H\Mm\to i_{0}^{\ast}P^{0}H\Mm$
is an equivalence.  Thus $\Mm'(G\times S)=0$ (this also follows from
Proposition~\ref{thm:64}).  The claim now follows from a simple
diagram chase.
\end{pf}

\begin{rem}
\label{rem:37} The second assertion of Proposition~\ref{thm:185} can
also be deduced directly from Corollary~\ref{thm:188}.
\end{rem}

\begin{cor}
\label{thm:186} If $X=H\Mm$ is a zero slice and $\piu_{0}X=0$ then $X$
is contractible.  \qed
\end{cor}

\begin{cor}
\label{thm:187} The $(-1)$-slice of $S^{-1}$ is $\Sigma^{-1}H\Am$.  The zero
slice of $S^{0}$ is $H\Zm$.
\end{cor}

\begin{pf}
The first assertion follows easily from Part \thmListItem{1} of
Proposition~\ref{thm:185}.  For the second assertion note that the
$S^{0}\to H\Am$ is a $P^{0}$-equivalence, so the zero slice of $S^{0}$
is $P^{0}H\Am$.  Consider the fibration sequence
\[
H\underline{I}\to H\Am\to H\Zm,
\]
in which $\underline{I}=\ker \Am\to \Zm$ is the augmentation ideal.
The leftmost term is slice positive by Proposition~\ref{thm:64}, and
the rightmost term is in $\leslice{G}{0}$ by Proposition~\ref{thm:185}.  The claim
now follows from Lemma~\ref{thm:7}.
\end{pf}

\begin{cor}
\label{thm:22}
For $K\subset G$, the $m|K|$-slice of $\slicecell(m,K)$ is 
\[
H\Zm\wedge \slicecell(m,K)
\]
and the $(m|K|-1)$-slice of $\Sigma^{-1}\slicecell(m,K)$ is 
\[
H\Am \wedge \Sigma^{-1}\slicecell(m,K).
\]
\end{cor}

\begin{pf}
Using the fact that $G_{+}\smashove{K}(\slot)$ commutes with the
formation of the slice tower (Proposition~\ref{thm:26}) it suffices to
consider the case $K=G$.  The result then follows from
Corollaries~\ref{thm:70} and ~\ref{thm:187}.
\end{pf}

\subsubsection{Cellular slices, isotropic and pure spectra}
\label{sec:even-spectra}

\begin{defin}
\label{def:13} A $d$-slice is {\em cellular} if it is of the form
$H\Zm\wedge \refine$, where $\refine$ is a wedge of slice cells of
dimension $d$.  A cellular slice is {\em isotropic} if $\refine$ can
be written as a wedge of slice cells, none of which is free (i.e., of
the form $G_{+}\wedge S^{n}$).  A cellular slice is {\em pure} if
$\refine$ can be written as a wedge of regular slice cells (those of
the form $\slicecell(m,K)$, and not $\Sigma^{-1}\slicecell(m,K)$).
\end{defin}

\begin{defin}
\label{def:31} A $G$-spectrum $X$ {\em has cellular slices} if
$\slc{n}{X}$ is cellular for all $d$, and is {\em isotropic} or {\em
pure} if its slices are isotropic or pure.
\end{defin}

\begin{lem}
\label{thm:118}
Suppose that $f:X\to Y$ is a map of cellular $d$-slices and
$\piu_{d}f$ is an isomorphism.  Then $f$ is a weak equivalence.
\end{lem}

\begin{pf}
The proof is by induction on $|G|$.   If $G$ is the trivial group,
the result is obvious since $X$ and $Y$ are Eilenberg-MacLane
spectra.   Now suppose we know the result for all proper $H\subset
G$, and consider the map of isotropy separation sequences 
\[
\xymatrix{
\epfam\wedge X  \ar[r]\ar[d]  & X  \ar[d]\ar[r] &
 \tepfam \wedge X\ar[d] \\
\epfam\wedge Y \ar[r] & Y \ar[r] & \tepfam \wedge Y\mathrlap{\ .}
}
\]
By the induction hypothesis, the left vertical map is a weak
equivalence.  If $d$ is not congruent to $0$ or $-1$ modulo $|G|$ then
the rightmost terms are contractible, since every slice cell of
dimension $d$ is induced.  Smashing with $S^{m\rho_{G}}$ for suitable
$m$, we may therefore assume $d=0$ or $d=-1$.  Smashing with $S^{1}$
in case $d=-1$ we reduce to the case $d=0$ and therefore assume that
$X=H\Mm_{0}$ and $Y=H\Mm_{1}$ with $\Mm_{0}$ and $\Mm_{1}$ permutation Mackey
functors.  The result then follows from part~\thmItemref{part:2} of
Lemma~\ref{thm:109}.
\end{pf}

\begin{prop}
\label{thm:189}
Suppose that $X$ and $Y$ have cellular slices.  If $f:X\to Y$ has the
property that $\piu_{\ast}f$ is an isomorphism.  Then $f$ is a weak
equivalence.
\end{prop}

\begin{pf}
It suffices to show that for each $d$ the induced map of slices 
\begin{equation}
\label{eq:106}
\slc{d}X\to \slc{d}Y
\end{equation}
is a weak equivalence.  Since the map of ordinary spectra underlying
the slice tower is the Postnikov tower, the map satisfies the
conditions of Lemma~\ref{thm:118}, and the result follows.
\end{pf}

For certain slices, the condition on $Y$ in Proposition~\ref{thm:189} 
can be dropped.  

\begin{lem}
\label{thm:157} Suppose that $f:X\to Y$ is a map of $0$-slices and $X$
is cellular.  If $\piu_{0}f$ is an isomorphism then $f$ is an
equivalence.
\end{lem}

\begin{pf}
Write $X=H\Mm$ and $Y=H\Mm'$, and let $S$ be a finite $G$-set.  Consider
the diagram
\[
\xymatrix{
\Mm(S)  \ar[r]^{\therefore\sim}\ar[d]_{\sim}  & \Mm'(S) 
  \ar[d]^{\text{mono}} \\
\Mm(G\times S)^{G} \ar[r]_{\sim} & \Mm'(G\times S)^{G} }
\]
in which the vertical maps come from the action mapping $G\times S\to
S$ (see the discussion preceding Lemma~\ref{thm:109}).  The bottom
arrow is an isomorphism by assumption.  The vertical maps are
monomorphisms by Proposition~\ref{thm:185}.  The left vertical map is
an isomorphism since $\Mm$ is a permutation Mackey functor
(part~\thmItemref{part:3} of Lemma~\ref{thm:109}).  The result follows.
\end{pf}

\begin{prop}
\label{thm:66} Suppose that $f:X\to Y$ is a map of $d$-slices, $X$ is
cellular and $d\not\equiv -1\mod p$ for any prime $p$ dividing $|G|$.  If
$\piu_{d}X\to \piu_{d}Y$ is an isomorphism then $f$ is a weak
equivalence.
\end{prop}

\begin{pf}
Let $C$ be the mapping cone of $f$.  We know that $C\ge d$.  We will
show that 
\[
[\slicecell,C ]^{G} = 0
\]
for all slice cells $\slicecell$ with $\dim\slicecell \ge d$.  This
will show (Lemma~\ref{thm:2}) that $C<d$ and hence must be
contractible since its identity map is null.  The assertion is obvious
when $G$ is the trivial group.  By induction on $|G|$ we may assume
$\slicecell$ is not induced.  If $d$ is divisible by $|G|$ we may
smash with $S^{-d/|G|\,\rho_{G}}$ and reduce to the case $d=0$ which
is Lemma~\ref{thm:157}.  It remains to show that $\pi^{G}_{m\rho_{G}}
C=0$ when $m|G|\ge d$ and that $\pi^{G}_{m\rho_{G}-1}C=0$ when
$m|G|-1\ge d$.  Since
\[
d\not\equiv 0, -1\mod |G|,
\]
the conditions in fact implies $m|G|-1>d$.  So we are in the situation
$m|G|-1>d$ and we need to show that both $\pi^{G}_{m\rho_{G}} C$ and
$\pi^{G}_{m\rho_{G}-1} C$ are zero.  The exact sequence
\[
\pi^{G}_{m\rho_{G}}Y \to \pi^{G}_{m\rho_{G}}C \to  \pi^{G}_{m\rho_{G}-1}X 
\]
gives the vanishing of $\pi^{G}_{m\rho_{G}}C$.  For the remaining case
consider the exact sequence
\[
\pi^{G}_{m\rho_{G}-1}Y \to \pi^{G}_{m\rho_{G}-1}C \to
\pi^{G}_{m\rho_{G}-2}X\to \pi^{G}_{m\rho_{G}-2}Y.
\]
As above, the left group vanishes since $Y$ is a $d$-slice and
$S^{m\rho_{G}-1}>d$.  Lemma~\ref{thm:156} below implies that the left
vertical map in
\[
\xymatrix{
\pi^{G}_{m\rho_{G}-2}X  \ar[r]\ar[d]  &   \pi^{G}_{m\rho_{G}-2}Y\ar[d] \\
\piu_{mg-2}X  \ar[r]_{\approx}        & \piu_{mg-2}Y
}	
\]
is monomorphism, and therefore so is the top horizontal map.
Thus $\pi^{G}_{m\rho_{G}-1}C=0$ by exactness.
\end{pf}

\begin{lem}
\label{thm:156}
Suppose $\slicecell$ is a slice cell of dimension $d$.  If $m|G|-1>
d$ then the restriction mapping
\[
\pi^{G}_{m\rho_{G}-2}H\Zm\wedge \slicecell\to 
\piu_{mg-2}H\Zm\wedge \slicecell
\]
is a monomorphism.
\end{lem} 

\begin{pf}
When $G$ is trivial the map is an isomorphism.  By induction on $|G|$
we may therefore assume $G$ is not the trivial group and that
$\slicecell$ is not induced, in which case $\slicecell=S^{k\rho_{G}}$
or $\slicecell= S^{k\rho_{G}-1}$.  Note that
\[
S^{m\rho_{G}-2} = S^{(m-1)\rho_{G}-1}\wedge S^{\rho_{G}-1}\ge (m-1)|G|-1 > (m-2)|G|
\]
so that both $\pi^{G}_{m\rho_{G}-2}H\Zm\wedge S^{k\rho_{G}}$ and
$\pi^{G}_{m\rho_{G}-2}H\Zm\wedge S^{k\rho_{G}-1}$ are zero unless $k=m-1$.
The group $\pi^{G}_{m\rho_{G}-2}H\Zm \wedge S^{(m-1)\rho_{G}-1}$ is
zero since it is isomorphic to
\[
\pi^{G}_{m\rho_{G}-1}H\Zm\wedge S^{(m-1)\rho_{G}}
\]
and $S^{m\rho_{G}-1}\ge m|G|-1 > (m-1)|G|$.   This leaves the group
\[
\pi^{G}_{m\rho_{G}-2}H\Zm\wedge S^{(m-1)\rho_{g}} \approx
\pi^{G}_{\rho_{G}-2}H\Zm
\]
whose triviality was established in Example~\ref{eg:8}.
\end{pf}

\subsubsection{The special case in which $G$ is a finite $2$-group}
\label{sec:special-case-g=ztn}

In this section we record some results which are special to the case
in which $G$ has order a power of $2$.  The results about even slices
are used in the proof of the Reduction Theorem in
\S\ref{sec:cert-auxil-spectra}.  The results on odd slices were used
in an earlier approach to the main results of this paper, but are no
longer.  We include them here because they provide useful tools for
investigating slices of various spectra.  Throughout this section the
group $G$ will be a finite $2$-group.

Suppose that $X$ is a $G$-spectrum with the property that $\piu_{d}X$ is a
free abelian group.   In \S\ref{sec:refinement-homotopy} we will define 
a {\em refinement of $\piu_{d}X$} to be a map
\[
c:\refine\to X
\]
in which $\refine$ is a wedge of slice cells of dimension $d$, with
the property that the map $\piu_{d}\refine\to\piu_{d}X$ is an
isomorphism.

\begin{prop}
\label{thm:81} 
If $\refine\to X$ is a refinement of $\piu_{2k}X$ then
the canonical map
\[
H\Zm\wedge \refine\to \slc{2k}{X}
\]
is an equivalence.
\end{prop}

\begin{pf}
By Corollary~\ref{thm:22} (and the fact that the formation of slices
commutes with the formation of wedges), the map
\[
\refine\to H\Zm\wedge\refine
\]
induces an equivalence
\[
\slc{2k}\refine\to H\Zm\wedge\refine.
\]
Applying $\slc{2k}$ to $\refine\to X$ then leads to a map
\[
H\Zm\wedge \refine \to\slc{2k}{X}
\]
which, since the slice tower refines the Postnikov tower, is an
equivalence of underlying non-equivariant spectra.  The result now
follows from Proposition~\ref{thm:66}, since the only prime dividing
$|G|$ is $2$.
\end{pf}

Proposition~\ref{thm:81} gives some control over the even slices of a
$G$-spectrum $X$ when $G$ is a $2$-group.  The odd slices are
something of a different story, and getting at them requires some
knowledge of the equivariant homotopy type of $X$.  Note that by
Proposition~\ref{thm:185} any Mackey functor can occur in an odd
slice.  On the other hand, only special ones can occur in even slices.

\begin{cor}\label{thm:86}
If $\slicecell$ is a slice cell of odd dimension $d$, then  for any $X$,
\[
[\slicecell,X]^{G} \approx [\slicecell,\slc{d}{X}]^{G}.
\]
\end{cor}

\begin{pf}
Since the formation of $\slc{d}{X}$ commutes with the functors
$i_{H}^{\ast}$, induction on $|G|$ reduces us to the case when
$\slicecell$ is not an induced slice cell.  So we may assume
$\slicecell=S^{m\rho_{G}-1}$.  Smashing $\slicecell$ and $X$ with
$S^{-m\rho_{G}}$, and using Corollary~\ref{thm:70} reduces to the case
$m=0$, which is given by Proposition~\ref{thm:33}.
\end{pf}

The situation most of interest to us in this paper is when the odd
slices are contractible.  Proposition~\ref{thm:85} below gives a
useful criterion.

\begin{prop}
\label{thm:85} For a $G$ spectrum $X$ and an odd integer $d$, the
following are equivalent:

\begin{thmList}
\item The $d$-slice of $X$ is contractible;
\item For every slice cell $\slicecell$ of dimension $d$, $[\slicecell,X]^{G}=0$.
\end{thmList}
\end{prop}

\begin{pf}
By Corollary~\ref{thm:86} (which requires $d$ to be odd), there is an
isomorphism
\[
[\slicecell,X]^{G}=[\slicecell,\slc{d}{X}]^{G}.
\]
By Lemma~\ref{thm:2}, the vanishing of this group implies that
$\slc{d}{X}<d$ and hence must be contractible, since it is also $\ge
d$.
\end{pf}

\begin{cor}
\label{thm:87} Suppose that $d$ is odd.  If $X\to Y\to Z$ is a
cofibration sequence, and the $d$-slices of $X$ and $Z$ are
contractible, then the $d$-slice of $Y$ is contractible.
\end{cor}

\begin{pf}
This is immediate from Proposition~\ref{thm:85} and the long exact
sequence of homotopy classes of maps.
\end{pf}

\begin{rem}
\label{rem:48} Using the slice spectral sequence one can easily show
that a pure spectrum always admits a refinement of homotopy groups.
Thus the results above say that a spectrum $X$ is {\em pure} if and
only if the even homotopy groups admit an equivariant refinement, and
the ``slice homotopy groups'' $\pi^{H}_{m\rho_{H}-1}X$ are all zero
whenever $H\subset G$ is non-trivial.
\end{rem}

\subsection{Further multiplicative properties of the slice filtration}
\label{sec:furth-mult-prop}

In this section we show that the slice filtration has the expected
multiplicative properties for pure spectra.  Our main result is
Proposition~\ref{thm:179} below.  It has the consequence that if $X$
and $Y$ are pure spectra, and $E_{r}^{s,t}(\slot)$ is the slice
spectral sequence, then there is a map of spectral sequences
\[
E_{r}^{s,t}(X)\otimes E_{r}^{s',t'}(Y) \to E_{r}^{s+s', t+t'}(X\wedge Y)
\]
representing the pairing $\pim_{\ast}X\wedge \pim_{\ast}Y\to
\pim_{\ast}(X\wedge Y)$.  In other words, multiplication in the slice
spectral sequence of pure spectra behaves in the expected manner.
We leave the deduction of this property from Proposition~\ref{thm:179}
to the reader.

\begin{prop}
\label{thm:179}
If $X\ge n$ is pure and $Y\ge m$ has cellular slices, then $X\wedge
Y\ge n+m$.
\end{prop}

\begin{pf}
We need to show that $P^{n+m-1}(X\wedge Y)$ is contractible.   By
Lemma~\ref{thm:76} the map
\[
X\wedge Y \to P^{n+m-1}X\wedge P^{n+m-1}Y
\]
is a $P^{n+m-1}$-equivalence, so we may reduce to the case in which
the slice filtrations of $X$ and $Y$ are finite.  That case in turn
reduces to the situation in which
\begin{align*}
X &= H\Zm\wedge \slicecell(m,K)\\
Y &= H\Zm\wedge \slicecell'
\end{align*}
in which $\slicecell'$ is any slice cell.  By induction on $|G|$ the
assertion further reduces to the case in which neither $\slicecell$
nor $\slicecell'$ is induced.  Thus we may assume
\begin{align*}
X &= H\Zm\wedge S^{k\rho_{G}} \\
Y &= H\Zm\wedge S^{\ell\rho_{G}}\quad \text{or}\quad H\Zm\wedge
\Sigma^{-1}S^{\ell\rho_{G}},
\end{align*}
in which case the result follows from Proposition~\ref{thm:6}.
\end{pf}

\section{The complex cobordism spectrum}
\label{sec:case-mu}

From here forward we specialize to the case $G=\ztn$, and for
convenience localize all spectra at the prime $2$.  Write
\[
g=|G|,
\]
and let $\gamma\in G$ be a fixed generator.   

\subsection{\texorpdfstring{The spectrum $\mutn{G}$}{The spectrum $\ssmutn{G}$}}
\label{sec:spectrum-mutng}

We now introduce our equivariant variation on the complex cobordism
spectrum by defining
\[
\mutn{G} = \norm_{\zt}^{G}\mur,
\]
where $\mur$ is the $\zt$-equivariant {\em real bordism} spectrum of
Landweber~\cite{MR0222890} and Fujii~\cite{MR0420597} (and further
studied by Araki~\cite{MR614829} and Hu-Kriz~\cite{MR1808224}).  In
\S\ref{sec:real-bordism} we will give a construction of $\mur$ as a
commutative algebra in $\ugspectra{\zt}$.  The norm is taken along the
unique inclusion $\zt\subset G$.  Since the norm is symmetric
monoidal, and its left derived functor may be computed on the spectra
underlying cofibrant commutative rings (Proposition~\ref{thm:16}), the
spectrum $\mutn{G}$ is an equivariant commutative ring spectrum.  For
$H\subset G$ the unit of the restriction-norm adjunction
(Proposition~\ref{thm:119}) gives a canonical commutative algebra map
\begin{equation}\label{eq:72}
\mutn{H} \to i_{H}^{\ast}\mutn{G}.
\end{equation}

By analogy with the shorthand $i_{0}^{\ast}$ for restriction along the
inclusion of the trivial group, we will employ the shorthand notation
\[
i_{1}^{\ast} = i_{\zt}^{\ast}
\]
for the restriction map $\ugspectra{G}\to\ugspectra{\zt}$ induced by
the unique inclusion $\zt\subset G$.  Restricting, one has a
$\zt$-equivariant smash product decomposition
\begin{equation}\label{eq:73}
i_{1}^{\ast}\mutn{G} = \bigwedge_{j=0}^{g/2-1} \gamma^{j}\mur.
\end{equation}

\subsection{Real bordism, real orientations and formal groups}
\label{sec:real-bordism-real}

We begin by reviewing work of Araki~\cite{MR614829} and
Hu-Kriz~\cite{MR1808224} on real bordism.

\subsubsection{The formal group}
\label{sec:formal-group}

Consider $\cp^{n}$ and $\cp^{\infty}$ as pointed $\zt$-spaces under
the action of complex conjugation, with $\cp^{0}$ as the base point.
The fixed point spaces are $\rp^{n}$ and $\rp^{\infty}$.  There are
homeomorphisms
\begin{equation}\label{eq:28}
\cp^{n}/\cp^{n-1}\equiv S^{n\rho_{2}},
\end{equation}
and in particular an identification $\cp^{1}\equiv S^{\rho_{2}}$.

\begin{defin}[Araki~\cite{MR614829}]
\label{def:4} Let $E$ be a $\zt$-equivariant homotopy commutative ring
spectrum.  A {\em real orientation} of $E$ is a class $\gbar{x}\in
\tilde E_{\zt}^{\rho_{2}}(\cp^{\infty})$ whose restriction to
\[
\tilde E_{\zt}^{\rho_{2}}(\cp^{1})
=\tilde E_{\zt}^{\rho_{2}}(S^{\rho_{2}})
\approx E_{\zt}^{0}(\text{pt})
\]
is the unit.  A {\em real oriented spectrum} is a
$\zt$-equivariant ring spectrum $E$ equipped with a real orientation.
\end{defin}

If $(E,\gbar{x})$ is a real oriented spectrum and $f:E\to E'$ is an equivariant
multiplicative map, then 
\[
f_{\ast}(\gbar{x})\in (E')^{\rho_{2}}(\cp^{\infty})
\]
is a real orientation of $E'$.  We will often not distinguish in
notation between $\gbar{x}$ and $f_{\ast}\gbar{x}$.

\begin{eg}
\label{eg:16}
The zero section $\cp^{\infty}\to MU(1)$ is an equivariant
equivalence, and defines a real orientation
\[
\gbar{x}\in \mur^{\rho_{2}}(\cp^{\infty}),
\]
making $\mur$ into a real oriented spectrum.
\end{eg}

\begin{eg}
\label{eg:18}
From the map
\[
\mur\to i_{1}^{\ast}\mutn{G}
\]
provided by~\eqref{eq:72}, the spectrum $i_{1}^{\ast}\mutn{G}$ gets a
real orientation which we'll also denote
\[
\gbar{x} \in ({\mutn{G}})^{\rho_{2}}(\cp^{\infty}).
\]
\end{eg}

\begin{eg}
\label{eg:19}
If
$(H,\gbar{x}_{H})$ and $(E,\gbar{x}_{E})$ are two real oriented
spectra then $H\wedge E$ has two real orientations given by
\[
\gbar{x}_{H} = \gbar{x}_{H}\otimes 1\text{ and } \gbar{x}_{E} =
1\otimes \gbar{x}_{E}.
\]
\end{eg}

The following result of Araki follows easily from the
homeomorphisms~\eqref{eq:28}.

\begin{thm}[Araki~\cite{MR614829}]
\label{thm:107}
Let $E$ be a real oriented cohomology theory.  There are isomorphisms
\begin{align*}
E^{\star}(\cp^{\infty}) &\approx E^{\star}\LL \gbar{x}\RR \\
E^{\star}(\cp^{\infty}\times \cp^{\infty}) &\approx E^{\star}\LL \gbar{x}\otimes 1, 1\otimes \gbar{x}\RR \\
\end{align*}
\end{thm}

Because of Theorem~\ref{thm:107}, the map
$\cp^{\infty}\times\cp^{\infty}\to\cp^{\infty}$ classifying the tensor
product of the two tautological line bundles defines a formal group
law over $\pi^{G}_{\star}E$.  Using this, much of the theory relating
formal groups, complex cobordism, and complex oriented cohomology
theories works for $\zt$-equivariant spectra, with $\mur$ playing the
role of $MU$.  For information beyond the discussion below,
see~\cite{MR614829,MR1808224}.

\begin{rem}
\label{rem:58} 
A real orientation $\gbar{x}$ corresponds to a {\em
coordinate} on the corresponding formal group.  Because of this we
will use the terms interchangeably, preferring ``coordinate'' when the
discussion predominantly concerns the formal group, and ``real
orientation'' when it concerns spectra.
\end{rem}

The standard formulae from the theory of formal groups give elements
in the $RO(\zt)$-graded homotopy groups $\pi^{\zt}_{\star}E$ of real
oriented $E$.  For example, there is a map from the Lazard ring to
$\pi^{\zt}_{\star}E$ classifying the formal group law.  Using
Quillen's theorem to identify the Lazard ring with the complex
cobordism ring this map can be written as
\[
MU_{\ast}\to \pi^{\zt}_{\star}E.
\]
It sends $MU_{2n}$ to $\pi^{\zt}_{n\rho_{2}}E$.  When $E=\mur$ this
splits the forgetful map
\begin{equation}\label{eq:71}
\pi^{\zt}_{n\rho_{2}}\mur\to \piu_{2n}\mur = \pi_{2n}MU,
\end{equation}
which is therefore surjective.   A similar discussion applies to
iterated smash products of $\mur$ giving

\begin{prop}
\label{thm:171} For every $m>0$, the above construction gives a ring
homomorphism
\begin{equation}\label{eq:86}
\piu_{\ast}\bigwedge^{m}\mur\to \bigoplus_{j}\pi^{\zt}_{j\rho_{2}}\bigwedge^{m}\mur
\end{equation}
splitting the forgetful map
\begin{equation}\label{eq:75}
\bigoplus_{j}\pi^{\zt}_{j\rho_{2}}\bigwedge^{m}\mur \to \piu_{\ast}\bigwedge^{m}\mur.
\end{equation}
In particular, \eqref{eq:75} is a split surjection.  \qed
\end{prop}

It is a result of Hu-Kriz\cite{MR1808224} that~\eqref{eq:75} is in fact
an isomorphism.  This result, and a generalization to $\mutn{G}$ can
be recovered from the slice spectral sequence.

The class
\[
\gbar{x}_{H}\in H^{\rho_{2}}_{\zt}(\cp^{\infty};\zltm)
\]
corresponding to $1\in H^{0}_{\zt}(\text{pt},\zltm)$ under the isomorphism
\[
H^{\rho_{2}}_{\zt}(\cp^{\infty};\zltm) \approx
H^{\rho_{2}}_{\zt}(\cp^{2};\zltm) \approx
H^{0}_{\zt}(\text{pt},\zltm)
\]
defines a real orientation of $H\zltm$.  As in Example~\ref{eg:19}, the
classes $\gbar{x}$ and $\gbar{x}_{H}$ give two orientations of
$E=H\zltm\wedge\mur$.  By Theorem~\ref{thm:107} these are related by a
power series
\begin{align*}
\gbar{x}_{H}
&= \log_{F}(\gbar{x}) \\
&=\gbar{x}+\sum_{i>0}\gbar{\mm}_{i}\gbar{x}^{i+1},
\end{align*}
with
\[
\gbar{\mm}_{i}\in \pi^{\zt}_{i\rho_{2}}H\zltm\wedge \mur.
\]
This power series is the {\em logarithm} of $F$.  Similarly, the
invariant differential on $F$ gives classes $(n+1)\gbar{\mm}_{n}\in
\pi^{\zt}_{n\rho_{2}}\mur$.  The coefficients of the formal sum give
\[
\gbar{a}_{ij}\in \pi^{\zt}_{(i+j-1)\rho_{2}}\mur.
\]

\begin{rem}
\label{rem:12} Since the generator of $\zt$ acts by $(-1)^{n}$ on
\[
H_{2n}i_{0}^{\ast}S^{n\rho_{2}}=\piu_{2n}H\Zm\wedge S^{n\rho_{2}},
\]
it acts also acts by $(-1)^{n}$ on the non-equivariant class $m_{n}$
underlying $\gbar{\mm}_{n}$ and by $(-1)^{n}$ on
$\piu_{2n}\bigwedge^{m}\mur=\pi_{2n}\bigwedge^{m}MU$.
\end{rem}

If $(E, \gbar{x}_{E})$ is a real oriented spectrum then $E\wedge \mur$
has two orientations
\begin{align*}
\gbar{x}_{E} &= \gbar{x}_{E}\otimes 1 \\
\gbar{x}_{R} &= 1\otimes\gbar{x}.
\end{align*}
These two orientations are related by a power series
\begin{equation}\label{eq:74}
\gbar{x}_{R} = \sum \gbar{b}_{i}x_{E}^{i+1}
\end{equation}
defining classes
\[
\gbar{b}_{i}=\gbar{b}^{E}_{i}\in \pi^{\zt}_{i\rho_{2}}E\wedge \mur.
\]
The power series~\eqref{eq:74} is an isomorphism over
$\pi^{\zt}_{\star}E\wedge\mur$
\[
F_{E}\to F_{R}
\]
of the formal group law for $(E, \gbar{x}_{E})$ with the formal group
law for $(\mur,\gbar{x})$.

\begin{thm}[Araki~\cite{MR614829}]\label{thm:28}
The map
\[
E_{\star}[\gbar{b}_{1},\gbar{b}_{2},\dots] \to \pi^{\zt}_{\star}E\wedge \mur
\]
is an isomorphism.
\qed
\end{thm}

Araki's theorem has an evident geometric counterpart.  For each $j$
choose a map
\[
S^{j\rho_{2}}\to E\wedge \mur
\]
representing $\gbar{b}_{j}$.  As in \S\ref{sec:meth-poly-algebr}, let
\[
S[\gbar{b}_{j}]=\bigvee_{k\ge 0}S^{k\cdot j\rho_{2}}
\]
be the free associative algebra on $S^{j\rho_{2}}$ and
\[
S[\gbar{b}_{j}]\to E\wedge \mur
\]
the homotopy associative algebra map extending~\eqref{eq:79}.  Using
the multiplication map, smash these together to form a map of spectra
\begin{equation}\label{eq:82}
E[\gbar{b}_{1},\gbar{b}_{2},\dots]\to E\wedge\mutn{G},
\end{equation}
where
\[
E[\gbar{b}_{1},\gbar{b}_{2},\dots] = E\wedge \hocolim_{k}S[\gbar{b}_{1}]\wedge
S[\gbar{b}_{2}]\wedge\cdots\wedge S[\gbar{b}_{k}].
\]
The map on $RO(\zt)$-graded homotopy groups induced by~\eqref{eq:82}
is the isomorphism of Araki's theorem.   This proves
\begin{cor}
\label{thm:175}
If $E$ is a real oriented spectrum then there is a weak equivalence
\[
E\wedge\mur \approx E[\gbar{b}_{1},\gbar{b}_{2},\dots].
\]
\qed
\end{cor}

\begin{rem}
\label{rem:27} If $E$ is strictly associative then~\eqref{eq:82} is
a map of associative algebras, and the above identifies $E\wedge \mur$
as a twisted monoid ring over $E$.
\end{rem}

As in \S\ref{sec:meth-poly-algebr}, write 
\[
S^{0}[\gbar{b}_{1},\gbar{b}_{2},\dots] = \hocolim_{k}S^{0}[\gbar{b}_{1}]\wedge
S^{0}[\gbar{b}_{2}]\wedge\cdots\wedge S^{0}[\gbar{b}_{k}],
\]
and
\[
S^{0}[G\cdot\gbar{b}_{1},G\cdot\gbar{b}_{1},\dots ] = \norm_{\zt}^{G}
S^{0}[\gbar{b}_{1},\gbar{b}_{2},\dots].
\]
Using Proposition~\ref{thm:77} one can easily check that
$S^{0}[G\cdot\gbar{b}_{1},G\cdot\gbar{b}_{2},\dots]$ is a wedge of
isotropic regular slice cells.  Finally, let
\[
\mutn{G}[G\cdot \gbar{b}_{1},G\cdot\gbar{b}_{2},\dots] =
\mutn{G}\wedge S^{0}[G\cdot\gbar{b}_{1},G\cdot\gbar{b}_{1},\dots] 
\]

\begin{cor}
\label{thm:176} For $H\subset G$ of index $2$, there is an equivalence
of $H$-equivariant associative algebras
\[
i_{H}^{\ast}\mutn{G} \approx \mutn{H}[H\cdot
\gbar{b}_{1},H\cdot\gbar{b}_{2},\dots].
\]
\end{cor}

\begin{pf}
Apply $\norm_{\zt}^{H}$ to the decomposition of Corollary~\ref{thm:175}
with $E=\mur$.
\end{pf}

\subsubsection{The unoriented cobordism ring}
\label{sec:unor-cobord-ring}
Passing to geometric fixed points from
\[
\bar x:\cp^{\infty}\to \Sigma^{\rho_{2}}\mur
\]
gives the canonical inclusion
\[
a:\rp^{\infty} = MO(1)\to \Sigma MO,
\]
defining the $MO$ Euler class of the tautological line bundle.  There
are isomorphisms
\begin{align*}
MO^{\ast}(\rp^{\infty})&\approx MO^{\ast}\LL a\RR \\
MO^{\ast}(\rp^{\infty}\times\rp^{\infty})&\approx MO^{\ast}\LL a\otimes
1, 1\otimes a\RR
\end{align*}
and the multiplication map
$\rp^{\infty}\times\rp^{\infty}\to\rp^{\infty}$ gives a formal group
law over $MO_{\ast}$.  By Quillen~\cite{Qui:FGL}, it is the universal
formal group law $F$ over a ring of characteristic $2$ for which
$F(a,a)=0$.

As described by Quillen~\cite[Page~53]{MR0290382}, the formal group can
be used to give convenient generators for the unoriented cobordism
ring.  Let 
\[
e\in H^{1}(\rp^{\infty};\Z/2)
\]
be the $H\Z/2$ Euler class of the tautological line bundle.  Over
$\pi_{\ast}H\Z/2\wedge MO$ there is a power series relating $e$ and
the image of the class $a$
\[
e = \ell(a) = a+ \sum \alpha_{n} a^{n+1}.
\]

\begin{lem}
\label{thm:178}
The composite series
\begin{equation}\label{eq:50}
(a+\sum \alpha_{2^{j}-1} a^{2^{j}})^{-1}\circ \ell(a) = a +
\sum_{j>0} h_{j}a^{j+1}
\end{equation}
has coefficients in $\pi_{\ast}MO$.  The classes $h_{j}$ with
$j+1=2^{k}$ are zero.  The remaining $h_{j}$ are polynomial generators
for the unoriented cobordism ring
\begin{equation}\label{eq:51}
\pi_{\ast} MO =\Z/2[h_{j}, j\ne 2^{k}-1].
\end{equation}
\end{lem}

\begin{pf}
The assertion that $h_{j}=0$ for $j+1=2^{k}$ is straightforward.
Since the sequence
\begin{equation}
\label{eq:84}
\pi_{\ast}MO\to \pi_{\ast}H\Z/2\wedge MO \rightrightarrows
\pi_{\ast}H\Z/2\wedge H\Z/2\wedge MO
\end{equation}
is a split equalizer, to show that the remaining $h_{j}$ are in
$\pi_{\ast}MO$ it suffices to show that they are equalized by the
parallel maps in~\eqref{eq:84}.  This works out to showing that the
series~\eqref{eq:50} is invariant under substitutions of the form
\begin{equation}\label{eq:85}
e \mapsto e+ \sum \zeta_{m} e^{2^{m}},
\end{equation}
The series~\eqref{eq:50} is characterized as the unique isomorphism of
the formal group law for unoriented cobordism with the additive group,
having the additional property that the coefficients of $a^{2^{k}}$
are zero.  This condition is stable under the
substitutions~\eqref{eq:85}.  The last assertion follows from
Quillen's characterization of $\pi_{\ast}MO$.
\end{pf}

\begin{rem}
\label{rem:21} Recall the real orientation $\gbar{x}$ of
$i_{1}^{\ast}\mutn{G}$ of Example~\ref{eg:18}.  Applying the
$RO(G)$-graded cohomology norm (\S\ref{sec:other-uses-norm}) to
$\gbar{x}$, and then passing to geometric fixed points, gives a class
\[
\phig\norm (\gbar{x}) \in MO^{1}(\rp^{\infty}).
\]
One can easily check that $\phig\norm(\gbar{x})$ coincides with the $MO$
Euler class $a$ defined at the beginning of this section.  Similarly
one has
\[
\phig\norm (\gbar{x}_{H})= e.
\]
Applying $\phig\norm$ to $\log_{\gbar{F}}$ and using the fact that it is
a ring homomorphism (Proposition~\ref{thm:125}) gives
\[
e = a+ \sum \phig\norm(\gbar{m}_{k}) a^{k+1}.
\]
It follows that
\[
\phig\norm(\gbar{m}_{k}) = \alpha_{k}.
\]
\end{rem}

\subsection{Refinement of homotopy groups}
\label{sec:refinement-homotopy}

We begin by focusing on a simple consequence of Proposition~\ref{thm:171}.

\begin{prop}
\label{thm:80} For every $m>1$, every element of
\[
\pi_{2k}\left(\bigwedge^{m}MU\right)
\]
can be refined to an equivariant map
\[
S^{k\rho_{2}}\to  \bigwedge^{m}\mur.
\]
\qed
\end{prop}

This result expresses an important property of the $\zt$-spectra given
by iterated smash products of $\mur$.  Our goal in this section is to
formulate a generalization to the case $G=\ztn$.

\begin{defin}
\label{def:2} Suppose $X$ is a $G$-spectrum.  A {\em refinement of
$\piu_{k}X$} is a map
\[
c:\refine\to X
\]
in which $\refine$ is a wedge of slice cells of dimension $k$, 
inducing an isomorphism
\[
\piu_{k}\refine\to\piu_{k}X.
\]
A {\em refinement of the homotopy groups of $X$} (or a {\em refinement
of homotopy of $X$}) is a map
\[
\refine=\bigvee \refine_{k}\to X
\]
whose restriction to each $\refine_{k}$ is a refinement of $\piu_{k}$.
\end{defin}

\begin{rem}
\label{rem:22}   Let $\sigma_{G}(\Z)$ be the sign representation of
$G$ on $\Z$.    There is an $G$-module isomorphism
$\piu_{|G|}S^{\rho_{G}}\approx\sigma_{G}(\Z)$, and more generally
\[
\piu_{n|H|}\big(G_{+}\underset{H}{\wedge}S^{n\rho_{H}}\big) \approx
\ind_{H}^{G}\sigma_{H}(\Z)^{\otimes n}.
\]
This implies that when $k$ is even, a necessary condition for
$\piu_{k}X$ to admit a refinement is that it be isomorphic as a
$G$-module to a sum
\[
\bigoplus_{H\subset G }M_{H,k}
\]
where $M_{H,k}$ is zero unless $|H|$ divides $k$ and is a sum of
copies of $\ind_{H}^{G}\big(\sigma_{H}(\Z)^{\otimes \ell}\big)$ when
$k=\ell |H|$.  Adding the further condition that for every $H\subset
G$, with $k=\ell |H|$, every element in $\piu_{k}X$ transforming in
$\sigma_{H}(\Z)^{\otimes\ell}$ refines to an element of
$\pi^{H}_{\ell\rho_{H}}X$ makes it sufficient.  A similar analysis
describes the case in which $k$ is odd.
\end{rem}

\begin{rem}
\label{rem:28}
Using Remark~\ref{rem:22} one can check that a refinement of
$\piu_{k}X$ consists of isotropic slice cells if and only if
$\piu_{k}X$ does not contain a free $G$-module as a summand.   
\end{rem}

The splitting~\eqref{eq:86} used to prove Proposition~\ref{thm:80} is
multiplicative.  This too has an important analogue.

\begin{defin}
\label{def:14}
Suppose that $R$ is an equivariant associative algebra.  A {\em
multiplicative refinement of homotopy} is an associative algebra map
$\refine\to R$ which, when regarded as a map of $G$-spectra is a refinement
of homotopy.
\end{defin}

\begin{prop}%
\label{thm:126}%
For every $m\ge 1$ there exists a multiplicative refinement of homotopy
\[
\refine\to \bigwedge^{m}\mutn{G},
\]
with $\refine$ a wedge of regular isotropic slice cells.
\end{prop}

Two ingredients form the proof of Proposition~\ref{thm:126}.  The
first, Lemma~\ref{thm:127} below, is a description of
$\piu_{\ast}\mutn{G}$ as a $G$-module.  The computation is of interest
in its own right, and is used elsewhere in this paper.  It is proved
in \S\ref{sec:homotopy-result}.  The second is the classical
description of $\piu_{\ast}(\bigwedge^{m}\mutn{G})$, $m>1$, as a
$\piu_{\ast}\mutn{G}$-module.

\begin{lem}
\label{thm:127}
There is a sequence of elements $r_{i}\in \piu_{2i}\mutn{G}$ with the
property that
\begin{equation}
\label{eq:41}
\piu_{\ast}\mutn{G} = \Z_{{2}}[G\cdot r_{1}, G\cdot r_{2},\dots],
\end{equation}
in which $G\cdot r_{i}$ stands for the sequence
\[
(r_{i},\dots \gamma^{\frac{g}{2}-1}r_{i})
\]
of length $g/2$.
\end{lem}

We refer to the condition~\eqref{eq:41} by saying that the elements
$r_{i}\in \piu_{2i}\mutn{G}$ form a set of {\em $G$-algebra generators
for $\piu_{\ast}\mutn{G}$.}

\begin{rem}
\label{rem:35} Lemma~\ref{thm:127} completely describes
$\piu_{\ast}\mutn{G}$ as a representation of $G$.  To spell it out,
recall from Remark~\ref{rem:12} that the action of the generator of
$\zt$ on $\piu_{2i}\mutn{G}$ is by $(-1)^{i}$.  The elements
$r_{i}\in\piu_{2i}\mutn{G}$ therefore satisfy
$\gamma^{\frac{g}{2}}r_{i}=(-1)^{i}r_{i}$ and transform in the
representation induced from the sign representation of $\zt$ if $i$ is
odd and in the representation induced from the trivial representation
of $\zt$ if $i$ is even.  Lemma~\ref{thm:127} implies that the map
from the symmetric algebras on the sum of these representations to
$\piu_{\ast}\mutn{G}$ is an isomorphism.
\end{rem}

\begin{rem}
\label{rem:41} The fact that the action of $\zt$ on
$\piu_{2i}\mutn{G}$ is either a sum of sign or trivial representations
means that it cannot contain a summand which is free.  The same is
therefore true of the $G$-action.  By Remark~\ref{rem:28} this implies
that only isotropic slice cells may occur in a refinement of
$\piu_{2i}\mutn{G}$.
\end{rem}

Over $\piu_{\ast}\mutn{G}\wedge \mutn{G}$, there are two formal group
laws, $F_{L}$ and $F_{R}$ coming from the canonical orientations of
the left and right factors.   There is also a canonical isomorphism
between them, which can be written as
\[
x_{R} = \sum b_{j}x_{L}^{j+1}.
\]
Write
\[
G\cdot b_{i}
\]
for the sequence
\[
b_{i}, \gamma b_{i},\dots,\gamma^{g/2-1}b_{i}.
\]

The following result is a standard computation in complex cobordism.

\begin{lem}
\label{thm:173}
The ring $\piu_{\ast}\mutn{G}\wedge \mutn{G}$ is given by
\[
\piu_{\ast}\mutn{G}\wedge \mutn{G} = \piu_{\ast}\mutn{G}[G\cdot
b_{1},G\cdot b_{2},\dots].
\]
For $m>1$,
\[
\piu_{\ast}\bigwedge^{m}\mutn{G}
=
\piu_{\ast}\mutn{G}\wedge \bigwedge^{m-1}\mutn{G}
\]
is the polynomial ring
\[
\piu_{\ast}\mutn{G}[G\cdot  b_{i}^{(j)} ],
\]
with
\begin{align*}
i&=1,2,\dots,\quad\text{and} \\ j&=1,\dots,m-1.
\end{align*}
The element $b_{i}^{(j)}$ is the class $b_{i}$ arising from the
$j^{\text{th}}$ factor of $\mutn{G}$ in $\bigwedge^{m-1}\mutn{G}$.
\end{lem}

\begin{pf}
The second assertion follows from the first and the K\"unneth formula.
If not for the fact that $G$ acts on both factors of
$i_{0}^{\ast}\mutn{G}$, the first assertion would also follow
immediately from the K\"unneth formula and the usual description of
$\mustarmu$.  The quickest way to deduce it from the apparatus we have
describe so far is to let $G\subset G'$ be an embedding of index $2$
into a cyclic group, write
\[
\mutn{G}\wedge \mutn{G} \approx i_{G}^{\ast}\mutn{G'}
\]
and use Corollary~\ref{thm:176}.   
\end{pf}

\begin{rem}
\label{rem:36} 
As with Lemma~\ref{thm:127}, the lemma above actually
determines the structure of $\piu_{\ast}\mutn{G}\wedge \mutn{G}$ as a
$G$-equivariant $\piu_{\ast}\mutn{G}$-algebra.  See
Remark~\ref{rem:35}.
\end{rem}

\begin{pf*}{Proof of Proposition~\ref{thm:126}, assuming
Lemma~\ref{thm:127}} This is a straightforward application of the
method of twisted monoid rings of \S\ref{sec:meth-poly-algebr}.  To
keep the notation simple we begin with the case $m=1$.  Choose a
sequence $r_{i}\in\piu_{2i}\mutn{G}$ with the property described in
Lemma~\ref{thm:127}.  Let
\begin{equation}\label{eq:79}
\rr_{i}: S^{i\rho_{2}}\to i_{1}^{\ast}\mutn{G},
\end{equation}
be a representative of the image of $r_{i}$ under the
splitting~\eqref{eq:86}.   Since $\mutn{G}$ is a commutative algebra,
the method of twisted monoid rings can be used to construct
an associative algebra map
\begin{equation}\label{eq:81}
S^{0}[G\cdot\rr_{1},G\cdot\rr_{2},\dots]\to \mutn{G},
\end{equation}
Using Proposition~\ref{thm:77} one can easily check that
$S^{0}[G\cdot\rr_{1},G\cdot\rr_{2},\dots]$ is a wedge of regular isotropic
$G$-slice cells.  Using Lemma~\ref{thm:127} one then easily checks
that~\eqref{eq:81} is multiplicative refinement of homotopy.  The case
$m\ge 1$ is similar, using in addition Lemma~\ref{thm:173} and the
collection $\{r_{i}, b_{i}(j) \}$.
\end{pf*}

\subsection{\texorpdfstring{Algebra generators for $\piu_{\ast}\mutn{G}$}{Algebra generators for $\piu_{\ast}\ssmutn{G}$}}
\label{sec:homotopy-result}

In this section we will describe convenient algebra generators for
$\piu_{\ast}\mutn{G}$.  Our main results are Proposition~\ref{thm:31}
(giving a criterion for a sequence of elements $r_{i}$ to ``generate''
$\piu_{\ast}\mutn{G}$ as a $G$-algebra, as in Lemma~\ref{thm:127}) and
Corollary~\ref{thm:96} (specifying a particular sequence of $r_{i}$).
Proposition~\ref{thm:31} directly gives Lemma~\ref{thm:127}.

We remind the reader that the notation $\hu_{\ast}X$ refers to the
homology groups $H_{\ast}(i_{0}^{\ast}X)$ of the non-equivariant
spectrum underlying $X$.

\subsubsection{A criterion for a generating set}
\label{sec:criterion}

Let 
\[
\mm_{i}\in H_{2i}MU = \piu_{2i}H\Zm\wedge \mur
\]
be the coefficient of the universal
logarithm.  Using the
identification~\eqref{eq:73}
\[
i_{1}^{\ast}\mutn{G} = \bigwedge_{j=0}^{g/2-1} \gamma^{j}\mur
\]
and the K\"unneth formula, one has
\begin{gather*}
\hu_{\ast}\mutn{G} = \zlt[\gamma^{j}m_{k}],
\\
\intertext{where}
k=1,2,\dots, \\
j = 0,\dots, g/2-1.
\end{gather*}
By the definition of the $\gamma^{j}m_{k}$ and Remark~\ref{rem:12},
the action of $G$ on $\hu_{\ast}\mutn{G}$ is given by
\begin{equation}\label{eq:26}
\gamma\cdot \gamma^{j}\mm_{k}=\begin{cases}
\gamma^{j+1}\mm_{k} &\qquad j < g/2-1 \\
(-1)^{k}\mm_{k} &\qquad j = g/2-1\mathrlap{\ .}
\end{cases}
\end{equation}

Let
\begin{align*}
I &=\ker \piu_{\ast}\mutn{G}\to \zlt \\
I_{H} &=\ker \hu_{\ast}\mutn{G}\to\zlt
\end{align*}
denote the augmentation ideals, and
\begin{align*}
Q_{\ast} &= I/I^{2} \\
QH_{\ast} &= I_{H}/I_{H}^{2}
\end{align*}
the modules of indecomposable, with $Q_{2m}$ and $QH_{2m}$ indicating
the homogeneous parts of degree $2m$ (the odd degree parts are zero).
The module $QH_{\ast}$ is the free abelian group with basis
$\{\gamma^{j}\mm_{k} \}$, and from Milnor~\cite{Mil:MU}, one knows
that the Hurewicz homomorphism gives an isomorphism
\[
Q_{2k}\to QH_{2k}
\]
if $2k$ is not of the form $2(2^{\ell}-1)$, and an exact sequence
\begin{equation}\label{eq:27}
Q_{2(2^{\ell}-1)}\rightarrowtail QH_{2(2^{\ell}-1)} \twoheadrightarrow \Z/2
\end{equation}
in which the rightmost map is the one sending each $\gamma^{j}\mm_{k}$ to $1$.

Formula~\eqref{eq:26} implies that the $G$-module $QH_{2k}$ is the
module induced from the sign representation of $\zt$ if $k$ is odd
and from the trivial representation if $k$ is even.

\begin{lem}
\label{thm:106} Let $r = \sum a_{j}\gamma^{j}\mm_{k}\in QH_{2k}$.  The unique
$G$-module map
\begin{align*}
\Z_{(2)}[G] &\to QH_{2k}  \\
1 &\mapsto r
\end{align*}
factors through a map
\[
\Z_{(2)}[G]/(\gamma^{g/2}-(-1)^{k}) \to QH_{2k}
\]
which is an isomorphism if and only if $\sum a_{j}\equiv 1 \mod 2$.
\end{lem}

\begin{pf}
The factorization is clear, since $\gamma^{g/2}$ acts with eigenvalue
$(-1)^{k}$ on $QH_{2k}$.  Use the unique map $\Z_{(2)}[G]\to QH_{2k}$
sending $1$ to $\mm_{k}$ to identify $QH_{2k}$ with
$A=\Z_{(2)}[G]/(\gamma^{g/2}-(-1)^{k})$.  The main assertion is then
that an element $r=\sum a_{j}\gamma^{j}\in A$ is a unit if and only if
$\sum a_{j}\equiv 1\mod 2$.  Since $A$ is a finitely generated free
module over the Noetherian local ring $\Z_{(2)}$, Nakayama's lemma
implies that the map $A\to A$ given by multiplication by $r$ is an
isomorphism if and only if it is after reduction modulo $2$.  So $r$
is a unit if and only if it is after reduction modulo $2$.  But
$A/(2)=\Z/2[\gamma]/(\gamma^{g/2}-1)$ is a local ring with nilpotent
maximal ideal $(\gamma-1)$.  The residue map
\[
A/(2)\to A/(2,\gamma-1)=\Z/2
\]
sends $\sum a_{j}\gamma^{j}\mm_{k}$ to $\sum a_{j}$.  The result follows.
\end{pf}

\begin{lem}
The $G$-module $Q_{2(2^{\ell}-1)}$ is isomorphic to the module induced
from the sign representation of $\zt$.  For $y\in QH_{2(2^{\ell}-1)}$,
the unique $G$-map
\begin{align*}
\zlt[G] &\to QH_{2(2^{\ell}-1)} \\
1 &\mapsto y
\end{align*}
factors through a map
\[
A=\zlt[G]/(\gamma^{g/2}+1)\to Q_{2(2^{\ell}-1)}
\]
which is an isomorphism if and only if $y=(1-\gamma)r$ where $r\in
QH_{2(2^{\ell}-1)}$ satisfies the condition $\sum a_{j}=1\mod 2$ of
Lemma~\ref{thm:106}.
\end{lem}

\begin{pf}
Identify $QH_{2(2^{\ell}-1)}$ with $A$ by the map sending $1$ to
$\mm_{2^{\ell}-1}$.  In this case $A$ is isomorphic to
$\Z_{(2)}[\zeta]$, with $\zeta$ a primitive $g^{\text{th}}$ root
of unity, and in particular is an integral domain.  Under this
identification, the rightmost map in~\eqref{eq:27} is the quotient of
$A$ by the principal ideal $(\zeta-1)$.  Since $A$ is an integral
domain, this ideal is a rank $1$ free module generated by any element
of the form $(1-\gamma)r$ with $r\in A$ a unit.  The result follows.
\end{pf}

This discussion proves
\begin{prop}
\label{thm:31}
Let
\[
\{r_{1},r_{2},\dots \}\subset\piu_{\ast}\mutn{G}
\]
be any sequence of elements whose images 
\[
s_{k}\in QH_{2k}
\]
have the property that for $k\ne 2^{\ell}-1$, $s_{k}=  \sum a_{j}\gamma^{j}\mm_{k}$
with
\[
\sum a_{j}\equiv 1\mod 2,
\]
and $s_{2^{\ell}-1}= (1-\gamma)\left(\sum a_{j}\gamma^{j}\mm_{2^{\ell}-1} \right)$, with
\[
\sum a_{j}\equiv 1\mod 2.
\]
Then the sequence
\[
\{r_{1},\dots \gamma^{\frac{g}{2}-1}r_{1},r_{2},\dots,\gamma^{\frac{g}{2}-1}r_{2},\dots \}
\]
generates the ideal $I$, and so
\[
\zlt[r_{1},\dots
\gamma^{\frac{g}{2}-1}r_{1},r_{2},\dots,\gamma^{\frac{g}{2}-1}r_{2},\dots]\to
\piu_{\ast}\mutn{G}
\]
is an isomorphism.
\qed
\end{prop}

\subsubsection{Specific generators}
\label{sec:specific-generators}

We now use the action of $G$ on $i_{1}^{\ast}\mutn{G}$ to define
specific elements $\bar{r}_{i}\in\pi^{\zt}_{i\rho_{2}}\mutn{G}$ refining
a sequence satisfying the condition of Proposition~\ref{thm:31}.

Write
\[
\gbar{F}(\gbar{x},\gbar{y})
\]
for the formal group law over $\pi^{\zt}_{\star}\mutn{G}$, and
\[
\log_{\gbar{F}}(\gbar{x})=\gbar{x} + \sum_{i>0}\mbm_{k}\gbar{x}^{k+1}
\]
for its logarithm.  This defines elements
\[
\mbm_{k}\in \pi^{\zt}_{k\rho_{2}} H\zltm\wedge \mutn{G}.
\]
We define the elements
\begin{equation}
\label{eq:107}
\rr_{k}\in \pi^{\zt}_{k\rho_{2}} \mutn{G}
\end{equation}
to be the coefficients of the unique strict isomorphism of $\gbar{F}$
with the $2$-typification of $\gbar{F}^{\gamma}$.  The Hurewicz images
\[
\rr_{k}\in\pi^{\zt}_{k\rho_{2}}H\zltm\wedge \mutn{G}
\]
are given by the power series identity
\begin{equation}\label{eq:48}
\sum \rr_{k}\gbar{x}^{k+1} =
\left(\gbar{x}+ \sum
\gamma(\gbar{m}_{2^{\ell}-1})\gbar{x}^{2^{\ell}}\right)^{-1}
\circ \log_{\gbar{F}}(\gbar{x}).
\end{equation}
Modulo decomposables this becomes
\begin{equation}\label{eq:49}
\rr_{k} \equiv
\begin{cases}
\gbar{m}_{k}-\gamma\gbar{m}_{k} &\qquad k = 2^{\ell}-1 \\
\gbar{m}_{k} &\qquad \text{otherwise.}
\end{cases}
\end{equation}
This shows that the elements $\rr_{k}$ satisfy the condition of
Proposition~\ref{thm:31}, hence

\begin{cor}
\label{thm:96} The classes $r_{k}=i_{0}^{\ast}\rr_{k}$ form a set of
$G$-algebra generators for $\piu_{\ast}\mutn{G}.$  \qed
\end{cor}

These are the specific generators with which we shall work.  Though it
does not appear in the notation, the classes $\rr_{i}$ depend on the
group $G$.  In \S\ref{sec:periodicity-theorem} we will need to
consider the classes $\rr_{i}$ for a group $G$ and for a subgroup
$H\subset G$.  We will then use the notation
\[
\rr^{H}_{i}\text{ and }\rr^{G}_{i}
\]
to distinguish them.

The following result establishes an important property of these
specific $\rr_{k}$.  In the statement below, the symbol $N$ is the
norm map on the $RO(G)$-graded homotopy groups of commutative rings.

\begin{prop}
\label{thm:124} For all $k$
\[
\phig\norm(\rr_{k}) = h_{k}\in \pi_{k}MO,
\]
where the $h_{k}$ are the classes defined in
\S\ref{sec:unor-cobord-ring}.  In particular, the set
\[
\{\phig\norm(\gbar{r_{k}})\mid k\ne 2^{\ell}-1 \}
\]
is a set of polynomial algebra generators of $\pi_{\ast}MO$, and for
all $\ell$
\[
\phig\norm(\rr_{2^{\ell}-1})=\uc_{2^{\ell}-1} = 0.
\]
\end{prop}

\begin{pf}
From Remark~\ref{rem:21} we know that
\begin{align*}
\phig\norm \gbar{x} & = a\\
\phig \norm\gbar{x}_{H} & = e\\
\phig\norm\gbar{m}_{n} &= \alpha_{n}.
\end{align*}
Corollary~\ref{thm:123} implies that
\[
\phig\norm\gamma\gbar{m}_{n} =
\phig\norm\gbar{m}_{n},
\]
so we also know that
\[
\phig\norm\gamma\gbar{m}_{n} = \alpha_{n}.
\]
Since the Hurewicz homomorphism
\[
\xymatrix{
\pi_{\ast}\phig\mutn{G} \ar[d]_{\approx}\ar[r] &  \pi_{\ast}\phig (H\zltm\wedge \mutn{G})\ar[d]^{\approx}\\
\pi_{\ast}MO \ar[r] & \pi_{\ast}H\Z/2[b]\wedge MO
}
\]
is a monomorphism, we can calculate $\phig\norm\rr_{k}$
using~\eqref{eq:48}.  Applying $\phig\norm$ to~\eqref{eq:48}, and
using the fact that it is a ring homomorphism gives
\begin{align*}
a + & \sum (\phig\norm\rr_{k})a^{k+1}  \\
& = \left(a+\sum (\phig\norm \gamma\gbar{m}_{2^{\ell}-1})a^{2^{\ell}}
\right)^{-1}\circ \left(a+\sum
(\phig\norm\gbar{m}_{k})a^{k+1}\right) \\
& = \left(a+\sum \alpha_{2^{\ell}-1}a^{2^{\ell}}
\right)^{-1}\circ \left(a+\sum \alpha_{k}a^{k+1}\right).
\end{align*}
But this is the identity defining the classes $h_{k}$.
\end{pf}

In addition to
\[
\uc_{k}=\phig \norm(\rr_{k}) \in \pi_{k}\phig \mutn{G}=\pi_{k}MO
\]
there are some important classes $f_{k}$ attached to these specific
$\rr_{k}$.

\begin{defin}
\label{def:10}
Set
\[
f_{k}=a_{\bar\rho_{G}}^{k} \norm\rr_{k}  \in \pi_{k}^{G}\mutn{G},
\]
where $\bar{\rho}_{G}=\rho_{G}-1$ is the reduced regular
representation.
\end{defin}

The relationship between these classes is displayed in the following
commutative diagram.
\begin{displaymath}
\xymatrix
@R=8mm
@C=8mm
{
   &S^{k}\ar[dl]_(.5){a_{\bar\rho_{G}}^{k}}\ar[d]^(.5){f_{k}}\ar[dr]^(.3){h_{k}}\\
S^{k\rho_{G}}\ar[r]^(.5){\norm\rr_{k}}
   &\mutn{G}\ar[r]^(.5){}
       &\tepfam\wedge \mutn{G}\\
}
\end{displaymath}

\section{The Slice Theorem and the Reduction Theorem}
\label{sec:slice-theorem-i}

Using the method of twisted monoid rings one can
show the Slice Theorem and the Reduction Theorem to be equivalent.  In
\S\ref{sec:from-reduct-theor} we formally state the Reduction Theorem,
and assuming it, prove the Slice Theorem.  In \S\ref{sec:converse} we
establish a converse, for associative algebras $R$ which are pure and
which admit a multiplicative refinement of homotopy by a polynomial
algebra.  Both assertions are used in the proof of the Reduction
Theorem in \S\ref{sec:reduction-theorem}.

\subsection{From the Reduction Theorem to the Slice Theorem}
\label{sec:from-reduct-theor}

We now state the Slice Theorem, using the language of
\S\ref{sec:even-spectra}.

\begin{thm}[Slice Theorem]
\label{thm:89} The spectrum $\mutn{G}$ is an isotropic pure spectrum.
\end{thm}

For the proof of the slice theorem, let 
\[
A=S^{0}[G\cdot\rr_{1}, G\cdot \rr_{2},\dots] \to \mutn{G} 
\]
be the multiplicative refinement of homotopy constructed in
\S\ref{sec:refinement-homotopy} using the method of
twisted monoid rings, and the specific generators of
\S\ref{sec:specific-generators}.  Let $J$ be the left $G$-set defined
by
\[
J=\coprod_{i} G/\zt.
\]
As described in \S\ref{sec:meth-poly-algebr}, the spectrum $A$ is the indexed
wedge
\[
A=\bigvee_{f\in \nat^{J}} S^{\rho_{f}},
\]
in which $\rho_{f}$ is the unique multiple of the regular
representation of the stabilizer group of $f$ having dimension
\[
\dim f= 2\sum_{j\in J}j\, f(j).
\]
As in Example~\ref{eg:9}, let
\[
M_{d} \subset A
\]
be the monomial ideal consisting of the indexed wedge of the
$S^{\rho_{f}}$ with $\dim f\ge d$.  Then $M_{2d-1}=M_{2d}$, and the
$M_{2d}$ fit into a sequence
\[
\cdots \hookrightarrow M_{2d+2}\hookrightarrow M_{2d}\hookrightarrow
M_{2d-2}\hookrightarrow\cdots.
\]
The quotient 
\[
M_{2d}/M_{2d+2}
\]
is the indexed wedge
\begin{equation}
\label{eq:61}
\refine_{2d} = \bigvee_{\dim f=2d} S^{\rho_{f}}
\end{equation}
on which $A$ is acting through the multiplicative map $A\to S^{0}$
(Examples~\ref{eg:9} and~\ref{eg:6}).  The $G$-spectrum~\eqref{eq:61} is a
wedge of regular isotropic slice cells of dimension $2d$.

Replace $\mutn{G}$ with a cofibrant $A$-module, and form
\[
K_{2d} = \mutn{G}\underset A\wedge M_{2d}.
\]
The $K_{2d}$ fit into a sequence
\[
K_{2d+2}\hookrightarrow K_{2d}\hookrightarrow \cdots.
\]

\begin{lem}
\label{thm:172}  The sequences 
\begin{gather*}
K_{2d+2}\to K_{2d} \to K_{2d}/K_{2d+2} \\
K_{2d}/K_{2d+2}\to \mutn{G}/K_{2d+2}\to \mutn{G}/K_{2d}
\end{gather*}
are weakly equivalent to cofibration sequences.  There is an
equivalence
\begin{equation}\label{eq:31}
K_{2d}/K_{2d+2}\approx R(\infty)\wedge \refine_{2d}
\end{equation}
in which 
\[
R(\infty) = \mutn{G}\underset{A}{\wedge} S^{0}.
\]
\end{lem}

\begin{pf}
Since the map $K_{2d+2}\to K_{2d}$ is the inclusion of a wedge summand
it is an $h$-cofibration of spectra, and the first assertion follows
from Proposition~\ref{thm:116} and Corollary~\ref{thm:102}.  The
second assertion follows from the associativity of the smash product
\[
\mutn{G}\smashove{A} (M_{2d}/M_{2d+1}) \approx
(\mutn{G}\smashove{A}S^{0})\wedge \refine_{2d} \approx R(\infty)\wedge
\refine_{2d}.
\]
This completes the proof.
\end{pf}

The Thom map
\[
\mutn{G}\to H\zltm
\]
factors uniquely through an $\mutn{G}$-module map
\[
R(\infty)\to H\zltm.
\]
The following important result will be proved in
\S\ref{sec:proof-theor-refthm:21}.

\begin{thm}[The Reduction Theorem]
\label{thm:21}
The map
\[
R(\infty)\to H\zltm
\]
is a weak equivalence.
\end{thm}

The case $G=\zt$ of the Reduction Theorem is Proposition~4.9 of
Hu-Kriz\cite{MR1808224}.  Its analogue in motivic homotopy theory
appears in unpublished work of the second author and Morel.

To deduce the Slice Theorem from Theorem~\ref{thm:21} we need two
simple lemmas.

\begin{lem}
\label{thm:92}
The spectrum $K_{2d+2}$ is slice $2d$-positive.  
\end{lem}

\begin{pf}
The class of left $A$-modules $M$ for which $M\smashove{A}M_{2d+2}>2d$
is closed under homotopy colimits and extensions.   It contains every
module of the form $\Sigma^{k}G/H_{+}\wedge A$, with $k\ge 0$.   Since
$A$ is $(-1)$-connected this means it contains every $(-1)$-connected
cofibrant $A$-module.   In particular it contains the cofibrant
replacement of $\mutn{G}$.
\end{pf}

\begin{lem}
\label{thm:93}
If Theorem~\ref{thm:21} holds then $\mutn{G}/K_{2d+2}\le 2d$.
\end{lem}

\begin{pf}
This is easily proved by induction on $d$, using the fact that
\[
R(\infty)\wedge 
\refine_{2d}\to \mutn{G}/K_{2d+2} \to \mutn{G}/K_{2d}.
\]
is weakly equivalent to a cofibration sequence (Lemma~\ref{thm:172}).
\end{pf}

\begin{pf*}{Proof of the Slice Theorem assuming the Reduction Theorem}
It follows from the fibration sequence
\[
K_{2d+2}\to\mutn{G}\to \mutn{G}/K_{2d+2}, 
\]
Lemmas~\ref{thm:92} and~\ref{thm:93} above, and Lemma~\ref{thm:7} that 
\[
P^{2d+1}\mutn{G} \approx P^{2d}\mutn{G} \approx \mutn{G}/K_{2d+2}.
\]
Thus the odd slices of $\mutn{G}$ are contractible and the $2d$-slice
is weakly equivalent to 
\[
R(\infty)\wedge \refine_{2d}\approx
H\zltm\wedge \refine_{2d}.
\]
This completes the proof.
\end{pf*}

\subsection{A converse}
\label{sec:converse}

The arguments of the previous section can be reversed.  Suppose that
$R$ is a $(-1)$-connected associative algebra which we know in advance
to be pure, and that $A\to R$ is a multiplicative refinement of
homotopy, with 
\[
A=S^{0}[G\cdot\gbar{x}_{1},\dots]
\]
a twisted monoid ring having the property that $|\gbar{x}_{i}|>0$ for
all $i$.  Note that this implies that $\piu_{0}R=\Z$ and that
$P_{0}^{0}R=H\Zm$.  Let $M_{d+1}\subset A$ be the monomial ideal
consisting of the slice cells in $A$ of dimension $> d$, write
\[
\tilde P_{d+1}R = M_{d+1}\smashove{A}R 
\]
and
\[
\tilde P^{d}R = R/\tilde P_{d+1}R \approx (A/M_{d+1})\smashove{A}R.
\]
Then the $\tilde P^{d}R$ form a tower. Since $M_{d+1}> d$ and $R\ge 0$
(Proposition~\ref{thm:33}), the spectrum $\tilde P_{d+1}R$ is slice
$d$-positive.  There is therefore a map
\begin{equation}
\label{eq:11}
\tilde P^{d}R\to P^{d}R,
\end{equation}
compatible with variation in $d$.  

\begin{prop}
\label{thm:112}
The map~\eqref{eq:11} is a weak equivalence.   The tower $\{\tilde
P^{d}R \}$ is the slice tower for $R$.
\end{prop}

By analogy with the slice tower, write  $\tilde P_{d'}^{d}R$
for the homotopy fiber of the map
\[
\tilde P^{d}R\to \tilde P^{d'-1}R,
\]
when $d'\le d$.

We start with a lemma concerning the case $d=0$.

\begin{lem}
\label{thm:131}
Let $n\ge 0$.   If the map 
\[
\tilde P^{0}R\to P^{0}R
\]
becomes an equivalence after applying $P^{n}$, then for every
$d\ge 0$ the map
\[
\tilde P_{d}^{d}R \to P_{d}^{d}R
\]
becomes an equivalence after applying $P^{d+n}$.
\end{lem}

\begin{pf}
Write $\refine_{d}= M_{d}/M_{d+1}$.   Then there are equivalences
\[
\tilde P_{d}^{d}R \approx \refine_{d}\smashove{A}R
\approx \refine_{d}\wedge (S^{0}\smashove{A}R) \approx
\refine_{d}\wedge \tilde P_{0}^{0}R.
\]
Since $A\to R$ is a refinement of homotopy and $R$ is pure, the
analogous map 
\[
\refine_{d}\wedge P_{0}^{0}R\to P_{d}^{d}R
\]
is also a weak equivalence.   Now consider the following diagram
\[
\xymatrix{
*+{\refine_{d}\wedge P^{n} (\tilde P_{0}^{0}R)}  \ar[r]^{\sim}\ar[d]  &  *+{\refine_{d}\wedge P^{n} (P_{0}^{0}R)}  \ar[d] \\
*+{P^{d+n} (\refine_{d}\wedge \tilde P_{0}^{0}R)}  \ar[r]\ar[d]_{\sim}        & 
*+{P^{d+n} (\refine_{d}\wedge P_{0}^{0}R)}  \ar[d]^{\sim}\\
*+{P^{d+n} (\tilde P_{d}^{d}R)} \ar[r] & *+{P^{d+n} (P_{d}^{d}R)} 
}
\]
The top map is an equivalence by assumption.  The bottom vertical maps
are the result of applying $P^{d+n}$ to the weak equivalences just
described.     Since $\refine_{d}$ is a wedge of regular slice cells of dimension $d$,
Corollary~\ref{thm:70} implies that the upper vertical maps are weak
equivalences.   It follows that the bottom horizontal map is a weak
equivalence as well.
\end{pf}

\begin{pf*}{Proof of Proposition~\ref{thm:112}}
 We will show by induction on $k$ that
for all $d$, the map 
\[
P^{d+k}(\tilde P^{d}R) \to 
P^{d+k}(P^{d}R)  
\]
is a weak equivalence.  By the strong convergence of the slice tower
(Theorem~\ref{thm:34}) this will give the result.  The induction
starts with $k=0$ since $\tilde P_{d+1}R>d$ and so $R\to \tilde
P^{d}R$ is a $P^{d}$-equivalence.  For the induction step, suppose we
know the result for some $k>0$, and consider
\[
\xymatrix{
P^{d+k}\tilde P_{d}^{d}R  \ar[r]\ar[d]_{\sim}  
&   P^{d+k}(\tilde P^{d}R)  \ar[r]\ar[d]_{\sim}   
 & P^{d+k}(\tilde P^{d-1}R)   \ar[d]\\
P^{d+k} (P_{d}^{d}R)  \ar[r]
&   P^{d+k}( P^{d}R)  \ar[r]
 & P^{d+k}( P^{d-1}R)  
}
\]
The bottom row is a cofibration sequence since it can be identified
with 
\[
P_{d}^{d}R\to P^{d}R\to P^{d-1}R.
\]
The middle vertical map is a weak equivalence by the induction
hypothesis, and the left vertical map is a weak equivalence by the
induction hypothesis and Lemma~\ref{thm:131}.  It follows that the
cofiber of the upper left map is weakly equivalent to
$P^{d+k}(P^{d-1}R)$ and hence is $(d+k+1)$-slice null (in fact $d$
slice null).  The top row is therefore a cofibration sequence by
Corollary~\ref{thm:197}, and so the rightmost vertical map is a weak
equivalence.  This completes the inductive step, and the proof.
\end{pf*}

\section{The Reduction Theorem}
\label{sec:reduction-theorem}

We will prove the Reduction Theorem by induction on $g=|G|$.  The case
in which $G$ is the trivial group follows from Quillen's results.  We
may therefore assume that we are working with a non-trivial group $G$
and that the Reduction Theorem is known for all proper subgroups of
$G$.  In the first subsection below we collect some consequences of
this induction hypothesis.  The proof of the induction step is in
\S\ref{sec:proof-theor-refthm:21}.

\subsection{Consequences of the induction hypothesis}
\label{sec:cons-induct-hypoth}

This next result holds for general $G$.

\begin{lem}
\label{thm:192} 
Suppose that $X$ is pure spectrum and $\refine$ is a
wedge of regular slice cells.  Then $\refine\wedge X$ is pure.  If $X$
is pure and isotropic and $\refine$ is regular isotropic, then
$\refine\wedge X$ is pure and isotropic.
\end{lem}

\begin{pf}
Using Proposition~\ref{thm:26} one reduces to the case in which
$\refine=S^{m\rho_{G}}$.  In that case the claim follows from
Corollary~\ref{thm:70}.
\end{pf}

\begin{prop}
\label{thm:162} Suppose $H\subset G$ has index $2$.  If the Slice
Theorem holds for $H$ then the spectrum $i_{H}^{\ast}\mutn{G}$ is an
isotropic pure spectrum.
\end{prop}

\begin{pf}
This is an easy consequence of Corollary~\ref{thm:176}, which gives an
associative algebra equivalence
\[
i_{H}^{\ast}\mutn{G} \approx \mutn{H}[H\cdot \gbar{b}_{1},H\cdot\gbar{b}_{2},\dots].
\]
This shows that $i_{H}^{\ast}\mutn{G}$ is a wedge of smash products of
even dimensional isotropic slice cells with $\mutn{H}$, and hence (by
Lemma~\ref{thm:192}) an isotropic pure spectrum since $\mutn{H}$
is.
\end{pf}

\begin{prop}
\label{thm:170} Suppose $H\subset G$ has index $2$.  If the Slice
Theorem holds for $H$ then the map 
\[
i_{H}^{\ast}R(\infty)\to i_{H}^{\ast}H\zlt
\]
is an equivalence.
\end{prop}

\begin{pf}
By Proposition~\ref{thm:162} we know that $i_{H}^{\ast}\mutn{G}$ is
pure.  The claim then follows from Proposition~\ref{thm:112}.
\end{pf}

\subsection{Certain auxiliary spectra}
\label{sec:cert-auxil-spectra}

Our proof of the Reduction Theorem will require certain 
auxiliary spectra.    For an integer $k>0$ we define 
\[
R(k) = \mutn{G}/(G\cdot\rr_{1},\dots, G\cdot\rr_{k}) =
\mutn{G}\smashove{A} A'
\]
where 
\begin{align*}
A &=S^{0}[G\cdot\rr_{1},G\cdot\rr_{2},\dots] \\
A' &=S^{0}[G\cdot\rr_{k+1},G\cdot\rr_{k+2},\dots].
\end{align*}
The spectrum $R(k)$ is a right $A'$-module, and as in the case of
$\mutn{G}$ described in \S\ref{sec:slice-theorem-i}, the filtration of
$A'$ by the ``dimension'' monomial ideals leads to a filtration of
$R(k)$ whose associated graded spectrum is
\[
R(\infty)\wedge A'.
\]
Thus the reduction theorem also implies that $R(k)$ is a pure
isotropic spectrum.  By the results of the previous section, the
induction hypothesis implies that $i_{H}^{\ast}R(k)$ is pure and
isotropic.

We know from Proposition~\ref{thm:81} that when $m$ is even, the slice
$\slc{m}{R(k)}$ is given by
\[
\slc{m}{R(k)} \approx H\zlt\wedge \refine_{m} 
\]
where $\refine\subset A'$ is the summand consisting of the wedge of
slice cells of dimension $m$.  When $m$ is odd the above discussion
implies that $T\wedge \slc{m}{R(k)}$ is contractible for any $G$-CW
spectrum $T$ built entirely from induced $G$-cells.  In particular, the
equivariant homotopy groups of $\epfam\wedge R(k)$ may be investigated
by smashing the slice tower of $R(k)$ with $\epfam$, and we will do so
in \S\ref{sec:proof-theor-refthm:21}, where we will exploit some very
elementary aspects of the situation.

\subsection{Proof of the Reduction Theorem}
\label{sec:proof-theor-refthm:21}

As mentioned at the beginning of the section, our proof of the
Reduction Theorem is by induction on $|G|$, the case of the trivial
group being a result of Quillen.  We may therefore assume that $G$ is
non-trivial, and that the result is known for all proper subgroups
$H\subset G$.  By Proposition~\ref{thm:170} this implies that the map
\[
R(\infty)\to H\zltm
\]
becomes a weak equivalence after applying $i_{H}^{\ast}$.   

For the induction step we smash the map in question with the isotropy
separation sequence \eqref{eq:3}
\[
\xymatrix@C=2ex{
\epfam\wedge R(\infty) \ar[r]\ar[d]_{f} & R(\infty) \ar[r]\ar[d]^{g} &
{\tepfam}\wedge R(\infty) \ar[d]^{h}\\
\epfam\wedge H\zltm  \ar[r]& H\zltm\ar[r] & {\tepfam}\wedge H\zltm\mathrlap{\ .} \\
}
\]
By the induction hypothesis, the map $f$ is an equivalence.  It
therefore suffices to show that the map $h$ is, and that, as
discussed in Remark~\ref{rem:8}, is equivalent to showing that
\begin{equation}\label{eq:9}
\pi^{G}_{\ast}h:\pi_{\ast}\phig R(\infty)\to\pi_{\ast}\phig H\zltm
\end{equation}
is an isomorphism.

We first show that the two groups are abstractly isomorphic.

\begin{prop}
\label{thm:146}
The ring $\pi_{\ast}\phig H\zltm$ is given by
\[
\pi_{\ast}\phig H\zltm = \Z/2[b],
\]
with
\[
b=u_{2\sigma}a_{\sigma}^{-2}\in\pi_{2}\phig H\zltm\subset
a_{\sigma}^{-1} \pi^{G}_{\star} H\zltm.
\]
The groups $\pi_{\ast}\phig R(\infty)$ are given by
\[
\pi_{n}\phig R(\infty)=
\begin{cases}
\Z/2 &\quad n\ge 0 \text{ even} \\
0 &\quad \text{otherwise.}
\end{cases}
\]

\end{prop}

\begin{pf}
The first assertion is a restatement of Proposition~\ref{thm:94}.  For
the second we will make use of the monoidal geometric fixed point
functor $\phigm$.  The main technical issue is to take care that at
key points in the argument we are working with spectra $X$ for which
$\phig X$ and $\phigm X$ are weakly equivalent.

Recall the definition
\[
R(\infty)=\mutn{G}_{c}\smashove{A}S^{0},
\]
where for emphasis we've written $\mutn{G}_{c}$ as a reminder that
$\mutn{G}$ has been replaced by a cofibrant $A$-module (see
\S\ref{sec:meth-poly-algebr}).  Proposition~\ref{thm:191} implies that
$R(\infty)$ is cofibrant, so there is an isomorphism
\[
\pi_{\ast}\phig R(\infty)\approx \pi_{\ast}\phigm R(\infty)
\]
(Proposition~\ref{thm:14}).  For the monoidal geometric fixed point
functor, Proposition~\ref{thm:191} gives an isomorphism
\[
\phigm(R(\infty))=\phigm(\mutn{G}_{c}\smashove{A}S^{0})\approx \phigm\mutn{G}_{c}
\smashove{\phigm A} S^{0}.
\]

We next claim that there are associative algebra {\em isomorphisms}
\[
\phigm A\approx
S^{0}[\phig \norm \gbar{r}_{1},\phig \norm\gbar{r}_{2},\dots]\approx
S^{0}[\phizt \gbar r_{1},\phizt \gbar r_{2},\dots].  
\]
For the first, decompose $A$ into an indexed wedge, and use
Proposition~\ref{thm:73}.  For the second use the fact that the
monoidal geometric fixed point functor distributes over wedges, and
for $V$ and $W$ representations of $\zt$, can be computed in terms of
the isomorphisms
\[
\phigm(\norm_{\zt}^{G}(S^{-W}\wedge S^{V})) \approx \phigm
(S^{-\ind_{\zt}^{G}W}\wedge S^{\ind_{\zt}^{G}V}) \approx \phiztm
(S^{-W}\wedge S^{V}).
\]

By Proposition~\ref{thm:15}, $\phigm\mutn{G}_{c}$ is a cofibrant
$\phigm A$-module, and so
\[
\phigm\mutn{G}_{c}
\smashove{\phigm A} S^{0}\approx 
\phigm\mutn{G}/(\phigm\norm \gbar r_{1},\phigm\norm \gbar r_{2},\dots).
\]
Since $\mutn{G}_{c}$ is a cofibrant $A$-module, and the polynomial
algebra $A$ has the property that $S^{-1}\wedge A$ is cofibrant, the
spectrum underlying $\mutn{G}_{c}$ is cofibrant
(Corollary~\ref{thm:57}).  This means that
\[
\phigm \mutn{G}_{c} 
\]
and
\[
\phig\mutn{G}_{c}\sim \phig\mutn{G}\sim MO
\]
are related by a functorial zig-zag of weak equivalences
(Proposition~\ref{thm:14}).  Putting all of this together, we arrive
at the equivalence
\[
\phig R(\infty) \sim MO/(\phizt \gbar r_{1},\phizt \gbar r_{2},\dots).
\]
By Proposition~\ref{thm:124}
\[
\phig \rr_{i} = \begin{cases}
\uc_{i} & i\ne 2^{k}-1 \\
0 & i=2^{k}-1.
\end{cases}
\]
From this is an easy matter to compute $\pi_{\ast}MO/(\phig \gbar
r_{1},\phig \gbar r_{2},\dots)$ using the cofibration sequences
described at the end of \S\ref{sec:quotient-rings}.  The outcome is as
asserted.
\end{pf}

Before going further we record a simple consequence of the above
discussion which will be used in \S\ref{sec:rg-graded-slice}.

\begin{prop}
\label{thm:72}
The map
\[
\pi_{\ast}\phig\mutn{G} =\pi_{\ast}MO\to
\pi_{\ast}\phig H\zltm
\]
is zero for $\ast>0$.  
\qed
\end{prop}

A simple multiplicative property reduces the problem of showing
that~\eqref{eq:9} is an isomorphism to showing that it is surjective
in dimensions which are a power of $2$.

\begin{lem}
\label{thm:181}
If for every $k\ge 1$, the class $b^{2^{k-1}}$ is in the image of
\begin{equation}\label{eq:90}
\pi_{2^{k}}\phig \mutn{G}/(G\cdot\rr_{2^{k}-1}) \to \pi_{2^{k}}\phig H\zltm,
\end{equation}
then~\eqref{eq:9} is surjective, hence an isomorphism.
\end{lem}

\begin{pf}
By writing
\[
R(\infty) = \mutn{G}/(G\cdot\rr_{1})\smashove{\mutn{G}}
\mutn{G}/(G\cdot\rr_{2})\smashove{\mutn{G}}\cdots
\]
we see that if for every $k\ge 1$, $b^{2^{k-1}}$ is in the image
of~\eqref{eq:90}, then all products of the $b^{2^{k-1}}$ are in the
image of
\begin{equation}\label{eq:91}
\pi_{\ast}\phig R(\infty)\to \pi_{\ast}\phig H\zltm.
\end{equation}
Hence every power of $b$ is in the image of~\eqref{eq:91}.
\end{pf}

In view of Lemma~\ref{thm:181}, the Reduction Theorem follows from
\begin{prop}
\label{thm:182}
For every $k\ge 1$, the class $b^{2^{k-1}}$ is in the image of
\[
\pi_{2^{k}}\phig(\mutn{G}/(G\cdot\rr_{2^{k}-1})) \to \pi_{2^{k}}\phig(H\zltm).
\]
\end{prop}

To simplify some of the notation, write
\[
\ck = 2^{k}-1
\]
and
\[
M_{k}=\mutn{G}/(G\cdot \rr_{\ck}).
\]
Since $S^{\ck\rho_{G}}$ is obtained from $S^{\ck}$ by attaching
induced $G$-cells, the restriction map
\[
\pi^{G}_{\ck\rho_{G}+1}\tepfam\wedge M_{k}\to 
\pi^{G}_{\ck+1}\tepfam\wedge M_{k}
\]
is an isomorphism (Remark~\ref{rem:20}).  The element of interest in
this group (the one hitting $b^{2^{k-1}}$) arises from the class
\[
\norm\rr_{\ck}\in \pi^{G}_{\ck\rho_{G}}\mutn{G}
\]
and the fact that it is zero for two reasons in
$\pi^{G}_{\ck\rho_{G}}\tepfam \wedge M_{k}$ (it has been coned off in the
formation of $M_{k}$, and it is zero in
$\pi^{G}_{\ck\rho_{G}}\tepfam\wedge \mutn{G}=\pi_{\ck}MO$ by
Proposition~\ref{thm:124}).  We make this more precise and prove
Proposition~\ref{thm:182} by chasing the class $\norm\rr_{\ck}$ around
the sequences of equivariant homotopy groups arising from the diagram
\begin{equation}
\label{eq:117}
\xymatrix@C=2ex{
\epfam\wedge \mutn{G} \ar[r]\ar[d] & \mutn{G} \ar[r]\ar[d] &
{\tepfam}\wedge \mutn{G}\ar[d] \\
\epfam\wedge M_{k} \ar[r]\ar[d] & M_{k} \ar[r]\ar[d] &
{\tepfam}\wedge M_{k} \ar[d] \\
\epfam\wedge H\zltm  \ar[r]& H\zltm\ar[r] & {\tepfam}\wedge H\zltm\mathrlap{\ .}
}
\end{equation}

We start with the top row.  By Proposition~\ref{thm:124} the image of
$\norm\rr_{\ck}$ in 
\[
\pi^{G}_{\ck\rho_{G}}\tepfam\wedge \mutn{G}\approx
\pi^{G}_{\ck}\tepfam\wedge \mutn{G}\approx
\pi_{\ck}MO
\]
is zero.  
There is therefore a class
\[
y_{k}\in \pi^{G}_{\ckp\rho_{G}} \epfam\wedge \mutn{G}
\]
lifting $\norm\rr_{\ck}$.    The key computation, from which everything
follows is

\begin{prop}
\label{thm:183}
The image under 
\[
\pi^{G}_{\ckp\rho_{G}}\epfam\wedge\mutn{G}  \to
\pi^{G}_{\ckp\rho_{G}}\epfam\wedge H\zltm,
\]
of any choice of $y_{k}$ above, is non-zero.
\end{prop}

\begin{pf*}{Proof of Proposition~\ref{thm:182} assuming Proposition~\ref{thm:183}}
We continue chasing around the diagram~\eqref{eq:117}.  By
construction the image of $y_{k}$ in $\pi^{G}_{\ckp\rho_{G}}\epfam\wedge
M_{k}$ maps to zero in $\pi^{G}_{\ckp\rho_{G}} M_{k}$.  It therefore
comes from a class
\[
\tilde y_{k}\in \pi^{G}_{\ckp\rho_{G}+1} \tepfam\wedge M_{k}.
\]
The image of $\tilde y_{k}$ in $\pi^{G}_{\ckp\rho_{G}+1} \tepfam\wedge
H\zltm$ is non-zero since it has a non-zero image in
\[
\pi^{G}_{\ckp\rho_{G}} \epfam\wedge H\zltm
\]
by Proposition~\ref{thm:183}.  Now consider the commutative square
below, in which the horizontal maps are the isomorphisms
(Remark~\ref{rem:20}) given by restriction along the fixed point inclusion
$S^{2^{k}}\subset S^{\ck\rho_{G}+1}$:
\[
\xymatrix{
{\pi^{G}_{\ckp\rho_{G}+1} \tepfam\wedge M_{k}}\ar[r]^-{\approx} \ar[d] &
 {\pi^{G}_{2^{k}} \tepfam\wedge M_{k}} \ar[d] \\
{\pi^{G}_{\ckp\rho_{G}+1} \tepfam\wedge H\zltm}
\ar[r]_-{\approx} & {\pi^{G}_{2^{k}}
\tepfam\wedge H\zltm\mathrlap{\ .}
}
}
\]
The group on the bottom right is cyclic of order $2$, generated by
$b^{2^{k-1}}$.  We've just shown that the image of $\tilde y_{k}$
under the left vertical map is non-zero.  It follows that the right
vertical map is non-zero and hence that $b^{2^{k-1}}$ is in its image.
\end{pf*}

The remainder of this section is devoted to the proof of
Proposition~\ref{thm:183}.

The advantage of Proposition~\ref{thm:183} is that it entirely
involves $G$-spectra which have been smashed with $\epfam$, and which
(as discussed in \S\ref{sec:cert-auxil-spectra}) therefore fall under
the jurisdiction of the induction hypothesis.  In particular, the map
\begin{equation}\label{eq:60}
\epfam\wedge \mutn{G}\to \epfam\wedge H\zltm
\end{equation}
can be studied by smashing the slice tower of $\mutn{G}$ with $\epfam$.

We can cut down some the size of things by making use of the spectra
introduced in \S\ref{sec:cert-auxil-spectra}.  Factor~\eqref{eq:60}
as
\[
\epfam\wedge \mutn{G}\to \epfam\wedge R(\ck-1)\to \epfam\wedge H\zltm,
\]
and replace $y_{k}$ with its image 
\[
y_{k}\in \pi^{G}_{\ckp\rho_{G}}\epfam\wedge R(\ck-1).
\]

\begin{lem}
\label{thm:148}
For $0<m<\ckp g$, 
\[
\pi_{\ckp\rho_{G}}\epfam\wedge \slc{m}{R(\ck-1)} = 0.
\]
There is an exact sequence
\[
\xymatrix{
\pi^{G}_{\ckp\rho_{G}}\epfam\wedge P_{\ckp g}R(\ck-1) \ar[r] & 
\pi^{G}_{\ckp\rho_{G}}\epfam\wedge R(\ck-1) \ar[d] \\
&\pi^{G}_{\ckp\rho_{G}}\epfam\wedge H\zltm=\Z/2\mathrlap{\ .}
}
\]
\end{lem}

\begin{pf}
Because of the induction hypothesis, we know that the spectrum 
\[
\epfam\wedge \slc{m}{R(\ck-1)}
\]
is contractible when $m$ is odd, and that when $m$ is even it is
equivalent to
\[
\epfam\wedge H\Z\wedge \refine_{m},
\]
where $\refine\subset S^{0}[G\cdot\rr_{\ck},\dots]$ is the summand
consisting of the wedge of slice cells of dimension $m$.  Since
$1<m<\ck g$ all of these cells are induced.  This implies that the map
\[
\epfam\wedge H\Z\wedge \refine_{m}\to H\Z\wedge \refine_{m}
\]
is an equivalence, since 
\[
\epfam\to S^{0}
\]
is an equivalence after restricting to any proper subgroup of $G$.
But 
\[
\pi^{G}_{\ck\rho_{G}} H\Z\wedge \refine_{m} = \pi^{G}_{0}H\Z\wedge
S^{-\ck\rho_{G}}\wedge\refine_{m}=0
\]
since 
\[
H\Z\wedge S^{-\ck\rho_{G}}\wedge\refine_{m}
\]
is an $(m-\ck g)$-slice and $m-\ck g<0$.  This proves the first
assertion.  It implies that the map
\[
\pi^{G}_{\ckp\rho_{G}}\epfam\wedge P_{\ckp g}R(\ck-1) \to
\pi^{G}_{\ckp\rho_{G}}\epfam\wedge P_{1}R(\ck-1) 
\]
is surjective.  As mentioned in \S\ref{sec:cert-auxil-spectra},
Proposition~\ref{thm:81} implies that $\slc{0}{R(\ck-1)}=H\zlt$, and
so the second assertion follows from the exact sequence of the
fibration
\[
\epfam\wedge P_{1}R(\ck-1) \to \epfam\wedge R(\ck-1)\to \epfam\wedge \slc{0}{R(\ck-1)}.
\]
\end{pf}

The exact sequence in Lemma~\ref{thm:148} converts the problem of
showing that $y_{k}$ has non-zero image in
$\pi^{G}_{\ckp\rho_{G}}\epfam\wedge H\zltm$ to showing that it is not in
the image of
\[
\pi^{G}_{\ckp\rho_{G}}\epfam\wedge P_{\ckp g}R(\ck-1).
\]
We now isolate a property of this image 
that is not shared by $y_{k}$.  Recall that $\gamma$ is a fixed
generator of $G$.
\begin{prop}
\label{thm:165}
The image of
\[
\pi^{G}_{\ckp\rho_{G}}\epfam\wedge P_{\ckp g} R(\ck-1) \to
\pi^{G}_{\ckp\rho_{G}}
R(\ck-1)
\xrightarrow{i_{0}^{\ast}} \piu_{\ckp g}
R(\ck-1)
\]
is contained in the image of $(1-\gamma)$.
\end{prop}

The class $y_{k}$ does not have the property described in
Proposition~\ref{thm:165}.  Its image in $\piu_{\ckp g}R(\ck-1)$ is
$i_{0}^{\ast}\norm\rr_{\ck}$ which generates a sign representation of
$G$ occurring as a summand of $\piu_{\ckp g}R(\ck-1)$.  Thus once
Proposition~\ref{thm:165} is proved the proof of the Reduction Theorem
is complete.

The proof of Proposition~\ref{thm:165} makes use of the $RO(G)$-graded
Mackey functor 
\[
\pim_{\ckp\rho_{G}}(X)
\]
and the transfer map 
\begin{equation}
\label{eq:97}
\pim_{\ckp\rho_{G}}(X)(\zt)\to 
\pim_{\ckp\rho_{G}}(X)(\text{pt}),
\end{equation}
in which $\zt$ is regarded as a finite $G$-set through the unique
surjective map $G\to\zt$.  By definition
(\S\ref{sec:mackey-functors-1}) of the covariant part
${\pim_{\ck\rho_{G}}}_{\ast}$ of the Mackey functor, the
map~\eqref{eq:97} is given by the map of equivariant homotopy groups
\[
\pi^{G}_{\ckp\rho_{G}}(X\wedge {\zt}_{+}) \to 
\pi^{G}_{\ckp\rho_{G}} (X)
\]
induced by the unique surjective map $\zt\to\text{pt}$.

There are two steps in the proof of Proposition~\ref{thm:165}.  First
it is shown (Corollary~\ref{thm:164}) that the image of
\[
\pi^{G}_{\ckp\rho_{G}}\epfam\wedge P_{\ckp g}R(\ck-1)
\to \pi^{G}_{\ckp\rho_{G}}
R(\ck-1)
\]
is contained in the image of the transfer map just described.
We then show (Lemma~\ref{thm:163}) that the image of
the transfer map in $\piu_{\ckp g}R(\ck-1)$ is in the image of
$(1-\gamma)$.  

\begin{lem}
\label{thm:44} Let $M\ge 0$ be a $G$-spectrum, and regard $\zt$ as a
finite $G$-set using the unique surjective map $G\to\zt$.  The image of
\[
\pi^{G}_{0}\epfam\wedge
M\to\pi^{G}_{0}M
\]
is the image of the transfer map
\[
\pi^{G}_{0}M\wedge{\zt}_{+} \to
\pi^{G}_{0}M.
\]
\end{lem}

\begin{pf}
As mentioned in Remark~\ref{rem:42}, the space $\epfam$ can be taken
to be the space $S^{\infty}_{+}$ on which $\gamma$ acts through the
antipodal action.  The standard cell decomposition in this model has
$0$-skeleton ${\zt}_{+}$.  Since $M$ is $(-1)$-connected
(Proposition~\ref{thm:65}) this implies that
$\pi^{G}_{0}{\zt}_{+}\wedge M\to \pi^{G}_{0}\epfam\wedge M$ is
surjective, and the claim follows.
\end{pf}

\begin{cor}
\label{thm:43} The image of
\[
\pi^{G}_{\ckp\rho_{G}}\epfam\wedge P_{\ckp g}R(\ck-1) \to
\pi^{G}_{\ckp\rho_{G}}P_{\ckp g}R(\ck-1)
\]
is contained in the image of the transfer map.
\end{cor}

\begin{pf}
This follows from Lemma~\ref{thm:44} above, after the identification
\[
\pi^{G}_{\ckp\rho_{G}} P_{\ckp g}R(\ck-1) \approx
\pi^{G}_{0}S^{-\ckp\rho_{G}} \wedge P_{\ckp g}R(\ck-1)
\]
and the observation that
\[
S^{-\ckp\rho_{G}} \wedge P_{\ckp g}R(\ck-1) \approx
P_{0}(S^{-\ckp\rho_{G}} \wedge R(\ck-1))
\]
is $\ge 0$.
\end{pf}

\begin{cor}
\label{thm:164}
The image of
\[
\pi^{G}_{\ckp\rho_{G}}\epfam\wedge P_{\ckp g}R(\ck-1)
\to
\pi^{G}_{\ckp\rho_{G}}R(\ck-1)
\]
is contained in the image of the transfer map.
\end{cor}

\begin{pf}
Immediate from Corollary~\ref{thm:43} and the naturality of the transfer.
\end{pf}

The remaining step is the special case $X=P_{\ckp g}R(\ck-1)$,
$V=\ckp\rho_{G}$ of the next result.

\begin{lem}
\label{thm:163} Let $X$ be a $G$-spectrum, $V$ a virtual
representation of $G$ of virtual dimension $d$, and regard $\zt$ as a
finite $G$-set through the unique surjective map $G\to\zt$.  Write
$\epsilon \in\{\pm{1}\}$ for the degree of
\[
\gamma: i_{0}^{\ast}S^{V}\to i_{0}^{\ast}S^{V}.
\]
The image of
\[
\pi^{G}_{V}(X\wedge{\zt}_{+})\to\pi^{G}_{V}X \to
\piu_{d}X
\]
is contained in the image of
\[
(1+\epsilon\gamma):\piu_{d}X\to \piu_{d}X.
\]
\end{lem}

\begin{pf}
Consider the diagram
\[
\xymatrix{
\pi^{G}_{V} (X\wedge {\zt}_{+}) \ar[r]\ar[d] &
\pi^{G}_{V}X \ar[d] \\
{\piu_{d} (X\wedge {\zt}_{+})} \ar[r] &
\piu_{d}X\mathrlap{\ .}
}
\]
The non-equivariant identification
\[
{\zt}_{+}\approx S^{0}\vee S^{0}
\]
gives an isomorphism of groups of non-equivariant stable maps
\[
[S^{V}, X\wedge {\zt}_{+}] \approx
[S^{V},X] \oplus [S^{V},X],
\]
and so an isomorphism of the group in the lower left hand corner
with
\[
\piu_{d}X\oplus \piu_{d}X
\]
under which the generator $\gamma\in G$ acts as
\[
(a,b)\mapsto (\epsilon\gamma b, \epsilon\gamma a).
\]
The map along the bottom is $(a,b)\mapsto a+b$.  Now the image of the
left vertical map is contained in the set of elements invariant under
$\gamma$ which, in turn, is contained in the set of elements of the
form
\[
(a, \epsilon\gamma a).
\]
The result follows.
\end{pf}

\begin{pf*}{Proof of Proposition~\ref{thm:165}}
As described after its statement, Proposition~\ref{thm:165} is a
consequence of Corollary~\ref{thm:164} and Lemma~\ref{thm:163}.
\end{pf*}

\section{The Gap Theorem}
\label{sec:gap-theorem}

The proof of the Gap Theorem was sketched in the introduction, and
various supporting details were scattered throughout the paper.  We
collect the story here for convenient reference.

Given the Slice Theorem, the Gap Theorem is a consequence of the
following special case of Proposition~\ref{thm:108}

\begin{prop}
\label{thm:152}
Suppose that $G=\ztn$ is a non-trivial group, and $m\ge 0$.  Then
\[
H^{i}_{G}(S^{m\rho_{G}};\zltm) = 0 \qquad \text{ for }0< i<4.
\]
\qed
\end{prop}

\begin{lem}[The Cell Lemma] \label{thm:151} Let $G=\ztn$ for some
$n>0$.  If $\slicecell$ is an isotropic slice cell of even dimension,
then the groups $\pi^{G}_{k}H\zltm\wedge \slicecell$ are zero for
$-4<k<0$.
\end{lem}

\begin{pf}
Suppose that
\[
\slicecell=G_{+}\smashove{H}S^{m\rho_{H}}
\]
with $H\subset G$ non-trivial.  By the Wirthm\"uller isomorphism
\[
\pi^{G}_{k} H\zltm\wedge \slicecell \approx
\pi^{H}_{k}H\zltm\wedge S^{m\rho_{H}},
\]
so the assertion is reduced to the case $\slicecell=S^{m\rho_{G}}$ with $G$
non-trivial.  If $m\ge 0$ then $\pi^{G}_{k}H\zltm\wedge \slicecell=0$ for $k<0$.
For the case $m<0$ write $i=-k$,  $m'=-m>0$, and
\[
\pi^{G}_{k}H\zltm\wedge \slicecell = H^{i}_{G}(S^{m'\rho_{G}};\zltm).
\]
The result then follows from Proposition~\ref{thm:152}.
\end{pf}
 
\begin{thm}
\label{thm:153}
If $X$ is pure and isotropic, then
\[
\pi^{G}_{i}X=0\qquad -4<i<0.
\]
\end{thm}

\begin{pf}
Immediate from the slice spectral sequence for $X$ and the Cell Lemma.
\end{pf}

\begin{cor}
\label{thm:154} If $Y$ can be written as a directed homotopy colimit
of isotropic pure spectra, then
\[
\pi^{G}_{i}X=0\qquad -4<i<0.
\]
\qed
\end{cor}

\begin{thm}[The Gap Theorem]
\label{thm:155}
Let $G=\ztn$ with $n>0$ and $\D\in \pi_{\ell\rho_{G}}\mutn{G}$ be
any class.  Then for $-4<i<0$
\[
\pi^{G}_{i}\D^{-1}\mutn{G} = 0.
\]
\end{thm}

\begin{pf}
The spectrum $\D^{-1}\mutn{G}$ is the homotopy colimit
\[
\hocolim_{j} \Sigma^{-j\,\ell \rho_{G}} \mutn{G}.
\]
By the Slice Theorem, $\mutn{G}$ is pure and isotropic.  But then
the spectrum
\[
\Sigma^{-j\,\ell \rho_{G}} \mutn{G}
\]
is also pure and isotropic, since for any $X$
\[
\slc{m}{\Sigma^{\rho_{G}}X} \approx \Sigma^{\rho_{G}}\slc{m-g}{X}
\]
by Corollary~\ref{thm:70}.  The result then follows from
Corollary~\ref{thm:154}.
\end{pf}

\section{The Periodicity Theorem}
\label{sec:periodicity-theorem}

In this section we will describe a general method for producing
periodicity results for spectra obtained from $\mutn{G}$ by inverting
suitable elements of $\pi_{\star}^{G}\mutn{G}$.  The Periodicity
Theorem (Theorem~\ref{thm:51}) used in the proof of
Theorem~\ref{thm:20} is a special case.  The proof relies on a small
amount of computation of $\pi^{G}_{\star}\mutn{G}$.

\subsection{\texorpdfstring{The $RO(G)$-graded slice spectral sequence for $\mutn{G}$}{The $RO(G)$-graded slice spectral sequence for $\ssmutn{G}$}}
\label{sec:rg-graded-slice}

\begin{figure}
\includegraphics[]{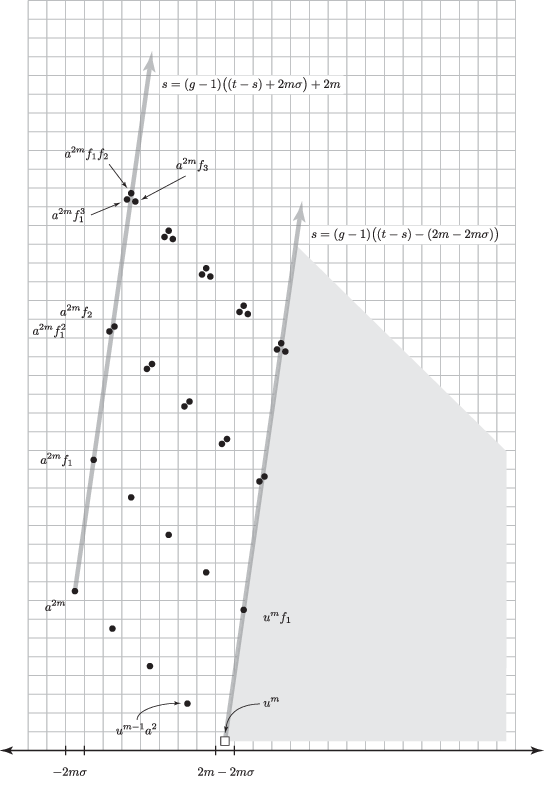}
\caption{The slice spectral sequence for $\pi^{G}_{-2m\sigma+\ast}\mutn{G}$}
\label{fig:4}
\end{figure}

\begin{figure}
\includegraphics[]{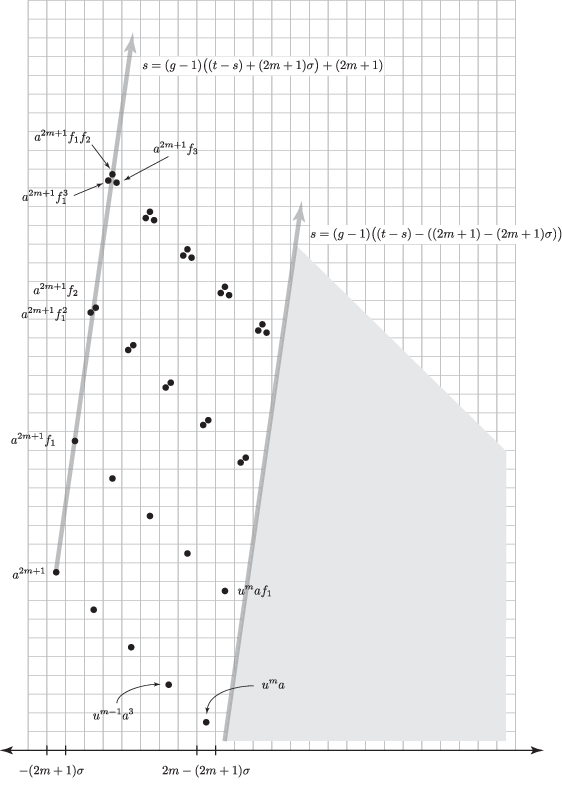}
\caption{The slice spectral sequence for $\pi^{G}_{-(2m+1)\sigma+\ast}\mutn{G}$}
\label{fig:5}
\end{figure}

\begin{figure}
\includegraphics[]{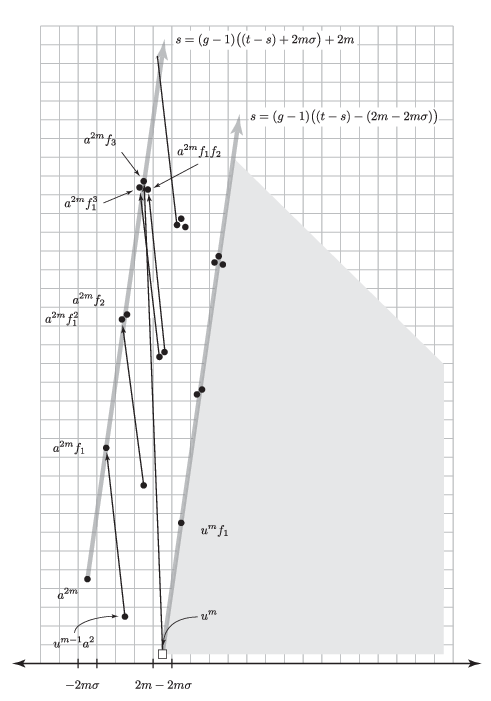}
\caption{Differentials on $u^{m}$}
\label{fig:6}
\end{figure}

Let $\sigma=\sigma_{G}$ be the real sign
representation of $G$, and
\[
u=u_{2\sigma} \in \pi_{2-2\sigma}^{G}
H\zltm
\]
the element defined in Definition~\ref{def:54}.  Since
\[
\slc{0}{\mutn{G}} = H\Zm
\]
the powers $u^{m}$ define elements
\[
u^{m}\in E_{2}^{0,2m-2m\sigma} = \pi^{G}_{2m-2m\sigma}\slc{0}{\mutn{G}}
\]
in the $E_{2}$-term of the $RO(G)$-graded slice spectral sequence 
\[
E_{2}^{s,t} = \pi^{G}_{t-s}\slc{\dim t}{\mutn{G}} \implies
\pi_{t-s}\mutn{G},
\]
with $t \in -2m\sigma+\Z$.  Our periodicity theorems depend on the
fate of these elements.  To study them it is convenient to consider
odd negative multiples of $\sigma$ as well, and to investigate the
slice spectral sequences for $\pi_{\ast -k\sigma}$ for $k\ge 0$.

It turns out to be enough to investigate the groups $E_{2}^{s,t}$ with
$s\ge (g-1)((t-s)-(k-k\sigma))$, where $g=|G|$.  The situation is
depicted in Figures~\ref{fig:4}-\ref{fig:6}.  We have, in fact,
already described all of the groups in this range.  To see this write
$t'=\dim t$ so that $t=t'+ (k-k\sigma)$, and
\[
E_{2}^{s,t} = \pi^{G}_{t'-s+k}S^{k\sigma}\wedge \slc{t'}{\mutn{G}}.
\]
Since $S^{k\sigma}\wedge \slc{t'}{\mutn{G}}\ge t'$, part \thmListItem{3} of
Proposition~\ref{thm:40} tells us that this group vanishes if 
\[
t'-s+k < \lfloor t'/g \rfloor, 
\]
and hence if
\[
s > (g-1)\big((t-s)+k\sigma ) + k.
\]
This gives the vanishing line depicted in
Figures~\ref{fig:4}-\ref{fig:6}.  Now $\slc{t'}{\mutn{G}}$ is
contractible unless $t'$ is even, in which case it a wedge of
$G$-spectra of the form $H\Zm\wedge \slicecell$ where $\slicecell$ is
a slice cell of dimension $t'$.  Since the restriction of $\sigma$ to
any proper subgroup is trivial, when
$\slicecell=G_{+}\underset{H}{\wedge}S^{\ell'\rho_{H}}$ is an induced
slice cell, there are isomorphisms
\[
S^{k\sigma}\wedge H\Zm\wedge\slicecell \approx
G_{+}\underset{H}{\wedge}\big( S^{k\sigma}\wedge
H\Zm\wedge S^{\ell'\rho_{H}}\big) \approx
G_{+}\underset{H}{\wedge} \big( S^{k}\wedge
H\Zm\wedge S^{\ell'\rho_{H}}\big)
\]
an so $\pi^{G}_{t'-s+k}S^{k\sigma}\wedge H\Zm\wedge \slicecell$ is
isomorphic to 
\[
\pi^{H}_{t'-s} H\Zm\wedge S^{\ell'\rho_{H}}.
\]
This group vanishes if 
\[
t'-s  < \ell' =  t'/h  \qquad (h=|H|),
\]
so certainly when 
\[
t'-s \le t'/g,
\]
or, equivalently when
\[
s\ge  (g-1)\big((t-s) - (k - k\sigma)\big).
\]
Thus in this range only the non-induced slice cells
contribute.     

The only even dimensional slice cells which are not induced are those
of the form $S^{\ell\rho_{G}}$.   We are therefore studying the groups 
\[
\pi^{G}_{j} H\Zm\wedge S^{k\sigma}\wedge S^{\ell\rho_{G}}
\]
with $j \le \ell+k$ and $k,\ell\ge 0$.  

\begin{lem}
\label{thm:84} 
For $k,\ell\ge 0$ and $j\le \ell+k$ the group
\[
\pi^{G}_{j} H\Zm\wedge S^{k\sigma}\wedge S^{\ell\rho_{G}} 
\]
is given by
\[
\pi^{G}_{j} H\Zm\wedge S^{k\sigma}\wedge S^{\ell\rho_{G}}  \approx
\begin{cases}
0 &\quad \text{if }(j-\ell)<0\text{ or }(j-\ell)\text{ is odd}    \\
\Z/2\cdot \{a_{\bar{\rho}}^{\ell}a_{\sigma}^{k-2m}u_{2\sigma}^{m}\} & \quad \text{if }(j-\ell) = 2m
\ge 0 \text{ and }\ell >0  \\
\zlt\cdot\{u_{2\sigma}^{m}\} &\quad \text{if }(j-\ell) = 2m
\ge 0 \text{ and }\ell  = 0.
\end{cases}
\]
\end{lem}

\begin{pf}
This computation reduces to the one described in Example~\ref{eg:20}.
To see this, write 
\[
S^{k\sigma}\wedge S^{\ell\rho_{G}} = 
S^{(k+\ell)\sigma}\wedge S^{\ell}\wedge S^{\ell(\rho_{G}-\sigma-1)}, 
\]
and consider the map
\begin{equation}
\label{eq:115}
a_{\bar{\rho}-\sigma}^{\ell}:\pi^{G}_{j} H\Zm\wedge S^{(k+\ell)\sigma}\wedge S^{\ell}
\to 
\pi^{G}_{j} H\Zm\wedge S^{k\sigma}\wedge S^{\ell\rho_{G}}.
\end{equation}
given by multiplication by $a_{\bar{\rho}-\sigma}^{\ell}$.  When
$\ell=0$ this map it is an isomorphism.  When $\ell>0$, the space
$S^{\ell(\rho_{G}-\sigma-1)}$ has the structure of a $G$-CW complex
with one $0$-cell and all other $G$-cells of positive dimension and
induced from proper subgroups.  Smashing with $S^{(k+\ell)\sigma}$ and
using the fact that the restriction of $\sigma$ to every proper
subgroup of $G$ is trivial, one finds that $S^{k\sigma}\wedge
S^{\ell\rho_{G}}$ is obtained from $S^{(k+\ell)\sigma}\wedge S^{\ell}$
by attaching induced $G$-cells of dimension greater than $(k+2\ell)$.
This implies that the map $a_{\bar{\rho}-\sigma}^{\ell}$ is an
isomorphism for $j<k+2\ell$, and so certainly for $j\le k+\ell$ since
$\ell>0$.  Thus in the range of interest, multiplication by
$a_{\bar{\rho}-\sigma}^{\ell}$ is isomorphism, and the computation
reduces to the evaluation of
\[
\pi^{G}_{j}
H\Zm\wedge S^{(k+\ell)\sigma}\wedge S^{\ell}.
\]
These groups were described in Example~\ref{eg:20}.
\end{pf}

To complete the description of the $E_{2}$-term of the $RO(G)$-graded
slice spectral sequence in this range we need to identify the summand
of non-induced slices of $\mutn{G}$.  From the associative algebra
equivalence
\[
\bigvee_{k\in\Z} \slc{k}{\mutn{G}} \sim H\Z\wedge S^{0}[G.\rr_{1},\dots]
\]
this is equivalent to identifying the summand of non-induced slice
cells in the twisted monoid ring
\[
S^{0}[G\cdot\rr_{1},\dots].
\]
Since the smash product of an induced spectrum with any spectrum is
induced, we can do this by identifying the summand of non-induced
slice cells in each
\[
S^{0}[G\cdot\rr_{i}]
\]
and smashing them together.   

Take the generating inclusion 
\[
\rr_{i}:S^{i\rho_{\zt}} \to S^{0}[\rr_{i}],
\]
apply $\norm_{\zt}^{G}$ to obtain 
\[
\norm\rr_{i}:S^{i\rho_{G}}\to S^{0}[G\cdot\rr_{i}],
\]
and extend it to an associative algebra map 
\begin{equation}
\label{eq:101}
S^{0}[\norm\rr_{i} ]\to S^{0}[G\cdot\rr_{i}].
\end{equation}

\begin{lem}
\label{thm:68}
The map~\eqref{eq:101} is the inclusion of the summand of non-induced
slice cells.
\end{lem}

\begin{pf}
The distributive law expresses $S^{0}[G\cdot\rr_{i}]=
\norm_{\zt}^{G}S^{0}[\rr_{i}]$ as an indexed wedge (see \S\ref{sec:meth-poly-algebr})
\[
S^{0}[G\cdot\rr_{i}] \approx \bigvee_{f:G/\zt\to\nat} S^{V_{f}},
\]
and $V_{f}=\bigoplus_{i=1}^{g/2}\gamma^{i}f(\gamma^{i})\rho_{\zt}$.
We now decompose the right hand side into an ordinary wedge over the
$G$-orbits.  Since an indexed wedge over a $G$-orbit is induced from
the stabilizer of any element of the orbit, the summand of non-induced
slice cells consists of those $f$ which are constant.   If $f:G/\zt$ is the
constant function with value $n$, then $V_{f}=n\rho_{G}$, so the
summand of non-induced slice cells is
\[
\bigvee_{\underline{n}} S^{n\rho_{G}}.
\]
The result follows easily from this.
\end{pf}

Smashing these together gives
\begin{cor}
\label{thm:82}
The associative algebra map 
\[
S^{0}[\norm\rr_{1},\dots]\to 
S^{0}[G\cdot\rr_{1},\dots]
\]
is the inclusion of the summand of non-induced slice cells.   
\end{cor}

To put this all together, consider the $\Z\times RO(G)$-graded ring
\[
\Z_{(2)}[a,f_{i}, u]/(2a, 2 f_{i})  \\
\]
with 
\begin{align*}
|a| &=(1,1-\sigma) \\
|f_{i}| &= (i(g-1),i g) \\
|u| &= (0,2-2\sigma).
\end{align*}
Define a map
\begin{equation}
\Z_{(2)}[a,f_{i}, u]/(2a, 2 f_{i}) \to \bigoplus_{\substack{s,k\ge 0 \\
t\in \ast-k\sigma}}E_{2}^{s,t}
\end{equation}
by 
\begin{align*}
f_{i} &\mapsto a_{\bar{\rho}}^{i}\norm\rr_{i}\in
E_{2}^{i(g-1), ig}=\pi_{i}^{G}\slc{ig}{\mutn{G}} \\
a &\mapsto a_{\sigma}\in E_{2}^{1,1-\sigma} = \pi_{-\sigma}
\slc{0}{0}\mutn{G} \\
\end{align*}
and by sending $u$ to the element $u\in E_{2}^{0,2-2\sigma}$ described
at the beginning of this section.  The combination of
Lemmas~\ref{thm:84} and~\ref{thm:68} gives
\begin{prop}
\label{thm:133}
The map
\begin{equation}
\label{eq:13}
\Z_{(2)}[a,f_{i}, u]/(2a, 2 f_{i}) \to \bigoplus_{\substack{s,k\ge 0 \\
t\in \ast-k\sigma}}E_{2}^{s,t}
\end{equation}
is an isomorphism in the range
\[
s\ge (g-1)((t-s)-(k-k\sigma)).  
\]
\qed
\end{prop}

We now turn to the differentials.  By construction, the $f_{i}$ are
the representatives at the $E_{2}$-term of the slice spectral sequence
of the elements defined in Definition~\ref{def:10} (and also called
$f_{i}$).  They are therefore permanent cycles.  Similarly, the
element $a$ is the representative of $a_{\sigma}$ and also a permanent
cycle.  This leaves the powers of $u$.  The case $G=C_{2}$ of the
following result appears in unpublished work of Araki and in
Hu-Kriz~\cite{MR1808224}.

\begin{thm}[Slice Differentials Theorem]
\label{thm:45}
In the slice spectral sequence for $\pi^{G}_{\star}\mutn{G}$ the
differentials $d_{i}u^{2^{k-1}}$ are zero for $i<r=1+(2^{k}-1)g$, and
\[
d_{r}u^{2^{k-1}} = a^{2^{k}}f_{2^{k}-1}.
\]
\end{thm}

\begin{rem}
\label{rem:5} It follows from Proposition~\ref{thm:133} that what lies
on the ``vanishing line''
\[
s = (g-1)\big((t-s)+k\sigma\big)+k
\]
is the algebra
\[
\Z_{(2)}[a,f_{i}]/(2a, 2f_{i}).
\]
In Proposition~\ref{thm:124} it was shown that the kernel of the map
\[
\Z_{(2)}[a_{\sigma},f_{i}]/(2a, 2f_{i})\to
\pi^{G}_{\star}\mutn{G}\to \pi^{G}_{\star}\phig\mutn{G}=\pi_{\ast}MO[a_{\sigma}^{\pm1}]
\]
is the ideal $(2, f_{1}, f_{3}, f_{7},\dots)$.  The only possible
non-trivial differentials into the vanishing line must therefore land
in this ideal.
\end{rem}

For the proof of Theorem~\ref{thm:45} the reader may find it helpful
to consult Figure~\ref{fig:6}.

\begin{pf*}{Proof of Theorem~\ref{thm:45}}
We establish the differential by induction on $k$.  Assume the result
for $k'<k$.  Then what's left in the range $s\ge (g-1)(t-s-k)$ after
the differentials assumed by induction is the sum of two modules over
$\Z_{(2)}[f_{i}]/(2 f_{i})$.  One is generated by $a^{2^{k}}$ and is
free over the quotient ring
\[
\Z/2[f_{i}]/(f_{1},f_{3},\dots, f_{2^{k-1}-1}).
\]
The other is generated by $u^{2^{k-1}}$.  Since the differential must
take its value in the ideal $(2,a,f_{1},f_{3},\dots)$, the next (and
only) possible differential on $u^{2^{k-1}}$ is the one asserted in
the theorem.  So all we need do is show that the classes $u^{2^{k-1}}$
do not survive the spectral sequence.  For this it suffices to do so
after inverting $a$.  Consider the map
\[
a_{\sigma}^{-1}\pi^{G}_{\star}\mutn{G}
\to a_{\sigma}^{-1}\pi^{G}_{\star}H\zltm.
\]
We know the $\Z$-graded homotopy groups of both sides, since they can
be identified with the homotopy groups of the geometric fixed point
spectrum.  If $u^{2^{k-1}}$ is a permanent cycle, then the class
$a^{-2^{k}}u^{2^{k-1}}$ is as well, and represents a class with
non-zero image in $\pi^{G}_{\ast}\phig H\zltm$.  This contradicts
Proposition~\ref{thm:72}.
\end{pf*}

\begin{rem}
\label{rem:3}
After inverting $a_{\sigma}$, the differentials
described in Theorem~\ref{thm:45} describe completely the
$RO(G)$-graded slice spectral sequence.  The spectral sequence starts
from
\[
\Z/2[f_{i}, a^{\pm1}, u].
\]
The class $u^{2^{k-1}}$ hits a unit multiple of $f_{2^{k}-1}$, and so
the $E_{\infty}$-term is
\[
\Z/2[f_{i},i\ne 2^{k}-1][a^{\pm1}]=MO_{\ast}[a^{\pm1}]
\]
which we know to be the correct answer since $\phig \mutn{G}=MO$.
This also shows that the class $u^{2^{k-1}}$ is a
permanent cycle modulo $(\rr_{2^{k}-1})$.  This fact corresponds to
the main computation in the proof of Theorem~\ref{thm:21} (which, of
course we used in the above proof).  The logic can be reversed, and
in~\cite{MR1808224} the results are established in the reverse order
(for the group $G=C_{2}$).
\end{rem}

Write
\[
\normrbar_{k} = N \rr_{2^{k}-1} \in
\pi^{G}_{(2^{k}-1)\rho_{G}}\mutn{G},
\]
and note that with this notation
\[
f_{2^{k}-1}= a_{\bar{\rho}}^{2^{k}-1}\normrbar_{k}.
\]
In the proof of the corollary below we will make use of the identity
\begin{equation}
\label{eq:116}
f_{2^{k+1}-1}\normrbar_{k}
=a_{\bar{\rho}}^{2^{k+1}-1}\normrbar_{k+1}\normrbar_{k}
= f_{2^{k}-1} a_{\bar{\rho}}^{2^{k}} \normrbar_{k+1}.
\end{equation}

The map $\normrbar_{k}:S^{(2^{k}-1)\rho_{G}}\to\mutn{G}$ is
represented at the $E_{2}$-term of the $RO(G)$-graded slice spectral
sequence by a map $S^{(2^{k}-1)\rho_{G}}\to\slc{(2^{k-1})g}{\mutn{G}}$
which we will also call $\normrbar_{k}$.  Multiplying, this defines
elements $\normrbar_{k}u^{2^{k}}$ in the $E_{2}$-term of the
$RO(G)$-graded slice spectral sequence.
\begin{cor}
\label{thm:50}
In the $RO(G)$-graded slice spectral sequence for
$\mutn{G}$, the class $\normrbar_{k}u^{2^{k}}$
is a permanent cycle.
\end{cor}

\begin{pf}
Write 
\[
r=1+(2^{k+1}-1)g.
\]
Theorem~\ref{thm:45} implies that differentials
$d_{i}(\normrbar_{k}u^{2^{k}})=\normrbar_{k} d_{i}(u^{2^{k}})$ are
zero for $i<r$, and
\[
d_{r}(\normrbar_{k} u^{2^{k}}) = \normrbar_{k}
a^{2^{k+1}}f_{2^{k+1}-1} = 
a^{2^{k+1}} f_{2^{k}-1} a_{\bar{\rho}}^{2^{k}} \normrbar_{k+1},
\]
the second equality coming from~\eqref{eq:116} above.  But from the
earlier differential
\[
d_{r'}u^{2^{k-1}} = a^{2^{k}} f_{2^{k}-1}
\]
where $r'=1+(2^{k}-1)g<r$, we also have
\[
d_{r'} (u^{2^{k-1}} a^{2^{k}} a_{\bar{\rho}}^{2^{k}}
\normrbar_{k+1}) = a^{2^{k+1}} f_{2^{k}-1} a_{\bar{\rho}}^{2^{k}} \normrbar_{k+1}
\]
so that in fact $d_{r}(\normrbar_{k} u^{2^{k}})=0$.  The target of the
remaining differentials work out to be in a region of the spectral
sequence which is already zero at the $E_{2}$-term.  So once we check
this, the proof is complete.

To check the claim about the vanishing region first note that with our
conventions, differential $d_{i+1}$ of the $RO(G)$-graded slice spectral
sequence maps a sub-quotient of
\[
\pi^{G}_{m}\slc{n}{X}
\]
to a sub-quotient of 
\[
\pi^{G}_{m-1}\slc{n+i}{X}.
\]
The class in question starts out at the $E_{2}$-term as
\[
\normrbar_{k}u^{2^{k}}\in \pi^{G}_{2^{k}(2-2\sigma)+(2^{k}-1)\rho_{G}}\slc{(2^{k}-1)g}{\mutn{G}}
\]
so we are interested in the groups
\[
\pi^{G}_{2^{k}(2-2\sigma)+(2^{k}-1)\rho_{G}-1}\slc{(2^{k}-1)g+i}{\mutn{G}} 
\]
or, equivalently 
\[
\pi^{G}_{2^{k+1}-1}\big(S^{2^{k+1}\sigma}\wedge
S^{-(2^{k}-1)\rho_{G}}\wedge \slc{(2^{k}-1)g+i}{\mutn{G}} \big)
\]
with $i+1 > r = 1+(2^{k+1}-1)g$.    To simplify the notation, write 
\[
X_{i} = S^{-(2^{k}-1)\rho_{G}}\wedge \slc{(2^{k}-1)g+i}{\mutn{G}},
\]
so that the group we are interested in is 
\begin{equation}
\label{eq:47}
\pi^{G}_{2^{k+1}-1}\big(S^{2^{k+1}\sigma}\wedge X_{i}\big).
\end{equation}
Now
\[
X_{i}\ge i.
\]
so Proposition~\ref{thm:40} implies that 
\[
\pi^{G}_{j}X_{i}=0
\]
for $j<\lfloor i/g \rfloor $.  Since $S^{2^{k+1}\sigma}$ is
$(-1)$-connected this means that if $i\ge 2^{k+1}g$ the
group~\eqref{eq:47} is trivial.  The remaining values of $i$ are
strictly between $(2^{k+1}-1)g$ and $(2^{k+1})g$, and hence not
divisible by $g$.  But since $\mutn{G}$ is pure, when $i$ is not
divisible by $g$ the spectrum $\slc{(2^{k}-1)g+i}{\mutn{G}}$ is
induced from a proper subgroup of $G$, hence so is $X_{i}$. There is
therefore an equivalence
\[
S^{2^{k+1}\sigma}\wedge X_{i} \approx
S^{2^{k+1}}\wedge X_{i},
\]
and so 
\[
\pi^{G}_{2^{k+1}-1}\big(S^{2^{k+1}\sigma}\wedge X_{i}\big) = 
\pi^{G}_{2^{k+1}-1}\big(S^{2^{k+1}}\wedge X_{i}\big) = 0
\]
since $X_{i}\ge 0$.
\end{pf}

\subsection{Periodicity theorems}
\label{sec:periodicity-theorem-1}

We now turn to our main periodicity theorem.  As will be apparent to
the reader, the technique can be used to get a much more general
result.  We have chosen to focus on a case which contains what is
needed for the proof of Theorem~\ref{thm:20}, and yet can be stated
for general $G=\ztn$.

Our motivating example is the spectrum $\kr$ of ``real''
$K$-theory~\cite{atiyah66:_k}.  Multiplication by the real Bott class
$\rr_{1}\in\pi_{\rho_{2}}\kr$ is an isomorphism, giving $\kr$ an
$S^{\rho_{2}}$-periodicity.  On the other hand, the representation
$4\rho_{2}$ admits a $\spin$ structure, and the construction of the
$\ko$-orientation of $\spin$ bundles leads to a ``Thom'' class
$u\in\pi^{\zt}_{8}\kr\wedge S^{4\rho_{2}}$.  This class is represented
at the $E_{2}$-term of the slice spectral sequence by $u_{4\rho_{2}}$.
Multiplication by $\rr_{1}^{4} u$ is then an equivariant map
$S^{8}\wedge\kr\to\kr$ whose underlying map of non-equivariant spectra
is an equivalence.  It therefore gives an equivalence $S^{8}\wedge
\kr^{h\zt}\approx \kr^{h\zt}$.  Since the map $\ko\to \kr^{h \zt}$ is
an equivalence, this gives the $8$-fold periodicity of $\ko$.

In our situation we begin with an equivariant commutative ring $R$, a
representation $V$ of $G$, and an element $D\in\pi^{G}_{V}R$.  We
manually create a spectrum with $S^{V}$-periodicity by working with
the homotopy colimit, $D^{-1}R$, of the sequence
\[
R \xrightarrow{D}{} S^{-V}\wedge R \xrightarrow{D}{} S^{-2V}\wedge
R \to\cdots.
\]
The unit inclusion
\[
S^{0}\to D^{-1}R
\]
gives a map
\[
H\Zm=\slc{0}{S^{0}}\to \slc{0}{D^{-1}R}
\]
and hence defines, for every oriented representation $W$ of $G$,
elements
\[
u_{W}\in \pi_{\dim W-W}\slc{0}{R}=E_{2}^{0,\dim W-W}
\]
in the $E_{2}$-term of the $RO(G)$-graded slice spectral sequence for
$\pi^{G}_{\star}D^{-1}R$.   We will show, under certain hypotheses on
$D$, that there is an integer $k>0$ with the property that
$u_{kV}$ is a permanent cycle.   Let
$u\in\pi^{G}_{\star}D^{-1}R$  be any element representing
$u_{k V}$.   Then the equivariant map 
\[
S^{k \dim V} \wedge D^{-1}R\xrightarrow{u}{} S^{k V}\wedge 
D^{-1}R\xrightarrow{D^{k}}{} D^{-1}R
\]
induces an equivalence of underlying, non-equivariant spectra, and
hence an equivalence of homotopy fixed point spectra
\[
\big(S^{k\dim V}\wedge D^{-1}R\big)^{hG} \to 
\big(D^{-1}R\big)^{hG}. 
\]
This establishes a periodicity theorem for the homotopy fixed point
spectrum $(D^{-1}R)^{h G}$.

The exposition is cleanest when one exploits multiplicative properties
of the spectrum $D^{-1}R$.  There are some easy general things to say
at first.  The spectrum $D^{-1}R$ is certainly an $R$-module, and
inherits a homotopy commutative multiplication (over $R$) from $R$.
The technique of \cite[\S{}VIII.4]{elmendorf97:_rings} can be used to
show that the non-equivariant spectrum underlying $D^{-1}R$ has a
unique commutative algebra structure for which the map
$i_{0}^{\ast}R\to i_{0}^{\ast}D^{-1}R$ is a map of commutative rings.

With an additional assumption on $D$, one can go further.  Let
$H\subset G$ be a subgroup, and suppose that there is an $m>0$ for
which the norm $\norm_{H}^{G}\big(i_{H}^{\ast}D\big)$ divides $D^{m}$.
Write $D^{m}=D'\cdot\norm_{H}^{G}\big(i_{H}^{\ast}D\big)$, and to keep
the notation compact abbreviate $\norm_{H}^{G}\big(i_{H}^{\ast}D\big)$
to $\norm_{H}^{G}D$.  Then there is a commutative diagram
\[
\xymatrix@!C@!R{
\norm_{H}^{G}R  \ar[r]^-{\norm_{H}^{G}(D)}\ar[d]  &
\norm_{H}^{G}\big(S^{-V}\wedge R\big) \ar[r]^-{\norm_{H}^{G}(D)}\ar[d]
&\norm_{H}^{G}\left(S^{-2V}\wedge R \right)   \ar[d] \cdots \\
R  \ar[r]^-{\norm_{H}^{G}(D)}\ar[d]^-{1}  &  S^{-V'}\wedge R \ar[r]^-{\norm_{H}^{G}(D)}\ar[d]^-{D'}   &  S^{-2V'}\wedge R \ar[d]^-{{D'}^{2}}\cdots \\
R  \ar[r]^-{D^{m}}        &  S^{-m V}\wedge R \ar[r]^-{D^{m}}
& S^{-2m V}\wedge R \cdots
}
\]
in which $V'=\ind_{H}^{G}V$.  Passing to the colimit gives a map
\[
\norm_{H}^{G}i_{H}^{\ast}\big(D^{-1}R\big)\to D^{-1}R
\]
extending the iterated multiplication.  This allows one to form norms
of elements in $\pi_{\star}^{H}D^{-1}R$ as if $D^{-1}R$ were an {\em equivariant}
commutative ring.

A necessary condition for $D^{-1}R$ to actually be an equivariant
commutative ring, is that for {\em every} $H\subset G$ the norm
$\norm_{H}^{G}i_{H}^{\ast}D$ divides a power of $D$.  In fact the
condition is also sufficient.  The proof of the result below is
described in~\cite{hill13:_equiv_multip_closur}.

\begin{prop}
\label{thm:61} Let $R$ be an equivariant commutative ring and
$D\in\pi^{G}_{\star}R$.  If $D$ has the property that for every
$H\subset G$, the element $\norm_{H}^{G}i_{H}^{\ast}D$ divides a power
of $D$, then the spectrum $D^{-1}R$ has a unique equivariant
commutative algebra structure for which the map $R\to D^{-1}R$ is a
map of commutative rings.  
\qed
\end{prop}

We will not make use of Proposition~\ref{thm:61}, as the ad hoc
formation of norms from the non-trivial subgroups of $G$ is sufficient
for our purpose.

Suppose that $t\in RO(G)$ and $u\in\pi^{H}_{t}D^{-1}R$ is represented
at the $E_{2}$-term of the $RO(H)$-graded slice spectral sequence by
the image of $u'\in\pi^{H}_{t}H\Z$ under the map $\pi^{H}_{t}H\Zm\to
\pi^{H}_{t}\slc{0}{D^{-1}R}$ induced by the unit.  We then have an
$H$-equivariant commutative diagram
\[
\xymatrix{
 &   \ar[dl]_{u} S^{t}  \ar[d]\ar[drr]^{u'}  &    & \\
D^{-1}R        &  P_{0}D^{-1}R \ar[r]\ar[l] &
\slc{0}{D^{-1}R} & H\Zm. \ar[l]
}
\]
The maps in the bottom row are maps of homotopy commutative ring
spectra.  Since the formation of slice sections commutes with filtered
colimits, if $\norm_{H}^{G}D$ divides a power of $D$ then the spectra
along the bottom row also come equipped with maps
$\norm_{H}^{G}(\slot)\to(\slot)$ extending the iterated
multiplication, and compatible with the maps between them.  This means
we may apply the norm to the whole diagram to produce
\[
\xymatrix{
 &   \ar[dl]_{\norm_{H}^{G}u} S^{\ind_{H}^{G}t}  \ar[d]\ar[drr]^{\norm_{H}^{G}u'}  &    & \\
D^{-1}R        &  P_{0}D^{-1}R \ar[r]\ar[l] &
\slc{0}{D^{-1}R} & H\Zm. \ar[l]
}
\]
showing that $\norm_{H}^{G}u'$ is a permanent cycle representing the
class $\norm_{H}^{G}u\in\pi_{\ind_{H}^{G}t}D^{-1}R$.

We will take $R$ to be the spectrum $\mutn{G}$.  In order to specify
the element $D$ we need to consider all of the spectra $\mutn{H}$ for
$H\subset G$, and we'll need to distinguish some of the important
elements of the homotopy groups we've specified.  We use~\eqref{eq:72}
to map
\[
\pi^{H}_{\star}\mutn{H} \to \pi^{H}_{\star}\mutn{G},
\]
and  make all of our computations in $\pi^{H}_{\star}\mutn{G}$.  Let
\[
\rr_{i}^{H} \in \pi^{\zt}_{i\rho_{2}} \mutn{H}\subset
\pi^{\zt}_{i\rho_{2}} \mutn{G}
\]
be the element defined in \S\ref{sec:specific-generators}, and
\[
\normrbar_{k}^{H} = \norm_{\zt}^{H} (\rr_{2^{k}-1}^{H})\in
\pi^{H}_{(2^{k}-1)\rho_{H}}\mutn{G}.
\]
Finally, in addition to $g=|G|$ we will write $h=|H|$ for
$H\subset G$.  

\begin{thm}
\label{thm:47} Let $\D\in\pi^{G}_{\ell \rho_{G}}\mutn{G}$ be any class
having the property that for every nontrivial $H\subset G$, the
element $\norm_{H}^{G}i_{H}^{\ast}D$ divides a power of $D$, and whose
image in $\pi^{H}_{\star}\mutn{G}$ is divisible by
$\normrbar_{g/h}^{H}$.  The class $u_{2\rho_{G}}^{2^{g/2}}$ is a
permanent cycle in the $RO(G)$-graded slice spectral sequence for
$\pi^{G}_{\star}D^{-1}\mutn{G}$.
\end{thm}

\begin{pf}
By Corollary~\ref{thm:50}, for each nontrivial subgroup $H\subset G$,
the class $\normrbar_{g/h}^{H} u_{2\sigma_{H}}^{2^{g/h}}$ is a
permanent cycle in the $RO(H)$-graded slice spectral sequence for
$\pi^{H}_{\star}\mutn{G}$.  Since $i_{H}^{\ast}D$ is divisible by
$\normrbar_{g/h}^{H}$, the class $u_{2\sigma_{h}}^{2^{g/h}}$ is then a
permanent cycle in the $RO(G)$-graded slice spectral sequence for
$\pi^{H}_{\star}D^{-1}\mutn{G}$.  From this inventory of permanent
cycles, and the ad hoc norm described above, we will show that
$u_{2\rho_{G}}^{2^{g/2}}$ is also a permanent cycle.

To begin, note that if $H\subset G$ has index $2$, then
$\ind_{H}^{G}1=1+\sigma_{G}$.  It follows from Lemma~\ref{thm:138}
that 
\[
u_{2\rho_{G}}=u_{2 \sigma_{G}}^{g/2} \norm_{H}^{G}
u_{2\rho_{H}}. 
\]
This leads to the formula
\[
u_{2\rho_{G}}^{k} = \prod_{0\ne H\subset G} \norm_{H}^{G}(u_{2\sigma_{H}}^{k
h/2}).
\]
When $k=2^{g/2}$ we have $kh/2=2^{g/2}h/2\ge 2^{g/h}$ for every $h\ne
1$ dividing $g$, so every term in the product is a permanent cycle
(the inequality is an equality only when $h=2$).  This completes the
proof.
\end{pf}

Write $\normr^{G} = u_{2\rho_{G}}\normrbar_{1}^{G}$.

\begin{cor}
\label{thm:139}
In the situation of Theorem~\ref{thm:47} the class
\begin{equation}\label{eq:44}
(\normr^{G})^{2^{g/2}} =
u_{2\rho_{G}}^{2^{g/2}}(\normrbar_{1}^{G})^{2\cdot 2^{g/2}}
\end{equation}
is a permanent cycle.  Any class in $\pi^{G}_{2\cdot g\cdot
2^{g/2}}\D^{-1}\mutn{G}$ represented by~\eqref{eq:44} restricts to a
unit in $\piu_{\ast}\D^{-1}\mutn{G}$.
\end{cor}

\begin{pf}
The fact that~\eqref{eq:44} is a permanent cycle is immediate from
Theorem~\ref{thm:47}.  Since the slice tower refines the Postnikov
tower, the restriction of an element in the $RO(G)$-graded group
$\pi^{G}_{\star}\D^{-1}\mutn{G}$ to $\piu_{\ast}\D^{-1}\mutn{G}$ is
determined entirely by any representative at the $E_{2}$-term of the
slice spectral sequence.  Since $u_{2\rho_{G}}$ restricts to $1$, the
restriction of any representative of~\eqref{eq:44} is equal to the
restriction of $(\normrbar_{1}^{G})^{2\cdot 2^{g/2}}$, which is a unit
since $\normrbar_{1}^{G}$ divides $\D$.
\end{pf}

This gives 

\begin{thm}
\label{thm:51} With the notation of Theorem~\ref{thm:47}, if $M$
is any equivariant $\D^{-1}\mutn{G}$-module, then multiplication by
$(\normr^{G})^{2^{g/2}}$ is a weak equivalence
\[
\Sigma^{2\cdot g\cdot 2^{g/2}}i_{0}^{\ast}M\to i_{0}^{\ast}M
\]
and hence an isomorphism
\[
(\normr^{G})^{2^{g/2}}:\pi_{\ast}M^{hG}\to \pi_{\ast+2\cdot g\cdot
2^{g/2}}M^{hG}.
\]
\qed
\end{thm}

For example, in the case of $G=\zt$ the groups
$\pi_{\ast}(D^{-1}\mutn{G})^{hG}$ are periodic with period $2*2*2=8$
and for $G=\zf$ there is a periodicity of $2*4*2^{2}=32$.  For $G=\ze$
we have a period of $2*8*2^{4}=256$.

\begin{rem}
\label{rem:46}
Suppose that $D\in\pi^{G}_{\star}R$ is of the form
\[
D=\norm_{\zt}^{G}x.   
\]
Then for $\zt\subset H\subset G$ one has
\[
\norm_{H}^{G}i_{H}^{\ast}D = D^{g/h}.
\]
Indeed,
\[
\norm_{H}^{G}i_{H}^{\ast}D = 
\norm_{H}^{G}i_{H}^{\ast}\norm_{\zt}^{G}x = 
\norm_{H}^{G}(\norm_{\zt}^{H})^{g/h}= \norm_{\zt}^{G}x^{g/h} = D^{g/h}.
\]
Since each $\normrbar_{k}^{H}$ has this form, any class $D$ which is a
product of $\norm_{H}^{G}\normrbar_{k}^{H}$ has the property required
for Theorems~\ref{thm:47} and~\ref{thm:51}.
\end{rem}

\begin{cor}[The Periodicity Theorem]
\label{thm:140}
Let $G=\ze$, and
\[
\D=(\norm_{\zt}^{\ze}\normrbar_{4}^{\zt})\,
(\norm_{\zf}^{\ze}\normrbar_{2}^{\zf}) \,
(\normrbar_{1}^{\ze}) \in \pi^{G}_{19\rho_{G}}\mutn{G}.
\]
Then multiplication by $(\normr^{G})^{16}$ gives an isomorphism
\[
\pi_{\ast}(\D^{-1}\mutn{G})^{hG}\to
\pi_{\ast+256}(\D^{-1}\mutn{G})^{hG}.
\]
\qed
\end{cor}

\begin{rem}
\label{rem:45} 
For a periodicity theorem, one gets a sufficient
inventory of powers of $u_{2\sigma_{H}}$ as permanent cycles as long
as for each $H$, some $\normrbar^{H}_{j}$ is inverted.  This is also
enough to prove the Homotopy Fixed Point Theorem.  Our particular
choice of $\normrbar^{H}_{g/h}$ is dictated by the requirements of the
Detection Theorem.
\end{rem}

\section{The Homotopy Fixed Point Theorem}
\label{sec:fixed-points-homot}

Until now we haven't had occasion to refer to the function
$G$-spectrum of maps from a pointed $G$-space $S$ to a $G$-spectrum
$X$, which exists as part of the completeness of $\gspectra{G}$ as a
topological $G$-category.    We will write $X^{S}$ for this object, so
that 
\[
\ugspectra{G}(Z,X^{S}) = \ugspectra{G}(Z\wedge S, X).
\]

\begin{defin}
\label{def:11}
A $G$-spectrum $X$ is {\em cofree} if the map
\begin{equation}
\label{eq:55}
X\to X^{EG_{+}}
\end{equation}
adjoint to the projection map $EG_{+}\wedge X\to X$ is a weak
equivalence.
\end{defin}

If $X$ is cofree then the map
\[
\pi^{G}_{\ast}X \to \pi^{G}_{\ast}X^{EG_{+}} = \pi_{\ast}X^{hG}
\]
is an isomorphism.  The main result of this section
(Theorem~\ref{thm:52}) asserts that any module over $D^{-1}\mutn{G}$
is cofree.

The map~\eqref{eq:55} is an equivalence of underlying spectra, and
hence becomes an equivalence after smashing with any $G$-CW spectrum 
built entirely out of free $G$-cells.   In particular, the map
\begin{equation}
\label{eq:98}
EG_{+}\wedge X\xrightarrow{\sim}{} EG_{+}\wedge (X^{EG_{+}})
\end{equation}
is an equivariant equivalence.  One exploits this, as
in~\cite{MR763905}, by making use of the pointed $G$-space $\tilde EG$
defined by the cofibration sequence
\begin{equation}
\label{eq:12}
EG_{+}\to S^{0}\to \tilde EG.
\end{equation}

\begin{lem}
\label{thm:166}
For a $G$-spectrum $X$, the following are equivalent:
\begin{thmList}
\item For all non-trivial $H\subset G$, the spectrum $\phih X$ is
contractible.
\item The map $EG_{+}\wedge X\to X$ is a weak equivalence.
\item The $G$-spectrum $\tilde EG\wedge X$ is contractible.
\end{thmList}
\end{lem}

\begin{pf}
The equivalence of the second and third conditions is immediate from
the cofibration sequence defining $\tilde EG$.    Since $EG_{+}$ is
built from free $G$-cells, condition \thmListItem{2}  implies condition \thmListItem{1}.
For $H\subset G$ non-trivial, we have 
\[
\phih(\tilde EG\wedge X)\approx
\phih(\tilde EG)\wedge \phih(X)\approx
S^{0}\wedge \phih(X).
\]
Since the non-equivariant spectrum underlying $\tilde EG$ is
contractible, condition \thmListItem{1} therefore implies that $\phih(\tilde
EG\wedge X)$ is contractible for {\em all} $H\subset G$.  But this
means that $\tilde EG\wedge X$ is contractible
(Proposition~\ref{thm:229}).
\end{pf}

\begin{cor}
\label{thm:228}
If $R$ is an equivariant ring spectrum satisfying the equivalent
conditions of Lemma~\ref{thm:166} then any module over $R$ is cofree.
\end{cor}

The condition of Corollary~\ref{thm:228} requires $R$ to be an
equivariant ring spectrum in the weakest possible sense, namely that $R$
possesses a unital multiplication (not necessarily associative) in
$\ho\ugspectra{G}$.  Similarly, the ``module'' condition is also one
taking place in the homotopy category.

\begin{pf}
Let $M$ be an $R$-module.  Consider the diagram
\begin{equation}
\label{eq:43}
\xymatrix{
EG_{+}\wedge M  \ar[r]\ar[d]  & M  \ar[r]\ar[d]   &  \tilde EG\wedge M \ar[d]\\
EG_{+}\wedge M^{EG_+}  \ar[r]        &  M^{EG_+}  \ar[r]
& \tilde EG\wedge M^{EG_+} 
}
\end{equation}
obtained by smashing $M\to M^{EG_+}$ with the sequence~\eqref{eq:12}.
The fact that $R$ satisfies the condition \thmListItem{1} of Lemma~\ref{thm:166}
implies that any $R$-module $M'$ does since $\phih(M)$ is a retract of
$\phih(R\wedge M)\approx \phih(R)\wedge \phih(M)$.  Thus both $M$ and
$M^{EG_+}$ satisfy the conditions of Lemma~\ref{thm:166}, and the
terms on the right in~\eqref{eq:43} are contractible.  The left
vertical arrow is the weak equivalence~\eqref{eq:98}.  It follows that
the middle vertical arrow is a weak equivalence.
\end{pf}

Turning to our main purpose, we now consider a situation similar to
the one in \S\ref{sec:periodicity-theorem-1}, and fix a class
\[
\D\in\pi^{G}_{\ell\rho_{G}}\mutn{G}
\]
with the property that for all non-trivial $H\subset G$ the restriction of
$\D$ to $\pi^{H}_{\ast}\mutn{G}$ is divisible by $\normrbar^{H}_{k}$
for some $k$ which may depend on $H$.

\begin{thm}[Homotopy Fixed Point Theorem]
\label{thm:52}
Any module $M$ over $\D^{-1}\mutn{G}$ is cofree, and so
\[
\pi^{G}_{\ast}M \to \pi_{\ast}M^{hG}
\]
is an isomorphism.
\end{thm}

\begin{pf}
We will show that $D^{-1}\mutn{G}$ satisfies condition \thmListItem{1} of
Lemma~\ref{thm:166}.   The result will then follow from
Corollary~\ref{thm:228}.  Suppose that $H\subset G$ is non-trivial.
Then
\[
\phih(D^{-1}\mutn{G}) \approx
\phih(D)^{-1}\phih(\mutn{G}).
\]
But $D$ is divisible by $\normrbar^{H}_{k}$, and so $\phih(D)$ is
divisible by 
\[
\phih(\normrbar^{H}_{k}) =
\phih(\norm_{\zt}^{H}(\rr^{H}_{2^{k}-1})) y
= \Phi^{\zt}(\rr^{H}_{2^{k}-1})
\]
which is zero by Proposition~\ref{thm:124}.  This completes the proof.
\end{pf}

\begin{cor}
\label{thm:141}
In the situation of Corollary~\ref{thm:140}, the map ``multiplication
by $\normr^{G}$'' gives an isomorphism
\[
\pi^{G}_{\ast}(\D^{-1}\mutn{G})\to
\pi^{G}_{\ast+256}(\D^{-1}\mutn{G}).
\]
\end{cor}

\begin{pf}
In the diagram
\[
\xymatrix{
{\pi^{G}_{\ast}(\D^{-1}\mutn{G})}\ar[r]\ar[d] &
{\pi^{G}_{\ast+256}(\D^{-1}\mutn{G})} \ar[d] \\
{\pi_{\ast}(\D^{-1}\mutn{G})^{hG}}\ar[r] &
{\pi^{G}_{\ast+256}(\D^{-1}\mutn{G})^{hG}}
}
\]
the vertical maps are isomorphisms by Theorem~\ref{thm:52}, and the
bottom horizontal map is an isomorphism by Corollary~\ref{thm:140}.
\end{pf}

\section{The Detection Theorem}\label{sec-detect}

\subsection{Outline of the proof}
\label{sec:outline-proof}

We now turn to the proof of the Detection Theorem.   For the
convenience of the reader, we restate the result.

\begin{thm} [The Detection Theorem] \label{thm:142} If
$\theta_{j}\in\pi_{2^{j+1}-2}S^{0}$ is an element of Kervaire
invariant $1$, and $j>2$, then the image of $\theta_{j}$ in
$\pi_{2^{j+1}-2}\magic$ is non-zero.
\end{thm}

To recapitulate, we are working with the group $G=\ze$, and the
spectrum $\Omega$ is the spectrum of $G$-fixed points in
$\premagic=\D^{-1}\mutn{G}$, with $\D\in\pi_{19\rho}\mutn{G}$ the
element specified in Corollary~\ref{thm:140}.

\begin{thm}[Algebraic Detection Theorem]
\label{thm:78}
If 
\[
x\in \ext_{MU_{*} (MU)}^{2, 2^{j+1}}\left(MU_{*}, MU_{*}\right)
\]
is any element mapping to 
\[
h_{j}^{2}\in \ext_{\sta}^{2,2^{j+1}}\left(\Z/2,\Z/2 \right)
\]
in the $E_{2}$-term of the classical Adams spectral sequence, and
$j>2$, then the image of $x$ in $H^{2} (\ze; \piu_{2^{j+1}}\premagic)$
is nonzero.
\end{thm}

We will prove the Algebraic Detection Theorem by establishing the
following.

\begin{prop}
\label{thm:230} 
For $j>2$, there is a map
\begin{equation}
\label{eq:151}
H^{2}(\ze;\pi_{2^{j+1}}\premagic) \to
\Q/\Z
\end{equation}
making the diagram
\begin{equation}
\label{eq:52}
\xymatrix{
*++{\ext^{2,2^{j+1}}_{\mustarmu}(\mustar,\mustar)}\ar[r]\ar[d]  & *++{H^{2}(\ze;\pi_{2^{j+1}}\premagic)}  \ar[d] \\
*++{\ext^{2,2^{j+1}}_{\sta}(\Z/2, \Z/2)}  \xyhookar[r]  &  *++{\Q/\Z}
}
\end{equation}
commute.
\end{prop}

In~\eqref{eq:52}, the bottom row is the Kervaire invariant
homomorphism sending $h_{j}^{2}$ to $1/2$.  Since the vector space
\[
\ext^{2,2^{j+1}}_{\sta}(\Z/2, \Z/2)
\]
has dimension $1$, with basis $h_{j}^{2}$
(Adams~\cite[Theorem~2.5.1]{adams60:_hopf}), the Kervaire invariant
homomorphism is completely specified by this property, and is a
monomorphism.  In plain language, Proposition~\ref{thm:230} asserts
that the Kervaire invariant homomorphism, thought of as a map
\[
\ext^{2,2^{j+1}}_{\mustarmu}(\mustar,\mustar) \to\Q/\Z,
\]
factors through $H^{2}(\ze,\pi_{2^{j+1}}\premagic)$.  This directly
implies Theorem~\ref{thm:78}.

\begin{rem}
\label{rem:24} All three of these results (Theorems~\ref{thm:142}
and~\ref{thm:78}, and Proposition~\ref{thm:230}) remain true without
the restriction $j>2$.  The other cases $j\le 2$ require separate
arguments, and are not needed for the proof of Theorem~\ref{thm:20},
so we do not include them.
\end{rem}

We now describe the proof of Proposition~\ref{thm:230}, deferring the
details to later subsections.  In order to construct the
map~\eqref{eq:151} we use of the theory of formal $A$-modules to
construct a $\ze$-equivariant ring homomorphism from
$\piu_{\ast}\premagic$ to much smaller ring.  Let $A=\Z_{2}[\zeta]$ be
the $2$-adic completion of the ring obtained by adjoining an
$8^{\text{th}}$ root of unity to the ring of integers, and $\ltgroup$
the Lubin-Tate formal $A$-module over $A$ associated to any  choice
of power series $f(x)\in A\LL x\RR$ satisfying (see~\S\ref{sec:formal-a-modules})
\begin{align*}
f(x) &\equiv \pi x \mod (x^{2}) \\
f(x) &\equiv x^{2} \mod (\pi),
\end{align*}
with uniformizer $\pi=\zeta-1$.  By construction, there is an isomorphism
\begin{align*}
A &\xrightarrow{\approx}{} \Endo(\ltgroup) \\
a &\mapsto [a](x)
\end{align*}
satisfying $[a]'(0)=a$.  Using the map $\gamma\mapsto \zeta$ to
identify the group of $8^{\text{th}}$ roots of unity with $\ze$ gives
an action of $\ze$ on $\ltgroup$ extending the canonical action of
$\zt$ by formal multiplication by $-1$.  As described in
$\S\ref{sec:conjugation-action}$ below, this data is classify by a
$\ze$-equivariant map of graded rings
\begin{equation}
\label{eq:153}
\piu_{\ast}\mutn{\ze}\to A_{\ast}
\end{equation}
in which $A_{\ast}=A[u^{\pm1}]$, $|u|=2$, and in which the action of
the chosen generator $\gamma\in\ze$ is the $A$-algebra map sending $u$
to $\zeta u$.  The first thing to check about this map is
\begin{prop}
\label{thm:239}
The image of $D\in\pi_{19\rho}\mutn{\ze}$ under
\[
\pi_{19\rho}\mutn{\ze} \to \piu_{152} \mutn{\ze} \to A_{152}
\]
is a unit, hence~\eqref{eq:153} factors uniquely through a
$\ze$-equivariant map
\begin{equation}
\label{eq:155}
\piu_{\ast}\premagic\to A_{\ast}.
\end{equation}
\end{prop}

Let
\[
\chi:H^{2}(\ze;\piu_{2^{j+1}}\premagic) \to 
H^{2}(\ze;A_{2^{j+1}})
\]
be the map of cohomology groups induced by~\eqref{eq:155}.  Using $\chi$,
form the rightmost arrow in the diagram below 
\begin{equation}
\label{eq:152}
\xymatrix{
*++{\ext^{2,2^{j+1}}_{\mustarmu}(\mustar,\mustar)}\ar[r]\ar[d]  & *++{H^{2}(\ze;\piu_{2^{j+1}}\premagic)}  \ar[d]^{\tfrac{2}{\pi}\cdot\chi} \\
*++{\ext^{2,2^{j+1}}_{\sta}(\F_{2}, \F_{2})}  \ar@{^{(}->}@<-.3ex>[r]
&  *++{H^{2}(\ze;A_{2^{j+1}})\mathrlap{\ .}}
}
\end{equation}
For the bottom arrow, note that both $\ext^{2,2^{j+1}}_{\sta}(\Z/2,
\Z/2)$ and $H^{1}(\ze;A_{2^{j+1}}/(\pi))$ are cyclic of order $2$,
and hence isomorphic by a unique isomorphism.  The bottom arrow
in~\eqref{eq:152} is defined to be the map corresponding to the
connecting homomorphism
\[
H^{1}(\ze;A_{2^{j+1}}/(\pi)) \to H^{2}(\ze;A_{2^{j+1}})
\]
under this isomorphism.  For $j>2$, the action of $\gamma$ on $u^{2^{j}}$ is trivial, and so
$H^{2}(\ze,A_{2^{j+1}})\approx A_{2^{j+1}}/(8)$, and one easily checks
that this map is a monomorphism.  

The main point is the commutativity of the diagram.   Once that is
established, the map~\eqref{eq:151} can be taken to be the composition
of the right vertical arrow in~\eqref{eq:152} with any map (dashed arrow)
\[
\xymatrix{
      &  *++{H^{2}(\ze;A_{2^{j+1}})} \ar@{-->}[d] \\
*++{H^{1}(\ze;A_{2^{j+1}}/(\pi))}  \xyhookar[ur]\xyhookar[r]        & *++{\Q/\Z.}
}
\]
factoring the inclusion through the connecting homomorphism.

Checking the commutativity of~\eqref{eq:152} involves some technical
details about the groups
$\ext^{2,2^{j+1}}_{\mustarmu}(\mustar,\mustar)$.  The following lemma
can be read off
from~\cite[Theorem~1.5]{shimomura81:_novik_rm_ext_sp2_at_prime} (see
\S\ref{sec:addendum}).

\begin{lem}
\label{thm:236}
For $j>1$, the map
\[
\ext^{1,2^{j+1}}_{\mustarmu}(\mustar,\mustar/(2)) \to 
\ext^{2,2^{j+1}}_{\mustarmu}(\mustar,\mustar)
\]
is surjective after localizing at $2$. \qed
\end{lem}

Lemma~\ref{thm:236} enables us to replace the upper left corner
of~\eqref{eq:152} with the group
$\ext^{1,2^{j+1}}_{\mustarmu}(\mustar,\mustar/(2))$, and verify the
commutativity of
\begin{equation}
\label{eq:156}
\xymatrix{
*++{\ext^{1,2^{j+1}}_{\mustarmu}(\mustar,\mustar/(2))}\ar[r]\ar[d]  & *++{H^{2}(\ze;\pi_{2^{j+1}}\premagic)}  \ar[d]^{\tfrac{2}{\pi}\cdot\chi} \\
*++{\ext^{2,2^{j+1}}_{\sta}(\Z/2, \Z/2)}  \ar@{^{(}->}@<-.3ex>[r]  &  *++{H^{2}(\ze;A_{2^{j+1}}).}
}
\end{equation}

The key technical point in doing this is

\begin{prop}
\label{thm:237}
The maps 
\begin{gather*}
\ext^{1,2^{j+1}}_{\mustarmu}(\mustar,\mustar/(2)) \to
\ext_{\sta}^{2,2^{j+1}}(\Z/2,\Z/2)  \\
\ext^{1,2^{j+1}}_{\mustarmu}(\mustar,\mustar/(2)) \to H^{1}(\ze;A_{2^{j+1}}/(\pi))
\end{gather*}
are surjective and have the same kernel.
\end{prop}

Proposition~\ref{thm:237} gives the commutativity of the left square in
\[
\xymatrix@C=1em{  
*++{\ext^{1,2^{j+1}}_{\mustarmu}(\mustar,\mustar/(2))}
\ar[r]\ar@{->>}[d] &
*++{H^{1}(\ze;\piu_{2^{j+1}}\premagic/(2))}   \ar[d]\ar[r]
&
*++{H^{2}(\ze;\piu_{2^{j+1}}\premagic)} \ar[d]^-{\tfrac{2}{\pi}\cdot\chi} 
\\
*++{\ext^{2,2^{j+1}}_{\sta}(\Z/2,\Z/2)} 
\ar[r]^{\approx}&
*++{H^{1}(\ze;A_{2^{j+1}}/(\pi))}   
\xyhookar[r] 
&
*++{H^{2}(\ze;A_{2^{j+1}})}. \\
}
\]
The commutativity of the right hand square follows from the naturality
of the connecting homomorphism.  The outer square is~\eqref{eq:156}.
This completes our summary of the proof of the Proposition~\ref{thm:230}
and the Detection Theorem.

\begin{rem}
\label{rem:59}
The argument of this section can be easily adapted to prove a
detection theorem for $\mutn{C_{2^{n}}}$ as long as $n\ge 3$.   The
result does not hold in the cases $n<3$.    What fails is the
assertion in Proposition~\ref{thm:237} that the two maps have the same
kernel.   This assertion makes essential use of the fact that the
reduction of the Lubin-Tate group over $A/(\pi)$ has height greater
than $2$.   
\end{rem}

The remainder of this section is devoted to filling in the details of
this outline.  We begin in \S\ref{sec:formal-a-modules} by recalling
the Lubin-Tate formal $A$-module~\cite{LT} and some simple but useful
results relating the power series $[a](x)$ to the $\pi$-adic valuation
of $a$.  We turn in \S\ref{sec:group-acti-homog} to the ideas
connecting the Adams-Novikov $E_{2}$-term to group cohomology.  In
\S\ref{sec:conjugation-action} we describe the ``conjugation action''
and prove Proposition~\ref{thm:240} which describes the functor
co-represented by $\piu_{\ast}\mutn{G}$ on the category of
$G$-equivariant graded commutative rings.  Setting all of this up
brings us as far as the statement of Proposition~\ref{thm:239} which
is proved in \S\ref{sec:fund-repr}.  Proposition~\ref{thm:237} is
proved in \S\ref{sec:technical-results}.  The proof relies heavily on
the computations in~\cite{MRW}
and~\cite{shimomura81:_novik_rm_ext_sp2_at_prime}, in the form of
Proposition~\ref{thm:231}.  An addendum to this section discusses how
these computations are made, and how they lead to Lemmas~\ref{thm:236}
and Proposition~\ref{thm:231}.

The reader may also wish to consult~\cite{MR3050712} for another
presentation of these ideas.

\subsection{Formal $A$-modules and the Lubin-Tate group}
\label{sec:formal-a-modules} 

Let $A$ and $R$ be commutative rings, and $e:A\to R$ a ring
homomorphism.   A ($1$-dimensional) {\em formal $A$-module over $R$}
is a formal group law $F$ over $R$, equipped with a ring homomorphism
\begin{align*}
A &\to \Endo(F) \\
a &\mapsto [a](x)
\end{align*}
with the property that $[a]'(0)=e(a)$.  In the case of interest to us,
$e$ is a monomorphism (in fact the identity map), and we will not
distinguish in notation between $a$ and $e(a)$.

Formal $A$-modules were introduced by Lubin and Tate in their
work~\cite{LT} on local class field theory.  For $A$ the ring of
integers in a local field with finite residue field, they constructed
a formal $A$-module over $A$ itself, unique up to isomorphism.  Their
construction starts with a choice of uniformizer $\pi\in A$ and a
power series
\[
f(x) \in A\LL x\RR
\]
intended to be the endomorphism $[\pi](x)$.  Writing $q$ for the order
of the residue field, the power series $f$ is required to satisfy
\begin{align*}
f(x) &\equiv \pi x \mod (x^{2}) \\
f(x) &\equiv x^{q} \mod (\pi).
\end{align*}
For example, $f(x)$ could be taken to be $\pi x+ x^{q}$.
Lubin and Tate showed that such an $f$ determines a formal $A$-module
in which the formal sum is the unique power series $\ltgroup(x,y)\in A\LL x\RR$
satisfying
\begin{align*}
\ltgroup(x,y) &\equiv x + y \mod (x,y)^{2} \\
\ltgroup(f(x),f(y)) &= f(\ltgroup(x,y)), 
\end{align*}
and for $a\in A$ the power series $[a](x)$ is the unique power series
satisfying
\begin{align*}
[a](x) &\equiv a x \mod (x)^{2} \\
[a](f(x)) &= f([a](x)).
\end{align*}
In particular, one does indeed have $[\pi](x)=f(x)$.

Continuing with the Lubin-Tate formal $A$-module, for $a\in A$ write
\[
[a](x) \equiv \alpha_{d}x^{d} + \cdots \mod (\pi)
\]
with $0\ne \alpha_{d}\in A/(\pi)$.  One easily checks that the
function $\nu(a)=\log_{q}(d)$ defines a valuation on $A$.  The the
fact that $[\pi](x)=f(x)$ implies that $\nu$ is the unique valuation
for which $\nu(\pi)=1$.

We are interested in the case
\[
A=\Z_{2}[\zeta],
\]
with $\pi = \zeta-1$, and any fixed choice of $f(x)$.  Since
$\nu(\zeta-1) = 1$, $\nu(\zeta^{2}-1) = 2$ and $\nu(\zeta^{4}-1)=4$,
and since any unit in $A$ is congruent to $1$ modulo $\pi$, this means
that modulo $(\pi)$
\begin{align*}
[\zeta-1](x) &\equiv x^{2} +\cdots \\
[\zeta^{2}-1](x) &\equiv x^{4} +\cdots \\
[\zeta^{4}-1](x) &\equiv x^{16} +\cdots, 
\end{align*}
and so 
\begin{equation}
\label{eq:110}
\begin{aligned}{}
[\zeta](x) &\equiv x \underset{\ltgroup}{+} x^{2} +\cdots \\
& \equiv x + x^{2} +\cdots \\
[\zeta^{2}](x) &\equiv x \underset{\ltgroup}{+} x^{4} +\cdots \\
& \equiv x + x^{4} +\cdots \\
[\zeta^{4}](x) &\equiv x \underset{\ltgroup}{+} x^{16} +\cdots  \\
&\equiv x + x^{16} +\cdots.
\end{aligned}
\end{equation}

These congruences play an important role in the proof of
Proposition~\ref{thm:239}.

\begin{rem}
\label{rem:60} The formulae~\eqref{eq:110} are independent of the
choice of $f(x)$.  In particular they hold for a choice of $f$ leading
to a $2$-typical formal group law.
\end{rem}

\subsection{Group actions and homogeneous formal group laws}
\label{sec:group-acti-homog} 
We now turn to the relationship between
group actions on formal group laws and group cohomology.  Our eventual
goal involves some explicit formulas, so we begin with a relatively
detailed summary.

\subsubsection{Homogeneous formal group laws}
\label{sec:homog-form-group}

Suppose that $R_{\ast}=\bigoplus R_{n}$ is a graded commutative ring.
By Quillen's work~\cite{Qui:FGL,Ad:SHGH} the set of graded ring
homomorphisms
\[
\mustar \to R_{\ast}
\]
is in one to one correspondence with the set of formal group laws $F$
over $R$ which are {\em homogeneous of degree $-2$} in the sense that the
formal sum
\[
F(x,y)
\]
is homogeneous of degree $-2$, when $x$ and $y$ are given degree
$-2$.    In terms of the power series 
\[
F(x,y) = x + y + \sum a_{ij}x^{i}y^{j}
\]
this means that $a_{ij}$ has degree $2(i+j-1)$.   

The graded ring $\mustarmu=\pi_{\ast}MU\wedge MU$ co-represents the
functor associating to a graded ring $R_{\ast}$ the set of pairs
$(F_{1}$, $F_{2})$ of homogeneous formal group laws and an isomorphism
$g:F_{1}\to F_{2}$ between them which is {\em strict} in the sense
that it is given by a power series of the form
\[
[g](x) = x + O(x^{2}).
\]
More generally, the ring $\pi_{\ast}MU^{(n)}$ ($n$-fold smash
product) co-represents functor associating to a graded ring $R_{\ast}$
the set of chains 
\[
F_{1}\to\dots\to F_{n}
\]
of homogeneous formal group laws and strict isomorphisms over $R$.
The standard convention is that the homogeneous formal group law
$F_{i}$ is the one classified by the map
\[
\pi_{\ast}MU\to \pi_{\ast}MU\wedge\cdots\wedge MU
\]
induced from the inclusion of the $i^{\text{th}}$ factor 
\[
S^{0}\wedge\dots\wedge MU\wedge \dots \wedge S^{0}\to MU\wedge \cdots\wedge MU.
\]

Taken together, the pair $(\mustar,\mustarmu)$ forms the {\em Hopf
algebroid} which co-represents the functor associating to a graded
commutative ring $R_{\ast}$ the groupoid of homogeneous formal group
laws over $R_{\ast}$ and strict isomorphisms.  For the definition of
Hopf algebroid, the reader is referred
to~\cite[Definition~A1.1.1]{Rav:MUnew}.

\subsubsection{Group actions}
\label{sec:group-actions-formal}

Let $\mfg$ be the category of pairs $(R,F)$, with $F$ a formal
group law over a commutative ring $R$, and in which a morphism
\[
(f,\psi):(R_{1},F_{1})\to (R_{2},F_{2})
\]
consists of a ring homomorphism $f:R_{1}\to R_{2}$, and an isomorphism
of formal group laws $\psi:F_{2}\xrightarrow{\approx}{}
f^{\ast}F_{1}$.  Morphisms can also be described as ring homomorphisms
\begin{align*}
h:R_{1}\LL x\RR &\to R_{2}\LL x \RR \\
h(r) &= f(r) \\
h(x) &= \psi(x)
\end{align*}
which are compatible with the formal sum in the sense that the diagram 
\[
\xymatrix{
*++{R_{1}\LL x\RR}  \ar[r]^{h}\ar[d]_{F_{1}}  &  *++{R_{2}\LL x\RR} \ar[d]^{F_{2}} \\
*++{R_{1}\LL x,y\RR}  \ar[r]_{h}   & *++{R_{2}\LL x,y\RR}
}
\]
commutes.  Let $\mfgh$ be the analogous category of {\em homogeneous}
formal group laws over graded rings, and strict isomorphisms.

The categories $\mfg$ and $\mfgh$ are related by the {\em
strictification} functor
\begin{align*}
\mfg & \to\mfgh \\
(R,F) &\mapsto (R_{\ast},F^{h}).
\end{align*}
The ring $R_{\ast}=R[u]$ is obtained from $R$ by adjoining a
polynomial variable $u$ with $|u|=2$,  and $F^{h}$ is the unique
formal power series satisfying
\begin{equation}
\label{eq:160}
u\, F^{h}(x,y) = F(u x, u x).
\end{equation}
The strictification of a map $(f,\psi):(R,F)\to (R',F')$ is the pair
\[
(f^{h}, \psi^{h}):(R_{\ast},F^{h})\to (R'_{\ast}, (F')^{h})
\]
with $f^{h}(u)= \psi'(0)u$ and $\psi^{h}(x)$ the unique power series
satisfying $u\, \psi^{h}(x) = \psi(u x)$.  

A {\em (left) action} of a group $G$ on a pair $(R,F)\in\mfg$ is a
map of monoids
\[
G\to \mfg((R,F), (R,F)),
\]
and a {\em strict (left) action} of a group on $(R_{\ast},F)\in\mfgh$ is
a map
\[
G\to \mfgh((R_{\ast},F),(R_{\ast},F)).
\]
The strictification functor converts a group action into a strict one.

A left action of $G$ on $(R,F)$ corresponds to a left action of $G$ on
$R\LL x\RR$.  We will use the notation
\begin{align*}
r &\mapsto g\,r \quad r\in R \\
x & \mapsto [g](x)
\end{align*}
for this action. 

\begin{eg}
\label{eg:12} Suppose that $E$ is a complex oriented, homotopy
commutative ring spectrum, and that a finite group $G$ acts on $E$ by
homotopy multiplicative maps.  Let $F$ denote the corresponding
homogeneous formal group law over $\pi_{\ast}E$.  Then the action of
$G$ on $E^{\ast}(\cp^{\infty})$ gives a strict action of $G$ on
$(\pi_{\ast}E,F)$.
\end{eg}

\begin{eg}
\label{eg:21} The group $\zt$ acts on any $(R,F)\in\mfg$ as the
identity map on $R$ and formal multiplication by $(-1)$ on $F$.
\end{eg}

\subsubsection{Group cohomology and the Adams-Novikov spectral
sequence}
\label{sec:group-cohom-adams}
When $(R_{\ast},F)$ is a homogeneous formal group law equipped with a
strict action of a group $G$, there is a map
\[
\ext^{s,t}_{\mustarmu}(\mustar,\mustar) \to H^{s}(G;R_{t}).
\]
Conceptually, it arises from the inclusion functor of the subcategory
of $\mfgh$ whose only object is $(R_{\ast},F)$ and whose monoid of
self maps is given by the action of $G$.  For the purposes of explicit
computation, it is  conveniently described as derived from a map of
Hopf algebroids
\begin{equation}
\label{eq:119}
(\mustar,\mustarmu) \to (R_{\ast}, C(G;R_{\ast})),
\end{equation}
in which $C(G;R_{\ast})$ is the ring of (set-theoretic) functions from
$G$ to $R_{\ast}$. 

The Hopf algebroid $(R_{\ast}, C(G;R_{\ast}))$ is the one expressing
the group action of $G$ on $R_{\ast}$.  The ``source'' map
\[
\eta_{L}:R_{\ast}\to C(G;R_{\ast})
\]
sends $r\in R_{\ast}$ to the constant function with value $r$, and the
``target'' map $\eta_{R}:R_{\ast}\to C(G;R_{\ast})$ is the transpose
of the action mapping.  It associates to $r\in R_{\ast}$ the function
sending $g\in G$ to $g\cdot r$.  The coproduct
\[
\Delta:C(G;R_{\ast}) \to C(G;R_{\ast})\underset{R_{\ast}}{\otimes}
C(G;R_{\ast})
\]
is the composition of the map
\[
C(G;R_{\ast})\to C(G\times G,R_{\ast})
\]
dual to multiplication in $G$, and the isomorphism
\[
C(G;R_{\ast})\underset{R_{\ast}}{\otimes} C(G;R_{\ast}) 
\xrightarrow{\approx}{} C(G\times G,R_{\ast}) 
\]
given by setting
\[
(f_{1}\otimes f_{2})(g_{1},g_{2}) =f_{1}(g_{1})\cdot g_{1}f_{2}(g_{2}).
\]
The map~\eqref{eq:119} consists of the map $\mustar\to R_{\ast}$
classifying the homogeneous formal group law $F$, and the map
$\mustarmu\to C(G,R_{\ast})$, defined by declaring the composition
\[
\mustarmu\to C(G,R_{\ast})\xrightarrow{\ev_{g}}{} R_{\ast}
\]
to be the map classifying the strict isomorphism 
\[
[g](x):F\to g^{\ast}F.
\]

When the $G$-action on $(R_{\ast}, F)$ arises, as in
Example~\ref{eg:12}, from an action of $G$ on a complex oriented
homotopy commutative ring spectrum $E$, the map~\eqref{eq:119} is the
$E_{2}$-term of a map of spectral sequences abutting to the
homomorphism $\pi_{\ast}S^{0}\to \pi_{\ast}E^{hG}$.  We couldn't quite
find this result in the literature
(though~\cite[Proposition~6.7]{devinatz04:_homot_morav} is close).  To
see it, let
\[
C^{\bullet}(G;E) = \map^{G}(EG_{\bullet},E)
\]
be the cosimplicial spectrum of $G$-maps from the bar construction
model for $EG_{\bullet}$ into $E$.  Thus
\[
C^{n}(G;E) = \prod_{G^{n}}E
\]
and $\tot C^{\bullet}(G;E)$ is the homotopy fixed point spectrum
$E^{hG}$.  The cosimplical ring $[n]\mapsto \pi_{\ast}C(G^{n},E)$
is the nerve of the Hopf-algebroid $(\pi_{\ast}E,
C(G,\pi_{\ast}E))$ and forms the cobar complex for calculating
$H^{\ast}(G,\pi_{\ast}E)$.   The homotopy fixed point spectral
sequence is the homotopy spectral sequence of this cosimplicial
spectrum.   

Choose a complex orientation for $E$, and for every $n\ge 0$ let
\begin{equation}
\label{eq:54}
MU^{(n+1)}\to C^{n}(G;E) \in \ho\spectra
\end{equation}
be a representative of the unique homotopy class of homotopy
multiplicative maps whose restriction to the $i^{\text{th}}$ smash
factor of $MU^{(n+1)}$ is the composition of the complex orientation
\[
MU\to E = C^{0}(G;E)
\]
with the cosimplicial structure map $C^{0}(G;E)\to C^{n}(G;E)$
corresponding to the inclusion of the $i^{\text{th}}$-vertex of
$\Delta[n]$.  The maps~\eqref{eq:54} fit into a homotopy commutative
map of cosimplicial resolutions
\begin{equation}
\label{eq:159}
\xymatrix@C=2em{
S^{0}  \ar[r]\ar[d]  &  MU \ar@<.5ex>[r]\ar@<-.5ex>[r]\ar[d] & MU\wedge MU \ar@<1ex>[r]\ar@<0ex>[r]\ar@<-1ex>[r]\ar[d] & MU\wedge MU\wedge MU\cdots\ar[d] \\
E^{hG}  \ar[r]        & C^{0}(G;E) \ar@<.5ex>[r] \ar@<-.5ex>[r]
& C^{1}(G;E)\ar@<1ex>[r]\ar@<0ex>[r]\ar@<-1ex>[r] & C^{2}(G;E)\cdots.
}
\end{equation}
If this were an actually commutative diagram, the desired spectral
sequence would be the one derived from the induced map of ``Tot
towers.''  Even though it is only homotopy commutative, the fact that
the top row is an $MU$ Adams resolution and the spectra in the bottom
row are complex oriented means that it can be still be refined to a
map of towers.

This result doesn't quite seem to appear in the literature, though an
assertion along these lines is made in~\cite[pp 289-90]{Mlr:Rel}, and
the case in which $MU$ is replaced by $E_{n}$
is~\cite[Proposition~6.2]{devinatz04:_homot_morav}.  To spell out the
details, we begin with some generalities about the Tot tower of a
cosimplical spectrum $X^{\bullet}$.  Let $NX^{n}$ be the iterated
mapping cone of the coface maps
\[
d^{i}:X^{n-1}\to X^{n}, \quad i =1,\dots, n.
\]
The spectrum $NX^{n}$ is a retract of $X^{n}$ (it is split by the
inclusion of the iterated homotopy fiber of the co-degeneracy maps).
By construction it depends, as a retract of $X^{n}$, functorially on
$X^{\bullet}$.  The spectra $NX^{n}$ fit into a sequence
\[
NX^{0}\xrightarrow{d^{0}}{} NX^{1}\xrightarrow{d^{0}}{}\cdots
\]
which is a ``complex'' in the sense that the composition of any two
maps is null homotopic.  

The homotopy spectral sequence of the cosimplicial spectrum
$X^{\bullet}$ is derived from the tower $\{\tot_{n}X^{\bullet} \}$.
For our purposes it is easier to work with the fibers of the map from
$\tot X^{\bullet}$.  Write $FX^{0}=\tot X^{\bullet}$, and define
$FX^{n}$ to be the homotopy fiber of the map $\tot X^{\bullet}\to
\tot_{n-1}X^{\bullet}$.  Then there is a functorial fibration sequence
\begin{equation}
\label{eq:89}
FX^{n}\to FX^{n-1}\to \Sigma^{-(n-1)} NX^{(n-1)}.
\end{equation}
Of course the homotopy spectral sequence can also be derived from the
tower $\{FX^{n} \}$, for example by using it to reconstruct the
Tot-tower.

To simplify the notation, write $X^{\bullet}=MU^{\bullet+1}$ and
$Y^{\bullet}=C^{\bullet}(G;E)$ for the cosimplicial spectra occurring
in the top and bottom rows of~\eqref{eq:159}.  The complex
$NX^{\bullet}$ is the standard $MU$-Adams resolution for $S^{0}$, and
in this case the fibration sequence~\eqref{eq:89} is equivalent to
\[
FX^{n}\to FX^{n-1}\to MU\wedge FX^{n-1}.
\]
The consequence we need of this is the characterizing property of an
Adams tower: if $R$ is any $MU$-module, then the connecting homomorphism
\begin{equation}
\label{eq:99}
[FX^{n},R] \hookrightarrow [\Sigma^{-n-1}NX^{n},R]
\end{equation}
is a monomorphism.

Our aim is to construct a map of towers
\begin{equation}
\label{eq:87}
\{FX^{n} \}\to \{FY^{n} \}.
\end{equation}
Suppose by induction we have produced a homotopy commutative diagram
\[
\xymatrix{
FX^{n-1} \ar[r]\ar[d]   &  \Sigma^{-(n-1)}NX^{n-1} \ar[d]\\
FY^{n-1} \ar[r]         & \Sigma^{-(n-1)}NY^{n-1}, 
}
\]
and choose any map  $FX^{n}\to FY^{n}$ making
\[
\xymatrix{
\Sigma^{-n}NX^{n-1} \ar[r]\ar[d] & FX^{n}  \ar[r]\ar[d]  &  FX^{n-1} \ar[d] \\
\Sigma^{-n}NY^{n-1} \ar[r]       & FY^{n}  \ar[r]        &  FY^{n-1} 
}
\]
commute up to homotopy.  We claim that the diagram
\[
\xymatrix{
FX^{n}  \ar[r]\ar[d]  &  \Sigma^{-n}NX^{n} \ar[d] \\
FY^{n}  \ar[r]        &  \Sigma^{-n}NY^{n}.
}
\]
also commutes up to homotopy.    The claim completes the induction step,
and gives~\eqref{eq:87}.

To verify the claim, consider
\begin{equation}
\label{eq:114}
\xymatrix{
\Sigma^{-n}NX^{n-1}\ar[r]\ar[d] & FX^{n}  \ar[r]\ar[d]  &  \Sigma^{-n}NX^{n} \ar[d] \\
\Sigma^{-n}NY^{n-1}\ar[r] &FY^{n}  \ar[r]        &
\Sigma^{-n)}NY^{n}\mathrlap{\ .}
}
\end{equation}
The outermost square commutes since it is a retract of a suspension of
one of the squares in~\eqref{eq:159}.  The spectrum
$\Sigma^{-(n+1)}NY^{n+1}$ admits the structure of an $MU$-module since
it is a retract of a suspension of $C^{n+1}(G;E)$ which is complex
orientable.  Taking it for $R$ in the monomorphism~\eqref{eq:99} shows
that the commutativity of the outer square in~\eqref{eq:114} implies
the commutativity of the right hand square.  This verifies the claim.

\subsubsection{The conjugation action}
\label{sec:conjugation-action}

Applying the strictification functor to Example~\ref{eg:21}, one is
led to the ``conjugation action'' of $\zt$ on homogeneous formal group
laws over graded rings.

Let $F$ be a homogeneous formal group law over a graded
commutative ring $R_{\ast}$, and $c:R_{\ast}\to R_{\ast}$ any ring
homomorphism with the property that $c:R_{2n}\to R_{2n}$ the map given
by multiplication by $(-1)^{n}$.  The homogeneous formal group law
$F^{c}:= c^{\ast}F$ is given by
\[
F^{c}(x,y) = - F(-x, -y).
\]
The power series
\[
c(x) = - [-1]_{F}(x)
\]
has the property that $c\circ c(x)=x$ and provides both a strict
isomorphism $F\to F_{c}$ and its inverse $F^{c}\to F$.  These combine
to give an action of $\zt$ on $(R_{\ast},F)$ which we call the {\em
conjugation action} associated to $c:R_{\ast}\to R_{\ast}$.

The map $c$ is completely specified on the even degree elements in
$R_{\ast}$ and in general there are as many conjugation actions as
there are ways of extending $c$ to all of $R_{\ast}$.  In the examples
of interest to us, $R_{\ast}$ will be evenly graded, and so there is
exactly one conjugation action.

\begin{eg}
\label{eg:22} 
If $E$ is a real-oriented spectrum, then the underlying
$\zt$ action on $(i_{0}^{\ast}E)^{\ast}\LL \cp^{\infty}\RR$ is the
conjugation action.  The case $E=\mur$ is universal in the sense that
the map $\mustar\to R_{\ast}$ classifying a homogeneous formal group
law is equivariant for any choice of conjugation action.
\end{eg}

We now generalize Example~\ref{eg:22}.  Let $G=\ztn$, and give
$i_{1}^{\ast}\mutn{G}$ the real orientation coming from the unit
\[
\mur\to i_{1}^{\ast}\mutn{G}
\]
of the norm-restriction adjunction on equivariant commutative
algebras (Example~\ref{eg:18}).   Examples~\ref{eg:12}
and~\ref{eg:22} then equip  $(\piu_{\ast}\mutn{G}, F)$ with a
$G$-action extending the conjugation action of $\zt$.   

\begin{prop}
\label{thm:240} The pair $(\piu_{\ast}\mutn{G},F)$ equipped with its
$G$-action is {\em universal} in the sense that map associating to a
$G$-equivariant
\begin{equation}
\label{eq:53}
f:\piu_{\ast}\mutn{G} \to R_{\ast},
\end{equation}
the pair $(R\ast,f^{\ast}F)$ with its induced $G$-action, is a
bijection between the set of equivariant maps~\eqref{eq:53} and the set of
$(R_{\ast},F)$ equipped with a $G$-action extending the one on
$R_{\ast}$, and the conjugation action of $\zt\subset G$.  
\end{prop}

\begin{pf}
Suppose that $(R_{\ast}, F)$ is a homogeneous formal group law over a
graded ring, equipped with a $G$-action extending the conjugation
action of $\zt$.   Choose a generator $\gamma\in G$.   This data is
equivalent to an isomorphism 
\[
\tau:F_{1}\to F_{0}
\]
having the property that the composite of the chain of isomorphisms 
\begin{gather*}
F_{g/2}\xrightarrow{\tau_{g/2-1}}{} F_{g/2-1}\xrightarrow{}{}\cdots
\xrightarrow{\tau_{1}}{}F_{1}\xrightarrow{\tau_{0}}{} F_{0} \\
\tau_{i} = (\gamma^{i})^{\ast}\tau
\end{gather*}
is the conjugation isomorphism $c$.  The claim then follows from the
decomposition~\eqref{eq:73} and the description of
$\pi_{\ast}MU\wedge\dots\wedge MU$ in terms of chains of composable
strict isomorphisms.
\end{pf}

\subsection{The fundamental representation}
\label{sec:fund-repr}

As described in \S\ref{sec:outline-proof} these ideas can be used to
construct a $\ze$-equivariant ring homomorphism from
$\piu_{\ast}\premagic$ to much smaller ring.  Let $A=\Z_{2}[\zeta]$ be
as in \S\ref{sec-detect}, and $\ltgroup$ the Lubin-Tate formal
$A$-module over $A$ associated to a power series series $f(x)$ leading
to a $2$-typical formal group law.  Using $\zeta$ to identify the
group of $8^{\text{th}}$-roots of unity with $\ze$, we get a $\ze$
action on $(A,\ltgroup)$.  From this apply the strictification functor
to get a strict action of $\ze$ on $(A_{\ast}, \ltgroup^{h})$.  With
an eye toward Proposition~\ref{thm:239}, we invert the class $u$ and
re-define $A_{\ast}$ to be $A[u,u^{-1}]$.  The underlying $\zt$-action
is the conjugation action, so Proposition~\ref{thm:240} provides a
$\ze$-equivariant map
\begin{equation}
\label{eq:46}
\piu_{\ast}\mutn{\ze}\to A_{\ast},
\end{equation}
classifying $(A_{\ast},\ltgroup^{h})$ with its $\ze$-action.

\begin{prop}
The image of $D$ under
\[
D\in\pi_{19\rho}\mutn{\ze} \to \piu_{152} \mutn{\ze} \to A_{152}
\]
is a unit, hence~\eqref{eq:46} factors through a $\ze$-equivariant map 
\[
\piu_{\ast}\premagic\to A_{\ast}.
\]
\end{prop}

\begin{pf}
We must show that the classes $r_{1}^{\ze}$, $r_{3}^{\zf}$, and
$r_{15}^{\zt}$ all map to units in $A_{\ast}$.  It suffices to show
that they do so in $A_{\ast}/(\pi)$.  By definition
(\S\ref{sec:specific-generators}) the image of $r_{1}^{\ze}$ in
$A_{2}/(\pi)$ is the coefficient of $x^{2}$ in the isomorphism of
$\gamma^{\ast}\ltgroup^{h}$ with the $2$-typification of $\ltgroup^{h}$.  Since
$\ltgroup^{h}$ is already $2$-typical, this is just the coefficient of
$x^{2}$ in the power series $[\zeta]_{\ltgroup^{h}}(x)$ in the
homogeneous formal group law.  By~\eqref{eq:110} this coefficient is
congruent to $u$ modulo $\pi$, hence a unit.  Equation~\eqref{eq:110}
similarly shows that $r_{3}^{\zf}$ maps to $u^{3}$, and $r_{15}^{\zt}$
to $u^{15}$ modulo $(\pi)$.  This completes the proof.
\end{pf}

\subsection{Technical results}
\label{sec:technical-results}

In this section we describe explicitly the maps 
\begin{equation}
\label{eq:157}
\begin{gathered}
\ext^{1,2^{j+1}}_{\mustarmu}(\mustar,\mustar/(2)) \to
\ext_{\sta}^{2,2^{j+1}}(\Z/2,\Z/2)  \\
\ext^{1,2^{j+1}}_{\mustarmu}(\mustar,\mustar/(2)) \to H^{1}(\ze;A_{2^{j+1}}/(\pi))
\end{gathered}
\end{equation}
occurring in the statement of Proposition~\ref{thm:237}.  The results are
Propositions~\ref{thm:144} and~\ref{thm:234} below.  Combined with
Lemma~\ref{thm:236} and Proposition~\ref{thm:231} they directly imply
Proposition~\ref{thm:237}.

\subsubsection{Preliminaries}
\label{sec:preliminaries}

We remind the reader that everything has been localized at the prime
$p=2$.  The $2$-typification $\funiv^{(2)}$ of the universal formal
group law $\funiv$ is classified by a map $\bpstar\to\mustar$.  This
map extends to an equivalence of Hopf algebroids
\[
(\bpstar, \bpstarbp) \xrightarrow{\approx}{} (\mustar,\mustarmu),
\]
giving an isomorphism
\begin{equation}
\label{eq:88}
\ext^{s,t}_{\bpstarbp}(\bpstar,\bpstar)\xrightarrow{\approx}{}
\ext^{s,t}_{\mustarmu}(\mustar,\mustar).
\end{equation}
Our proofs depend heavily on the computations in~\cite{MRW}
and~\cite{shimomura81:_novik_rm_ext_sp2_at_prime}, which are stated
for $BP$.  Because of the above isomorphism they apply equally well to
$MU$.

In order to describe explicit computations, we fix the identification
\[
\bpstarbp= \Z_{(2)}[v_{1}, v_{2},\dots, t_{1},t_{2},\dots]
\]
in which the $v_{i}$ are the Hazewinkel generators, and the elements
$t_{i}\in \bpstarbp$ are the coefficients of the universal isomorphism
\begin{align*}
\eta_{R}^{\ast}\funiv^{(2)} &\to \funiv^{(2)}=
\eta_{L}^{\ast}\funiv^{(2)} \\
x &\mapsto \sum\nolimits^{\funiv^{(2)}} t_{n} x^{2^{n}}.
\end{align*}
We will not distinguish in notation between the $v_{i}$ and $t_{i}$ in
$\bpstarbp$, and their images in $\mustarmu$.

An important role in the proof of~\ref{thm:230} is played by the
element $t_{1}$.  Since any coordinate $x$ is $2$-typical modulo
$x^{3}$, the class $t_{1}$ is also given by the coefficient of the
universal isomorphism of $\eta_{R}^{\ast}{\funiv}$ with $\eta_{L}^{\ast}\funiv$
\[
x \mapsto x + t_{1}x^{2} + \cdots.
\]
With the standard conventions this is the inverse of the universal
strict isomorphism over $\pi_{\ast}MU\wedge MU$, which goes from
$\eta_{L}^{\ast}\funiv$ to $\eta_{R}^{\ast}\funiv.$

\subsubsection{The Adams-Novikov $2$-line}

Let
\begin{align*}
\delta_{2} &: \ext_{\mustarmu}^{1, 2^{j+1}}(\mustar, \mustar/(2)) \xrightarrow{}{}
\ext_{\mustarmu}^{2, 2^{j+1}}(\mustar, \mustar) \quad\text{and}\\
\delta_{1} &:\ext_{\mustarmu}^{0, 2^{j+1}}(\mustar, \mustar/(2,
v_{1}^{\infty}))\xrightarrow{} \ext_{\mustarmu}^{1,
2^{j+1}}(\mustar, \mustar/(2))
\end{align*}
be the connecting homomorphisms associated to the short exact
sequences of $\mustarmu$ co-modules
\begin{gather*}
0\to \mustar\xrightarrow{2}{} \mustar\to \mustar/(2) \to 0\\
0\to \mustar/2\to v_{1}^{-1}\mustar\to \mustar/(2,v_{1}^{\infty}) \to 0.
\end{gather*}
Our description of the maps~\eqref{eq:157} relies on the following
computation, whose proof is discussed in \S\ref{sec:addendum}.  We
employ the standard ``cobar construction'' notation for elements (see,
for example,~\cite[Definition A1.2.11]{Rav:MUnew}).

\begin{prop}[\cite{MRW,shimomura81:_novik_rm_ext_sp2_at_prime}]
\label{thm:231} For $j>1$, the $\Z/2$-vector space
$\ext_{\mustarmu}^{1, 2^{j+1}}(\mustar, \mustar/(2))$ has a basis
consisting of the elements,
\[
v_{1}^{2^{j}-1}[t_{1}], \text{ } v_{1}^{2^{j}-2}[t_{1}^{2}]
\]
and the image under $\delta_{1}$ of certain elements of the form
\[
v_{2}^{s 2^{k}}/v_{1}^{2^{k}} \in 
\ext_{\mustarmu}^{0, 2^{j+1}}(\mustar, \mustar/(2,
v_{1}^{\infty})),
\]
with $s$ odd.  \qed
\end{prop}

We will also need

\begin{lem}
\label{thm:235} For $k\ge 2$
the connecting homomorphisms $\delta_{1}$ and $\delta_{2}$ satisfy
the following congruences modulo the ideal $(2, v_{1}^{2})$:
\begin{align*}
\delta_{1}(v_{2}^{s 2^{k}}/v_{1}^{2^{k}}) &\equiv
v_{2}^{(s-1)2^{k}}[t_{1}^{2^{k+1}}]  \\
\delta_{2}v_{1}^{2^{k}}[t_{1}] &\equiv 0 \\
\delta_{2}v_{1}^{2^{k}}[t_{1}^{2}] &\equiv 0 \\
\delta_{2}\delta_{1}(v_{2}^{s 2^{k}}/v_{1}^{2^{k}})  &\equiv
v_{2}^{(s-1)2^{k}}[t_{1}^{2^{k}}\vert t_{1}^{2^{k}}].
\end{align*}
\end{lem}

\begin{pf}
This is a straightforward (and long-known) computation using the
structure formulae
\begin{align*}
\eta_{R}(v_{1}) &= v_{1} + 2 t_{1} \\
\eta_{R}(v_{2}) &\equiv v_{2} +v_{1}t_{1}^{2} + v_{1}^{2}t_{1}\mod 2.
\end{align*}
The assertion about $\delta_{1}$ is easy to check.  The structure
formulae imply that $\eta_{R}(v_{1}^{2})\equiv v_{1}^{2}$ modulo $4$, so
one may work modulo $(4,v_{1}^{2})$ when computing $\delta_{2}$.  The
terms $v_{1}^{2^{k}}[t_{1}]$ and $v_{1}^{2^{k}}[t_{1}^{2}]$ are
already in this ideal, giving the first two assertions about
$\delta_{2}$.  The last makes use of the congruences
\begin{align*}
\delta_{2}[t_{1}^{2^{k+1}}] &\equiv [t_{1}^{2^{k}}\vert t_{1}^{2^{k}}]
\mod 2 \\
\eta_{R}v_{2}^{4i} &\equiv v_{2}^{4i} \mod (4, v_{1}^{2}).
\end{align*}
Since $s$ is odd and $k\ge 1$, $(s-1)2^{k}$ is divisible by $4$.  This
means that
\[
\delta_{2}(v_{2}^{(s-1)2^{k}}[t_{1}^{2^{k+1}}]) \equiv
v_{2}^{(s-1)2^{k}}\delta_{2}([t_{1}^{2^{k+1}}]) \equiv
v_{2}^{(s-1)2^{k}}[t_{1}^{2^{k}}\vert t_{1}^{2^{k}}]\mod (2, v_{1}^{2}).
\]
This completes the proof.
\end{pf}

\subsubsection{The proof of Proposition~\ref{thm:237}}
\label{sec:map-classical-adams}

Given Lemma~\ref{thm:236} and Proposition~\ref{thm:231},
Proposition~\ref{thm:237} is an immediate consequence of the two
results below.

\begin{prop}
\label{thm:144} For $j>1$, the map
\[
\ext^{1,2^{j+1}}_{\mustarmu}(\mustar,\mustar) \to \ext^{2,2^{j+1}}_{\sta}(\Z/2,\Z/2).
\]
is given by
\begin{align*}
v_{1}^{2^{j}-1}[t_{1}] & \mapsto 0 \\
v_{1}^{2^{j}-2}[t_{1}^{2}] & \mapsto 0 \\
\delta_{1}(v_{2}^{s 2^{k}}/v_{1}^{2^{k}}) & \mapsto 0 \qquad (s > 1)\\
\delta_{1}(v_{2}^{2^{j-1}}/v_{1}^{2^{j-1}}) & \mapsto
h_{j}^{2}.
\end{align*}
\end{prop}

\begin{pf}
This follows directly from Lemma~\ref{thm:235} and the fact that the
map from $\mustarmu$ to the dual Steenrod algebra 
given by
\begin{align*}
v_{i} &\mapsto 0 \\
t_{i} &\mapsto \chi(\xi_{i})^{2}.
\end{align*}
\end{pf}

\begin{rem}
\label{rem:54}
The map from the Adams-Novikov $E_{2}$-term to the classical Adams
$E_{2}$-term has been completely determined for $s\le 2$ and all $t$.
For a comprehensive discussion, the reader is referred
to~\cite[Chapter~5]{Rav:MUnew}.   
\end{rem}

We next turn to the second map in Proposition~\ref{thm:237}.  When
$j>2$ the action of $\ze$ on $A_{2^{j+1}}$ is trivial, and the group
$H^{1}(\ze;A_{2^{j+1}}/(\pi))$ is cyclic of order two, generated by
the cohomology class of the cocycle whose value on $\gamma$ is
$u^{2^{j}}$.  Let us denote this class $\theClass$.

\begin{prop}
\label{thm:234} For $j>2$, the map
\[
\ext^{1,2^{j+1}}_{\mustarmu}(\mustar,\mustar/(2))\to
H^{1}(\ze;A_{2^{j+1}}/(\pi)).
\]
is given by
\begin{align*}
v_{1}^{2^{j}-1}[t_{1}] &\mapsto 0 \\
v_{1}^{2^{j}-2}[t_{1}^{2}] &\mapsto 0 \\
\delta_{1}(v_{2}^{s 2^{k}}/v_{1}^{2^{k}}) & \mapsto 0 \qquad (s > 1)\\
\delta_{1}(v_{2}^{2^{j-1}}/v_{1}^{2^{j-1}}) & \mapsto \theClass.
\end{align*}
\end{prop}

\begin{pf}
Let $\nu$ be the valuation on $A$ normalized so that $\nu(\pi)=1$.
Since $\nu(2)=4$, 
\[
[2]_{\ltgroup^{h}}(x) \equiv x^{2^{4}} + \cdots  \mod \pi,
\]
and $v_{1}$ and $v_{2}$ both map to zero in $A_{\ast}/(\pi)$.
This gives the first line, and makes the second a consequence of
Lemma~\ref{thm:235}.   Lemma~\ref{thm:235} also gives the identity
\[
\delta_{1}(v_{2}^{2^{j}}/v_{1}^{2^{j}}) = [t_{1}^{2^{j}}],
\]
so to determine the image of $\delta_{1}(v_{2}^{2^{j}}/v_{1}^{2^{j}})$
we need to work out the image of $t_{1}$ under the map of
Hopf-algebroids
\begin{equation}
\label{eq:158}
(\mustar,\mustarmu) \to (A_{\ast}/(\pi), C(\ze,A_{\ast}/(\pi))).  
\end{equation}
As explained at the end of~\ref{sec:preliminaries}, the
element $t_{1}$ occurs as the coefficient of $x^{2}$ the isomorphism
$\eta_{R}^{\ast}\funiv^{(2)}\to \funiv^{(2)}$, inverse to the
universal strict isomorphism.  Since we have chosen a $2$-typical
coordinate on $\ltgroup^{h}$, the element $t_{1}$ is therefore sent,
under the map of Hopf-algebroids~\eqref{eq:158}, to the $1$-cocycle on
$\ze$ whose value on $\gamma$ is the coefficient of $x^{2}$ in the
inverse of the power series $[\zeta](x)$.  By~\eqref{eq:110} and the
formula for strictification~\eqref{eq:160} this is $-u^{2^{j}}\equiv
u^{2^{j}}$ modulo $(\pi)$.
\end{pf}

\subsection{Addendum} 
\label{sec:addendum} 

Lemma~\ref{thm:236} and Proposition~\ref{thm:231} do not quite appear
literature in a readily accessible form, and the purpose of this
addendum is to outline their proofs, explaining how the key points
can be read off from the results of~\cite{MRW}
and~\cite{shimomura81:_novik_rm_ext_sp2_at_prime}.  To conform with
the notation of these references, we will use $BP$ rather that $MU$ in
this section.

The paper~\cite{MRW} introduced the chromatic approach to studying
the groups 
\[
\ext_{\bpstarbp}^{s,\ast}(\bpstar,\bpstar).
\]
The computation begins with the fact that for $s>0$ one has
\[
\ext_{\bpstarbp}^{s,\ast}(\bpstar,\bpstar)\otimes \Q =0.
\]
This means that the connecting homomorphism
\[
\ext^{1,2^{j+1}}_{\bpstarbp}(\bpstar,\bpstar/2^{\infty})\to 
\ext^{2,2^{j+1}}_{\bpstarbp}(\bpstar, \bpstar)
\]
is an isomorphism.    The assertion of Proposition~\ref{thm:231} is
that the map 
\[
\ext^{1,2^{j+1}}_{\bpstarbp}(\bpstar, \bpstar/2) \to 
\ext^{1,2^{j+1}}_{\bpstarbp}(\bpstar, \bpstar/2^{\infty})
\]
induced by the inclusion 
\[
\bpstar/2 \xrightarrow{1/2}{} \bpstar/2^{\infty}
\]
is surjective, and that the left had group is spanned by the elements
listed.  Continuing with the chromatic approach, one is led to the
following diagram (in which, to manage the size, we have abbreviated
$\ext^{s,t}_{\bpstarbp}(\bpstar,M)$ to $\ext^{s,t}(M)$)
\begin{equation}
\label{eq:149}
\xymatrix@C=1em{
*++{\ext^{0,2^{j+1}}(\bpstar/(2,v_{1}^{\infty}))}  \ar[r]^{\delta_{1}}\xyhookar[d]  & *++{\ext^{1,2^{j+1}}(\bpstar/(2))}   \ar[r]\ar[d]   & *++{\ext^{1,2^{j+1}}(v_{1}^{-1}\bpstar/(2))}    \ar[d]\\
*++{\ext^{0,2^{j+1}}(\bpstar/(2^{\infty},v_{1}^{\infty}))}   \ar[r]        &  *++{\ext^{1,2^{j+1}}(\bpstar/(2^{\infty}))}     \ar[r]         &*++{\ext^{1,2^{j+1}}(v_{1}^{-1}\bpstar/(2^{\infty}))\mathrlap{\ .}  }
}
\end{equation}
The rightmost column is analyzed using the Miller-Ravenel change of
rings theorem (see~\cite{MilRav} or~\cite[Chapter 6, \S1]{Rav:MUnew})
which identifies it with the map
\[
H^{1}(\Z_{2}^{\times};\zlt[v_{1}^{\pm1}]/(2)) \to 
H^{1}(\Z_{2}^{\times};\zlt[v_{1}^{\pm1}]/(2^{\infty})) 
\]
in which $\lambda\in\Z_{2}^{\times}$ acts on $v_{1}$ with eigenvalue
$\lambda$.  This is easily calculated, and one finds that the map is
indeed surjective, and that
$\ext^{1,2^{j+1}}_{\bpstarbp}(\bpstar,v_{1}^{-1}\bpstar/(2))$ has
dimension $2$, with basis the image of $v_{1}^{2^{j-1}}[t_{1}]$ and
$v_{1}^{2^{j-2}}[t_{1}^{2}]$.  This reduces Proposition~\ref{thm:231}
to the assertion that the left vertical arrow is surjective (hence an
isomorphism), and that the upper left group has a basis consisting of
the elements of the form $v_{2}^{s 2^{k}}/v_{1}^{2^{k}}$.  For this
one first appeals to the invariant prime ideal theorem
(\cite{Morava,Land:Hom}, or see~\cite[Theorem 4.3.2]{Rav:MUnew}) for
the fact that
\[
\ext^{0,\ast}_{\bpstarbp}(\bpstar,\bpstar/(2,v_{1})) = \Z/2[v_{2}].
\]
It follows that any invariant element in $\bpstarbp/(2^{\infty},
v_{1}^{\infty})$ has the form 
\begin{equation}
\label{eq:148}
\frac{v_{2}^{s 2^{k}} + r}{v_{1}^{\ell}2^{i+1}}
\end{equation}
with $r\in (2,v_{1})\subset \bpstar$.  We now come to the key point.
It turns out that a necessary (but not sufficient) condition that such
an element to be invariant is that the indices satisfy the inequality
\begin{equation}
\label{eq:150}
\ell \le 2^{k-i}+ 2^{k-i-1}.
\end{equation}
This can be extracted from the stronger conditions of
\cite[Theorem~3.3]{shimomura81:_novik_rm_ext_sp2_at_prime}, in which
the symbol $x_{n}$ is an explicitly defined element, congruent to
$v_{2}^{2^{n}}$ modulo $(2,v_{1})$ and $y_{i}$ is an explicitly
defined element congruent to $v_{1}^{2^{i}}$ modulo $2$.
From~\eqref{eq:150} it follows that for an element of the
form~\eqref{eq:148} to be invariant and have degree $2^{j+1}$, the
numbers $i$, $j$, $k$, and $\ell$ must satisfy
\[
6 s 2^{k} -2(2^{k-i}+2^{k-i-1})\le 2^{j+1} \le  6 s 2^{k}.
\]
Expanding, and dividing both sides by $2^{k+1}$ gives
\[
3 s -2^{-i} - 2^{-i-1} \le 2^{j-k} \le  3s.
\]
Since $s\ge 1$ and $i\ge 0$, one has 
\[
3 s -2^{-i} - 2^{-i-1} \ge 3s - 3/2 > 1,
\]
and so $k< j$.  This implies $2^{j-k}$ is even, and so must equal
$3s-1$.  This in turn means that $2^{-i} + 2^{-i-1}>1$, and so $i$
must be $0$.

It thus follows from the inequality~\eqref{eq:150} that the invariant
elements of degree $2^{j+1}$ in $\bpstar/(2^{\infty},v_{1}^{\infty})$
have the form
\[
\frac{v_{2}^{s 2^{k}} + r}{v_{1}^{2^{k}}2}.
\]
Thus the left vertical map in~\eqref{eq:149} is surjective.   Since 
\[
\frac{v_{2}^{s 2^{k}}}{v_{1}^{2^{k}}2}
\]
is already invariant, a simple induction shows that the elements
stated form a basis.

\appendix

\section{The category of equivariant orthogonal spectra}
\label{sec:equiv-orth-spectra}

In this appendix we recall the definition and some basic properties of
the theory of equivariant orthogonal spectra.   For further details
and references the reader is referred to Mandell-May~\cite{MR1922205}
and to Mandell-May-Schwede-Shipley~\cite{MR1806878}.

One of the reasons we have chosen to use equivariant orthogonal
spectra is that it has many convenient category theoretic properties.
These are independent of the homotopy theory of equivariant
orthogonal spectra, and so we make two passes through the theory, one
focused on the category theory, and the other on the homotopy theory.

Our main new innovation is the theory of the norm (\S\ref{sec:norm}).
Most of the category theoretic aspects apply to any symmetric monoidal
category, and things work out much cleaner at that level of
generality.

\subsection{Category theory preliminaries}
\label{sec:categ-theory-prel}

\subsubsection{Symmetric monoidal categories}
\label{sec:symm-mono-categ}

\begin{defin}
\label{def:36}
A {\em symmetric monoidal category} is a category $\cat V$ equipped
with a functor 
\[
\otimes:\cat V\times\cat V\to \cat V,
\]
a unit object $\one\in \cat V$, a natural associativity isomorphism 
\[
a_{ABC}:(A\otimes B)\otimes C \approx A\otimes (B\otimes C)
\]
a natural commutativity isomorphism 
\[
s_{AB}:A\otimes B\approx B\otimes A
\]
and a unit isomorphism 
\[
\iota_{A}:\one\otimes A\approx A.
\]
This data is required to satisfy the associative and commutative
coherence axioms, as well as the strict symmetry axiom.
\end{defin}

The two coherence axioms express that all of the ways one might get
from one iterated tensor product to another using the associativity
and commutativity transformations coincide.  The strict symmetry axiom
is that the square of the commutativity transformation is the identity
map.  See~\cite{MR1712872}, or Borceux~\cite[\S6.1]{borceux2}.

Even though it requires six pieces of data to specify a symmetric
monoidal category we will usually indicate one with a triple $\cat
V=(\cat V_{0},\otimes,\one)$.

A symmetric monoidal category is {\em closed} if for each $A$, the
functor $A\otimes (\slot)$ has a right adjoint $(\slot)^{A}$, which
one can think of as an ``internal hom.''  Note that
\[
\cat V(\one, X^{A}) \approx \cat V(A,X)
\]
so that one can recover the usual hom from the internal hom.

\subsubsection{Sifted colimits, commutative and associative algebras}
\label{sec:commutative-algebras} 

In a closed symmetric monoidal category, the monoidal product commutes
with colimits in each variable.  It follows easily that the iterated
monoidal product
\[
X\mapsto X^{\otimes n}
\]
commutes with all colimits over indexing categories $I$ for which the
diagonal $I\to I^{n}$ is {\em final} in the sense of \cite[\S
IX.3]{MR1712872}.  If $I\to I\times I$ is final, then for all $n\ge
2$, $I\to I^{n}$ is also final.

\begin{defin}
A category $I$ is {\em sifted} the diagonal embedding $I\to I\times I$
is final.  
\end{defin}

Equivalently (see~\cite[15.2 (c)]{MR0327863}, or \cite[Theorem~2.15]{MR2757312}), a small category $I$ is sifted if
and only if the formation of colimits over $I$ in sets commutes with
finite products.

\begin{defin}
A {\em sifted colimit} is a colimit over a sifted category.
\end{defin}

Examples of sifted colimits include reflexive coequalizers and
directed colimits.  In fact the class of sifted colimits is
essentially the smallest class of colimits containing reflexive
coequalizers and directed colimits.  See, for
example~\cite{adamek:_what_are_sifted_colim, MR0327863}.

Let $\cat V=(\cat V_{0},\otimes,\one)$ be a closed symmetric
monoidal category.  

\begin{defin}
\label{def:30} An {\em associative algebra} in $\cat V$ is an object
$A$ equipped with a multiplication map $A\otimes A\to A$ which is
unital and associative.  A {\em commutative algebra} is an associative
algebra for which the multiplication map is commutative.
\end{defin}

The categories of associative and commutative algebras (and algebra
maps) in $\cat V$ are denoted $\asscat V$ and $\commcat V$,
respectively.

The following straightforward result holds more generally for algebras
over any operad.  The existence of colimits in the algebra categories
is proved by expressing any algebra as a reflexive coequalizer of a
diagram of free algebras.  There is an even more general result for
algebras over a triple~\cite[Proposition~4.3.1]{borceux2}

\begin{prop}
\label{thm:1}
Suppose that $\cat V$ is a closed symmetric monoidal category.  The
forgetful functors 
\begin{align*}
\asscat V &\to \cat V \\
\commcat V &\to \cat V \\
\end{align*}
create limits.  If $\cat V$ is cocomplete these functors have left
adjoints
\begin{align*}
X &\mapsto T(X) = \coprod_{n\ge 0} X^{\otimes n} \\
X &\mapsto \sym (X) = \coprod_{n\ge 0} X^{\otimes n}/\Sigma_{n},
\end{align*}
the categories $\asscat V$ and $\commcat V$ are cocomplete, and the
``free'' functors above commute with all sifted colimits.   \qed
\end{prop}

A {\em left module} over an associative algebra $A$ is an object $M$
equipped with a unital and associative left multiplication 
\[
A\otimes M\to M.
\]
Similarly a {\em right module} is an object $N$ equipped with a
unital, associative right multiplication $N\otimes A\to N$.   Given a
left $A$-module $M$ and a right $A$-module $N$ one defines
$N\tensove{A} M$ by the (reflexive) coequalizer
\[
N\otimes A\otimes M\rightrightarrows N\otimes M\rightarrow N\tensove{A}M.
\]
When $A$ is commutative, a left $A$-module can be regarded as a right
$A$-module by the action 
\[
M\otimes A\xrightarrow{\text{flip}}{} A\otimes M \to M.
\]
Using this, the formation $M\tensove{A}N$ makes the category of left
$A$-modules into a symmetric monoidal category.

\subsubsection{Enriched categories}
\label{sec:enriched-categories}

In this section we briefly describe the basic notions of enriched
categories..  The reader is referred to~\cite{Kelly:Enr}
or~\cite[Ch.~6]{borceux2} for further details.

Suppose that $\cat V = (\cat V_{0},\otimes,\one)$ is a symmetric monoidal
category.  

\begin{defin}
\label{def:33} A $\cat V$-category $\cat C$ consists of a collection
$\ob \cat C$ called the {\em objects} of $\cat C$, for each pair
$X,Y\in\ob\cat C$ a {\em morphism object} $\cat C(X,Y)\in \ob \cat
V_{0}$, for each $X$ an {\em identity morphism} $\one \to \cat C(X,X)$
and for each triple $X,Y,Z$ of objects of $\cat C$ a {\em composition
law}
\[
\cat C(Y,Z)\otimes \cat C(X,Y)\to \cat C(X,Z).
\]
This data is required to satisfy the evident unit and associativity
properties.
\end{defin}

As is customary, we write $X\in \cat C$ rather than $X\in\ob\cat C$.
Most of the notions of ordinary category theory carry through in the
context of enriched categories, once formulated without reference to
``elements'' of mapping objects.  For example a functor $F:\cat C\to
\cat D$ of $\cat V$-categories consists of a function
\[
F:\ob\cat C\to\ob \cat D
\]
and for each pair of objects $X,Y\in\cat C$ a $\cat V$-morphism 
\[
F:\cat C(X,Y)\to \cat D(FX,FY)
\]
compatible with the unit and composition.  A natural transformation
between two functors $F$ and $G$ is a function assigning to each
$X\in\cat C$ a map $T_{X}:\one\to \cat D(FX,GX)$ which for every $X,Y$
makes the diagram 
\begin{equation}
\label{eq:17}
\xymatrix{
\cat C(X,Y) \ar[r]^-{T_{Y}\otimes F} \ar[d]^{G\otimes T_{X}} & \cat
D(FY,GY)\otimes
\cat D(FX,FY) \ar[d] \\
\cat D(GX, GY)\otimes \cat D(FX,GX) \ar[r] & \cat D(FX, GY)
}
\end{equation}
commute.  

There is an ordinary category $\cat C_{0}$ underlying the enriched
category $\cat C$.  The objects of $\cat C_{0}$ are the objects of
$\cat C$, and one defines
\[
\cat C_{0}(X,Y) = \cat V_{0}(\one, \cat C(X,Y)).
\]
If $\cat V$ itself underlies a $\cat W$-enriched category, then any
$\cat V$-category $\cat C$ has an underlying $\cat W$-category, whose
underlying ordinary category is $\cat C_{0}$.

When $\cat V$ is a closed symmetric monoidal category, the internal
hom defines an enrichment of $\cat V$ over itself, with underlying
category $\cat V_{0}$.

When $\cat V$ is closed, a natural transformation $F\to G$ can be
described as a map
\[
\one\to \prod_{X\in \cat C} \cat D(FX,GX)
\]
which equalizes the two arrows 
\begin{equation}
\label{eq:1}
\prod_{X\in\cat C}\cat D(FX,GX)\rightrightarrows \prod_{X,Y\in\cat
C}\cat D(FX,GY)^{\cat C(X,Y)}.
\end{equation}
describing the two ways of going around~\eqref{eq:17}.

We will write $\vcat V$ for the $2$-category of $\cat V$-categories,
and denote the category of enriched functors $\cat C\to \cat D$ as
$\vcat V(\cat C,\cat D)_{0}$.  When $\cat V$ is closed and contains
products indexed by the collection of pairs of objects of $\cat C$,
the category $\vcat V(\cat C,\cat D)_{0}$ underlies an enriched
category $\vcat V(\cat C,\cat D)$ in which the object of natural
transformations from $F$ to $G$ is given by the equalizer
of~\eqref{eq:1}.

\subsection{Equivariant orthogonal spectra}

\subsubsection{Equivariant spaces}\label{sec:equivariant-spaces}

Let $\spaces$ be the category of pointed, compactly generated weak
Hausdorff spaces (in the sense of \cite{MR0251719}).  The category
$\spaces$ is symmetric monoidal under the smash product, with unit
$S^{0}$.  A {\em topological} category is a category enriched over
$(\spaces, \wedge, S^{0})$.

\begin{rem}
\label{rem:66} 
Working with compactly generated {\em weak Hausdorff
spaces} has many benefits, but it does create some technical issue.
Colimits are computed by forming the colimit in topological spaces,
replacing the topology by the compactly generated topology, and then
forming the universal quotient which is weak Hausdorff.  This last
step can alter the underlying point set.  It does not, however,  in
the case of pushouts along closed inclusion.   More precisely, given a
pushout diagram
\[
\xymatrix{
A  \ar[r]\ar[d]  & X
\ar[d] \\
B  \ar[r]        & Y
}
\]
of topological spaces in which $A\to X$ is a closed inclusion, if $A$,
$X$, and $B$ are compactly generated and weak Hausdorff then so is
$Y$.  This follows from~\cite[Proposition~2.5]{MR0251719} and the
remark about adjunction spaces immediately preceding its statement.
Among other things this means that the smash product of two compactly
generated weak Hausdorff spaces can be computed as the smash product
of the underlying compactly generated spaces.
\end{rem}

Now suppose that $G$ is a group.  Let $(\ugspaces G,
\wedge, S^{0})$ be the topological symmetric monoidal category of
pointed spaces with a left $G$-action and spaces of equivariant maps.
With this structure $\ugspaces{G}$ is a closed symmetric monoidal
category, with internal mapping spaces $\gspaces G(X,Y)=Y^{X}$ given
by the space of non-equivariant maps, with the conjugation action of
$G$.  

A word about notation.  The expression ``category of $G$-spaces'' can
reasonably refer to three objects, depending on what is meant by a
map.  It can be an ordinary category, a category enriched over
topological spaces, or a category in which the hom objects are the
$G$-spaces of non-equivariant maps.  As indicated above we will use
$\gspaces{G}$ to denote the category enriched over $G$-spaces, with
$\gspaces G(X,Y)$ denoting the $G$-space of non-equivariant maps, and
$\ugspaces{G}$ for the {\em topological} category of $G$-spaces, and
spaces of equivariant maps.

We will be making use of categories enriched over $\ugspaces{G}$.  As
in~\cite{MR1922205}, we will refer to them as {\em topological
$G$-categories} (or just {\em $G$-categories} for short).  Let $\gcat
G$ denote the collection of topological $G$-categories, and write
$\gcat G(\cat C, \cat D)$ for the enriched category of functors and
left $G$-spaces of natural transformations.  The symbol $\gcat G(\cat
C,\cat D)^{G}$ will denote the topological category of functors and
spaces of equivariant natural transformations.

\subsubsection{Change of group}
\label{sec:change-group-1}

Suppose that $H\subset G$ is a subgroup.  The restriction functor 
\[
\ugspaces{G}\to \ugspaces{H}
\]
has continuous left and a right adjoints given by
\begin{align*}
Y &\mapsto G_{+}\smashove{H} Y \\
Y &\mapsto \ugspaces{H}(G_{+},Y).
\end{align*}
These two constructions are basic examples of indexed monoidal
products (see \S\ref{sec:index-mono-prod}).  Because $\ugspaces{G}$ is
pointed there is a canonical equivariant map 
\[
G_{+}\smashove{H} Y \to \ugspaces{H}(G_{+},Y).
\]

\subsubsection{The basic indexing categories}
\label{sec:basic-index-categ}

For a real orthogonal representation $V$ of $G$ let $O(V)$ be the
orthogonal group of non-equivariant linear isometric maps $V\to V$.
The group $G$ acts on $O(V)$ by conjugation, and the group of fixed
points is the orthogonal group of equivariant maps.  Given orthogonal
representations $V$ and $W$, we define $O(V,W)$ to be the Stiefel
manifold of linear isometric embeddings of $V$ into $W$, with the
conjugation action of $G$.  The $G$ fixed points in $O(V,W)$ are the
equivariant orthogonal embeddings.  The group $O(W)$ acts transitively
on $O(V,W)$ on the left.  A choice of embedding $V\to W$ identifies
$O(V,W)$ with the homogeneous space $O(W)/O(W-V)$.

\begin{defin}
\label{def:37} The category $\igcat G$ is the topological $G$-category
whose objects are finite dimensional real orthogonal representations
of $G$,  and with morphism $G$-space $\igcat G(V,W)$ the Thom complex
\[
\igcat G(V,W)=\thom(O(V,W);W-V)
\]
of the ``complementary bundle'' $W-V$ over $O(V,W)$.   
\end{defin}

We will denote the topological category underlying $\igcat{G}$ with
the symbol $\uigcat{G}$.  Thus $\uigcat{G}(V,W)=\igcat{G}(V,W)^{G}$.

The $G$-space $\igcat G(V,W)$ can be thought of as the topologically
indexed wedge
\[
\bigvee_{V\to W} S^{W-V}.
\]
When $\dim W<\dim V$ it reduces to the one point space $\ast$.  When
$\dim W\ge \dim V$ one can get a convenient description by choosing an
orthogonal $G$-representation $U$ with $\dim U+\dim V=\dim W$ (for
example the trivial representation).  With this choice one has
\[
\igcat G(V,W) \approx O(V\oplus U,W)_{+}\smashove{O(U)}S^{U}.
\]

The fixed point space $\igcat{G}(V,W)^{G}$ is given by 
\begin{equation}
\label{eq:139}
\igcat{G}(V,W)^{G} = 
\icat(\fixg{V},\fixg{W})\wedge O(V^{\perp},W^{\perp})^{G}_{+},
\end{equation}
in which $V^{G}$ denotes space of invariant vectors in $V$, 
and $V^{\perp}$ its orthogonal complement.   The space
$O(V^{\perp},W^{\perp})^{G}$ in turn decomposes into the product 
\[
\prod_{\alpha}O(V_{\alpha},W_{\alpha})
\]
in which $\alpha$ is running through the set of non-trivial
irreducible representations of $G$,  and $V_{\alpha}\subset V$ and
$W_{\alpha}\subset W$  indicate the $\alpha$-isotypical parts.
 
When $G$ is the trivial group we will denote the category $\igcat{G}$
simply by $\icat$.  For any $G$ there is an inclusion $\icat\subset
\igcat{G}$ identifying $\icat$ with the full subcategory of objects
with trivial $G$-action.  There is also a forgetful functor
$\igcat{G}\to \icat$ which refines in the evident manner to a functor
from $\igcat{G}$ to the $G$-category of objects in $\icat$ equipped
with a $G$-action.  One can easily check that this is an equivalence.
For later reference, we single this statement out.  

\begin{prop}
\label{thm:3} 
The forgetful functor described above gives an
equivalence of $\igcat{G}$ with the topological $G$-category of
objects in $\icat$ equipped with a $G$-action.  Passage to fixed
points gives an equivalence of $\uigcat{G}$ with the topological
category of objects in $\icat$ equipped with a $G$-action.  \qed
\end{prop}

Proposition~\ref{thm:3} plays an important role in establishing one of
the basic properties of the norm (Proposition~\ref{thm:12}).  

\subsubsection{Orthogonal spectra}
\label{sec:orthogonal-spectra}

\begin{defin}
\label{def:5}
An {\em orthogonal $G$-spectrum} is a functor 
\[
\igcat G\to \gspaces{G}
\]
of topological $G$-categories.
\end{defin}

Informally, an orthogonal spectrum $X$ consists of a collection of
spaces $\{X_{V} \}$, and for each $V\to W$ a 
\[
S^{W-V}\wedge X_{V}\to X_{W}.
\]
These maps are required to be compatible with composition in $\igcat
G$, the action of $G$, and to vary continuously in $V\to W$.  More
formally, one has equivariant maps
\[
\thom(O(V,W);S^{W-V})\wedge X_{V}\to X_{W}
\]
compatible with composition.

\begin{defin}
\label{def:34} The topological $G$-category of {\em orthogonal
$G$-spectra} is the category
\[
\gspectra G = \gcat G(\igcat G,\gspaces G).
\]
The {\em (topological) category of $G$-spectra} is 
\[
\ugspectra G = \gcat G(\igcat G,\gspaces G)^{G}.
\]
\end{defin}

We will use the notation 
\[
\spectra=\gcat G(\icat,\spaces)
\]
to denote the category $\gspectra{G}$ for the case of the trivial group.

The ($G$-)category of orthogonal $G$-spectra is complete and
cocomplete (in the sense of enriched categories).  Both limits and
colimits in $\ugspectra{G}$ are computed objectwise:
\begin{align*}
(\varinjlim X^{\alpha})_{V} &= \varinjlim X^{\alpha}_{V}\\
(\varprojlim X^{\alpha})_{V} &= \varprojlim X^{\alpha}_{V}.
\end{align*}

Certain orthogonal $G$-spectra play a fundamental role.  For $V\in
\igcat G$ let
\[
S^{-V}:\igcat G\to \gspaces{G}
\]
be the functor co-represented by $V$.  By the Yoneda lemma
\[
\gspectra{G}(S^{-V},X) = X_{V}.
\]
For a pointed $G$-space $A$ let $S^{-V}\wedge A$ be the orthogonal
$G$-spectrum with
\[
\left(S^{-V}\wedge A \right)_{W}= \left(S^{-V} \right)_{W}\wedge A.
\]
Again, by Yoneda,
\[
\gspectra{G}(S^{-V}\wedge A,X) = \gspaces{G}(A, X_{V}).
\]
It also follows from the Yoneda lemma that every $X$ is functorially
expressed as a reflexive coequalizer
\begin{equation}\label{eq:32}
\bigvee_{V,W} S^{-W}\wedge
\igcat G(V,W)\wedge X_{V}\rightrightarrows
\bigvee_{V} S^{-V}\wedge X_{V} \to X.
\end{equation}
We call this the {\em tautological presentation} of $X$ and for ease
of typesetting, sometimes indicate it as
\begin{equation}\label{eq:33}
X=\varinjlim_{V} S^{-V}\wedge X_{V}.
\end{equation}

\subsubsection{Smash product}
\label{sec:smash-prod-orth}

The symmetric mon\-oid\-al structures on $\igcat G$ and $\gspaces G$
combine to give $\ugspectra G$ a symmetric monoidal structure (the Day
convolution), denoted $\wedge$.  The smash product of two orthogonal
$G$-spectra $X$ and $Y$ is defined to be the left Kan extension of
\[
(V,W)\mapsto X_{V}\wedge Y_{W}: \igcat{G}\times \igcat{G}\to
\gspaces G
\]
along the map 
\[
\igcat{G}\times\igcat{G}\to \igcat{G}
\]
sending $(V,W)$ to $V\oplus W$.  The smash product is thus
characterized by the fact that it commutes with enriched colimits in
both variables, and satisfies 
\[
S^{-V}\wedge S^{-W}= S^{-(V\oplus W)}.
\]
In terms of the tautological presentations
\begin{align*}
X=\varinjlim_{V}S^{-V}\wedge X_{V} \\
Y=\varinjlim_{W}S^{-W}\wedge Y_{W}
\end{align*}
one has
\begin{align*}
X\wedge Y &= \varinjlim_{V}S^{-V}\wedge X_{V}\wedge Y \\
          &= \varinjlim_{V}S^{-V}\wedge X_{V}\wedge \varinjlim_{W}S^{-W}\wedge Y_{W} \\
          &= \varinjlim_{V,W}S^{-V\oplus W}\wedge X_{V}\wedge Y_{W}.
\end{align*}
The above expression is, of course, an abbreviation for the
reflexive coequalizer diagram
\begin{multline*}
\bigvee_{\substack{V_{0},V_{1}, \\ W_{0},W_{1}}}
\igcat G(V_{0},V_{1})\wedge \igcat G(W_{0},W_{1})\wedge
S^{-V_{1}\oplus W_{1}}\wedge
X_{V_{0}}\wedge Y_{W_{0}}\\ \rightrightarrows
\bigvee_{V,W}S^{-V\oplus W}\wedge X_{V}\wedge Y_{W}.
\end{multline*}

\begin{prop}
\label{thm:63}
The category $\ugspectra G$ is a closed symmetric monoidal category
with respect to $\wedge$. \qed
\end{prop}

Smashing the tautological presentation of a general spectrum $X$ with
$S^{-V}$ gives a presentation of $S^{-V}\wedge X$ as a (reflexive)
coequalizer
\[
\bigvee_{W_{0},W_{1}} \igcat{G}(W_{0},W_{1})\wedge S^{-V\oplus W_{1}}\wedge
X_{W_{0}} \rightrightarrows
\bigvee_{W} S^{-V\oplus W}\wedge X_{W} \to S^{-V}\wedge X.
\]
This is {\em not} the tautological presentation of $S^{-V}\wedge X$, but
from it, one can read off the formula of the following lemma

\begin{lem}
\label{thm:111}
If $\dim W <\dim V$, then $(S^{-V}\wedge X)_{W}=\ast$.  If $\dim W\ge
\dim V$, then  there is a natural isomorphism of $G$-spaces
\[
(S^{-V}\wedge X)_{W} \approx O(V\oplus U,W)_{+}\smashove{O(U)} X_{U}
\]
where $U$ is any orthogonal $G$-representation with 
\[
\dim U+\dim V=\dim W.
\]
\qed
\end{lem}

\subsubsection{Variations on the indexing category}

There is a lot of flexibility in defining $\ugspectra G$.  In this
section we describe a variation which is especially convenient for
certain category theoretical properties, and will be used in our
construction of the norm.  We learned of the result below from Lars
Hesselholt and Mark Hovey.  It is due to Mandell-May
(\cite[Lemma~V.1.5]{MR1922205}).

\begin{prop}\label{thm:143}
Let $i:\icat\to \igcat G$ be the inclusion of the full subcategory of
trivial $G$-representations.  The functors 
\[
i^{\ast}:\gcat G(\igcat G,\gspaces G) \to \gcat G(\igcat, \gspaces G)
\]
and 
\[
i_{!}:\gcat G(\igcat,\gspaces G) \to \gcat G(\igcat G, \gspaces G)
\]
given by restriction and left Kan extension along $i$ are inverse
equivalences of enriched symmetric monoidal categories.
\end{prop}

In other words the symmetric monoidal (topological) category
$\ugspectra G$ can simply be regarded as the symmetric monoidal
(topological) category of objects in $\spectra$ equipped with a
$G$-action.

The proof of Proposition~\ref{thm:143} requires a simple technical
lemma (\cite[Lemma~V.1.1]{MR1922205}).

\begin{lem}
\label{thm:218}
Suppose that $V$ and $W$ are orthogonal $G$-represen\-ta\-tions with $\dim
V=\dim W$.  Then for any $U$
\[
O(V,U)\times O(V) \times O(W,V)\rightrightarrows
O(V,U) \times O(W,V)\to O(W,U)
\]
is a (reflexive) coequalizer in $\ugspaces G$.
\end{lem}

\begin{pf}
Since the forgetful functor $\ugspaces G\to \spaces$ preserves
colimits and reflects isomorphisms, it suffices to prove the result in
$\spaces$, where it is obvious, since the coequalizer diagram can be
split by choosing an orthogonal (non-equivariant) isomorphism of $V$
with $W$.
\end{pf}

\begin{pf*}{Proof of Proposition~\ref{thm:143}}
Since $i:\icat\to \igcat{G}$ is fully faithful, the left Kan extension
$i_{!}$ is fully faithful (see for example,
\cite[Corollary~X.3]{MR1712872}).  To show that it is essentially
surjective, let $W\in \igcat{G}$ be any object, and let $V\in \icat$
be a vector space of the same dimension as $W$.  Define $X$ by the
coequalizer
\[
(O(W,V) \times O(V))_{+}\wedge S^{-V}\rightrightarrows
O(W,V)_{+}\wedge S^{-V} \to X.
\]
Since $\igcat{G}(W,V)= O(W,V)$, $i_{\ast}X$ is given by the
coequalizer of
\[
(\igcat{G}(W,V) \times \igcat{G}(V,V))_{+}\wedge S^{-V}\rightrightarrows
\igcat{G}(W,V)_{+}\wedge S^{-V} \to i_{\ast}X.
\]
There is thus a natural map 
\begin{equation}
\label{eq:111}
i_{\ast}X\to S^{-W}.
\end{equation}
Evaluating at $U\in \igcat{G}$ and using Lemma~\ref{thm:218} shows
that~\eqref{eq:111} is an isomorphism.  Thus $S^{-W}$ is in the image
of $i_{\ast}$.  It then follows easily that $i_{\ast}$ is essentially
surjective.

Finally, the fact that $i_{\ast}$ is symmetric monoidal is immediate
from the fact that left Kan extensions commute.  It follows that
$i^{\ast}$ is as well, since it is the inverse equivalence.
\end{pf*}

\subsubsection{Equivariant commutative and associative algebras}
\label{sec:equiv-einfty-ainfty}

Using the notions described in \S\ref{sec:commutative-algebras} one
can transport many algebraic structures to $\ugspectra{G}$ using the
symmetric monoidal smash product.   

\begin{defin}
\label{def:35} A {\em $G$-equivariant commutative (associative)
algebra} is a commutative (associative) algebra with unit in
$\ugspectra G$.  
\end{defin}

The conventions of \S\ref{sec:commutative-algebras} would dictate that
we refer to the {\em topological} categories of $G$-equivariant
commutative and associative algebras as $\comm\ugspectra{G}$ and
$\ass\ugspectra{G}$.   To ease some of the typesetting it
will be convenient to employ the slightly abbreviated notation
\begin{align*}
\ugeinftycat{G} &= \comm\ugspectra{G} \\
\ugainftycat{G} &= \ass\ugspectra{G},
\end{align*}
and to write $\geinftycat{G}$ and $\gainftycat{G}$ for the
corresponding $G$-equivariant topological categories of not
necessarily equivariant algebra maps. 

Since $\ugspectra G$ is a closed symmetric monoidal category under
$\wedge$, Proposition~\ref{thm:1} implies that both $\ugeinftycat{G}$
and $\ugainftycat{G}$ are complete and cocomplete, and that the
forgetful functors
\begin{align*}
\ugeinftycat{G} &\to \ugspectra{G} \\
\ugainftycat{G} &\to \ugspectra{G} 
\end{align*}
create enriched limits, sifted colimits, and have left adjoints
\begin{align*}
\sym :\ugspectra{G} &\to \ugeinftycat{G} \\
T: \ugspectra{G}&\to \ugainftycat{G}.
\end{align*}

Similarly, there are categories of left and right modules over an
associative algebra $A$.  We will use the symbol $\rmod{A}$ for the
category of {\em left $A$-modules.}  As described in
\S\ref{sec:symm-mono-categ}, when $A$ is commutative, the category
$\rmod{A}$ inherits a symmetric monoidal product $M\smashove{A}N$
defined by the reflexive coequalizer diagram
\[
M\wedge A\wedge N \rightrightarrows M\wedge N \rightarrow
M\smashove{A} N.
\]

\subsection{Indexed monoidal products}
\label{sec:index-mono-prod}

\subsubsection{Covering categories and fiberwise constructions}

We begin with an example.  Suppose that $(\cat C,\otimes,
\one)$ is a symmetric monoidal category and that $I$ is a finite set.
Write $\cat C^{I}$ for the $I$-fold product of copies of $\cat C$.  For
notational purposes, and subsequent generalization it will be useful
to think of an object of $\cat C^{I}$ as a functor $X:I\to \cat C$,
with $I$ regarded as a category with no non-identity morphisms.  The
iterated monoidal product
\[
\otimes^{I} X = \bigotimes_{i\in I}Xi
\]
defines a functor 
\[
\otimes^{I}:\cat C^{I}\to \cat C.
\]
The functor $\otimes^{I}$ is natural in isomorphisms in $I$ (this is
just the {\em symmetry} of the symmetric monoidal structure).  In this
section we make use of the notion of a {\em covering category} to
exploit this naturality in a systematic way.

Let $\isosets$ be the groupoid of sets and isomorphisms.  Suppose that
$J$ is a category, and that $P:J\to \isosets$ is a functor with the
property that each $Pj$ is finite.  Then $P$ defines a $J$-diagram of
finite sets, and the iterated monoidal product defines for each $j$ a
functor
\begin{equation}
\label{eq:22}
\otimes^{Pj}:\cat C^{Pj}\to \cat C.
\end{equation}
These vary functorially in $j$.  This functoriality is expressed most
cleanly using the Grothendieck
construction~\cite[\S{}VI.8]{MR2017446} (see also~\cite[\S{}B.1]{MR1953060},
or \cite[p.~44]{MR1300636} where the special case in which
$\mathfrak{Cat}$ is replaced with $\mathfrak{Sets}$ is attributed to
Yoneda).

Suppose that $J$ is a category, and that $P:J\to \categories$ is a
functor.  The Grothendieck construction associates to $P$ the category 
\[
I=\int P
\]
of pairs $(j,s)$ with $j\in J$ and $s\in P(j)$.  The set of maps from
$(j,s)$ to $(j',s')$ is the set of pairs $(f,h)$ with $f:j\to j'$ a
map in $J$, and $h:Pf(s)\to s'$ a map in $Pj'$.  By regarding a finite
set as a category with no non-identity morphisms the Grothendieck
construction also applies to functors $P:J\to\isosets$.

A functor $p:I\to J$ arises from the Grothendieck construction of
$P:J\to\isosets$ if and only if it is satisfies the following two
conditions
\begin{itemize}
\item[i)] for every morphism $f:i\to j$ in $J$, and every $a\in I$ with
$pa=i$, there is a unique morphism $g$ with domain $a$, and with $pg=f$;
\item[ii)] for every morphism $f:i\to j$ in $J$, and every $b\in I$ with
$pb=j$, there is a unique morphism $g$ with range $b$, and with $pg=f$.
\end{itemize}
If $p:I\to J$ satisfies the above conditions, then $j\mapsto
p^{-1}(j)$ defines a functor
\[
p^{-1}:J\to\isosets.
\]
This structure is analogous
to the notion of a covering space, and we name it accordingly.

\begin{defin}
A functor $I\to J$ satisfying properties \thmListItem{1} and \thmListItem{2} above is
called a {\em covering category}.
\end{defin}

A covering category $p:I\to J$ in which each of the fibers $p^{-1}(j)$
is finite will be called a {\em finite} covering category.

The aggregate of the functors~\eqref{eq:22} is a functor 
\[
\tensorp p:\cat C^{I}\to \cat C^{J}
\]
given in terms of $p$ by
\[
\tensorp pX(j)= \bigotimes_{p(i)=j} X(i).
\]
We will have much more to say about this in the next few sections.
For now we focus on the general process that led to its construction.

Suppose we are given a formation of a category depending functorially
on a set $I$, or in other words a functor
\[
C:\isosets\to \categories.
\]
Given a covering category $p:I\to J$  let $C_{I}\to J$ be the category
obtained by applying the Grothendieck construction to the composite 
\[
J\to \isosets\xrightarrow{C}{} \categories
\]
in which the first functor is the one classifying $I\to J$.  Let
$C(p)$ be the category of sections of $C_{I}\to J$.  We will say that
$C(p)$ is constructed from $C$ by {\em working fiberwise}.  For
example the category constructed from $C(S)=\cat C^{S}$ by working
fiberwise is $\cat C^{I}$.  The category constructed from the constant
functor $C'(S)=\cat C$ is $\cat C^{J}$.

A natural transformation $C\to C'$ leads, via the same process, to a
functor $C(p)\to C'(p)$ which we will also describe as being
constructed by {\em working fiberwise}.

\subsubsection{Indexed monoidal products}
\label{sec:index-mono-prod-1}

When $(\cat C,\otimes,\one)$ is a symmetric monoidal category, the 
diagram category $\cat C^{I}$ can be regarded as a symmetric monoidal
category using the objectwise monoidal structure.

\begin{defin}
\label{def:38} Let $p:I\to J$ be a finite covering category and $(\cat
C,\otimes,\one)$ a symmetric monoidal category.  The {\em indexed
monoidal product} (along $p$) is the functor
\[
\tensorp p:\cat C^{I}\to \cat C^{J}
\]
constructed fiberwise from the iterated monoidal product.
\end{defin}

For some purposes the notation $X^{\otimes (I/J)}$ is preferable to
$p^{\otimes}_{\ast}X$.  When $J$ is the one point $G$-set this can be
further abbreviated to $X^{\otimes I}$.  We use this alternate
notation systematically when $\otimes$ is the smash product $\wedge$.

The properties of iterated monoidal products listed in the following
proposition are straightforward.

\begin{prop}
\label{thm:147}
The functor $\otimes^{I}:\cat C^{I}\to \cat C$ is symmetric monoidal.  If
\[
\otimes:\cat C^{2}\to \cat C
\]
commutes with colimits in each variable then so does $\otimes^{I}$.
In this situation $\otimes^{I}$ commutes with sifted colimits. \qed
\end{prop}

Applying Proposition~\ref{thm:147} fiberwise to a finite covering
category $p:I\to J$ gives 

\begin{prop}\label{thm:4}
The indexed monoidal product $\tensorp p:\cat C^{I}\to \cat C^{J}$ is
symmetric monoidal.  If
\[
\otimes:\cat C^{2}\to \cat C
\]
commutes with colimits in each variable then $\tensorp p$ commutes
with sifted colimits. \qed
\end{prop}

\begin{rem}
\label{rem:34}
Though it plays no role in this paper, it can be useful to observe
that the class of colimits preserved by $\tensorp{p}$ is slightly
larger than the class of sifted colimits.   For example $\tensorp{p}$
will commute with {\em objectwise reflexive coequalizers}, which are
diagrams of the form 
\[
X \rightrightarrows Y \to Z
\]
with the property that for each $j\in J$ there is a map $Yj\to Xj$
completing 
\[
Xj\rightrightarrows Yj
\]
to a reflexive coequalizer diagram.  The maps $Yj\to Xj$ are not
required to be natural in $j$.
\end{rem}

The following is also straightforward

\begin{prop}
\label{thm:5} Suppose that $p:I\to J$ and $q:J\to K$ are covering
categories.  Then $q\circ p$ is a covering category, which is finite
if $p$ and $q$ are.  In that case there is a natural isomorphism
\[
\tensorp q\circ \tensorp p \approx \tensorp {(q\circ p)}.  
\]
arising from the symmetric monoidal structure.   This natural
isomorphism is compatible with composition in the sense that if
\[
I\xrightarrow{p}{}J\xrightarrow{q}{}K\xrightarrow{r}{}L
\]
is a composition of finite covering categories categories, the diagram 
\[
\xymatrix{
\tensorp r\circ \tensorp q\circ \tensorp p  \ar[r]\ar[d]  & \tensorp r\circ \tensorp{(q\circ p)}
\ar[d] \\
\tensorp{(r\circ q)}\circ \tensorp p  \ar[r]        &\tensorp{(r\circ q\circ p)}
}
\]
(in which the associativity isomorphisms have been suppressed) commutes.
\qed
\end{prop}

The following results are also proved by working fiberwise.

\begin{prop}
\label{thm:10} Suppose that $(\cat C,\otimes, \one_{\cat C})$ and
$(\cat D, \wedge,\one_{\cat D})$ are symmetric monoidal
categories, and that
\begin{gather*}
F:\cat C\to \cat D \\
T:FX\wedge FY\to F(X\otimes Y) \\
\phi:\one_{\cat D} \to F\one_{\cat D}
\end{gather*}
form a lax monoidal functor.  If $p:I\to J$ is a finite covering
category then $T$ gives a natural transformation
\[
p^{T}_{\ast}:p_{\ast}^{\wedge}\circ F^{I} \to F^{J}\circ
\tensorp{p}
\]
between the two ways of going around
\[
\xymatrix{
\cat C^{I} \ar[r]^{F^{I}}\ar[d]_{\tensorp{p}} & \cat D^{I} \ar[d]^{p^{\wedge}_{\ast}} \\
\cat C^{J}\ar[r]_{F^{J}} & \cat D^{J}\mathrlap{\ .}
}
\]
If $T$ is a natural isomorphism, then so is $p^{T}$.  
\qed
\end{prop}

The association $p\mapsto p^{T}_{\ast}$ of
Proposition~\ref{thm:10} is compatible with the composition
isomorphism of Proposition~\ref{thm:5} in the evident sense.

Suppose that $p:I\to J$ is a covering category, and $f:\tilde J\to J$
is a functor.  Let $\tilde I$ be the ``rigid pullback'' category of
pairs $(j',i)\in \tilde J\times I$ with $f(i')= p(j)$, and in which a
morphism is a pair $(g,g')$ with $f(g)=p(g')$.  Then the functor
$\tilde p:\tilde I\to \tilde J$ defined by $(j',i)\mapsto j'$ is a
covering category.

\begin{prop}
\label{thm:46} In the situation described above, if $p:I\to J$ is
finite then the following commutes up to a natural isomorphism given
by the symmetric monoidal structure
\[
\xymatrix{
\cat C^{I} \ar[r] \ar[d]_{\tensorp{p}} & \cat C^{\tilde I} \ar[d]^{\tensorp{\tilde{p}}} \\
\cat C^{J} \ar[r]_{f^{\ast}} & \cat C^{\tilde J}\mathrlap{\ .} 
}
\]
\qed
\end{prop}

The categories $I$ and $J$ used in this paper arise from a left action
of a group $G$ on a finite set $A$.  Given such an $A$, let
$\cat{B}_{A}G$ be the category whose set of objects is $A$ and in
which a map $a\to a'$ is an element $g\in G$ with the property that
$g\,a=a'$.  When $A=\text{pt}$ we will abbreviate $\cat{B}_{A}G$ to
just $\cat{B}G$.  For any finite map $A\to B$ of $G$-sets, the
corresponding functor
\[
\cat{B}_{A}G\to \cat{B}_{B}G
\]
is a covering category.

In the following series of examples we suppose $H\subset G$ is a
subgroup, take $A=G/H$ to be the set of right $H$-cosets, and write
$p:A\to\text{pt}$ for the unique equivariant map.   In this case the
inclusion of the identity coset gives an equivalence 
\[
\cat{B} H\to \cat{B}_{A}G
\]
and hence an equivalence of functor categories
\[
\cat C^{\cat{B}_{A}G} \to \cat{C}^{\cat{B}H}.
\]
An inverse is provided by the left Kan extension.

\begin{eg}
\label{eg:10} 
Suppose $\cat C$ is the category of abelian groups, with
$\oplus$ as the symmetric monoidal structure.  Then
$\cat{C}^{\cat{B}_{A}G}$ is equivalent to the category of left
$H$-modules, and the functor $p^{\oplus}_{\ast}$ is left additive
induction.  If the symmetric monoidal structure is taken to be the
tensor product, then $\tensorp{p}$ is ``norm induction.''
\end{eg}

\begin{eg}
Now take $(\cat C,\otimes,\one)$ to be the category $(\spectra,\wedge,
S^{0})$ of orthogonal spectra.  From the above and
Proposition~\ref{thm:143}, the category $\spectra^{\cat{B}_{A}G}$ is
equivalent to the category of orthogonal $H$-spectra, and
$\spectra^{\cat{B}G}$ is equivalent to the category of orthogonal
$G$-spectra.  In this case $p^{\wedge}_{\ast}$ defines a
multiplicative transfer from orthogonal $H$-spectra to orthogonal
$G$-spectra.  This is the {\em norm}.  It is discussed more fully in
\S\ref{sec:norm} and \S\ref{sec:homot-prop-norm}.
\end{eg}

\begin{rem}
\label{rem:1} When $\cat C$ has all colimits and the tensor unit
$\one$ is the initial object one may form infinite ``weak'' monoidal
products, and the condition that $p:I\to J$ be finite may be dropped.
If $I$ is an infinite set and $\{Xi \}$ a collection of objects
indexed by $i\in I$ set
\[
\otimes^{I} Xi = \varinjlim_{I'\subset I\text{ finite}} \otimes^{I'}Xi
\]
in which the transition maps associated to $I'\subset I''$ are given
by tensoring with the unit
\[
\otimes^{I'}Xi \approx
\left(\otimes^{I'}Xi\right)\otimes\left(\otimes^{I''-
I'}\one\right) \to \otimes^{I''}Xi.
\]
The functor $\tensorp{p}$ is constructed by working fiberwise.
\end{rem}

\begin{rem}
The results of this section apply, with the obvious modifications, in
the setting of enriched categories.
\end{rem}

\subsubsection{Distributive laws}\label{sec:distributive-laws}

Continuing with the same notation, we now assume that the category
$\cat C$ comes equipped with two symmetric monoidal structures,
$\otimes$ and $\oplus$, and that $\otimes$ distributes over $\oplus$
in the sense that there is a natural isomorphism
\[
A\otimes (B\oplus C) \approx (A\otimes B)\oplus (A\otimes
C)
\]
compatible with all of the symmetries.  For a precise definition
see~\cite{MR0335598}, or the definition of {\em bipermutative
category} in~\cite[Chapter~VI]{may77:_e_e}.  In all of our examples,
$\oplus$ will be the categorical coproduct, and $A\otimes(\slot)$ will
commute with all colimits.  Given $p:I\to J$ and $q:J\to K$ we can
form
\[
q^{\otimes}_{\ast}\circ p_{\ast}^{\oplus}.
\]
Our goal is to express this in the form
\[
q^{\otimes}_{\ast}\circ p_{\ast}^{\oplus} = 
r_{\ast}^{\oplus} \circ \pi^{\otimes}_{\ast}.
\]

We start with the case in which $K$ is the trivial category, and
$p:I\to J$ is a map of finite sets.  Let $\Gamma=\Gamma(I/J)$ be the
set of sections $s:J\to I$ of $p$.  Write $\ev:J\times\Gamma \to I$
for the evaluation map, $\pi:J\times\Gamma\to \Gamma$ for the
projection, and with an eye toward generalization, $r:\Gamma\to
\{\text{pt}\}$ for the unique map.  The following lemma expresses the
usual distributivity expansion
\[
\bigotimes_{j\in J} \left(\bigoplus_{p(i)=j}Xi\right) \approx
\bigoplus_{s\in \Gamma} \left(\bigotimes_{j\in J}Xs(j)\right) 
\]
in functorial terms.

\begin{lem}
\label{thm:8}
The following diagram of functors commutes, up to a canonical natural
isomorphism given by the symmetries of the symmetric monoidal
structures
\[
\xymatrix@C=1em{
\cat C^{I} \ar[d]_{p^{\oplus}_{\ast}} \ar[rr]^{\ev^{\ast}} & & \cat C^{J\times \Gamma}
\ar[d]^{\pi^{\otimes}_{\ast}} \\
\cat C^{J} \ar[dr]_{q^{\otimes}_{\ast}} && \cat C^{\Gamma} \ar[dl]^{r^{\oplus}_{\ast}} \\
& \cat C &.
}
\]
\qed
\end{lem}

Working fiberwise, it is now a simple matter to deal with the more
general case in which $p:I\to J$ and $q:J\to K$ are covering categories.
Let $\Gamma$ be the category of pairs $(k,s)$, with $k\in K$ and $s$ a
section of $(q\circ p)^{-1}k\to q^{-1}k$.  A morphism $(k,s)\to
(k,s')$ in $\Gamma$ is a map $f:k\to k'$ making the following diagram
commute
\[
\xymatrix@C=1in{
(q\circ p)^{-1}k \ar[r]^{(q\circ p)^{-1}(f)} & (q\circ p)^{-1}k' \\
p^{-1}k \ar[u]^{s}\ar[r]_{p^{-1}(f)} & p^{-1}k'\ar[u]_{s'}.
}
\]
Write $\Gamma\underset{K}{\times}J$ for the fiber product, 
\[
\ev:\Gamma\underset{K}{\times}J\to I
\]
for the ``evaluation,'' and $\pi:\Gamma\underset{K}{\times}J\to J$ for
the projection.    By naturality in $I$ and $J$ in Lemma~\ref{thm:8} we have

\begin{prop}
\label{thm:9}
The following diagram of functors commutes, up to a canonical natural
isomorphism given by the symmetries of the symmetric monoidal
structures
\[
\xymatrix@C=1em{
\cat C^{I} \ar[d]_{p^{\oplus}_{\ast}} \ar[rr]^{\ev^{\ast}} & & \cat C^{J\underset{K}{\times} \Gamma}
\ar[d]^{\pi^{\otimes}_{\ast}} \\
\cat C^{J} \ar[dr]_{q^{\otimes}_{\ast}} && \cat C^{\Gamma} \ar[dl]^{r^{\oplus}_{\ast}} \\
& \cat C^{K} &
}
\]
\qed
\end{prop}

This formula is used in showing that the norm of a wedge of regular
slice cells is a wedge of regular slice cells
(Proposition~\ref{thm:77}), in the construction of monomial ideals
(\S\ref{sec:monomial-ideals}), and in describing the structure of
equivariant twisted monoid rings and their monomial ideals
(\S\ref{sec:meth-poly-algebr}).

\subsubsection{Indexed monoidal products and pushouts}
\label{sec:index-mono-prod-pushout}

The homotopy theoretic properties of the norm depend on a formula for
the indexed monoidal product of a pushout.  We describe here the
absolute case.  The fiberwise analogue is spelled out in
\S\ref{sec:very-flat-diagrams}.

Suppose that $(\cat C,\otimes)$ is a closed symmetric monoidal
category with finite colimits, and let $I$ be a finite set.  For $X\in
\cat C^{I}$ write $X^{\otimes I}$ for the iterated monoidal product.
Suppose we are given a pushout diagram
\begin{equation}
\label{eq:174}
\xymatrix{
A \ar[r]\ar[d]  &  B \ar[d] \\
X  \ar[r]        & Y
}
\end{equation}
in $\cat C^{I}$.  We wish to express $Y^{\otimes I}$ as an iterated
pushout.  Since the coequalizer diagram
\[
X\amalg A\amalg B\rightrightarrows X\amalg B\to Y
\]
can be completed to a reflexive coequalizer, the sequence 
\[
(X\amalg A\amalg B)^{\otimes I}\rightrightarrows (X\amalg B)^{\otimes
I}\to Y^{\otimes I}
\]
is a coequalizer (Proposition~\ref{thm:147}).  Using the
distributivity law of \S\ref{sec:distributive-laws} this can be
re-written as
\[
\coprod_{I=I_{0}\amalg I_{1}\amalg I_{2}} X^{\otimes I_{0}}\otimes
A^{\otimes I_{1}}\otimes B^{\otimes I_{2}}
\rightrightarrows
\coprod_{I=I_{0}\amalg I_{1}} X^{\otimes I_{0}}\otimes B^{\otimes
I_{1}}
\to Y^{\otimes I}.
\]
The horizontal arrows do not preserve the coproduct decompositions,
but the sequence can be filtered by the cardinality of the exponent of
$B$.  Define $\fil_{n}Y$ by the coequalizer diagram
\[
\coprod_{\substack{I=I_{0}\amalg I_{1}\amalg I_{2} \\
|I_{1}|+|I_{2}| \le n}} 
X^{\otimes I_{0}}\otimes
A^{\otimes I_{1}}\otimes B^{\otimes I_{2}}
\rightrightarrows
\coprod_{\substack{I=I_{0}\amalg I_{1} \\ |I_{1}|\le n}}
X^{\otimes I_{0}}\otimes B^{\otimes I_{1}}
\to \fil_{n}Y.
\]
Thus $\fil_{0}Y=X^{\otimes I}$ and $\fil_{|I|}Y = Y^{\otimes I}$.
There is an evident coequalizer diagram
\[
\coprod_{\substack{I=I_{0}\amalg I_{1}\amalg I_{2} \\
|I_{1}|+|I_{2}| = n}} 
X^{\otimes I_{0}}\otimes
A^{\otimes I_{1}}\otimes B^{\otimes I_{2}} 
\rightrightarrows
\fil_{n-1}Y \amalg \big(\coprod_{\substack{I=I_{0}\amalg I_{1} \\ |I_{1}| = n}}
X^{\otimes I_{0}}\otimes B^{\otimes I_{1}}\big)
\to \fil_{n}Y, 
\]
which can be re-written as  a pushout square
\[
\xymatrix{
{\displaystyle\coprod_{\substack{I=I_{2}\amalg I_{1}\amalg I_{0} \\
|I_{0}| = |I|- n}} X^{\otimes I_{0}}\otimes A^{\otimes I_{1}}\otimes
B^{\otimes I_{2}} } \ar[r]\ar[d] &
{\displaystyle\coprod_{\substack{I=I_{1}\amalg I_{0} \\ |I_{1}| = n}}
X^{\otimes I_{0}}\otimes B^{\otimes I_{1}}} \ar[d] \\
\fil_{n-1}Y  \ar[r]        & \fil_{n}Y\mathrlap{\ .}
}
\]
The upper left term may be replaced by its effective quotient
\[
\coprod_{|I_{1}|=n}
X^{\otimes I_{0}}\otimes \partial_A{B^{\otimes I_{1}}}
\]
in which $\partial_{A}B^{\otimes S}$ is defined by the coequalizer
diagram
\begin{equation}
\label{eq:34}
\coprod_{\substack{S=S_{0}\amalg S_{1}\amalg S_{2}\\ S_{0}\ne
\emptyset}} A^{\otimes S_{0}}\otimes A^{\otimes S_{1}}\otimes
B^{\otimes S_{2}} \rightrightarrows \coprod_{\substack{S=S_{0}\amalg
S_{1} \\ S_{0}\ne\emptyset}} A^{\otimes S_{0}}\otimes B^{\otimes
S_{1}} \to \partial_{A}B^{\otimes S}, 
\end{equation}
leading to a pushout square
\begin{equation}
\label{eq:59}
\entrymodifiers={+!!<0pt,\fontdimen22\textfont2>}
\xymatrix{
{\displaystyle 
\coprod_{\substack{I=I_{0}\amalg I_{1} \\ |I_{1}|=n}}
X^{\otimes I_{0}}\otimes \partial_A{B^{\otimes I_{1}}}}
  \ar[r]\ar[d]  &   
{\displaystyle 
\coprod_{\substack{I=I_{0}\amalg I_{1} \\ |I_{1}|=n}}
X^{\otimes I_{0}}\otimes B^{\otimes I_{1}}}\ar[d] \\
\fil_{n-1}Y^{\otimes I}  \ar[r]        & \fil_{n}Y^{\otimes I}\mathrlap{\ .}
}
\end{equation}
The object $\partial_{A}B^{\otimes S}$ can also be computed as the
coequalizer of 
\begin{equation}
\label{eq:175} \coprod_{\substack{S=S_{0}\amalg S_{1}\amalg S_{2}\\
|S_{0}|=|S_{1}|=1}} A^{\otimes S_{0}}\otimes A^{\otimes S_{1}}\otimes
B^{\otimes S_{2}} \rightrightarrows \coprod_{\substack{S=S_{0}\amalg
S_{1} \\ |S_{0}|=1}} A^{\otimes S_{0}}\otimes B^{\otimes S_{1}} \to
\partial_{A}B^{\otimes S},
\end{equation}

We call the map
\[
\partial_{A} B^{\otimes S}\to B^{\otimes S} 
\]
the {\em indexed corner map}, since in the absolute case with $|I|=2$
it reduces to the ``corner map'' in 
\[
\xymatrix{
A\otimes B\coprod B\otimes A  \ar[r]\ar[d]  &
B\otimes A\ar[d] \\
A\otimes B  \ar[r]        & B\otimes B
}
\]
from the pushout of the top and left arrows to the bottom right term.

\begin{rem}
\label{rem:47}
The category of arrows in $\cat C$ becomes a closed symmetric
monoidal category with 
\[
(A_{1}\to B_{1})\otimes (A_{2}\to B_{2})
\]
taken to be the corner map in 
\[
\xymatrix{
A_{1}\otimes A_{2}  \ar[r]\ar[d]  &  A_{1}\otimes B_{2} \ar[d] \\
B_{1}\otimes A_{2}  \ar[r]        & B_{1}\otimes B_{2}.
}
\]
If $A\to B$ is a map in $\cat C^{S}$ then $(A\to B)^{\otimes S}$ works
out to be $\partial_{A}B^{\otimes S}\to B^{\otimes S}$.  
\end{rem}

By working fiberwise, one obtains a similar iterated pushout
describing $p_{\ast}^{\otimes} Y$, involving the evident analogue
$\partial_{A}\tensorp{p}B\to \tensorp{p}B$ of
$\partial_{A}B^{\otimes S}\to B^{\otimes S}$.   

Taking $A=X$ and $B=Y$ in~\eqref{eq:174} gives a filtration of the
indexed monoidal product of any map.  In the case of a pushout
square~\eqref{eq:174} the two filtrations in fact coincide.  We
describe the absolute case.  The relative case follows easily by
working fiberwise.

\begin{prop}
\label{thm:265}
Let 
\begin{equation}
\label{eq:176}
\xymatrix{
A  \ar[r]\ar[d]  &  B \ar[d] \\
X  \ar[r]        & Y
}
\end{equation}
be a pushout square in $\cat C^{I}$.  
\begin{thmList}
\item The square
\[
\xymatrix{
\partial_{A}B^{\otimes I}  \ar[r]\ar[d]  & B^{\otimes I}  \ar[d] \\
\partial_{X}Y^{\otimes I}  \ar[r]        &  Y^{\otimes I}
}
\]
is a pushout square.
\item  The filtrations of $Y^{\otimes
I}$ arising from~\eqref{eq:176} and
\begin{equation}
\label{eq:83}
\xymatrix{
X  \ar[r]\ar[d]  &  X \ar[d] \\
Y  \ar[r]        & Y
}
\end{equation}
coincide.
\end{thmList}
\end{prop}

\begin{pf}
The proof is by induction on $n=|I|$, the case in which $n=1$ being
trivial.  Let $\fil_{m}Y^{\otimes I}$ be the filtration computed from
the pushout square~\eqref{eq:176}, and $\fil'_{m}Y^{\otimes I}$ be the
one computed from~\eqref{eq:83}.  The evident map of squares gives a
natural map $\fil_{m}Y^{\otimes I}\to\fil'_{m}Y^{\otimes I}$.
Consider the diagram
\[
\xymatrix{
{\displaystyle 
\coprod_{\substack{I=I_{0}\amalg I_{1} \\ |I_{1}|=m}}
X^{\otimes I_{0}}\otimes \partial_A{B^{\otimes I_{1}}}}
  \ar[r]\ar[d]  &   
{\displaystyle 
\coprod_{\substack{I=I_{0}\amalg I_{1} \\ |I_{1}|=m}}
X^{\otimes I_{0}}\otimes B^{\otimes I_{1}}}\ar[d] \\
{\displaystyle 
\coprod_{\substack{I=I_{0}\amalg I_{1} \\ |I_{1}|=m}}
X^{\otimes I_{0}}\otimes \partial_X{Y^{\otimes I_{1}}}}
  \ar[r]\ar[d]  &   
{\displaystyle 
\coprod_{\substack{I=I_{0}\amalg I_{1} \\ |I_{1}|=m}}
X^{\otimes I_{0}}\otimes Y^{\otimes I_{1}}}\ar[d] \\
\fil_{m-1}Y^{\otimes I}  \ar[r]        & \fil_{m}Y^{\otimes I}\mathrlap{\ .}
}
\]
If $m<n$ then the induction hypothesis and part \thmListItem{1} imply
that the upper square is a pushout square.  This shows that the map
$\fil_{m}Y^{\otimes I}\to \fil'_{m}Y^{\otimes I}$ is an isomorphism
$m<n$.  The case $m=n-1$ then gives an identification
\[
\fil_{n-1} Y^{\otimes I}= \partial_{X}Y^{\otimes I},
\]
which, when combined with the  pushout square 
\[
\xymatrix{
\partial_{A}B^{\otimes I}  \ar[r]\ar[d]  & B^{\otimes I}
\ar[d] \\
\fil_{n-1}Y^{\otimes I}  \ar[r]        & Y^{\otimes I}
}
\]
gives part \thmListItem{1} for $I$.
\end{pf}

By working in the category of arrows (as in Remark~\ref{rem:47}) one
can see that the formation of $\partial_A(\tensorp{p}B)$ is compatible
with the isomorphism coming from the distributive law.  
More explicitly, let $I\xrightarrow{p} J\xrightarrow{q} K$ be a
sequence of covering categories, and recall the basic diagram
\[
\xymatrix@C=1em{
\cat C^{I} \ar[d]_{p^{\oplus}_{\ast}} \ar[rr]^{\ev^{\ast}} & & \cat C^{J\underset{K}{\times} \Gamma}
\ar[d]^{\pi^{\otimes}_{\ast}} \\
\cat C^{J} \ar[dr]_{q^{\otimes}_{\ast}} && \cat C^{\Gamma} \ar[dl]^{r^{\oplus}_{\ast}} \\
& \cat C^{K} &\mathrlap{\phantom{MM}.}
}
\]
Suppose that $A\to B$ is a map in $\cat C^{I}$.  The distributivity
isomorphism in the arrow category is given by
\[
\xymatrix{
\partial_{(p^{\oplus}_{\ast}A)}(p^{\oplus}_{\ast} B)^{\otimes J}
\ar[r]\ar[d]_{\approx}  & (p^{\oplus}_{\ast} B)^{\otimes J} \ar[d]^{\approx}\\
r^{\oplus}_{\ast} \partial_{(\ev^{\ast}A)}
\pi^{\otimes}_{\ast}(ev^{\ast}B) \ar[r] & r^{\oplus}_{\ast}\pi^{\otimes}_{\ast}(ev^{\ast}B)\mathrlap{\ .}
}
\]
The fact that the left vertical arrow is an isomorphism is what
expresses the compatibility of $\partial_{A} q^{\otimes}_{\ast}B$ with
the distributive law.  

\subsubsection{Commutative algebras and indexed monoidal products}
\label{sec:comm-algebr-index}

By Proposition~\ref{thm:1}, if $\cat C$ is a co-complete closed
symmetric monoidal category, then $\commcat C$ is cocomplete.  The
restriction functor $p^{\ast}:\commcat C^{J}\to \commcat C^{I}$ then
has a left adjoint $p_{!}$ given by left Kan extension.

\begin{prop}
\label{thm:18} If $p:I\to J$ is a covering category, the following
diagram commutes
\[
\xymatrix{
\commcat C^{I}\ar[r]\ar[d]^{p_{!}} & \cat C^{I} \ar[d]^{\tensorp{p}}\\
\commcat C^{J}\ar[r] & \cat C^{J}\mathrlap{\ .}
}
\]
\end{prop}

\begin{pf}
For a commutative algebra $A\in \commcat C^{I}$, and $j\in J$ the
value of $p_{!}A$ at $j$ is calculated as the colimit over the
category $I/j$ of the restriction of $p$.  Since $p:I\to J$ is a
covering category, the category $I/j$ is equivalent to the discrete
category $p^{-1}j$, and so
\[
(p_{!}A)j = \otimes^{p^{-1}j}A,
\]
and the result follows.
\end{pf}

\subsubsection{Monomial ideals}
\label{sec:monomial-ideals}

Let $I$ be a set and consider the polynomial algebra 
\[
A=\Z[x_{i}], \quad i \in I.
\]
As an abelian group, it has a basis consisting of the monomials $x^{f}$, with 
\[
f:I\to \{0,1,2,\dots \}
\]
a function taking the value zero on all but finitely many elements, and
\[
x^{f} = \prod_{j\in J}x_{j}^{f(j)}.
\]
The collection of such $f$ is a monoid under addition, and we denote
it $\nat^{I}$.  If $D\subset \nat^{I}$ is a monoid ideal then the
subgroup $M_{D}\subset A$ with basis $\{x^{f}\mid f\in D \}$ is an
ideal.  These are the {\em monomial ideals} and they can be formed in
any monoidal product of free associative algebras in any closed
symmetric monoidal category.

Let $(\cat C,\otimes, \one)$ be a closed symmetric monoidal category.
Fix a set $I$ which we temporarily assume to be finite.  Given $X\in
\cat C^{I}$ let
\[
TX=\coprod_{n\ge 0}X^{\otimes n}
\]
be the free associative algebra generated by $X$.  Write
$A=\tensorp{p}TX\in \cat C$, where $p:I\to \text{pt}$ is the
unique map.  Then $A$ is an associative algebra in $\cat C$.  The
motivating example above occurs when $\cat C$ is the category of
abelian groups and $X$ is the constant diagram $Xi=\Z$.

Using Proposition~\ref{thm:9}, the object $A$ can be expressed as an
indexed coproduct
\[
A=\coprod_{f:I\to \nat} X^{\otimes f}
\]
where $\nat=\{0,1,2,\dots \}$ and 
\[
X^{\otimes f}=\bigotimes_{i\in I}X(i)^{\otimes f(i)}.  
\]
The set
\[
\nat^{I}=\{f:I\to \nat \}
\]
is a commutative monoid under addition of functions.  The
multiplication map in $A$ is the sum of the isomorphisms
\begin{equation}
\label{eq:29}
X^{\otimes f}\otimes X^{\otimes g}\approx X^{\otimes(f+g)}
\end{equation}
given by the symmetry of the monoidal product $\otimes$, and the
isomorphism
\[
X^{\otimes f(i)}\otimes X^{\otimes g(i)} \approx X^{\otimes(f(i)+g(i))}.
\]
For a monoid ideal $D\subset \nat^{I}$, set
\[
M_{D}=\coprod_{f\in D} X^{\otimes f}.
\]
The formula~\eqref{eq:29} for the multiplication in $A$ gives $M_{D}$
the structure of an ideal in $A$.  If $D\subset D'$ then the evident
inclusion $M_{D}\subset M_{D'}$ is an inclusion of ideals.

When $\cat C$ is pointed (in the sense that the initial object is
isomorphic to the terminal object), the map
\[
A\to A/M_{D}
\]
is a map of associative algebras, where $A/M_{D}$ is defined by the
pushout diagram 
\[
\xymatrix{
M_{D}   \ar[r] \ar[d]  & A \ar[d]\\
\ast   \ar[r] & A/M_{D}\mathrlap{\ ,}
}
\]
with $\ast$ denoting the terminal (and initial) object.

\begin{defin}
The ideal $M_{D}\subset A$ is the {\em  monomial ideal} associated to the
monoid ideal $D$.
\end{defin}

\begin{eg}
\label{eg:6}
Suppose that $\dim:\nat^{I}\to \nat$ is any homomorphism.   Given 
$d\in\nat$ the set 
\[
\{f\mid \dim f\ge d\}
\]
is a monoid ideal.   We denote the corresponding monomial ideal
$M_{d}$.     The $M_{d}$ form a
decreasing filtration
\[
\dots \subset M_{d+1}\subset M_{d}\subset\dots\subset M_{1}\subset M_{0}=A.
\]
When $\cat C$ is pointed, the quotient
\[
M_{d}/M_{d+1}
\]
is isomorphic as an $A$ bimodule to 
\[
A/M_{1}\otimes \coprod_{\dim f=d} X^{\otimes f},
\]
in which $A$ act through its action on the left factor. 
\end{eg}

\begin{rem}
\label{rem:4}
The quotient module is defined by the pushout square
\[
\xymatrix{
M_{d+1}\ar[r] \ar[d] & M_{d} \ar[d] \\
\ast \ar[r] & M_{d}/M_{d+1}\mathrlap{\ .}
}
\]
The pushout can be calculated in the category of left $A$-modules, $A$
bimodules, or just in $\cat C$.
\end{rem}

\begin{rem}
\label{rem:31} 
All of this discussion can be made to be covariant with
respect to inclusion in $I$.  Suppose that $I_{0}\subset I_{1}$ is an
inclusion of finite sets and $X_{1}:I_{1}\to \cat C$ is an
$I_{1}$-diagram.  Define $X_{0}:I_{1}\to \cat C$ by 
\[
X_{0}(i)= 
\begin{cases}
X_{1}(i) & i\in I_{0} \\
\ast &\text{otherwise}.
\end{cases}
\]
There is a natural map $X_{0}\to X_{1}$.  Let $A_{0}$ and $A_{1}$ be
the associative algebras constructed from the $X_{i}$ as described
above.  The algebra $A_{0}$ coincides with the one constructed
directly from the restriction of $X_{0}$ to $I_{0}$.  A monoid ideal
$D_{1}\subset\nat^{I_{1}}$ defines ideals $M_{D_{0}}\subset A_{0}$ and
$M_{D_{1}}\subset A_{1}$.  The monoid ideal $D_{0}$ is the same as
the one constructed from the intersection of $D_{0}$ with
$\nat^{I_{0}}$, where $\nat^{I_{0}}$ is regarded as a subset of
$\nat^{I_{1}}$ by extension by $0$.  There is a commutative diagram
\[
\xymatrix{
M_{D_{0}}  \ar[r]\ar[d] &M_{D_{1}}   \ar[d] \\
A_{0} \ar[r] & A_{1}\mathrlap{\ .} 
}
\]
Using this, the construction of monomial ideals can be extended to
the case of infinite sets $I$, by passing to the colimit over the
finite subsets.  As in the motivating example, when the set $I$ is
infinite, the indexing monoid $\nat^{I}$ is the set of finitely
supported functions.
\end{rem}

By working fiberwise, this entire discussion applies to the situation
of a (possibly infinite) covering category $p:I\to J$.  Associated to
$X:I\to \cat C$ is 
\[
A=\tensorp{p}TX\in
\asscat{C^{J}}=(\asscat{C})^{J}.
\]
In case $I/J$ is infinite, the algebra $A$ is formed fiberwise by
passing to the colimit from the finite monoidal products using the
unit map, as described in Remark~\ref{rem:1}.  As an object of
$\cat{C}^{J}$, the algebra $A$ decomposes into
\[
A=\coprod_{f\in\Gamma} X^{\otimes f}
\]
where $\Gamma$ is the set of sections of 
\[
\nat^{I/J}\to J
\]
with $\nat^{I/J}$ formed from the Grothendieck construction applied to 
\[
j\mapsto \nat^{I_{j}}\qquad (I_{j}= p^{-1}(j)).
\]
The category $\nat^{I/J}$  is a commutative monoid over $J$, and
associated to any monoid ideal $D\subset \nat^{I/J}$ over $J$, is a
monomial ideal $M_{D}\subset A$.    

The situation of interest in this paper (see \S\ref{sec:meth-poly-algebr})
is when $I\to J$ is of the form 
\[
\cat{B}_{K}G\to \cat{B}G
\]
associated to a $G$-set $K$, and the unique map $K\to\text{pt}$.  In
this case $\nat^{I/J}$ is the $G$-set $\nat^{K}$ of finitely supported
functions $K\to\nat$.  The relative monoid ideals are just the
$G$-stable monoid ideals.  A simple algebraic example arises in the
case of a polynomial algebra $\Z[x_{i}]$ in which a group $G$ is
acting on the set indexing the variables.

\subsection{The norm}
\label{sec:norm}

We now specialize the discussion of \S\ref{sec:index-mono-prod} to the
case $(\cat C,\otimes,\one)=(\spectra,\wedge, S^{0})$ and define the
{\em norm functor}.

Because of Proposition~\ref{thm:143} we may identify the category of
$G$-equivariant orthogonal spectra as the functor category
$\spectra^{\cat{B}G}$.  If $H\subset G$, then the functor
\[
i:\cat{B}H\to \cat{B}_{G/H}G
\]
sending the unique object to the coset $H/H$ is an equivalence of
categories.  This leads to an equivalence
\[
i^{\ast}:\spectra^{\cat{B}_{G/H}G}\to \spectra^{\cat{B}H},
\]
with inverse 
\[
i_{!}:\spectra^{\cat{B}H}\to \spectra^{\cat{B}_{G/H}G}
\]
given by the left Kan extension.  Write $p:\cat{B}_{G/H}G\to \cat{B}G$
for the functor corresponding to the $G$-map $G/H\to\text{pt}$.

\begin{defin}
The {\em norm functor} $\norm_{H}^{G}:\ugspectra{H}\to \ugspectra{G}$ is
the composite 
\[
\xymatrix{
\spectra^{\cat{B}H} \ar[r]^{i_{!}} \ar[dr]_{\norm_{H}^{G}} & \spectra^{\cat{B}_{G/H}G} \ar[d]^{p_{\ast}^{\wedge}} \\
&\spectra^{\cat{B}G}}
\]
\end{defin}

By Proposition~\ref{thm:4} we have

\begin{prop}
\label{thm:150} The functor $\norm_{H}^{G}$ is symmetric monoidal and
commutes with sifted colimits.\qed
\end{prop}

\begin{rem}
\label{rem:25} By Remark~\ref{rem:34}, the norm also commutes with the
formation of coequalizer diagrams in $\ugspectra{H}$ whose underlying
non-equivariant diagram in $\spectra$ extends to a reflexive
coequalizer.
\end{rem}

\begin{rem}
\label{rem:51}
We have defined the norm on the topological categories of equivariant
spectra.   Since it is symmetric monoidal it naturally extends to a
functor  of enriched categories
\[
\norm_{H}^{G}:\gspectra{H}\to \gspectra{G}
\]
compatible with the norm on spaces (and, in fact, spectra) in the
sense that it gives for every $X,Y\in\gspectra{H}$ a $G$-equivariant
map
\[
\norm_{H}^{G}(\gspectra{H}(X,Y)) \to
\gspectra{G}(\norm_{H}^{G}X,\norm_{H}^{G}Y).
\]
\end{rem}

By Proposition~\ref{thm:18}, on equivariant commutative algebras
the norm is the left adjoint of the restriction functor.

\begin{cor}
\label{thm:184}
The following diagram commutes up to a natural isomorphism given by
the symmetry of the smash product:
\[
\xymatrix{
\ugeinftycat{H} \ar[r] \ar[d] & \ugspectra{H} \ar[d]^{\norm_{H}^{G}} \\
\ugeinftycat{G} \ar[r] & \ugspectra{G}\mathrlap{\ .}  
}
\]
The left vertical arrow is the left adjoint to the restriction functor.
\end{cor}

\begin{rem}
\label{rem:11}
Because of Corollary~\ref{thm:184} we will refer to the left adjoint
to the restriction functor 
\[
\ugeinftycat{G}\to \ugeinftycat{H}
\]
as the {\em commutative algebra norm}, and denote it
\[
\norm_{H}^{G}:\ugeinftycat{H}\to \ugeinftycat{G}.
\]
\end{rem}

The Yoneda embedding gives a functor 
\begin{align*}
\icat^{op} &\to \spectra  \\
V &\mapsto S^{-V}.  
\end{align*}
By definition of $\wedge$ this is a symmetric monoidal functor, and we
are in the situation described in Proposition~\ref{thm:10}.  Thus if
$p:I\to J$ is a covering category, there is a natural isomorphism
between the two ways of going around
\begin{equation}
\label{eq:2}
\xymatrix{
\big(\icat^{op}\big)^{I} \ar[r] \ar[d]_{p^{\oplus}_{\ast}} & \spectra^{I} \ar[d]^{p^{\wedge}_{\ast}}\\
\big(\icat^{op}\big)^{J} \ar[r] & \spectra^{J}\mathrlap{\ .}
}
\end{equation}
Take $I=\cat{B}_{G/H}G$ and $J=\cat{B}G$.  Then the functor category
$\big(\icat^{op}\big)^{I}$ is equivalent to the category
$\big(\uigcat{H}\big)^{op}$ (Proposition~\ref{thm:3}), and
$\spectra^{I}$ is equivalent to $\ugspectra{H}$
(Proposition~\ref{thm:143}).  By naturality, the functor
\[
{\uigcat{H}}^{op} \to \ugspectra{H}
\]
corresponding to 
\[
\big(\icat^{op}\big)^{I}\to \spectra^{I}
\]
is just the Yoneda embedding, and so sends an orthogonal
$H$-representation $V$ to $S^{-V}$.  Similarly $\big(\icat^{op}\big)^{J}$ is
equivalent to ${\uigcat{G}}^{op}$, $\spectra^{J}$ is equivalent to the
category of orthogonal $G$-spectra, and the functor between them sends
an orthogonal $G$-representation $W$ to $W^{-W}$.  One easily checks
(as in Example~\ref{eg:10}) that the functor $p^{\oplus}_{\ast}$
corresponds to additive induction.  We therefore have a commutative
diagram
\[
\xymatrix{
\big(\uigcat{H}\big)^{op}  \ar[r]\ar[d]_{\ind_{H}^{G}}  & \ugspectra{H}  \ar[d]^{\norm_{H}^{G}}\\
\big(\uigcat{G}\big)^{op}  \ar[r] &  \ugspectra{G} 
}
\]
This proves

\begin{prop}
\label{thm:12} Let $V$ be a finite dimensional $H$-representation, and
set $W=\ind_{H}^{G}V$.  There is a natural isomorphism
\[
\norm_{H}^{G}\,S^{-V} \approx S^{-W}.
\]
\qed
\end{prop}

\subsection{$h$-cofibrations}
\label{sec:h-cofibrations}
Suppose that $\cat C$ is a complete topological category (and in
particular tensored and cotensored over $\spaces$).

\begin{defin}
\label{def:49}
A map $i:A\to X$ in $\cat C$ is an {\em $h$-cofibration} if it has the
homotopy extension property:  given $f:X\to Y$ and a homotopy
$h:A\otimes [0,1]\to Y$ with $h\vert_{A\otimes\{0 \}}=f\circ i$ there
is an extension of $h$ to $H:X\otimes[0,1]\to Y$.
\end{defin}

\begin{eg}
The mapping cylinder $A\to X \underset{A}{\cup}A\otimes [0,1]$ of
any map $A\to X$ is an $h$-cofibration.
\end{eg}

As is well-known, a map $i:A\to X$ is an $h$-cofibration if and only if
\[
\cyl i = X\otimes\{0 \}\underset{A\otimes\{0 \}}{\cup}A\otimes[0,1]
\to X\otimes [0,1]
\]
is the inclusion of a retract.

\begin{prop}
\label{thm:211} The class of $h$-cofibrations is stable under
composition, and the formation of coproducts and cobase change.  Given
a sequence 
\[
X_{1}\xrightarrow{f_{1}}{} \cdots\to X_{i}\xrightarrow{f_{i}}X_{i+1}\to\cdots
\]
in which each $f_{i}$ is an $h$-cofibration, the map 
\[
X_{1}\to \varinjlim_{i}X_{i}
\]
is an $h$-cofibration.  \qed
\end{prop}

\begin{prop}
\label{thm:216}
Any topological functor $L$ which is a continuous left adjoint
preserves the class of $h$-cofibrations.  \qed
\end{prop}

Now suppose that $\cat C$ has a symmetric monoidal structure $\otimes$
which is compatible with the cartesian product of spaces, in the sense
that for spaces $S$ and $T$, and objects $X,Y\in\cat C$ there is a natural isomorphism
\[
(X\otimes S)\otimes (Y\otimes T) \approx (X\otimes Y)\otimes (S\times T)
\]
compatible with the enrichment and the symmetric monoidal structures.
Then given $i:A\to X$ we may form
\[
i^{\otimes n}:A^{\otimes n}\to X^{\otimes n}
\]
and regard it as a map in the category $\cat C^{\cat B{\Sigma_{n}}}$
of objects in $\cat C$ equipped with a $\Sigma_{n}$-action.

\begin{prop}
\label{thm:212} If $A\to X$ is an $h$-cofibration in $\cat C$, then
for any $Z$, $A\otimes Z\to X\otimes Z$ is an $h$-cofibration. \qed
\end{prop}

\begin{prop}
\label{thm:204}
If $i:A\to X$ is an $h$-cofibration then $i^{\otimes n}$ is an
$h$-cofibration in $\cat C^{\cat B{\Sigma_{n}}}$.
\end{prop}

\begin{rem}
In the category of equivariant orthogonal spectra a version of this
result appears in \cite[Lemma 15.8]{MR1806878} (where the reader is
referred to \cite[Lemma XII.2.3]{elmendorf97:_rings}).
\end{rem}

\begin{pf}
The main point is to show that the diagonal inclusion
\begin{equation}
\label{eq:124}
\cyl(A^{\otimes n}\to X^{\otimes n}) \to
\cyl(A\to X)^{\otimes n}
\end{equation}
is the inclusion of a $\Sigma_{n}$-equivariant retract.  Granting this
for the moment, one constructs a $\Sigma_{n}$-equivariant retraction
of
\[
\cyl(A^{\otimes n}\to X^{\otimes n}) \to X^{\otimes n}\otimes [0,1]
\]
as the composition
\begin{multline*}
X^{\otimes n}\otimes [0,1]\xrightarrow{1\otimes
{\text{diag}}}{}X^{\otimes n}\otimes [0,1]^{n}\approx
(X\otimes[0,1])^{\otimes n} \\
\to \cyl(A\to X)^{\otimes n}\to
\cyl(A^{\otimes n}\to X^{\otimes n}).
\end{multline*}

For the retraction of~\eqref{eq:124} start with the pushout square
\[
\xymatrix{
A\otimes \{0 \}  \ar[r]\ar[d] & A\otimes[0,1]   \ar[d] \\
X  \ar[r]        & \cyl(A\to X)
}
\]
and consider the last stage of the filtration of $\cyl(A\to
X)^{\otimes n}$ constructed in
\S\ref{sec:index-mono-prod-pushout}
\begin{equation}
\label{eq:113}
\xymatrix{
\partial_{A}(A\otimes [0,1])^{\otimes n}
  \ar[r]\ar[d]  &   
(A\otimes [0,1])^{\otimes n} \ar[d] \\
\fil_{n-1}(\cyl(A\to X)^{\otimes n})  \ar[r]        & \cyl(A\to X)^{\otimes n}\mathrlap{\ .}
}
\end{equation}
Form the $\Sigma_{n}$-equivariant map
\[
\fil_{n-1}(\cyl(A\to X)^{\otimes n})\to X^{\otimes n} \to
\cyl(A^{\otimes n} \to X^{\otimes n})
\]
using the map $\cyl(A\to X)\to X$.  To extend it to
$\fil_{n}(\cyl(A\to X)^{\otimes n})=\cyl(A\to X)^{\otimes n}$ note
that the top row of~\eqref{eq:113} can be identified with the tensor
product of the identity map of $A^{\otimes n}$ with
\[
\partial_{\{0 \}}I^{n} \to I^{n}.
\]
This identification is compatible with the action of the symmetric
group.  The desired extension is then constructed using any
$\Sigma_{n}$-equivariant retraction of $I^{n}$ to the diagonal which
takes $\partial_{\{0 \}}I^{n}$ to $\{0 \}$
\end{pf}

Working fiberwise one concludes

\begin{prop}
\label{thm:210}
Suppose that $\cat C$ is as above, and $p:I\to J$ is a covering
category.   The indexed monoidal product 
\[
\tensorp{p}:\cat C^{I}\to \cat C^{J}
\]
preserves the class of $h$-cofibrations.  \qed
\end{prop}

We end with a technical result which is useful for establishing some
of the basic homotopy theoretic properties of equivariant orthogonal
spectra, especially in connection with the monoidal geometric fixed
point functor (\S\ref{sec:geom-fixed-points}).  Though it does not
appear explicitly in the literature in this form, it is a minor
variation of~\cite[Appendix, Proposition~3.9]{LMayS} and is proved in
the same manner.

\begin{lem}
\label{thm:105}
An $h$-cofibration in $\ugspectra{G}$ is an objectwise closed inclusion. 
\end{lem}

\begin{pf}
The assertion reduces immediately to showing that that
$h$-cofibrations in the category of compactly generated weak Hausdorff
spaces are closed inclusions.  For this latter fact suppose that
$A\subset X$ is an $h$-cofibration of compactly generated weak
Hausdorff spaces.  Then the mapping cylinder $\cyl(A\to X)$ is an
equalizer of two maps from $X\times [0,1]$ to itself.  Since $X\times
[0,1]$ is weak Hausdorff, it is closed.  But $A\subset X$ is the
inverse image of $\cyl (A\to X)$ under the inclusion $X\times\{1\}\to
X\times [0,1]$.  See~\cite[Pages~488-9]{LMayS}.  
\end{pf}

\section{Homotopy theory of equivariant orthogonal spectra}
\label{sec:homot-theory-equiv}

We now turn to the stable homotopy theory of equivariant orthogonal
spectra, the basis of which is the notion of {\em stable weak
equivalence} defined in~\S\ref{sec:homot-theory-ugsp}.  Our goal is to
set up the infrastructure needed for the proofs of properties
$\gswdes{1}$-$\gswdes{6}$, and for working with the formation of
indexed wedges, smash products, symmetric powers and their
compositions.  These latter are explicit constructions, and to work
with them in homotopy theory means determining, in each case, a full
subcategory of $\ugspectra{G}$ on which the construction preserves
weak equivalences, and which is {\em homotopically wide} in the sense
that it contains at least one object of each weak equivalence class,

The standard way of doing this is to complete the set of weak
equivalences to a Quillen model category structure, in such a way that
each of the constructions takes weak equivalences between cofibrant
objects to weak equivalences.  This can be done in this case, but a
problem arises when composing these operations.  For example, in all
of the standard model structures on $\ugspectra{G}$, the symmetric
powers of a cofibrant object are not cofibrant (or at least not known
to be).  The situation is reminiscent of the theory of unbounded
operators, in which a domain of definition needs to be specified, and
in which one can run into trouble trying to compose operators.  It
might be possible to find a model structure whose collection of
cofibrant objects is preserved by all of these constructions.  But
this is more than is really required.

This is a situation where the language of model categories tends to
obscure the basic task at hand.  What is needed is to determine, for a
given functor, a homotopically wide full subcategory on which the
functor preserves weak equivalences.  This problem depends only on the
weak equivalences, and is most naturally considered in the context of
homotopical categories.  With this in mind we begin our work using
homotopical categories, where the entire focus is on weak equivalences
and derived functors, and put off introducing a model category
structure until it is really needed.

Here is a summary of the contents of this appendix.  In
\S\ref{sec:categories-with-weak-2} we review the theory of homotopical
categories.  Section~\ref{sec:flat-maps} introduces various notions of
``flatness,'' which depend only on the class of weak equivalences, and
play an important role in determining the homotopical properties of
various functors.  In \S\ref{sec:equiv-stable-homot} we develop a
considerable amount of the stable homotopy theory of $\ugspectra{G}$
using only the language of homotopical categories.  This includes most
of the results used in \S\ref{sec:prop-gswd-gswd} to verify
$\gswdes{1}$-$\gswdes{5}$.  Our analysis is facilitated by an
approximation $\pist\ugspectra{G}$ to $\ho\ugspectra{G}$ which we
study as a homotopical category in its own right.  To go further it is
helpful to have a model structure around and in
\S\ref{sec:spectra-as-model} we define the {\em positive complete
model structure} on $\ugspectra{G}$.  This is a variation on the {\em
positive stable model structure} of~\cite{MR1922205} having the
convenient property that indexed wedges and smash products of
cofibrant objects are cofibrant.  Sections~\ref{sec:homot-prop-norm}
and~\ref{sec:symmetric-powers} describe the homotopy properties of
indexed smash and symmetric powers.  Section~\ref{sec:rings-modules}
contains a proof that the forgetful functor $\comm\ugspectra{G}\to
\ugspectra{G}$ creates a model structure.  The proofs of this that
appear in the literature are incomplete, and it does not seem possible
to give a complete proof without first analyzing the homotopy
properties of indexed smash products.
Section~\ref{sec:index-smash-prod-comm-rings} contains the important
result that the formation of indexed smash products is homotopical on
a subcategory of $\ugspectra{G}$ containing both the cofibrant objects
and the spectra underlying cofibrant commutative rings.  This result
is crucial for making use of the norm functor and is part of the
reason that we work outside of the framework of model categories.
Sections~\ref{sec:geom-fixed-points}
and~\ref{sec:geom-fixed-points-norm} contain results on the geometric
fixed point functor and its interaction with the constructions
described above.  Finally, Section~\ref{sec:real-bordism} contains a
construction of the real bordism spectrum $\mur$ on which all of the
results of this paper are based.  

\subsection{Homotopical categories and model categories}
\label{sec:categories-with-weak-2}

We begin by reviewing some notions from~\cite{dwyer04:_homot}.

\begin{defin}
\label{def:7}
A {\em homotopical category} is a category $\cat C$ equipped with a
class of morphisms called {\em weak equivalences} which contains all
identity maps, and satisfies the {\em
two out of six property} described below.
\end{defin}
The two out of six property asserts that in the situation 
\[
\bullet \xrightarrow{u} \bullet \xrightarrow{v} \bullet
\xrightarrow{w} \bullet
\]
if $v u$ and $w v$ are in $\mathcal W$ then so are $u$, $v$, $w$, and
$v w u$.  It implies the ``two out of three'' property (that two of
three maps in composition being weak equivalences implies the third
is), and that isomorphisms are weak equivalences.

\begin{rem}
\label{rem:62} If the weak equivalences have the property that a map
is a weak equivalence if and only if some functor applied to the map
becomes an isomorphism, then identity maps are weak equivalences, the
two out of six property automatically holds, and retracts of weak
equivalences are weak equivalences.
\end{rem}

Suppose that $\cat C$ is a homotopical category.
\begin{defin}
\label{def:55}
A {\em homotopy functor} is a functor $F:\cat C\to \cat D$ with the
property that $F(w)$ is an isomorphism whenever $w\in\weq$.   
\end{defin}

There is a universal homotopy functor $L:\cat C\to \ho\cat C$ called
the {\em the localization of $\cat C$ with respect to $\weq$.}  It is
characterized uniquely up to unique isomorphism by the following
universal property: for every category $\cat D$, and every homotopy
functor $F:\cat C\to \cat D$ there is a unique functor $\ho\cat C\to
D$ making the diagram
\[
\xymatrix{
\cat C \ar[dr]_{F}\ar[r]^-{L} & \ho\cat C \ar[d] \\
& \cat D 
}
\]
commute.  While this characterization may seem stronger than is
natural for characterizing an arrow in a $2$-category, it simplifies
the presentation.  The difference between this and the $2$-categorical
formulation amounts to the convention that the map $\cat C\to \ho\cat
C$ be the identity map on objects.  The category $\ho\cat C$ is the
{\em homotopy category of} $\cat C$.  Since the localization functor
$L$ is the identity map on objects, it tends not to appear in
notation.

Two issues emerge when working with homotopical categories.  One is to
find a description of $\ho\cat C(X,Y)$ and the other is to describe
conditions under which a functor $F:\cat C\to \cat D$ between
homotopical categories induces a functor $\ho F:\ho\cat C\to \ho\cat
D$.  For the first question the following can be helpful.

\begin{prop}
\label{thm-copied:99} The transformation $\cat
C(X,\slot)\to \ho\cat C(X,\slot)$ is the universal natural
transformation from $\cat C(X,\slot)$ to a homotopy functor.
\end{prop}

\begin{pf}
This is one situation where it is clearer to actually make use of the
notation $L:\cat C\to\ho\cat C$.  Spelled out, the assertion is that
if $F:\cat C\to\sets $ is a homotopy functor and $\cat C(X,\slot)\to
F$ a natural transformation, then there is a unique dotted arrow
making the diagram
\begin{equation}
\label{eq:144}
\xymatrix{
\cat C(X,\slot)  \ar[r]\ar[d]_{L}  & F \\
\ho\cat C(L X,L(\slot))  \ar@{-->}[ur]        & 
}
\end{equation}
commute.  Before describing the proof we make an observation about the
property characterizing the functor $L:\cat C\to \ho\cat C$.  For
homotopy functors $F$ and $G$ on $\cat C$, this property supplies unique
factorizations $F=\tilde F\circ L$ and $G=\tilde G\circ L$.  It also
implies that composition with $L$ gives a bijection between the set of
natural transformations $\tilde G\to \tilde F$ and $G\to F$.

With this in mind we now turn to the proof of the proposition.  By the
Yoneda Lemma, the horizontal arrow in~\eqref{eq:144} is given by an
element of $F(X)$.  By the remark above, the set of natural
transformations
\[
\ho\cat C(L X, L(\slot) ) \to F
\]
is in bijection with the set of natural transformations 
\[
\ho\cat C(L X,\slot)\to \tilde F
\]
which, again by Yoneda, is in one to one correspondence with the
elements of $\tilde F(LX)=F(X)$.  The map between these sets
corresponding to the two ways of going around~\eqref{eq:144} is the
identity.
\end{pf}

\begin{cor}
\label{thm:129} 
Suppose that $\cat C$ is a homotopical category, and
that $X\in\cat C$ has the property that $\cat C(X,\slot)$ is a
homotopy functor.  Then the natural transformation $\cat C(X,\slot)\to
\ho\cat C(X,\slot)$ is a bijection.
\end{cor}

\begin{pf}
Immediate from Proposition~\ref{thm-copied:99}.
\end{pf}

For the second question, there is an apparatus of definitions to
organize the situation.

\begin{defin}
\label{def-copied:56} A functor between homotopical categories is {\em
homotopical} if it sends weak equivalences to weak equivalences.
\end{defin}

By the universal property, a homotopical functor $F:\cat C\to \cat D$
induces a functor $\ho F:\ho\cat C\to \ho\cat D$.  Furthermore,
adjoint homotopical functors induce adjoint functors on the homotopy
categories.  But there are more general situations under which such a
functor is induced.  Suppose that $F:\cat C\to \cat D$ is a functor
between homotopical categories and that one can find a subcategory
$\cat C'\subset \cat C$ on which $F$ is homotopical (where the weak
equivalences in $\cat C'$ are taken to be those morphisms which are
weak equivalences in $\cat C$).  Then $F$ induces a functor
\[
\ho\cat C'\to \ho\cat D.
\]
If, in addition, $\ho\cat C'\to \ho \cat C$ is an equivalence of
categories, then one gets an induced functor $\ho\cat C\to \ho\cat D$
by composing with an inverse to this equivalence.

The situation becomes more manageable when there is a pair $(r,s)$
consisting of a functor $r:\cat C\to \cat C$ with the property that
$F\circ r$ is homotopical, and a natural weak equivalence $s:r\to\id$.
In that case $\cat C'$ can be taken to be the full subcategory
generated by the image of $r$, the induced functor $\lder F:\ho\cat
C\to \ho\cat D$ can be computed as
\[
\lder FX = F\circ r(X),
\]
and because of $s$, comes equipped with a natural transformation
between the two ways of going around the diagram
\[
\xymatrix{
\cat C  \ar[r]^{F}\ar[d]  &  \cat D \ar[d] \\
\ho\cat C  \ar[r]_{\lder F}\ar@{=>}[ur]|-*+{\scriptstyle{T}}        &
\ho\cat D\mathrlap{\ .}
}
\]
Together with this transformation, $\lder F$ is characterized by a
universal property.  It is most easily stated if we overload some of
the notation by using the symbol $F$ to denote the composite functor
\[
\cat C\xrightarrow{F} \cat D\to\ho\cat D
\]
and identify functors $\ho\cat C\to \ho\cat D$ with homotopy
functors $\cat C\to \ho\cat D$.     With these conventions we may
regard the transformation $T$ as going from $\lder F$ to $F$
\[
T:\lder F\to F.
\]
The universal property is that if $G:\cat C\to\ho\cat D$ is a homotopy
functor and $S:G\to F$ is a natural transformation, then there is a
unique natural transformation $G\to \lder F$ making 
\[
\xymatrix{
G  \ar[r]\ar[dr]_{S}  &  \lder F \ar[d]^{T} \\
       & F
}
\]
commute.  Put differently, $\lder F$ is the closest homotopy functor
to the left of $F$.  

The functor characterized by the above properties is the {\em left
derived functor} of $F$.  It is guaranteed to exist when $F$ is {\em
left deformable} in the sense that there is a pair $(r,s)$ as above,
and $F\circ r$ is homotopical.

A common situation arises when the weak equivalences on $\cat C$
refine to a model category structure, and $F$ takes weak equivalences
between cofibrant objects to weak equivalences.  In that case $F$ is
left deformable, and one may take $(r,s)$ to be a functorial cofibrant
replacement.

There are evident dual notions of a {\em right deformable functors $F$}
and a {\em right derived functor $\rder F$}.  For more on the
definition of derived functors the reader is referred to~\cite{Qui:HA}
for the case of model categories, and
to~\cite[Chapter~VII]{dwyer04:_homot} for the more general case of
homotopical categories.

When
\[
F:\cat C\leftrightarrows \cat D:G
\]
are adjoint functors between homotopical categories, and $F$ is left
deformable and $G$ is right deformable, then the derived functors 
\[
\lder F:\ho\cat C\leftrightarrows\ho\cat D:\rder G
\]
are adjoint.   See~\cite[Chapter~VII, \S44]{dwyer04:_homot}.

It is common, when there is no confusion likely, to drop the $\lder$
from $\lder F$ and not distinguish in notation between a functor and
its derived functors.   We follow this convention in the main body of
the paper, where the emphasis is on homotopy theory.

\subsection{Flat maps} \label{sec:flat-maps} 
The notion of a {\em flat map} and a {\em flat functor} was introduced
in the unpublished manuscript~\cite{HopGoerss:Mult} in order to
isolate useful classes of maps and objects on which left derived
functors can be computed.  Though the original context involved model
categories, the definitions involve only the weak equivalences and
belong most naturally to the theory of homotopical categories.  The
dual notion was coined a ``sharp map'' by Charles Rezk, and used for a
different purpose
in~\cite{rezk:_fibrat_and_homot_colim_of_simpl_sheav}.

\begin{defin}
\label{def:45} A functor $F:\cat C\to \cat D$ between categories with
weak equivalences is {\em flat} if it is homotopical and preserves
colimits.
\end{defin}

Typically the functor $F$ will be a left adjoint, and so will
automatically preserve colimits.

\begin{defin}
\label{def:46} Suppose that $\cat C$ is a homotopical category
possessing small colimits.  A map $f:A\to X$ in $\cat C$ is {\em flat}
if for every $A\to B$ and every weak equivalence $B\to B'$, the map
\[
X\underset{A}\cup B\to X\underset{A}\cup B'
\]
is a weak equivalence.
\end{defin}

In other words a morphism $f$ is flat if and only if ``cobase change
along $f$'' preserves weak equivalences.  Since cobase change is a
left adjoint this is equivalent to the flatness of the cobase change
functor.

\begin{eg}
\label{eg:13}
A model category is {\em left proper} if and only if every cofibration
is flat.
\end{eg}

\begin{prop}
\label{thm:83}
\begin{thmList}
\item Finite coproducts of flat maps are flat.
\item Composites of flat maps are flat
\item Any cobase change of a flat map is flat.
\item If a retract of a weak equivalence is a weak equivalence then a
retract of a flat map is flat.
\end{thmList}
\qed
\end{prop}

\begin{prop}
\label{thm:158}
Suppose that
\[
\xymatrix{
X_{1} \ar[d]_{\sim}  & A_{1} \ar[l]_{\flat}\ar[r]^{\flat}\ar[d]^{\sim} & Y_{1} \ar[d]^{\sim} \\
X_{2}  & A_{2} \ar[l]\ar[r]^{\flat} & Y_{2} \\
}
\]
is a diagram in which $A_2\to Y_2$ and both maps in the top row are
flat.  If the vertical maps are weak equivalences, then so is the map
\[
X_1\underset{A_1}\cup Y_1 \to
X_2\underset{A_2}\cup Y_2
\]
of pushouts.
\end{prop}

\begin{pf}
First suppose that $A_1 = A_2 = A$.  Then 
\[
X_1\underset A{\cup} Y_1 \to
X_1\underset A{\cup} Y_2
\]
is a weak equivalence since $A\to X_1$ is flat.  The map $X_1\to
X_1\underset A{\cup} Y_2$ is flat, since it is a cobase change of
$A\xrightarrow{\flat} Y_{2}$ along $A\to X_{1}$.  But this implies
that
\[
X_1\underset A{\cup} Y_2 \to
X_2\underset{X_1}\cup \left(X_1\underset A\cup Y_2\right) = 
X_2\underset A{\cup} Y_2.
\]
is a weak equivalence.  Putting these together gives the result in
this case.

For the general case, consider the following diagram
\[
\xymatrix{
X_1 \ar[d] & A_1 \ar[r]\ar[l]\ar[d] & Y_1 \ar[d] \\
X_1\underset{A_1}{\cup}A_2 \ar[d] & A_2 
   \ar[r]\ar[l]\ar@{=}[d] & A_2\underset{A_1}\cup Y_1 \ar[d] \\ 
X_2 & A_2 \ar[r]\ar[l] & Y_2\mathrlap{\ .}  
}
\]
The flatness of the maps $A_1\to X_1$ and $A_1\to Y_1$ implies that
the upper vertical maps (hence all the vertical maps) are weak
equivalences, and that the maps in the middle row are flat.  It also
implies that
\[
A_1\to X_1\underset{A_1}{\cup}Y_1
\]
is flat.  Since $A_1\to A_2$ is a weak equivalence, this means that
\[
X_1\underset{A_1}{\cup} Y_1\to A_2\underset{A_1}{\cup}\left(
X_1\underset{A_1}{\cup} Y_1\right)
\]
is a weak equivalence.  But this is the map from the pushout of the
top row to the pushout of the middle row.  By the case in which
$A_1=A_2$, the map from the pushout of the middle row to the pushout
of the bottom row is a also a weak equivalence.  This completes the
proof.
\end{pf}

\begin{rem}
\label{rem:52} 
If $\cat{C}$ has the property that every map can be
factored into a flat map followed by a weak equivalence, then the
above result holds with the assumption that only one of the maps in
the top row is a weak equivalence.  Suppose for instance that it is
the map $A_1\to X_1$, and factor $A_1\to Y_1$ into a flat map $A_1\to
Y'_1$ followed by a weak equivalence $Y'_1\to Y_1$ Now consider the
diagram
\[
\xymatrix{
X_1 \ar@{=}[d] & \ar[l]_{\flat} A_1 \ar[r]^{\flat}\ar@{=}[d]  & Y'_1\ar[d]^{\sim} \\
X_1\ar[d]_{\sim} & \ar[l]_{\flat} A_1 \ar[r]\ar[d]^{\sim} & Y_1 \ar[d]^{\sim}\\
X_2 & \ar[l] A_2 \ar[r]_{\flat} & Y_2 \mathrlap{\ .}  
}
\]
By Proposition~\ref{thm:158}, the map from the pushout of the top row
to the pushout of the middle row is a weak equivalence, as is the map
from the pushout of the top row to the pushout of the bottom row.  The
map from the pushout of the middle row to the pushout of the bottom
row is then a weak equivalence by the two out of three property of
weak equivalences.
\end{rem}

\begin{rem}
\label{rem:7} In the category $\ugspectra{G}$ equipped with the stable
weak equivalences (\ref{sec-copied:homot-theory-ugsp}), the
$h$-cofibrations will turn out to be flat.  The mapping cylinder
construction then factors every map into a flat map followed by a weak
equivalence, so Remark~\ref{rem:52} applies.
\end{rem}

Now suppose that $(\cat C,\otimes, \one)$ is a closed symmetric
monoidal category, equipped with a class $\mathcal W$ of weak
equivalences, making $\cat C$ into a homotopical category.

\begin{defin}
\label{def:59}
An object $X\in \cat C$ is {\em flat} if the functor $X\otimes
(\slot)$ is flat.   
\end{defin}

Showing that a symmetric monoidal structure on $\cat C$ induces one on
$\ho\cat C$ essentially comes down to exhibiting enough flat objects
in $\cat C$.  In \S\ref{sec:weak-equiv-smash} we will show that the
cellular objects of $\ugspectra{G}$ are flat.

\begin{rem}
\label{rem:67} Suppose that every object $Z\in\cat C$ admits a weak
equivalence equivalence $\tilde Z\to Z$ from a flat object $\tilde Z$.
If $X\to Y$ is a weak equivalence of flat objects, so is $X\otimes Z\to
Y\otimes Z$ for any $Z$.  This follows from the diagram
\[
\xymatrix{
X\otimes \tilde Z  \ar[r]^{\sim}\ar[d]_{\sim}  & X\otimes Z
\ar[d] \\
Y\otimes \tilde Z  \ar[r]_{\sim}        & Y\otimes Z\mathrlap{\ .}
}
\]
\end{rem}
\subsection{Equivariant stable homotopy theory}
\label{sec:equiv-stable-homot} 

The weak equivalences were defined
in~\S\ref{sec:homot-theory-ugsp} as the maps inducing isomorphisms of
stable homotopy groups.  Equipped with them $\ugspectra{G}$ becomes a
homotopical category, and the functor $\ugspectra{G}\to
\ho\ugspectra{G}$ is defined.  In this section we establish many of
the basic properties of $\ho\ugspectra{G}$, including most of the
results used in \S\ref{sec:prop-gswd-gswd} to verify
$\gswdes{1}$-$\gswdes{5}$ of \S\ref{sec:g-spectra}.  

\subsubsection{Stable weak equivalences and basic homotopical functors}
\label{sec-copied:homot-theory-ugsp} 

We begin with some basic homotopical functors.  
\begin{prop}
\label{thm:255}
The formation of filtered colimits along objectwise closed
inclusions is homotopical. 
\end{prop}

\begin{pf}
This is immediate from the fact that formation of homotopy groups
commutes with filtered colimits of closed inclusions.
\end{pf}

Since $h$-cofibrations are objectwise closed inclusions
(Lemma~\ref{thm:105}), Proposition~\ref{thm:255} applies to the
formation of filtered colimits along $h$-cofibrations.

The following three results, which are part of~\cite[Theorem
III.3.5]{MR1922205} (see also~\cite[Theorem 7.4~\thmListItem{4}]{MR1806878}),
imply that many basic functors are homotopical. 

\begin{prop}
\label{thm:137}
Suppose $f:X\to Y$ is a map and let $F\to X$ be the homotopy fiber,
defined by the pullback square
\[
\xymatrix{
F  \ar[r]\ar[d]  & PY
\ar[d] \\ 
X  \ar[r]_{f}        & Y
}
\]
in which $PY$ is the path spectrum of $Y$.  There is a long exact
sequence
\[
\cdots\to \pi^{H}_{k}F\to \pi^{H}_{k}X\to \pi^{H}_{k}Y \to
\pi^{H}_{k-1} F\to\cdots
\]
\end{prop}

\begin{pf*}{Sketch of proof}
This sequence is gotten by passing to the colimit from the exact sequence
\[
\cdots\to \pi^{H}_{k+V}F_{V}\to \pi^{H}_{k+V}X_{V}\to \pi^{H}_{k+V}Y_{V} \to
\pi^{H}_{k-1+V} F_{V}\to\cdots.
\]
\end{pf*}

\begin{prop}
\label{thm:53}
For any $X$, any $H\subset G$, and any $k\in\Z$ the suspension map 
\[
\pi^{H}_{k}X\to \pi^{H}_{k+1}S^{1}\wedge X
\]
is an isomorphism.
\end{prop}

\begin{pf*}{Sketch of proof}
Choose an exhausting sequence $\{V_{n} \}$ with the property that
$V_{n}\oplus \R\subset V_{n+1}$.  Then the map
$\pi^{H}_{k+V_{n}}X_{V_{n}}\to \pi^{H}_{k+V_{n+1}}X_{V_{n+1}}$ factors
through the suspension map $\pi^{H}_{k+1+V_{n}}S^{1}\wedge X_{V_{n}}$,
and so the sequence for computing $\pi^{H}_{k+1}S^{1}\wedge X$ threads
through the sequence for computing $\pi^{H}_{k}X$.
\end{pf*}

\begin{prop}
\label{thm:116} Let $X\to Y$ be an $h$-cofibration.
\begin{thmList}
\item  The map $Y\cup
CX\to Y/X$ is a weak equivalence.  
\item There is a natural long exact sequence of stable homotopy groups
\[
\dots\to \pi^{H}_{k}X\to \pi^{H}_{k}Y\to \pi^{H}_{k}(Y/X) \to \pi^{H}_{k-1}X\to\dots,
\]
in which the map $\pi_{k}Y\to \pi_{k}Y/X$ is induced by the evident
quotient map, and the connecting homomorphism $\pi^{H}_{k}Y/X\to
\pi^{H}_{k-1}X$ is induced by the maps
\[
Y/X  \leftarrow Y\cup CX \rightarrow \Sigma X.
\]
and the suspension isomorphism of Proposition~\ref{thm:53}.
\end{thmList}
\end{prop}

\begin{pf*}{Sketch of proof}
For the first part, since $A\to X$ is an $h$-cofibration, the map
$X\cup CA\to X/A$ is a homotopy equivalence, hence induces an
isomorphism of stable homotopy groups.  The result can then be deduced
from Proposition~\ref{thm:53} as in~\cite[III.2.1]{LMayS}.
\end{pf*}

\begin{cor}
\label{thm:225}
The $h$-cofibrations in $\ugspectra{G}$ are flat. \qed
\end{cor}

Proposition~\ref{thm:116} implies that the formation of mapping cones
is homotopical as is the formation of quotients of $h$-cofibrations.
It also gives parts \thmListItem{1} and \thmListItem{3} of the
Proposition below.  Part \thmListItem{2} follows from the fact that
the formation of unstable homotopy groups commutes with products and
the fact that filtered colimits commute with finite products.

\begin{prop}
\label{thm-copied:91}
\begin{thmList}
\item For any any  set of
spectra $\{X_{\alpha} \}$ the map 
\[
\bigoplus \pi^{G}_{\ast}X_{\alpha}\to \pi^{G}_{\ast}\bigvee X_{\alpha}
\]
is an isomorphism, hence the formation of wedges is homotopical.
\item For any any finite set of
spectra $\{X_{\alpha} \}$ the map 
\[
\pi^{G}_{\ast}\prod X_{\alpha} \to \prod \pi^{G}_{\ast}X_{\alpha}
\]
is an isomorphism, hence the formation of finite products is homotopical.
\item For any finite set of
spectra $\{X_{\alpha} \}$ the map 
\[
\bigvee X_{\alpha}\to \prod X_{\alpha}
\]
is a weak equivalence.
\end{thmList}
\qed
\end{prop}

\begin{cor}
\label{thm:90}
The category $\ho\ugspectra{G}$ is additive, and admits finite
products and arbitrary coproducts.  The coproducts are given by wedges
and the finite products by finite products.
\end{cor}

\begin{pf}
Let's begin with the case of coproducts.  Let $J$ be a set.  The
adjoint functors
\[
\bigvee:\big(\ugspectra{G}\big)^{J} \leftrightarrows \ugspectra{G}: \text{diag}
\]
are homotopical by Proposition~\ref{thm-copied:91}.  They therefore
induce adjoint functors 
\[
\bigvee:\big(\ho\ugspectra{G}\big)^{J} \leftrightarrows \ho\ugspectra{G}: \text{diag}
\]
on the homotopy categories.  This shows that arbitrary coproducts
exist in $\ho\ugspectra{G}$ and that they may be computed as wedges.
A similar argument shows that finite products exist, are computed as
products in $\ugspectra{G}$, and that the map from a finite coproduct
to a finite product is an isomorphism.  This endows the morphism sets
in $\ho\ugspectra{G}$ with the structure of commutative monoids.  That
they are in fact abelian groups can be seen by checking that for all
$X$, the ``shearing map'' $X\vee X\to X\times X$, with first component
the projection to the first summand and second component the coproduct
of the identity map with itself, is a weak equivalence.
\end{pf}

The ``indexed'' analogue of Proposition~\ref{thm-copied:91} is also
true, and appears as Proposition~\ref{thm-copied:246}.    It expresses
a kind of ``equivariant additivity'' on $\ho\ugspectra{G}$.

\subsubsection{Suspension and zero space}
\label{sec:susp-zero-space}

The suspension and zero space functors were defined in
Definition~\ref{def:6}.  Formation of the suspension spectrum is
nearly homotopical.

\begin{prop}
\label{thm:258} The suspension spectrum functor is homotopical on the
subcategory of non-degenerately based $G$-spaces.  The right derived
functor $\rder\Omega^{\infty}X$ may be computed as
\[
\rder\Omega^{\infty}X = \ho \varinjlim \Omega^{V_{n}}X_{V_{n}}
\]
where $\{V_{n}\}$ is any choice of exhausting sequence,
and $\Omega^{V_{n}}(\slot)$ is the $G$-space of non-equivariant maps.
\end{prop}

\begin{pf}
The assertion about suspension spectra follows from the fact that if
$K\to L$ is an equivariant weak equivalence of non-degenerately based
$G$-spaces, then so is
\[
S^{V}\wedge K\to S^{V}\wedge L
\] 
for any representation $V$.  This reduces to the statement that for
every $H\subset G$, the map
\[
S^{V^{H}}\wedge K^{H}\to
S^{V^{H}}\wedge L^{H}
\]
is a weak equivalence, assuming $K^{H}\to L^{H}$ is.  But this is a
standard fact.  The functor $\rder\Omega^{\infty}X = \hocolim
\Omega^{V_{n}}X_{V_{n}}$ is clearly homotopical, so what is needed for
the second assertion is to construct a functorial weak equivalence
$X\to X'$, in which $X'$ has the property that the map
\[
\rder \Omega^{\infty}X' \to \Omega^{\infty}X' 
\]
is a weak equivalence.  One way to do this is to define $X\to X'$ by
$X_{V}\to \hocolim_{n}\Omega^{V_{n}}X_{V\oplus V_{n}} = X'_{V}$.
(One can also take $X'$ to be the functorial fibrant replacement coming
from the small object construction in the positive complete
model structure of \S\ref{sec:strong-posit-stable}.)
\end{pf}

Adding a ``whisker'' provides a left deformation to $\Sigma^{\infty}$,
and the natural transformation $X\to X'$ appearing in the proof above
gives a right deformation of $\Omega^{\infty}$.  The derived
suspension spectrum and zero space functors therefore induce adjoint
functors on the homotopy categories
\[
\lder\Sigma^{\infty}:\ho\ugspaces{G} \leftrightarrows
\ho\ugspectra{G}:\rder\Omega^{\infty}.
\]

\subsubsection{An approximation to the homotopy category}
\label{sec-copied:relat-span-whit} 

Our further analysis of $\ho\ugspectra{G}$ is facilitated by an
approximation, $\pist\ugspectra{G}$.

Let
\begin{equation}
\label{eq-copied:166}
\varepsilon_{V}:S^{-V}\wedge S^{V}\to S^{0}
\end{equation}
be the map adjoint to the identity map of $S^{V}$.  Associated to a
linear isometric embedding $t:V\to W$ is a map
\begin{equation}
\label{eq:162}
S^{-W}\wedge S^{W} \to S^{-V}\wedge S^{V}
\end{equation}
One way to describe it is to note that the space of such maps is the
space of equivariant maps
\[
S^{W}\to \big(S^{-V}\wedge S^{V})_{W},
\]
and that
\[
\big(S^{-V}\wedge S^{V})_{W} \approx \thom(O(V,W);(W-V)\oplus V)
\approx O(V,W)_{+}\wedge S^{W}.
\]
The map~\eqref{eq:162} corresponds to smashing the identity map of
$S^{W}$ with the map $S^{0}\to O(V,W)^{G}_{+}$ sending the non-base
point to $t$.  The map~\eqref{eq:162} can also be expressed as
$\id\wedge \epsilon_{U}$ after 
rewriting the domain as
\[
S^{-V}\wedge S^{V}\wedge S^{-U}\wedge S^{U},
\]
with $U=W-t(V)$.  When $V<W$ the fixed point space $O(V,W)^{G}$ is
connected, and so the homotopy class~\eqref{eq:162} is independent of
the choice of $t$.

For $X, Y\in\ugspectra{G}$ let
\[
\pist\ugspectra{G}(X,Y) = \varinjlim_{V}\pi_{0}\ugspectra{G}(S^{-V}\wedge
S^{V}\wedge X, Y),
\]
in which the limit is taken over the partially ordered set of
representations of $G$ (\S\ref{sec:homot-theory-ugsp}).  We wish to
make $\pist\ugspectra{G}(X,Y)$ into the morphisms in a category.  For
this we need to define the composition law.  An element
$f\in\pist\ugspectra{G}$ is represented by a map $f_{V}:S^{-V}\wedge
S^{V}\wedge X\to Y$.  Given $f\in \pist\ugspectra{G}(X,Y)$ and
$g\in\pist\ugspectra{G}(Y,Z)$ represented by
\begin{align*}
f_{V}:S^{-V}\wedge S^{V}\wedge X &\to Y  \\
g_{W}:S^{-W}\wedge S^{W}\wedge Y &\to Z 
\end{align*}
the composition $g\circ f$ is defined to be the equivalence class of
the map 
\[
(g\circ f)_{W\oplus V}:S^{-W\oplus V}\wedge S^{W\oplus V}\wedge X\to Z
\]
constructed from the isomorphism 
\[
S^{-W\oplus V}\wedge S^{W\oplus V}\approx
S^{-W}\wedge S^{W}\wedge S^{-V}\wedge S^{V}
\]
and the composite
\[
S^{-W}\wedge S^{W}\wedge S^{-V}\wedge S^{V}\wedge X
\xrightarrow{\id\wedge f_{V}}
S^{-W}\wedge S^{W}\wedge Y \xrightarrow{g_{W}} Z.
\]
Associativity of the composition follows from the associativity of the
smash product.   

\begin{defin}
\label{def:9} The category $\pist\ugspectra{G}$ is the category whose
objects are those of $\ugspectra{G}$, with morphisms
$\pist\ugspectra{G}(X,Y)$, and the composition law described above.
\end{defin}

One thing that makes $\pist\ugspectra{G}$ so useful is that the hom
sets are easy to describe, and yet the functors $\pi^{H}_{k}$ factor
through it and are corepresentable.

\begin{prop}
\label{thm:88}
For all $k\in\Z$, there is a natural isomorphism
\begin{equation}
\label{eq:169}
\pist\ugspectra{G}(G/H_{+}\wedge S^{k}, Y) \approx \pi^{H}_{k}(Y).
\end{equation}
\end{prop}

\begin{pf}
Suppose $k\ge 0$.   Then
\begin{align*}
\pist\ugspectra{G}(G/H_{+}\wedge S^{k},Y) &= \varinjlim \pi_{0}\ugspectra{G}(S^{-V}\wedge
S^{V}\wedge G/H_{+}\wedge S^{k},Y) \\
& = \varinjlim \pi_{0}\ugspectra{H}(S^{-V}\wedge
S^{V}\wedge S^{k},Y) \\
&= \varinjlim \pi_{0}\ugspaces{H}(S^{V}\wedge S^{k},Y_{V}) \\
& = \varinjlim \pi^{H}_{k+V}Y_{V}
=\pi^{H}_{k}Y.
\end{align*}
Similarly, 
\begin{align*}
\pist\ugspectra{G}(G/H_{+}\wedge S^{-k},Y) &= \varinjlim \pi_{0}\ugspectra{G}(S^{-V}\wedge
S^{V}\wedge G/H_{+}\wedge S^{-k},Y) \\
& = \varinjlim \pi_{0}\ugspectra{H}(S^{-V}\wedge
S^{V}\wedge S^{-k},Y) \\
&= \varinjlim \pi_{0}\ugspaces{H}(S^{V},Y_{V+k}) \\
& = \varinjlim_{V} \pi^{H}_{V}Y_{V+k}
 = \varinjlim_{W>k} \pi^{H}_{W-k}Y_{W}
 = \pi^{H}_{-k}Y.
\end{align*}
\end{pf}

Proposition~\ref{thm:88} implies that a map $X\to Y\in\ugspectra{G}$
which becomes an isomorphism in $\pist\ugspectra{G}$ is a weak
equivalence.  An important example is

\begin{prop}
\label{thm:159}
Suppose that $V$ is a representation of $G$.   For every $X$, the map 
\begin{equation}
\label{eq:163}
S^{-V}\wedge S^{V}\wedge X\to X
\end{equation}
is an isomorphism in $\pist\ugspectra{G}$, hence a weak equivalence.
\end{prop}

\begin{pf}
We will show that for all $Y$, the map 
\[
\pist\ugspectra{G}(X,Y) \to
\pist\ugspectra{G}(S^{-V}\wedge S^{V}\wedge X,Y)
\]
is an isomorphism.   By definition,
\begin{equation}
\label{eq:173}
\pist\ugspectra{G}(X,Y) = \varinjlim_{W}(S^{-W}\wedge S^{W}\wedge X,Y),
\end{equation}
while 
\[
\pist\ugspectra{G}(S^{-V}\wedge S^{V}\wedge X,Y) =
\varinjlim_{U}\pi_{0}\,\ugspectra{G}(S^{-U}\wedge S^{U}\wedge S^{-V}\wedge S^{V}\wedge X,Y).
\]
Writing $W=U\oplus V$ and using the identification $S^{-W}\wedge
S^{W}\approx S^{-U}\wedge S^{U}\wedge S^{-V}\wedge S^{V}$, this last
colimit may be replaced by
\[
\varinjlim_{W>V}\pi_{0}\ugspectra{G}(S^{-W}\wedge S^{W}\wedge X,Y),
\]
since the set $\{U\mid U\oplus V >V \}$ is cofinal in the poset of all
representations.  But this clearly coincides with~\eqref{eq:173},
since $\{W\mid W>V \}$ is also cofinal in the poset of
representations.
\end{pf}

\begin{rem}
\label{rem:63}
The weak equivalence~\eqref{eq:163} is often written in the
form 
\[
S^{-V\oplus W}\wedge S^{W}\wedge X\to S^{-V}\wedge X.
\]
This is gotten from~\eqref{eq:163} by writing $S^{-V\oplus W}$
as $S^{-V}\wedge S^{-W}$ and writing the map as
\[
S^{-W}\wedge S^{W}\wedge \big(S^{-V}\wedge X\big) \to 
\big(S^{-V}\wedge X\big).
\]
\end{rem}

\begin{cor}
\label{thm:130}
Suppose that $V$ is a representation of $G$.   Smashing with $S^{V}$
and $S^{-V}$ are inverse equivalences of $\pist\ugspectra{G}$. \qed
\end{cor}

\begin{rem}
\label{rem:26} Corollary~\ref{thm:130} does not directly imply the
analogous statement for $\ho\ugspectra{G}$.  For that one needs to
know that smashing with $S^{V}$ and $S^{-V}$ are homotopical.  This
will be proved in \S\ref{sec:equiv-addit}.
\end{rem}

One consequence of Corollary~\ref{thm:130} is that
$\pist\ugspectra{G}$ is tensored over the equivariant
Spanier-Whitehead category $\swg{G}$ defined in
\S\ref{sec:g-spectra}.  The main point is to show that a map $K\to L$
in $\swg{G}$ gives a natural map $X\wedge K\to X\wedge L$ in
$\pist\ugspectra{G}$.  For this, suppose that the map $K\to L$ is
represented by a map of spaces $S^{V}\wedge K\to S^{V}\wedge L$.  This
latter map gives us an element of 
\[
\pist\ugspectra{G}(X\wedge
S^{V}\wedge K, X\wedge S^{V}\wedge L)
\]
and hence an element of
$\pist\ugspectra{G}(X\wedge K,Y\wedge L)$ under the isomorphism of
Corollary~\ref{thm:130}.

This fact leads to a form of Spanier-Whitehead duality in
$\pist\ugspectra{G}$.  Suppose that $K$ is a finite $G$-CW complex,
and that $L$ is a ``$V$-dual'' in the sense that there is a
representation $V$ of $G$ and maps in $\swg{G}$
\begin{align*}
K\wedge L&\to S^{V} \\
S^{V} &\to L\wedge K
\end{align*}
with the property that the composites
\begin{gather*}
S^{V}\wedge L\to L\wedge K\wedge L \to L\wedge S^{V} \\
K\wedge S^{V}\to K\wedge L\wedge S^{V}\to S^{V}\wedge K
\end{gather*}
are the symmetry isomorphism.   Then for $X, Y\in\pist\ugspectra{G}$
the composite
\begin{multline}
\label{eq:164}
\pist\ugspectra{G}(X,Y\wedge K) \to 
\pist\ugspectra{G}(X\wedge L,Y\wedge K\wedge L) \\ \to 
\pist\ugspectra{G}(X\wedge L, Y\wedge S^{V})	  \approx
\pist\ugspectra{G}(S^{-V}\wedge X\wedge L, Y)
\end{multline}
is an isomorphism, by the standard duality manipulation.   

Given $X\to Y\in\ugspectra{G}$, and any $Z$ there is a long exact
sequence
\begin{multline}
\label{eq:166}
\cdots\to \pist\ugspectra{G}(Z,S^{k}\wedge X)
\to \pist\ugspectra{G}(Z,S^{k}\wedge Y)
\to \pist\ugspectra{G}(Z,S^{k}\wedge (Y\cup CX)) \\
\to \pist\ugspectra{G}(Z,S^{k+1}\wedge X)\to\cdots.
\end{multline}
As in the proof of~\ref{thm:116}, this is proved with the argument
of~\cite[III.2.1]{LMayS}, using the analogue of
Proposition~\ref{thm:53} given as the special case of
Corollary~\ref{thm:130} in which $V$ is trivial.

There is also an easier long exact sequence in the other variable.
Let $A\to X$ be a map in $\ugspectra{G}$ and $Y$ any spectrum.  Then
there is a long exact sequence
\begin{equation}
\label{eq:168}
\cdots 
\to \pist\ugspectra{G}(S^{k}\wedge (X\cup CA),Y)
\to \pist\ugspectra{G}(S^{k}\wedge X,Y) 
\to \pist\ugspectra{G}(S^{k}\wedge A,Y) \to\cdots.
\end{equation}
Under the isomorphism given by Proposition~\ref{thm:88}, this is the long
exact sequence of Proposition~\ref{thm:137} associated to the
fibration sequence of function spectra
\[
Y^{X\cup CA} \to Y^{X}\to Y^{A}.
\]

\subsubsection{$\pist\ugspectra{G}$ as a homotopical category} 
\label{sec:pist-as-homot} 

We now study $\pist\ugspectra{G}$ as a homotopical category, and in
doing so establish the fact that the functor $\swg{G}\to
\ho\ugspectra{G}$ is fully faithful.

By Proposition~\ref{thm:88} the functors $\pi^{H}_{k}$ factor through
$\pist\ugspectra{G}$.  We make $\pist\ugspectra{G}$ into a homotopical
category by defining a map to be a weak equivalence if it induces an
isomorphism in $\pi^{H}_{k}$ for all $H\subset G$ and all $k\in\Z$.
Since a map in $\ugspectra{G}$ is a weak equivalence if and only if it
is so in $\pist\ugspectra{G}$, the canonical functor
\begin{equation}
\label{eq:199}
\ho\ugspectra{G} \to \ho\pist\ugspectra{G}
\end{equation}
is an isomorphism.  Corollary~\ref{thm:129} asserts that if $X\in\ugspectra{G}$ happens to
have the property that $\pist\ugspectra{G}(X,\slot)$ is a homotopy
functor, then $\pist\ugspectra{G}(X,\slot)\to
\ho\pist\ugspectra{G}(X,\slot)$ is an isomorphism.   Combining this
with the isomorphism~\eqref{eq:199} gives

\begin{lem}
\label{thm:249} 
If $X\in\ugspectra{G}$ has the property that
$\pist\ugspectra{G}(X,\slot)$ is a homotopy functor, then for all $Y$,
the maps
\begin{equation}
\label{eq-copied:165}
\pist\ugspectra{G}(X,Y) \xrightarrow{}
\ho\pist\ugspectra{G}(X,Y)\xleftarrow{\sim}\ho\ugspectra{G}(X,Y)
\end{equation}
are isomorphisms, and so $\ho\ugspectra{G}(X,Y)$ may be computed as
$\pist\ugspectra{G}(X,Y)$.
\end{lem}

\begin{prop}
\label{thm-copied:251}
For $k\in\Z$ the maps~\eqref{eq:169} and~\eqref{eq-copied:165} give isomorphisms
\[
\pi^{H}_{k}X\approx\pist\ugspectra{G}(G/H_{+}\wedge
S^{k},X)\approx\ho\ugspectra{G}(G/H_{+}\wedge S^{k},X).
\]
\end{prop}

\begin{pf}
The first isomorphism is given by Proposition~\ref{thm:88}, and it implies
that $\pist\ugspectra{G}(G/H_{+}\wedge S^{k},X)$ is a homotopy
functor of $X$.  Lemma~\ref{thm:249} then gives the second isomorphism. 
\end{pf}

\begin{cor}
\label{thm-copied:253}
A map $X\to Y$ in $\ugspectra{G}$ is a weak equivalence if and only if
it becomes an isomorphism in $\ho\ugspectra{G}$. \qed
\end{cor}

\begin{prop}
\label{thm-copied:244} When $X$ is of the form $X=S^{\ell}\wedge K$
with $K$ a finite $G$-CW complex, and $\ell\in\Z$, the functor
$\pist\ugspectra{G}(X,\slot)$ is a homotopy functor, and so for all $Y$
$\ho\ugspectra{G}(X,Y)$ may be computed as $\pist\ugspectra{G}(X,Y)$.
\end{prop}

\begin{pf}
Working through the skeletal filtration of $K$ and using the exact
sequence~\eqref{eq:168} reduces the claim to the case in which
$K=G/H_{+}\wedge S^{n}$.  But that case is Corollary~\ref{thm:88}.
\end{pf}

Note that 
\[
\pist\ugspectra{G}(S^{0}\wedge K,S^{0}\wedge L) = 
\varinjlim \pi_{0}\,\ugspaces{G}(S^{V}\wedge K,S^{V}\wedge
L).
\]
When $L$ is a finite $G$-CW complex, this is the definition of
$\swg{G}(K,L)$.  Thus Proposition~\ref{thm-copied:244} contains as a
special case

\begin{prop}
\label{thm:243}
The functor $\Sigma^{\infty}$ induces a fully faithful embedding
$\swg{G}\to \ho\ugspectra{G}$.\qed
\end{prop}

\subsubsection{Equivariant additivity}
\label{sec:equiv-addit}

Our next goal is to show that the formation of indexed wedges in
$\ugspectra{G}$ is homotopical.  We will do this, as
in~\cite{MR764596}, via a Spanier-Whitehead duality argument.  To make
this work we need to show that smashing with $S^{V}$ and $S^{-V}$ are
homotopical.  As mentioned in Remark~\ref{rem:26}, this implies that
they induce inverse functors on $\ho\ugspectra{G}$.  It also lays
the groundwork for our investigation of the homotopical properties of
the smash product in \S\ref{sec:weak-equiv-smash}.

\begin{lem}
\label{thm:136}
For a map $X\to Y$ in $\pist\ugspectra{G}$, the following are
equivalent
\begin{thmList}
\item The map $X\to Y$ is a weak equivalence.
\item For all $H\subset G$ and all $k\in\Z$ the map 
\[
\pist(G/H_{+}\wedge S^{k},X)\to 
\pist(G/H_{+}\wedge S^{k},Y)
\]
is an isomorphism.
\item For {\em some} representation $V$ of $G$, all $H\subset G$ and all $k\in\Z$ the map 
\[
\pist(G/H_{+}\wedge S^{k}\wedge S^{V},X)\to 
\pist(G/H_{+}\wedge S^{k}\wedge S^{V},Y)
\]
is an isomorphism.
\item For {\em all} representations $V$ of $G$, all $H\subset G$ and all $k\in\Z$ the map 
\[
\pist(G/H_{+}\wedge S^{k}\wedge S^{V},X)\to 
\pist(G/H_{+}\wedge S^{k}\wedge S^{V},Y)
\]
is an isomorphism.
\end{thmList}
\end{lem}

\begin{pf}
The equivalence of the first two statements in
Proposition~\ref{thm-copied:251}, and they imply the fourth by
Proposition~\ref{thm-copied:244}.  The fourth statement obviously
implies the third.  That the third statement implies the first two
is proved by induction on $|G|$, the assertion being trivial when $G$
is trivial.  We may therefore assume that part \thmListItem{3} holds,
and that part \thmListItem{2} holds for all proper $H\subset G$.  Let
$V_{0}\subset V$ be the subspace of invariant vectors.  Using the long
exact sequence~\eqref{eq:168}, and working by downward induction
through an equivariant cell decomposition of $S^{V}$, one sees that
for all $k\in\Z$ and all $H\subset G$, our assumptions imply that the
map
\[
\pist(G/H_{+}\wedge S^{k}\wedge S^{V_{0}},X)\to 
\pist(G/H_{+}\wedge S^{k}\wedge S^{V_{0}},Y)
\]
is an isomorphism.   But in $\pist\ugspectra{G}$ there is an
isomorphism  $S^{k}\wedge S^{V_{0}}\approx S^{k+\ell}$ with $\ell=\dim
V_{0}$, so this implies part \thmListItem{2}.
\end{pf}

We next show that both smashing with $S^{V}$ and smashing with $S^{-V}$
are homotopical functors.   Combined with Corollary~\ref{thm:130} this
implies that they induce inverse equivalences of $\ho\ugspectra{G}$.

\begin{prop}
\label{thm-copied:254}
Let $V$ be a representation of $G$.   The following conditions on a
map $X\to Y\in\pist\ugspectra{G}$ are equivalent
\begin{thmList}
\item  The map $X\to Y$ is a weak 
equivalence
\item  The map $S^{V}\wedge X\to S^{V}\wedge Y$ is a weak equivalence
\item The map $S^{-V}\wedge X\to S^{-V}\wedge Y$ is a weak equivalence.
\end{thmList}
\end{prop}

\begin{pf}
Since smashing with $S^{V}$ is the inverse equivalence of smashing
with $S^{-V}$ it suffices to establish the equivalence of the first
two assertions.  Now for any $X$, smashing with $S^{V}$ gives an
isomorphism
\[
\pist(G/H_{+}\wedge S^{k},S^{-V}\wedge X) \approx\pist(G/H_{+}\wedge
S^{k}\wedge S^{V},X),
\]
so the equivalence of the first two assertions is a consequence of Lemma~\ref{thm:136}.
\end{pf}

\begin{cor}
\label{thm:160}
Suppose that $V$ is a representation of $G$.   Smashing with $S^{V}$
and $S^{-V}$ are inverse equivalences of $\ho\ugspectra{G}$. \qed
\end{cor}

With Proposition~\ref{thm-copied:254} in place, we have the following
generalization of Proposition~\ref{thm-copied:244}.

\begin{prop}
\label{thm:161} When $X$ is of the form $X=S^{-V}\wedge K$, with $K$ a
finite $G$-CW complex, the functor $\pist\ugspectra{G}(X,\slot)$ is a
homotopy functor, hence 
\[
\pist\ugspectra{G}(X,\slot)\to
\ho\ugspectra{G}(X,\slot)
\]
is an isomorphism.
\end{prop}

\begin{pf}
By Corollary~\ref{thm:130} there is  an isomorphism.
\[
\pist\ugspectra{G}(S^{-V}\wedge K, (\slot)) \approx
\pist\ugspectra{G}(K, S^{V}\wedge (\slot)).
\]
But $S^{V}\wedge (\slot)$ is a homotopy functor by
Proposition~\ref{thm-copied:254}, and $\pist\ugspectra{G}(K,
(\slot))$ is a homotopy functor by~Proposition~\ref{thm-copied:244}. 
\end{pf}

Expanded out, Proposition~\ref{thm:161} gives the formula
\[
\ho\ugspectra{G}(S^{-V}\wedge K,Y) = \varinjlim_{W}[S^{W}\wedge
K,Y_{V\oplus W}]^{G}
\]
advertised in \S\ref{sec:homot-theory-ugsp} as~\eqref{eq:180}.   Taking
$S^{-V}\wedge K$ to be $S^{-V}\wedge S^{k}\wedge G/H_{+}$, $k\in
\Z$, this specializes to the isomorphism
\begin{equation}
\label{eq:179}
\ho\ugspectra{G}(S^{-V}\wedge S^{k}\wedge G/H_{+},X) \approx
\varinjlim_{W>-k}\pi^{H}_{W+k}X_{V\oplus W}.
\end{equation}
In particular, the expression
$\varinjlim_{W>-k}\pi^{H}_{W+k}X_{V\oplus W}$ is a homotopy functor of
$X$.  This fact is used in the proof of Proposition~\ref{thm:145},
which plays a fundamental role in establishing the positive complete
stable model category structure on $\ugspectra{G}$.

The fact that $\pist\ugspectra{G}$ is tensored over $\swg{G}$ also
gives control over homotopical properties of the smash product and of
indexed wedges.

\begin{cor}
\label{thm:132} If $X$ is of the form $S^{-W}\wedge K$, with $K$ a
$G$-CW complex and $W$ a representation of $G$, then the functor 
\[
(\slot)\wedge X:\ugspectra{G}\to \ugspectra{G}
\]
is homotopical.
\end{cor}

\begin{pf}
By Proposition~\ref{thm:255} we may assume $K$ to be finite.  In
addition, it suffices to show that smashing with $S^{-W}\wedge K$ is
homotopical as a functor from $\pist\ugspectra{G}$ to itself.  Suppose
that $Y\to Y'$ is a weak equivalence.  Let $L\in\swg{G}$ be a $V$-dual
of $K$.  By the isomorphism of Proposition~\ref{thm:88} it suffices to
show that for all $H\subset G$ and all $k\in\Z$, the map
\[
\pist\ugspectra{G}(G/H_{+}\wedge S^{k},Y\wedge X)\to
\pist\ugspectra{G}(G/H_{+}\wedge S^{k},Y'\wedge X)
\]
is an isomorphism.  Using the first part of the duality
isomorphism~\eqref{eq:164}, we can identify this map with
\[
\pist\ugspectra{G}(G/H_{+}\wedge S^{k}\wedge S^{W}\wedge L,S^{V}\wedge
Y)\to
\pist\ugspectra{G}(G/H_{+}\wedge S^{k}\wedge S^{W}\wedge L,S^{V}\wedge
Y'),
\]
and finally by Proposition~\ref{thm-copied:244}, with 
\[
\ho\ugspectra{G}(G/H_{+}\wedge S^{k}\wedge S^{W}\wedge L,S^{V}\wedge
Y)\to
\ho\ugspectra{G}(G/H_{+}\wedge S^{k}\wedge S^{W}\wedge L,S^{V}\wedge
Y').
\]
But this latter map is an isomorphism since $S^{V}\wedge Y\to
S^{V}\wedge Y'$ is a weak equivalence (Proposition~\ref{thm-copied:254}).
\end{pf}

\begin{prop}
\label{thm:134} Let $J$ be a finite $G$-set.  For any
$X\in\ugspectra{G}$, the canonical map $\bigvee_{j\in J} X\to
\prod_{j\in J}X$ is an isomorphism in $\pist\ugspectra{G}$, hence a
weak equivalence.
\end{prop}

\begin{pf}
The finite $G$-sets are self-dual in $\swg{G}$.  Since
\[
\bigvee_{j\in J}X \approx J_{+}\wedge X,
\]
the result follows from the duality isomorphism
\[
\pist\ugspectra{G}(Z,J_{+}\wedge X) \approx
\pist\ugspectra{G}(J_{+}\wedge Z, X) \approx
\pist\ugspectra{G}(Z, \prod_{j\in J}X) 
\]
once one checks that the composite map is the same as the one coming
from the canonical map from the (constant) finite indexed wedge to the
finite indexed product.  We leave this to the reader.
\end{pf}

\begin{cor}
\label{thm:128}
Let $J$ be a finite $G$-set and $X$ an equivariant $J$-diagram.   The
map
\[
\bigvee_{j\in J}X_{j}\to \prod_{j\in J}X_{j}
\]
is an isomorphism in $\pist\ugspectra{G}$, hence a weak equivalence.
\end{cor}

\begin{pf}
Let $U:\ugspectra{G}\to \spectra^{B_{J}G}$ be the pullback map
associated to the unique equivariant map $J\to \text{pt}$.   The
indexed wedge is the left adjoint to $U$ and the indexed product is
the right adjoint.   The natural transformation from the indexed wedge
to the indexed product is easily checked to satisfy the condition of
Lemma~\ref{thm:135} below.   This reduces us to checking the case in
which the $J$-diagram is constant at a $G$-spectrum $X$.  But that
case is covered by Proposition~\ref{thm:134}.
\end{pf}

We have used 

\begin{lem}
\label{thm:135} Suppose that $U:\cat D\to \cat C$ is a functor with a
left adjoint $L$ and right adjoint $R$, and that $L\to R$ is a natural
transformation.  If the composition
\begin{equation}
\label{eq:165}
\id\to UR \to \id
\end{equation}
of the adjoint to $L\to R$ with the counit of the adjunction is the
identity, then $L\to R$ is a retract of $LUR\to RUR$.
\end{lem}

\begin{pf}
Just apply $L\to R$ on the left to the composition~\eqref{eq:165} to get
\[
\xymatrix{
L  \ar[r]\ar[d]  & LUR  \ar[r]\ar[d]   & L  \ar[d]\\
R  \ar[r]        & RUR   \ar[r]         &   R\mathrlap{\ .}
}
\]
\end{pf}

Corollary~\ref{thm:128} implies the only non-trivial
part of the following ``indexed'' analogue of
Proposition~\ref{thm-copied:91}.

\begin{prop}
\label{thm-copied:246}
\begin{thmList}
\item\label{part:fip} The formation of finite indexed products is homotopical.
\item\label{part:11} Suppose that $J$ is a finite $G$-set, and $X:\cat{B}_{J}G\to \spectra$
is a functor.   The map 
\[
\bigvee_{j\in J} X_{j} \to \prod_{j\in J} X_{j}
\]
is a stable weak equivalence in $\ugspectra{G}$.   Hence the formation
of finite indexed wedges is homotopical.  
\item\label{part-copied:10} The formation of all indexed wedges is
homotopical. 
\end{thmList}
\qed
\end{prop}

\subsubsection{Change of group}

Let $H\subset G$ be a subgroup.  Specializing
Proposition~\ref{thm-copied:246} to the case $J=G/H$ gives the
homotopical properties of the ``change of group'' functors.  The
functor $i_{H}^{\ast}:\ugspectra{G}\to\ugspectra{H}$ is homotopical by
definition, and so induces a functor on the homotopy categories
\[
i_{H}^{\ast}:\ho\ugspectra{G} \to \ho\ugspectra{H}.
\]
Taking $J=G/H$ in Proposition~\ref{thm-copied:246} we see that the
left and right adjoints to $i_{H}^{\ast}$ are also homotopical, and
that the canonical natural transformation between them is a weak
equivalence.  They therefore induce left and right adjoints to the
restriction map on the homotopy categories, and the canonical map
between them is an isomorphism.  This is the {\em Wirthm\"uller
isomorphism}~\cite{MR0343260,MR764596}.

\subsubsection{Weak equivalences and the smash product}
\label{sec:weak-equiv-smash}

The smash product is not known to preserve weak equivalences, but it
does so in good cases.

\begin{defin}
\label{def-copied:58} An equivariant orthogonal spectrum is
{\em cellular} if it is in the smallest subcategory of $\ugspectra{G}$
containing the spectra of the form $G_{+}\underset{H}{\wedge}
S^{-V}\wedge S^{k}$ with $V$ a representation of $H$ and $k\ge 0$ and
which is closed under the formation of arbitrary coproducts, the
formation of mapping cones, and the formation of filtered colimits
along $h$-cofibrations.
\end{defin}

The small object argument shows that every $X$ receives, functorially,
a weak equivalence $\tilde X\to X$ from a cellular $\tilde X$.

\begin{prop}
\label{thm-copied:250} If $K$ is cellular then $K$ is {\em flat}: the
functor $X\mapsto X\wedge K$ preserves weak equivalences.
\end{prop}

\begin{pf}
By Corollary~\ref{thm:132} and the fact that the formation of indexed
wedges is homotopical (Proposition~\ref{thm-copied:246}) the result is true when
$K=G_{+}\underset{H}{\wedge} S^{-V}\wedge S^{k}$.  The functor
$X\wedge K$ is built from
\[
X\wedge G_{+}\underset{H}{\wedge} S^{-V}\wedge S^{k}
\]
by forming wedges, mapping cones, and filtered colimits along
$h$-cofibrations, all of which are homotopical by
Proposition~\ref{thm-copied:246}.
\end{pf}

Since every object is weakly equivalent to a cellular object, and
cellular objects are flat, Remark~\ref{rem:67} implies

\begin{prop}
\label{thm:263}
Suppose that $X\to Y$ is a weak equivalence of flat spectra.   Then
for any $Z$, the map $X\wedge Z\to Y\wedge Z$ is a weak equivalence. \qed
\end{prop}

Let $\flugspectra{G}\subset\ugspectra{G}$ be the full subcategory of
flat objects, considered as a homotopical category using the stable
weak equivalences.  Since every object of $\ugspectra{G}$ is weakly
equivalent to an object of $\flugspectra{G}$, the functor
\begin{equation}
\label{eq:170}
\ho\flugspectra{G}\to \ho\ugspectra{G}
\end{equation}
is an equivalence of categories.   The above results show

\begin{prop}
\label{thm:261} The smash product functor
\[
\flugspectra{G}\times \ugspectra{G} \to \ugspectra{G} 
\]
is homotopical. \qed
\end{prop}

The equivalence~\eqref{eq:170} and
Proposition~\ref{thm-copied:250} are enough to show that the smash
product descends to give $\ho\ugspectra{G}$ a symmetric monoidal
structure, and that the map $\swg{G}\to\ho\ugspectra{G}$ is symmetric
monoidal.  For a more refined statement, see \S\ref{sec:furth-homot-prop}.

\subsection{Spectra as a model category}
\label{sec:spectra-as-model}
\subsubsection{The positive complete model structure}
\label{sec:strong-posit-stable}
Let $\gencof$ be the set of maps
\begin{equation}
\label{eq-copied:168}
\gencof = \{G_{+}\smashove{H} S^{-V}\wedge S^{n-1}_{+}\to
G_{+}\smashove{H} S^{-V}\wedge D^{n}_{+} \}
\end{equation}
with $n\ge 0$, $H$ a subgroup of $G$ and $V$ a representation of $H$
{\em containing a non-zero invariant vector}.  We define the class
\[
\pscofib{G} \subset\ugspectra{G}
\]
of {\em positive complete cofibrations} to be the smallest collection
of maps in $\ugspectra{G}$ containing the maps
in~\eqref{eq-copied:168} and which is closed under coproducts, cobase
change along arbitrary maps, and filtered colimits.  A {\em positive
complete fibration} (or just {\em fibration}) is a map having the
right lifting property with respect to the class of maps in
$\pscofib{G}$ which are stable weak equivalences.

\begin{prop}
\label{thm-copied:248} The category $\ugspectra{G}$ equipped with the stable
weak equivalences, the positive complete cofibrations and the
positive complete fibrations forms a (cofibrantly generated)
Quillen model category.
\end{prop}

We will call this model structure the {\em positive complete
model structure}, and when we need to recruit a model structure for
some task, this will be the one we use.   Henceforth the terms
``cofibration,'' ``fibration'' and ``weak equivalence'' will refer to
``positive complete cofibration,'' ``positive complete
fibration,'' and ``stable weak equivalence.''

\begin{rem}
\label{rem:38} Since the maps in $\gencof$ are mapping cylinders
they are $h$-cofibrations.  This implies that the cofibrations in
$\ugspectra{G}$ are $h$-cofibrations (cf.~\cite[Lemma
III.2.5]{MR1922205}) and hence flat.  The cofibrant objects in
$\ugspectra{G}$ are cellular hence flat.
\end{rem}

The ``positive'' condition is needed for the study of commutative
algebras.  On the other hand, it creates some peculiarities in the
model structure.  For example, the zero sphere $S^{0}$ is not
cofibrant, nor is $S^{0}\wedge K$ when $K$ is a $G$-CW complex.  The
cofibrant replacements are given by
\[
S^{-1}\wedge S^{1}\wedge K\to S^{0}\wedge K.
\]
This means that  the adjunction
\[
\Sigma^{\infty}:\ugspaces{G} \leftrightarrows \ugspectra{G}: \Omega^{\infty}.
\]
is not a Quillen adjunction, even though the left adjoint preserves
all weak equivalences between non-degenerately based $G$-spaces, and
so barely needs to be derived.

The positive complete model structure does not quite appear in the
literature.  It is closely related to the {\em positive stable model
structure} of~\cite{MR1922205}.

The positive complete model structure is cofibrantly generated.
The set $\gencof$ is the set of generating cofibrations.  The set
$\genacyclic$ of generating acyclic cofibrations consists of the
analogous maps
\begin{equation}
\label{eq:172}
G_{+}\underset{H}{\wedge}S^{-V}\wedge I^{n-1}_{+}\to
G_{+}\underset{H}{\wedge}S^{-V}\wedge I^{n}
\end{equation}
together with the corner maps formed by smashing 
\begin{equation}
\label{eq:121}
G_{+}\underset{H}{\wedge}\big(S^{-V\oplus W}\wedge S^{W}\big)\to G_{+}\underset{H}{\wedge}\tilde S^{-V}
\end{equation}
with the maps $S^{n-1}_{+}\to D^{n}_{+}$.  The $H$-representation $V$
is assumed to have a non-zero invariant vector, while $W$ need not.
The map~\eqref{eq:121} is extracted from the factorization
\begin{equation}
\label{eq:147}
S^{-V\oplus W}\wedge S^{W}\to \tilde S^{-V} \to S^{-V}
\end{equation}
formed by applying the small object construction in $\ugspectra{H}$,
using the maps in $\gencof$.

A map $X\to Y$ has the right lifting property with respect to the
class of maps $\gencof$ if and only for each $H\subset G$ and each
representation $V$ of $H$ containing a non-zero invariant vector, the
map $X_{V}\to Y_{V}$ is an acyclic fibration in $\ugspaces{H}$.  Among
other things this implies that $X\to Y$ is a weak equivalence and that
the map $\tilde S^{-V}\to S^{-V}$ is a homotopy equivalence.  From
this one concludes that a map $X\to Y$ has the right lifting property
with respect to $\genacyclic$ if and only if for each subgroup
$H\subset G$ and each representation $V$ of $H$ containing a non-zero
invariant vector, the map $X_{V}\to Y_{V}$ is a fibration in
$\ugspectra{H}$, and for each representation $W$ of $H$ the square
\begin{equation}
\label{eq:130}
\xymatrix{
X_{V}  \ar[r]\ar[d]  &  \Omega^{W}X_{V\oplus W}  \ar[d] \\
Y_{V}  \ar[r]        &\Omega^{W}Y_{V\oplus W}
}
\end{equation}
is homotopy cartesian in $\ugspaces{H}$.   

\begin{prop}
\label{thm:145} If a map $X\to Y$ is a weak equivalence and has the
right lifting property with respect to $\genacyclic$ then it has the
right lifting property with respect to $\gencof$.
\end{prop}

\begin{pf}
We must show that the conditions imply that for each $H\subset G$ and
each representation $V$ of $H$ containing a non-zero invariant vector,
the map $X_{V}\to Y_{V}$ is an acyclic fibration in $\ugspaces{H}$.
Part of our assumption is that it is a fibration, so it remains to
show that it is a weak equivalence.  Choose an exhausting sequence
$\{V_{n} \}$.  Letting $W$ range through this sequence
in~\eqref{eq:130} leads to a homotopy cartesian square
\[
\xymatrix{
X_{V}  \ar[r]\ar[d]  &  \ho\varinjlim\Omega^{V_{n}}X_{V\oplus V_{n}}  \ar[d] \\
Y_{V}  \ar[r]        &  \ho\varinjlim\Omega^{V_{n}}Y_{V\oplus V_{n}}\mathrlap{\ .}
}
\]
Since $X\to Y$ is a weak equivalence, the rightmost vertical map is a
weak equivalence (by~\eqref{eq:179}), hence so is $X_{V}\to Y_{V}$.
\end{pf}

\begin{prop}
\label{thm:208}
Any cobase change along a map in $\genacyclic$ is a weak equivalence.
\end{prop}

\begin{pf}
Since the maps in $\genacyclic$ are flat, it suffices to check that the
maps in $\genacyclic$ are weak equivalences.   The only ones for which
this is not obvious are the corner maps.   Since they are flat, it
suffices to check that the quotients 
\[
G_{+}\underset{H}\wedge \left(\tilde S^{-W}/\big(S^{-V\oplus W}\wedge
S^{V}\big)\right)\wedge D^{m}/S^{m}
\]
are weakly contractible.  Since $D^{m}/S^{m}$ is flat, and
$G_{+}\underset{H}{\wedge}(\slot)$ is homotopical, it suffices to show
that
\[
\tilde S^{-W}/\big(S^{-V\oplus W}\wedge S^{V}\big)
\]
is weakly contractible in $\ugspectra{H}$, or, equivalently that the
leftmost map in~\eqref{eq:147} is a weak equivalence in
$\ugspectra{H}$.  But that is a consequence of
Proposition~\ref{thm:159} and the two out of three property.
\end{pf}

\begin{prop}
\label{thm:257}
A map $X\to Y$ is a fibration if and only if it has the right lifting
property with respect to $\genacyclic$
\end{prop}

\begin{pf}
Suppose that $A\to B$ is an acyclic cofibration.  Using the small
object construction with the maps in $\genacyclic$ factor it as
$A\to\tilde B\to B$ where $A\to \tilde B$ is a filtered colimit of
maps constructed by iterated cobase change along maps in $\genacyclic$
and $\tilde B\to B$ has the right lifting property with respect to
$\genacyclic$.  The map $A\to\tilde B$ is a weak equivalence by
Propositions~\ref{thm:208} and~\ref{thm:255}.  It follows that $\tilde
B\to B$ is a weak equivalence, and so by Proposition~\ref{thm:145},
has the right lifting property with respect to $\gencof$.  This means
that $A\to B$ is a retract of $A\to \tilde B$.  Since $X\to Y$ has the
right lifting property for $A\to\tilde B$ it also has this property
for $A\to B$.
\end{pf}

The verification of the model category axioms is now completely
straightforward and left to the reader. 

Let $H\subset G$ be a subgroup.
In the positive complete model category structures, the restriction functor
\[
i_{H}^{\ast}:\ugspectra{G}\to\ugspectra{H}
\]
preserves weak equivalences, fibrations and cofibrations.  This
implies

\begin{prop}
Let $H\subset G$ be a subgroup.   The restriction functor and its left
adjoint form a Quillen pair
\[
G_{+}\smashove{H}(\slot): \ugspectra{H} \leftrightarrows \ugspectra{G}:i_{H}^{\ast},
\]
as do the restriction functor and its right adjoint
\[
i_{H}^{\ast}: \ugspectra{G} \leftrightarrows \ugspectra{H}:\prod_{j\in
G/H}(\slot)_{j}.
\]
\qed
\end{prop}

\begin{cor}
\label{thm:168}
An indexed wedge of cofibrations is a cofibration. \qed
\end{cor}

Corollary~\ref{thm:168} is one of our reasons for introducing the
positive complete model structure.  The positive stable model
structure of~\cite{MR1922205} does not have this property.

Associated to any map $i:G'\to G$ of finite groups is a functor
$i^{\ast}:\ugspectra{G}\to \ugspectra{G'}$.  This functor has both a
left and right adjoint.  The functor $i^{\ast}$ sends the generating
cofibrations to indexed wedges of generating cofibrations, hence
cofibrations by Corollary~\ref{thm:168}.  Since it is a left adjoint
it therefore sends cofibrations to cofibrations.  It also sends the
generating acyclic cofibrations to weak equivalences.  To see this
note that the generators of the form $X\wedge (I^{n-1}_{+}\to
I^{n}_{+})$ are homotopy equivalences hence go to homotopy
equivalences.  To check that the corner maps go to weak equivalences,
it suffices to show that the maps~\eqref{eq:121} go to weak
equivalences.  Since $\tilde S^{-V} \to S^{-V}$ is a homotopy
equivalence, this is equivalent to showing that maps of the form
\[
G_{+}\underset{H}{\wedge}\big(S^{-V\oplus W}\wedge S^{W}\big)\to G_{+}\underset{H}{\wedge} S^{-V}
\]
go to weak equivalences.   But these maps go to an indexed wedge of
maps of the form 
\begin{equation}
\big(S^{-V'\oplus W'}\wedge S^{W'}\big)\to S^{-V'}
\end{equation}
which are weak equivalences.  Thus $i^{\ast}$ also sends acyclic
cofibrations to acyclic cofibrations.  This gives

\begin{prop}
\label{thm:241}
If $i:G'\to G$ is any homomorphism of finite groups, then the pullback functor 
\[
i^{\ast}:\ugspectra{G}\to \ugspectra{G'}
\]
is a left Quillen functor.
\end{prop}

For more along these lines see~\cite[Remark~V.3.13]{MR1922205}

\subsubsection{Smash product}
\label{sec:furth-homot-prop}
Equipped with the smash product and the positive complete model
category structure, $\ugspectra{G}$ is a {\em symmetric monoidal model
category} in the sense of Hovey~\cite[Definition
4.2.6]{hovey99:_model} and Schwede-Shipley~\cite{MR1997322}.  This
means that the analogue of Quillen's axiom SM7 holds (the {\em pushout
product axiom}), and for any cofibrant $X$, the map
\[
\tilde S^{0}\wedge X \to X
\]
is a weak equivalence, where $\tilde S^{0}\to S^{0}$ is a cofibrant
approximation.  As will be apparent to the reader the proof applies
equally well if ``cofibration'' is replaced by ``cellular.''

\begin{prop}
\label{thm:95} Equipped with the smash product, the positive complete
model structure is a symmetric monoidal model category.
\end{prop}

The positive complete model structure also satisfies the
monoid axiom~\cite[Definition~3.3]{MR1734325}.

\begin{prop}
\label{thm:259}
If $X\to Y$ is an acyclic cofibration in $\ugspectra{G}$, and $Z$ is
arbitrary then $X\wedge Z\to Y\wedge Z$ is a flat weak equivalence.
\end{prop}

We have stated these together to slightly streamline the proof.  When
cofibrations are flat, the monoid axiom implies the ``weak
equivalence'' part of the pushout product axiom once one knows the
``cofibration'' part.  Indeed suppose $A_{1}\to B_{1}$ is an acyclic
cofibration and $A_{2}\to B_{2}$ is a cofibration.  Then the vertical
arrows in the diagram
\[
\xymatrix{
A_{1}\wedge A_{2}  \ar[r]\ar[d]_{\sim}  & A_{1}\wedge B_{2}
\ar[d]^{\sim} \\
B_{1}\wedge A_{2}  \ar[r]        & B_{1}\wedge B_{2}
}
\]
are weak equivalences by the monoid axiom, and all of the arrows are
cofibrations by the ``cofibration'' part (Remark~\ref{rem:49}).  Since
cofibrations are flat, the map from $A_{1}\wedge B_{2}$ to the pushout
is a weak equivalence, and the desired weak equivalence assertion then
follows from two out of three.

\begin{pf*}{Proofs of Propositions~\ref{thm:95} and~\ref{thm:259}}
The unit axiom follows from Proposition~\ref{thm:263} since cofibrant
objects are cellular, hence flat (Remark~\ref{rem:38}).  The
pushout product axiom asserts that if $f_{1}:A_{1}\to B_{1}$ and
$f_{2}:A_{2}\to B_{2}$ are cofibrations, then the corner map from the
pushout of the left and top arrows in
\begin{equation}
\label{eq:171}
\xymatrix{
A_{1}\wedge A_{2}  \ar[r]\ar[d] & A_{1}\wedge B_{2}
\ar[d] \\
B_{1}\wedge A_{2} \ar[r]& B_{1}\wedge B_{2}
}
\end{equation}
to the bottom right term is a cofibration, and is acyclic if one of
$f_{1}$ or $f_{2}$ is.  It suffices to check the cofibration condition
when $f_{1}$, and $f_{2}$ are in $\gencof$ and so of the form
\begin{align*}
G\smashove{H_{1}}S^{-V_{1}}\wedge\big(S^{k-1} &\to D^{k}\big) \\
G\smashove{H_{2}}S^{-V_{2}}\wedge\big(S^{\ell-1} &\to D^{\ell}\big).
\end{align*}
But in that case the corner map is the smash product of
$G\smashove{H_{1}}S^{-V_{1}}$ with $G\smashove{H_{2}}S^{-V_{2}}$ with
the pushout product of $S^{k-1}\to D^{k}$ and $S^{\ell-1}\to
D^{\ell}$.  This is an indexed wedge of cofibrations hence a
cofibration.  As remarked above, once Proposition~\ref{thm:259} is
proved, we are done.  Since $X\to Y$ is a cofibration it is an
$h$-cofibration, so it suffices to show that $(Y/X)\wedge Z$ is weakly
contractible if $Y/X$ is.  But $Y/X$ is cofibrant, hence flat, so the
claim follows from Proposition~\ref{thm:263}.
\end{pf*}

\begin{rem}
\label{rem:49} The special case of the pushout product axiom for $\ast
\to A$ and $\ast \to B$ asserts that if $A$ and $B$ are cofibrant,
then so is $A\wedge B$.
\end{rem}

Hovey~\cite[Theorem~4.3.2]{hovey99:_model} now implies
\begin{cor}
\label{thm:113}
The left derived smash product makes $\ho\ugspectra{G}$ into a
complete symmetric monoidal category.
\end{cor}

\subsubsection{The canonical homotopy presentation}
\label{sec:canon-pres}
Let 
\[
\dots\subset V_{n}\subset V_{n+1}\subset\dots
\]
be an exhausting sequence of orthogonal $G$-representations, and
consider the transition diagram
\begin{equation}
\label{eq:10}
\xymatrix{
S^{-V_{n+1}}\wedge \igcat{G}(V_{n},V_{n+1})\wedge X_{n}\ar[r] \ar[d] &
S^{-V_{n+1}}\wedge X_{n+1} \\
S^{-V_{n}}\wedge X_{n}\mathrlap{\ .} &
}
\end{equation}
Write
\[
W_{n} = V_{n+1}-V_{n}
\]
for the orthogonal complement of $V_{n}$ in $V_{n+1}$.   The inclusion
$V_{n}\subset V_{n+1}$ gives an embedding 
\[
S^{W_{n}}\to \igcat{G}(V_{n},V_{n+1}), 
\]
and so from~\eqref{eq:10} a diagram
\[
\xymatrix{
S^{-(V_{n}\oplus W_{n})}\wedge S^{W_{n}}\wedge X_{V_{n}}\ar[r]\ar[d] &
S^{-V_{n+1}}\wedge X_{n+1} \\
S^{-V_{n}}\wedge X_{n} \mathrlap{\ .} &
}
\]
Putting these together as $n$ varies results in a system 
\begin{equation}
\label{eq:14}
\xymatrix@C=1em{
A_{0} & B_{0}\ar[l]_{\sim}\ar[r] 
&
A_{1} & B_{1}\ar[l]_{\sim}\ar[r] 
&
A_{2} & B_{2}\ar[l]_{\sim}\ar[r] 
&
A_{3} & B_{3}\ar[l]_{\sim} \ar[r]
&\dots.
}
\end{equation}
The system~\eqref{eq:14} maps to $X$ and a simple check of equivariant
stable homotopy groups shows that the map from its homotopy colimit to
$X$ is a weak equivalence.   Now for each $n$ let $C_{n}$ be the
homotopy colimit of the portion
\begin{equation}
\label{eq:30}
\xymatrix@C=1em{
A_{0} & B_{0}\ar[l]_{\sim}\ar[r] 
&
\dots\ar[r]
&
A_{n-1} & B_{n-1}\ar[l]_{\sim}\ar[r]
&A_{n}
}
\end{equation}
of~\eqref{eq:14}.   
Then $C_{n}$ is naturally weakly equivalent to
$A_{n}= S^{-V_{n}}\wedge X_{V_{n}}$, and the $C_{n}$ fit into a sequence
\begin{equation}
\label{eq:19}
C_{0} \to C_{1}\to C_{2}\to\dots
\end{equation}
whose homotopy colimit coincides with that of~\eqref{eq:14}.  This
gives the canonical homotopy presentation of $X$.  One can
functorially replace the sequence~\eqref{eq:19} with a weakly
equivalent sequence of cofibrations between cofibrant-fibrant objects.
The colimit of this sequence is naturally weakly equivalent to $X$.
It will be cofibrant automatically, and fibrant since the model
category $\ugspectra{G}$ is compactly generated.

We write the canonical homotopy presentation of $X$ as 
\[
X\approx\hocolim_{V_{n}}\,(S^{-V_{n}}\wedge X_{V_{n}})_{\text{cf}},
\]
or when more precision is needed, as a diagram 
\[
X\leftarrow \hocolim_{V_{n}}\,(S^{-V_{n}}\wedge X_{V_{n}})_{\text{c}}
\to \hocolim_{V_{n}}\,(S^{-V_{n}}\wedge X_{V_{n}})_{\text{cf}},
\]
with the subscript indicating cofibrant and cofibrant-fibrant
replacement.  

\subsection{Homotopy properties of the norm}
\label{sec:homot-prop-norm}
The purpose of this section is to establish Proposition~\ref{thm:199}
which asserts that indexed smash products have a left derived functor
which may be computed on cofibrant objects.  As will be apparent to
the reader, they can also be computed on cellular objects.  Many of
the technical results in this section are also required for our
analysis of symmetric powers and of commutative algebras.

Before formulating our main results, we generalize the situation slightly.

\subsubsection{Equivariant $J$-diagrams}
\label{sec:equiv-j-diagr}
Given a non-empty $G$-set $J$, consider the category
$\spectra^{\cat{B}_{J}G}$ of functors $\cat B_{J}G\to \spectra$.  A
choice of point $t$ in each $G$-orbit of $J$ gives an equivalence
\[
\spectra^{\cat{B}_{J}G} \approx \prod_{t}\ugspectra{G_{t}},
\]
where $G_{t}$ is the stabilizer of $t$.  We give
$\spectra^{\cat{B}_{J}G}$ the model structure corresponding to the
product of the positive complete model structures under this
equivalence.  The model structure is independent of the chosen points
in each orbit.  We will refer to the model category
$\spectra^{\cat{B}_{J}G}$ as the model category of {\em equivariant
$J$-diagrams} of spectra.

To be more explicit, a map of $J$-diagrams $X\to Y$ is a
{\em weak equivalence} if and only for each $j\in J$ the map $Xj\to
Yj$ is a weak equivalence in $\ugspectra{G_{j}}$.  The {\em generating
cofibrations} are the maps whose $j^{\text{th}}$ component has the
form
\[
{G_{j}}_{+}\underset{H_{j}}{\wedge}{S^{-V_{j}}}\wedge S^{m_{j}-1}_{+} \to
{G_{j}}_{+}\underset{H_{j}}{\wedge}{S^{-V_{j}}}\wedge D^{m_{j}}_{+}
\]
with $V_{j}$ a representation of $H_{j}$ having a non-zero invariant
vector.   They can be expressed without reference to points and
stabilizers as an indexed wedge
\begin{equation}
\label{eq:181}
p^{\vee}_{\ast}\left(S^{-V}\wedge \big(S^{n-1}_{+}\to  D^{n}_{+}\big)\right)
\end{equation}
with $p:J'\to J$ a finite surjective map of $G$-sets, and $V$ a
$G$-equivariant orthogonal vector bundle over $J'$ having a
nowhere-zero invariant section.  The {\em generating acyclic
cofibrations} are the maps of the form
\[
p^{\vee}_{\ast}S^{-V}\wedge\big(I^{n-1}_{+}\to  I^{n}_{+}\big)
\]
and those constructed as the corner map formed by smashing 
\begin{equation}
\label{eq:184}
p^{\vee}_{\ast}\big(S^{-V\oplus W}\wedge S^{W}\to \tilde S^{-V}\big)
\end{equation}
with the maps $S^{n-1}_{+}\to D^{n}_{+}$.  As in~\eqref{eq:121}, the
map~\eqref{eq:184} is extracted from the factorization
\begin{equation}
\label{eq:185}
S^{-V\oplus W}\wedge S^{W}\to \tilde S^{-V} \to S^{-V}
\end{equation}
by applying the small object construction in the category of
equivariant $J'$-diagrams using the generating cofibrations.  The map
$\tilde S^{-V}\to S^{V}$ is a homotopy equivalence.  

If $J\to K$ is a map of finite $G$-sets, the restriction
functor 
\[
\spectra^{B_{K}G} \to
\spectra^{B_{J}G} 
\]
has both a left and right adjoint, given by the two Kan extensions.
All three functors are homotopical, and the both the restriction
functor and its left adjoint send cofibrations to cofibrations.  This
means that the restriction functor is both a left and right Quillen
functor.

Let $p:J\to K$ be an equivariant map of finite $G$-sets.   The indexed
smash product gives a functor
\[
p^{\wedge}_{\ast}=(\slot)^{\wedge J/K}:\spectra^{\cat B_{J}G} \to \spectra^{\cat B_{K}G}.
\]
When $J\to K$ is the map $G/H\to\text{pt}$ this is the norm.  The
various homotopical properties of indexed and symmetric smash products
we require are most naturally expressed as properties of
$(\slot)^{\wedge J/K}$.  Working fiberwise, establishing these reduces
to the case $K=\text{pt}$.  To keep the discussion uncluttered we
focus on that case in this section, leaving the extension to the case
of more general $K$ to the reader.

\subsubsection{Indexed smash products and cofibrations}
\label{sec:index-smash-prod}
Let $p:J\to \text{pt}$ be the unique equivariant map and write the
indexed smash product as $(\slot)^{\wedge J}$.  Note that if $V$ is an
equivariant orthogonal vector bundle over $J$ then
\[
(S^{-V})^{\wedge J}=S^{-V'},
\]
where $V'$ is the orthogonal $G$-space of global sections of $V$.

\begin{lem}
\label{thm:247} Suppose that $A\to B$ is a generating cofibration in
$\spectra^{\cat{B}_{J}G}$.  The indexed corner map
$\partial_{A}B^{\wedge J}\to B^{\wedge J}$ is an indexed wedge
\[
\bigvee_{\Gamma} S^{-V}\wedge \big(S(W)_{+}\to D(W)_{+}\big)
\]
in which $\Gamma$ is a $G$-set, $V$ and $W$ are equivariant vector
bundles over $\Gamma$ and $V$ has a non-zero invariant section.  In
particular, $\partial_{A}B^{\wedge J}\to B^{\wedge J}$ is a
cofibration.
\end{lem}

\begin{pf}
This is a straightforward consequence of the distributive law
(Theorem~\ref{thm:9}) applied to~\eqref{eq:181}, and the compatibility
of the formation of $\partial_{A}B^{\wedge J}$ with indexed wedges, as
described at the end of \S\ref{sec:index-mono-prod-pushout}.
\end{pf}

\begin{prop}
\label{thm:195} Suppose that $J$ is a non-empty finite $G$-set.  If
$X\to Y$ is a cofibration of equivariant $J$-diagrams, the indexed
smash product
\[
X^{\wedge J} \to Y^{\wedge J}
\]
is an $h$-cofibration.  It is a cofibration between cofibrant objects
in $\ugspectra{G}$ if $X$ is cofibrant.
\end{prop}

\begin{pf}
The assertion that $X^{\wedge J}\to Y^{\wedge J}$ is an
$h$-cofibration is contained in Proposition~\ref{thm:210}.  For the
cofibration assertion we work by induction on $|J|$, and may therefore
assume the result to be known for any non-empty $J_{0}\subset J$ and
any $H\subset G$ stabilizing $J_{0}$ as a subset.  In particular, we
may assume that if $X$ is cofibrant, then $X^{\wedge J_{0}}$ is a
cofibrant $H$-spectrum for any non-empty proper $J_{0}\subset J$ and
any $H\subset G$ stabilizing $J_{0}$ as a subset.

We will establish the theorem in the case in which $X\to Y$ arises
from a pushout square of $J$-diagrams
\[
\xymatrix{
A \ar[r] \ar[d] & B \ar[d]\\
X \ar[r] & Y }
\]
in which $A\to B$ is a generating cofibration.  We will show in this
case that $X^{\wedge J}\to Y^{\wedge J}$ is an $h$-cofibration, and is
a cofibration if $X$ is cofibrant.  Since the formation of indexed
smash products commutes with directed colimits and retracts, the
proposition then follows from the small object argument.

Give $Y^{\wedge J}$ the filtration described in
\S\ref{sec:index-mono-prod-pushout}.  The successive terms are related
by the pushout square
\begin{equation}
\label{eq:122}
\entrymodifiers={+!!<0pt,\fontdimen22\textfont2>}
\xymatrix{
{\displaystyle 
\bigvee_{\substack{J=J_{0}\amalg J_{1} \\ |J_{1}|=n}}
X^{\wedge J_{0}}\wedge \partial_{A}B^{\wedge J_{1}}}
  \ar[r]\ar[d]  &   
{\displaystyle 
\bigvee_{\substack{J=J_{0}\amalg J_{1} \\ |J_{1}|=n}}
X^{\wedge J_{0}}\wedge B^{\wedge J_{1}}}\ar[d] \\
\fil_{n-1}Y^{\wedge J}  \ar[r]        & \fil_{n}Y^{\wedge J}\mathrlap{\ .}
}
\end{equation}
By Lemma~\ref{thm:247}, each of the maps 
\[
\partial_{A}B^{\wedge J_{1}} \to B^{\wedge J_{1}}
\]
is a cofibration.  If $X$ is cofibrant, then $X^{\wedge J_{0}}$ is
either $S^{0}$ or cofibrant by induction, hence
\[
X^{\wedge J_{0}} \wedge \partial_{A}B^{\wedge J_{1}} \to 
X^{\wedge J_{0}} \wedge B^{\wedge J_{1}} 
\]
is a cofibration by the pushout product axiom.  Since indexed
wedges preserve cofibrations, the top row of~\eqref{eq:122} is then a
cofibration and hence so is the bottom row.
\end{pf}

To show that the indexed smash product has a left derived functor we
need to augment Proposition~\ref{thm:195} and show that what when
$X\to Y$ is an acyclic cofibration, then $X^{\wedge J}\to Y^{\wedge
J}$ is a weak equivalence.  This can be proved with the above
argument once we know that the indexed corner maps
$\partial_{A}B^{\wedge J}\to B^{\wedge J}$ associated to the
generating acyclic cofibrations are weak equivalences.  But the
generating acyclic cofibrations contain the maps of the
form~\eqref{eq:184} so dealing with them requires understanding
something about indexed corner maps of fairly general cofibrations.
These can be studied as the indexed smash products of {\em maps} in a
different symmetric monoidal category.

\subsubsection{The category of arrows} 
\label{sec:indexed-corner-map}

Let $\augspectra{G}$ denote the category of maps $X=(X_{0}\to X_{1})$
in $\ugspectra{G}$, with morphisms the commutative diagrams.  As
mentioned in Remark~\ref{rem:47}, $\augspectra{G}$ can be made into a
closed symmetric monoidal category by defining
\[
(X_{1}\to X_{2})\wedge (Y_{1}\to Y_{2})
\]
to be the corner map, from the pushout of the top and left arrows in 
\[
\xymatrix{
X_{1}\wedge Y_{1}  \ar[r]\ar[d]  &  X_{2}\wedge Y_{1} \ar[d] \\
X_{1}\wedge Y_{2}  \ar[r]        & X_{2}\wedge Y_{2}
}
\]
to the bottom right corner.    The tensor unit is $\ast\to S^{0}$.

We give $\augspectra{G}$ the projective model structure in which a map
\begin{equation}
\label{eq:182}
(X_{1}\to X_{2}) \to (Y_{1}\to Y_{2})
\end{equation}
is a weak equivalence or fibration if each of $X_{i}\to Y_{i}$ is, and
is a cofibration if both $X_{1}\to Y_{1}$ and the corner map
\begin{equation}
\label{eq:183}
X_{2}\underset{X_{1}}{\cup} Y_{1}\to Y_{2}
\end{equation}
are cofibrations.  An object $X_{1}\to X_{2}$ is therefore cofibrant
if $X_{1}$ is cofibrant and $X_{1}\to X_{2}$ is a cofibration.  

The model structure on $\augspectra{G}$ is compactly generated.  The
generating (acyclic) cofibrations in $\augspectra{G}$ are of two
types.  Type I are the maps
\[
(K\to K) \to (L\to L)
\]
and type II are the maps 
\[
(\ast\to K) \to (\ast \to L)
\]
were $K\to L$ is running through the set $\gencof$ defined
in~\eqref{eq-copied:168} (respectively $\genacyclic$). 
\begin{prop}
\label{thm:194} 
Equipped with the structure just described,
$\augspectra{G}$ is a symmetric monoidal model category satisfying 
the monoid axiom.
\end{prop}

\begin{pf}
The proof follows the proof of Propositions~\ref{thm:95}
and~\ref{thm:259}, and, because of the special nature of the
generators, essentially reduces to it.  It suffices to check the
``cofibration'' assertion on generators.  In each of the three cases
(type I and type I, type II and type II, and mixed type) the result
reduces to the case of $\ugspectra{G}$.  Since the cofibrations are
$h$-cofibration, the monoid axiom reduces showing that if
$(\ast\to\ast)\to (X_{1}\to X_{2})$ is an acyclic cofibration and
$(Z_{1}\to Z_{2})$ is arbitrary, then both the domain and range in the
corner map of
\[
\xymatrix{
X_{1}\wedge Z_{1}  \ar[r]\ar[d]  & X_{1}\wedge Z_{2}
\ar[d] \\
X_{2}\wedge Z_{1}  \ar[r]        & X_{2}\wedge Z_{2}
}
\]
are weakly contractible.  But by the monoid axiom for $\ugspectra{G}$,
every term in the diagram is weakly contractible.  The claim then
follows since the left vertical arrow is an $h$-cofibration, hence
flat.  As pointed out before the statement of
Proposition~\ref{thm:95}, this implies the ``weak equivalence'' part
of the pushout product axiom.  The unit axiom is also straightforward
and left to the reader.
\end{pf}

The proof of Proposition~\ref{thm:194} is more or less completely
formal, and can be rewritten to apply to the arrow category of any
symmetric monoidal model category.  This is done in the recent
paper~\cite{hovey:_smith} of Hovey.

\subsubsection{Indexed corner maps and cofibrations}
\label{sec:indexed-corner-maps-1}

Proposition~\ref{thm:194} addresses the homotopy properties of
ordinary smash products in $\augspectra{G}$.  For the indexed smash
products we work in the arrow category $\aspectra{\cat{B}_{J}G}$ of
maps of equivariant $J$-diagrams, in the projective model structure.
Our aim is to establish Proposition~\ref{thm:267} which gives control
over the indexed corner maps in $\ugspectra{G}$
(Proposition~\ref{thm:266}).  It is the analogue in
$\aspectra{\cat{B}_{J}G}$ of Proposition~\ref{thm:195}.  In
preparation, we need to identify the generating (acyclic)
cofibrations.  As mentioned in the previous section, those in
$\augspectra{G}$ are of two types.  Type I are the maps
\[
(K\to K) \to (L\to L)
\]
and type II are the maps 
\[
(\ast\to K) \to (\ast \to L)
\]
were $K\to L$ is running through the set $\gencof$ defined
in~\eqref{eq-copied:168} (respectively $\genacyclic$). 
The generating (acyclic) cofibrations in $\aspectra{\cat{B}_{J}G}$ can be
taken to be the equivariant $J$-diagrams consisting entirely of type I
or type II generators.

\begin{rem}
\label{rem:55} 
A map~\eqref{eq:182} is an $h$-cofibration if both
$X_{1}\to Y_{1}$ and the corner map~\eqref{eq:183} are.  Since the
cofibrations in $\ugspectra{G}$ are $h$-cofibrations the same is true
of the cofibrations in $\augspectra{G}$.
\end{rem}

\begin{lem}
\label{thm:268} If $X\to Y$ is a generating cofibration in the
category of equivariant $J$-diagrams in $\augspectra{G}$, then the
indexed corner map
\[
\partial_{X}Y^{\wedge J}\to
Y^{\wedge J}
\]
is a cofibration between cofibrant objects in
$\augspectra{G}$.
\end{lem}

\begin{pf}
First note that for generating cofibrations of type I, the corner map is
\[
(\partial_{K}L^{\wedge J}\to \partial_{K}L^{\wedge J}) \to (L^{\wedge
J}\to L^{\wedge J})
\]
and in type II it is 
\[
(\ast \to \partial_{K}L^{\wedge J}) \to (\ast \to L^{\wedge J}).
\]
The assertion therefore reduces to Lemma~\ref{thm:247}.
\end{pf}

\begin{prop}
\label{thm:267} 
Suppose that $J$ is a finite $G$-set.
If
\[
X\to Y 
\]
is a cofibration in $\augspectra{G}$ and $X$ is cofibrant, then the
indexed smash product
\[
X^{\wedge J} \to Y^{\wedge J}
\]
is a cofibration between cofibrant objects.
\end{prop}

\begin{pf}
The proof proceeds exactly as in the case of
Proposition~\ref{thm:195}.  The filtration of
\S\ref{sec:index-mono-prod-pushout} and induction on $|J|$ reduce the
problem to showing that the indexed corner map (in $\aspectra{\cat{B}_{J}G}$)
\[
\partial_{ A}  Y^{\wedge J} \to
B^{\wedge J}
\]
is a cofibration between cofibrant objects, when $A\to B$ is a
cofibrant generator.  This is the content of Lemma~\ref{thm:268}.
\end{pf}

Specializing, we now have

\begin{prop}
\label{thm:266} If $X\to Y$ is a cofibration of equivariant
$J$-diagrams and $X$ is cofibrant, then the indexed corner map
$\partial_{X}Y^{\wedge J}\to Y^{\wedge J}$ is a cofibration between
cofibrant objects.
\end{prop}

\begin{pf}
If $X\to Y$ is a cofibration of cofibrant $J$-diagrams, then $(X\to
Y)$ is cofibrant $J$-diagram in $\augspectra{G}$, and so
\[
(X\to Y)^{\wedge J} = (\partial_{X}Y^{\wedge J}\to Y^{\wedge J})
\]
is cofibrant by Proposition~\ref{thm:267}.  
\end{pf}

The result below is not used elsewhere in this paper, but is useful in
other contexts.   Having come this far, we record it here.

\begin{prop}
\label{thm:270}
Suppose that $X\to Y\to Z$ is a sequence of cofibrations of
equivariant $J$-diagrams in $\ugspectra{G}$, and that $X$ is
cofibrant.   Then the map 
\[
\partial_{X}Z^{\wedge J} \to \partial_{Y}Z^{\wedge J}
\]
is a cofibration between cofibrant objects.   
\end{prop}

\begin{pf}
Define $Y\to Z'$ by the pushout square
\[
\xymatrix{
X  \ar[r]\ar[d]  & Y \ar@{-->}[dl]
\ar[d] \\
Z  \ar[r]        & Z'\mathrlap{\ .}
}
\]
As the dashed arrow indicates, the map $(X\to Z)\to (Y\to Z)$ is a
retract of $(X\to Z)\to (Y\to Z')$.  By Proposition~\ref{thm:267} the
map
\[
(\partial_{X}Z^{\wedge J}\to Z^{\wedge J}) \to 
(\partial_{Y}Z'^{\wedge J}\to Z'^{\wedge J}) 
\]
is a cofibration, hence so is the map 
\[
\partial_{X}Z^{\wedge J} \to 
\partial_{Y}Z'^{\wedge J},
\]
and therefore so is its retract
\[
\partial_{X}Z^{\wedge J} \to 
\partial_{Y}Z^{\wedge J}.
\]
\end{pf}

\subsubsection{Indexed smash products and  acyclic cofibrations}
\label{sec:indexed-corner-maps}

With the indexed corner maps of cofibrations under control we can now
turn to the acyclic cofibrations.   

\begin{lem}
\label{thm:233}
If $X\to Y$ is a generating acyclic cofibration in
$\spectra^{\cat{B}_{J}G}$, then the indexed corner map
\[
\partial_{X}Y^{\wedge J}\to
Y^{\wedge J}
\]
is an acyclic cofibration of cofibrant objects in $\ugspectra{G}$.
\end{lem}

\begin{pf}
We know from Proposition~\ref{thm:266} that the indexed corner maps
are cofibrations between cofibrant objects, so what remains is the
assertion that they are weak equivalences.  This can be reduced
further.  Suppose that $X\to Y$ is an acyclic cofibration in
$\spectra^{\cat{B}_{J}G}$ and we wish to show that the indexed corner
map $\partial_{X}Y^{\wedge J}\to Y^{\wedge J}$ is a weak equivalence.
Give $Y^{\wedge J}$ the filtration described in
\S\ref{sec:index-mono-prod-pushout}, in which the successive terms are
related by the pushout square
\[
\entrymodifiers={+!!<0pt,\fontdimen22\textfont2>}
\xymatrix{
{\displaystyle 
\bigvee_{\substack{J=J_{0}\amalg J_{1} \\ |J_{1}|=n}}
X^{\wedge J_{0}}\wedge \partial_{X}Y^{\wedge J_{1}}}
  \ar[r]\ar[d]  &   
{\displaystyle 
\bigvee_{\substack{J=J_{0}\amalg J_{1} \\ |J_{1}|=n}}
X^{\wedge J_{0}}\wedge Y^{\wedge J_{1}}}\ar[d] \\
\fil_{n-1}Y^{\wedge J}  \ar[r]        & \fil_{n}Y^{\wedge J}\mathrlap{\ .}
}
\]
By Proposition~\ref{thm:266} and the pushout product axiom, the upper
arrow is a cofibration, which, by induction on $|J|$, we may assume to
be acyclic when $n<|J|$.  Since the cofibrations are flat, this means
that the bottom arrow is an acyclic cofibration when $n<|J|$.  It
follows that in this case, the indexed corner map is a weak
equivalence if and only if the {\em absolute} map $X^{\wedge J}\to
Y^{\wedge J}$ is.

We now turn to the generating acyclic cofibrations.  
The generators of the form $X\wedge
\big(I^{n-1}_{+}\to I^{n}_{+}\big)$ are homotopy equivalences, hence
so are the absolute maps.     The other generators are of the form 
\begin{equation}
\label{eq:186}
\big(S^{n-1}_{+}\to D^{n_{+}}\big)\wedge
\big(p^{\vee}_{\ast}S^{-V\oplus W}\wedge S^{W}\to
p^{\vee}_{\ast}\tilde S^{-V} \big),
\end{equation}
where $p:J'\to J$ is a map of finite $G$-sets and $V$ and $W$ are
equivariant vector bundles over $J'$.  The fact that the norm is
multiplicative, together with the monoid axiom for $\augspectra{G}$,
reduces us to considering only the right hand factor
in~\eqref{eq:186}.  The distributive law further reduces us to the
case $J'=J$.  Finally, since the map $\tilde S^{-V}\to S^{V}$ is a
homotopy equivalence, we may replace $\tilde S^{V}$ with $S^{V}$.
Evaluating both sides using Proposition~\ref{thm:12} we see that the
assertion amounts to checking that
\[
S^{-V'\oplus W'}\wedge S^{W'}\to S^{-V'}
\]
is a weak equivalence, where $V'$ and $W'$ are
the $G$-spaces of global sections.   But this is
Proposition~\ref{thm:159} (see Remark~\ref{rem:63}).
\end{pf}

As with Lemma~\ref{thm:268}, the separate cases of type I and type II
generators reduce the result below to Lemma~\ref{thm:233}.

\begin{lem}
\label{thm:196}
If $X\to Y$ is a generating acyclic cofibration in
the category of equivariant $J$-diagrams in $\augspectra{G}$, then the
indexed corner map
\[
\partial_{X}Y^{\wedge J}\to
Y^{\wedge J}
\]
is an acyclic cofibration of cofibrant objects in $\augspectra{G}$. \qed
\end{lem}

\begin{prop}
\label{thm:221}
Suppose that $J$ is a finite $G$-set.   The functor 
\[
(\slot)^{\wedge J}:\aspectra{\cat B_{J}G}\to \augspectra{G}
\]
sends acyclic cofibrations between cofibrant objects to acyclic
cofibration between cofibrant objects, and hence weak
equivalences between cofibrant objects to weak equivalences between
cofibrant objects.
\end{prop}

\begin{pf}
The proof proceeds exactly as in the case of
Proposition~\ref{thm:195}.   That the second assertion follows from
the first is Ken Brown's Lemma (see, for example~\cite[Lemma
1.1.12]{hovey99:_model}).
\end{pf}

Specializing Proposition~\ref{thm:221}, we have

\begin{prop}
\label{thm:232} If $X\to Y$ is an acyclic cofibration in
$\spectra^{\cat{B}_{J}G}$ and $X$ is cofibrant, then both the indexed
corner map $\partial_{X}Y^{\wedge J}\to Y^{\wedge J}$ and the absolute
map $X^{\wedge J}\to Y^{\wedge J}$ are acyclic cofibrations between
cofibrant objects.  \qed
\end{prop}

\subsubsection{Homotopy properties of the norm}
\label{sec:homot-prop-norm-1}

With all this in hand we can now show that indexed smash products have
left derived functors.  From Proposition~\ref{thm:195},
Proposition~\ref{thm:232}, and Ken Brown's Lemma we have

\begin{prop}
\label{thm:11}
The indexed smash product
\[
(\slot)^{\wedge J}:\spectra^{\cat{B}_{J}G}\to \ugspectra{G}
\]
takes weak equivalences between cofibrant objects to weak
equivalences between cofibrant objects.  \qed
\end{prop}

This gives

\begin{prop}
\label{thm:199} 
The indexed smash product has a left derived functor 
\[
(\slot)^{\lwedge J}:\spectra^{\cat{B}_{J}G} \to \ho\ugspectra{G} 
\]
which may be computed as 
\[
X^{\lwedge J}= 
(X_{c})^{\wedge J}
\]
where $X_{c}\to X$ is a cofibrant approximation. 
\end{prop}

\subsection{Symmetric powers}
\label{sec:symmetric-powers}

We now turn to the homotopical properties of symmetric smash powers,
or just ``symmetric powers'' for short.

\subsubsection{Indexed symmetric powers}
\label{sec:index-symm-powers}

The $n^{\text{th}}$ symmetric (smash) power of a $G$-spectrum is the
orbit spectrum 
\[
\sym^{n}(X) = X^{\wedge n}/\Sigma_{n}.
\]
The homotopy properties of this functor are fundamental to
understanding the homotopy theory of equivariant commutative algebras.
For indexed smash products of commutative algebras the distributive
law leads one to consider indexed smash products of symmetric powers
\[
(\sym^{n}X)^{\wedge J}.
\]
These can be written as
\begin{equation}
\label{eq:198}
(\sym^{n}X)^{\wedge J} = (X^{\wedge n}/\Sigma_{n})^{\wedge J} \approx
X^{\wedge (\mathbf{n}\times J)}/\Sigma_{n}^{J}
\end{equation}
with $\mathbf n=\{1,\dots,n \}$.  This last expression is an {\em
indexed symmetric power}.  The definition and homotopy properties of
indexed symmetric powers are the subject of this section.

Before turning to the definition, we consider a more basic
situation.  Suppose that $i:\tilde G\to G$ is a surjective map of
groups with kernel $N$.  Then the functor $i^{\ast}:\ugspectra{G}\to
\ugspectra{\tilde G}$ has both a left and a right adjoint.  This is
most readily understood by thinking of $G$-spectra as objects of
$\spectra$ equipped with a $G$-action.  The left adjoint
$i_{!}:\ugspectra{\tilde G}\to \ugspectra{G}$ sends a spectrum $Y$ to
the {\em orbit spectrum} $Y/N$ equipped with its residual $G$-action.
The expression on the right of~\eqref{eq:198} is a special case of
this.  As in any diagram category, the orbit spectrum $Y/N$ is
computed objectwise: if $U$ is an orthogonal vector space then
$(Y/N)_{U}$ is the $G$-space $Y_{U}/N$.  For the homotopical
properties we need information about the $W$-space, for a
representation $W$ of $G$.  It is given by the formula
\[
(Y/N)_{W} = O(U,W)_{+}\underset{O(U)}{\wedge}(Y/N)_{U},
\]
where $U$ is any vector space of the same dimension as $W$ but with
trivial $G$-action.  Interchanging the colimits, this space can be
written as
\[
(O(U,W)_{+}\underset{O(U)}{\wedge}Y_{U})/N, 
\]
which, in turn is isomorphic to 
\[
Y_{W}/N,
\]
where now $W$ is regarded as a $\tilde G$ representation through the
map $\tilde G\to G$.   

We can now define indexed symmetric powers.   Let $I$ be a finite
$G$-set, and $\Sigma_{I}$ the group of (not necessarily equivariant)
automorphisms of $I$, with $G$ acting by conjugation.  Fix a
$G$-stable subgroup $\Sigma\subset \Sigma_{I}$ and regard $I$ as a
$\Sigma\rtimes G$-set through the projection map to $G$.  For a
$\Sigma\rtimes G$-equivariant $I$-diagram $X$ the {\em indexed
symmetric power} is the orbit $G$-spectrum
\[
\sym^{I}_{\Sigma}X = X^{\wedge I}/\Sigma.
\]
When the indexing set $I$ has a trivial $G$-action, $\Sigma$ is the
full symmetry group of $I$, and the equivariant $I$-diagram is the
constant diagram with value $X\in\ugspectra{G}$, then this construction is
the usual symmetric power $\sym^{|I|}X$ discussed above.   We will
usually not distinguish in notation between a $\Sigma\rtimes
G$-spectrum $X$ and the constant equivariant $I$-diagram with value $X$.

If $X\to Y$ is a map of $\Sigma\rtimes G$-equivariant $I$-diagrams,
the {\em indexed symmetric corner map} is the map of orbit $G$-spectra
\[
\partial_{X}\sym^{I}_{\Sigma}Y \to \sym^{I}_{\Sigma}Y
\]
obtained by passing to $\Sigma$ orbits from 
\[
\partial_{X}Y^{\wedge I}\to Y^{\wedge I}.
\]
It can also be regarded as the symmetric power $\sym^{I}_{\Sigma}(X\to
Y)$ of $(X\to Y)$ regarded as an object of the arrow category
$\augspectra{\cat{B}_{\Sigma\rtimes G}I}$.

\begin{rem}
\label{rem:73} Since the orbit spectrum functor is a continuous left
adjoint, it sends $h$-cofibrations to $h$-cofibrations.   For example, 
suppose that $X\to Y$ is a cofibration of cofibrant
$\Sigma\rtimes G$-equivariant $I$-diagrams.  By
Proposition~\ref{thm:195} and Proposition~\ref{thm:266} both the
indexed smash product
\[
X^{\wedge I}\to Y^{\wedge I}
\]
and the corner map
\[
\partial_{X}Y^{\wedge I}\to Y^{\wedge I}
\]
are cofibrations, and hence $h$-cofibrations, of $\Sigma\rtimes
G$-spectra.  This means that all four of the maps
\begin{align*}
\sym^{I}_{\Sigma}X&\to \sym^{I}_{\Sigma}Y \\
\partial_{X}\sym^{I}_{\Sigma} Y&\to \sym^{I}_{\Sigma}Y \\
(\egsigma{G})_{+}\smashove{\Sigma}X^{\wedge I} &\to
(\egsigma{G})_{+}\smashove{\Sigma}Y^{\wedge I} \\
(\egsigma{G})_{+}\smashove{\Sigma}\partial_{X}Y^{\wedge I} &\to
(\egsigma{G})_{+}\smashove{\Sigma}Y^{\wedge I} \\
\end{align*}
are $h$-cofibrations of $G$-spectra.
\end{rem}

Note that $X^{\wedge I}$ with its $\Sigma\rtimes G$-action is a
special case of an indexed monoidal product.  This means that the
distributive law applies to symmetric powers, and, given a pushout
square
\[
\xymatrix{
A  \ar[r]\ar[d]  & B  \ar[d] \\
X  \ar[r]        & Y\mathrlap{\ ,}
}
\]
there is a filtration of $\sym^{I}_{\Sigma}Y$ whose successive terms are
related by passing to $\Sigma$-orbits from the filtration described in
\S\ref{sec:index-mono-prod-pushout}.  

As described in~\cite{MR1922205}, the homotopy theoretic analysis of indexed
symmetric powers requires certain equivariant principal bundles.  For
the moment, let $\Sigma$ be any finite group with a $G$-action.

\begin{defin}
\label{def:48} An {\em equivariant universal $\Sigma$-space} is a
$\Sigma\rtimes G$-space $\egsigma{G}$ with the property that for each
finite $\Sigma\rtimes G$-set $S$, the space of $\Sigma\rtimes
G$-equivariant maps
\[
S\to \egsigma{G}
\]
is empty if some element of $S$ is fixed by a non-identity element of
$\Sigma$, and contractible otherwise.
\end{defin}

The defining property characterizes an equivariant universal
$\Sigma$-space up to $\Sigma\rtimes G$-equivariant weak homotopy
equivalence.  The space $\egsigma{G}$ is the total space of the
universal $G$-equivariant principal $\Sigma$-bundle.  It can be
constructed as a $\Sigma\rtimes G$-CW complex, with cells of the form
$S\times D^{m}$, where $S$ is a $\Sigma$-free $\Sigma\rtimes G$-set.
We will always assume our equivariant universal $\Sigma$-spaces are
$\Sigma\rtimes G$-CW complexes, in which case the characterization is
up to equivariant homotopy equivalence.

The symmetric powers of a cofibrant spectrum are
rarely cofibrant.  However they still have very good homotopy
theoretic properties.  Our main result is the following.

\begin{prop}
\label{thm:262} Suppose that $X\to Y$ is a cofibration between
cofibrant $\Sigma\rtimes G$-equivariant $I$-diagrams.  In the square
of $G$-spectra
\begin{equation}
\label{eq:177}
\xymatrix{
(\egsigma{G})_{+}\smashove{\Sigma}\partial_{X}Y^{\wedge I}
\ar[r]\ar[d]_{\sim}  &
(\egsigma{G})_{+}\smashove{\Sigma}Y^{\wedge I} 
\ar[d]^{\sim} \\
\partial_{X}\sym^{I}_{\Sigma}Y  \ar[r]        &
\sym^{I}_{\Sigma}Y\mathrlap{\ ,}
}
\end{equation}
every object is flat, the upper row is a cofibration between cofibrant
objects, the vertical maps are weak equivalences, and the bottom row
is an $h$-cofibration.  The horizontal maps are weak equivalences if
$X\to Y$ is.
\end{prop}

\begin{rem}
\label{rem:32}
By Proposition~\ref{thm:263} the maps in~\eqref{eq:177} asserted to be
weak equivalences remain so after smashing with any spectrum $Z$.
\end{rem}

\begin{rem}
\label{rem:72} The situation that comes up in studying the free
commutative algebra functor is that $X\to Y$ is a cofibration of
cofibrant $G$-spectra, regarded as a $\Sigma_{I}\rtimes G$-spectrum
through the map to $G$, and then regarded as a constant equivariant
$I$-diagram.  This map of equivariant $I$-diagrams is cofibrant by
Proposition~\ref{thm:241}, and so Proposition~\ref{thm:262} applies.
\end{rem}

Along the way to proving Proposition~\ref{thm:262} we will also show

\begin{prop}
\label{thm:198} 
The functors $(\egsigma{G})_{+}\smashove{\Sigma}(\slot)^{\wedge I}$
and  $\sym^{I}_{\Sigma}(\slot)$ take weak equivalences between
cofibrant objects to weak equivalences.
\end{prop}

\begin{rem}
\label{rem:65} 
Proposition~\ref{thm:262} is part the reason for the
{\em positive} condition in the model structure we have chosen.  The
result is not true for general {\em cellular} objects described in
\S\ref{sec:weak-equiv-smash}, though it is true for cellular object
built from cells of the form $G_{+}\smashove{H}S^{-V}\wedge D^{k}_{+}$
with $V$ non-zero.  The condition that $V$ is non-zero is used in the
proof of Proposition~\ref{thm:200}.
\end{rem}

We assertions about the top row in Proposition~\ref{thm:262} are most
easily analyzed in the arrow category $\aspectra{\cat B_{I}G}$.  

\begin{lem}
\label{thm:201} The functor
\[
\egsigma{G}_{+}\underset{\Sigma}{\wedge}(\slot)^{\wedge
I}:\augspectra{\cat{B}_{\Sigma\rtimes G}I}\to\augspectra{G}
\]
takes acyclic cofibrations between cofibrant objects to acyclic
cofibrations between cofibrant objects and hence 
weak equivalences between cofibrant objects to weak equivalences
between cofibrant objects.
\end{lem}

\begin{pf}
Let $X\to Y$ be an acyclic cofibration.  By working through an
equivariant cell decomposition of $\egsigma{G}$ and using SM7 for the
topological enrichment we reduce to showing that if $S$ is a
$\Sigma$-free $\Sigma\rtimes G$-set, then the map
\[
S_{+}\smashove{\Sigma}X^{\wedge I} \to
S_{+}\smashove{\Sigma}Y^{\wedge I} 
\]
is an acyclic cofibration between cofibrant objects.  This is an
indexed wedge of maps, indexed by the $\Sigma$-orbits $\mathcal
O\subset S$.  The summand corresponding to $\mathcal O$ is the map of
$G_{\mathcal O}$-spectra
\[
\mathcal O_{+}\smashove{\Sigma}X^{\wedge I}\to
\mathcal O_{+}\smashove{\Sigma}Y^{\wedge I},
\]
where $G_{\mathcal O}\subset G$ is the subgroup of $G$ preserving
$\mathcal O$.  Since $\mathcal O$ is a $\Sigma$-torsor, this is just 
the map of indexed smash products 
\[
X^{\wedge I'}  \to Y^{\wedge I'} 
\]
with $I' = \mathcal O\underset{\Sigma}{\times}I$, and is an acyclic
cofibration between cofibrant objects by Proposition~\ref{thm:221}.
The second assertion follows from the first by Ken Brown's Lemma.
\end{pf}

The vertical maps in~\eqref{eq:177} require a more detailed analysis.

\begin{defin}
\label{def:47} Suppose that $\Sigma$ is a group with an action of $G$,
and that $X$ is a $\Sigma\rtimes G$-spectrum.  We will say that $X$ is
{\em $\Sigma$-free as a $G$-spectrum} if for each orthogonal
$G$-representation $W$ the $\Sigma$-action on $X_{W}$ is free away
from the base point.
\end{defin}

\begin{prop}
\label{thm:200} If $X$ is a cofibrant $\Sigma\rtimes G$-equivariant
$I$-diagram, and $Z$ is any $\Sigma\rtimes G$-spectrum then $X^{\wedge
I}\wedge Z$ is a $\Sigma$-free $G$-spectrum.  The map
\begin{equation}
\label{eq:131}
(\egsigma{G})_{+}\smashove{\Sigma} (X^{\wedge I}\wedge Z)\to
(X^{\wedge I}\wedge Z)/\Sigma.
\end{equation}
is a weak equivalence in $\ugspectra{G}$.
\end{prop}

\begin{rem}
\label{rem:64}
We will mostly be interested in the case in which the $\Sigma$-action
on $Z$ is trivial.   In that case the equivalence~\eqref{eq:131}
takes the form
\[
\big((\egsigma{G})_{+}\smashove{\Sigma} X^{\wedge I}\big)\wedge Z\xrightarrow{\sim}{}
\sym^{I}_{\Sigma}(X)\wedge Z
\]
\end{rem}

\begin{rem}
\label{rem:53} The proof of Proposition~\ref{thm:200} is nearly
identical to that of~\cite[Lemma III.8.4]{MR1922205}.  We go through
the details because the statement is slightly more general, and in
order to correct a minor error in~\cite{MR1922205}.  The statements
of~\cite[Lemma III.8.4]{MR1922205}, and the related~\cite[Lemma
IV.4.5]{MR1922205} both use $E\Sigma_{i}$, whereas the object that
should really be used is $\eggroup{G}{\Sigma_{i}}$.  This makes the
proofs of~\cite[Theorem III.8.1]{MR1922205} and~\cite[Theorem
4]{MR2066508} on equivariant commutative rings incomplete.  The actual
homotopical analysis of commutative rings is more or less equivalent
to the homotopical analysis of the norm.  So it would seem that any
correct treatment needs to be built on the theory of the norm.
\end{rem}

\begin{pf*}{Proof of Proposition~\ref{thm:200}}
For the first assertion, it suffices to show that if $A\to B$
is a generating cofibration,
\[
\xymatrix{
A  \ar[r]\ar[d]  & B  \ar[d] \\
X_{0} \ar[r]        & X_{1}\mathrlap{\ ,}
}
\]
is a pushout square, and $X_{0}^{\wedge I}\wedge Z$ is $\Sigma$-free,
then $X_{1}^{\wedge I}\wedge Z$ is $\Sigma$-free.  We use the
filtration described in \S\ref{sec:index-mono-prod-pushout} and
consider the pushout square below
\begin{equation}
\label{eq:123}
\entrymodifiers={+!!<0pt,\fontdimen22\textfont2>}
\xymatrix{
{\displaystyle 
\bigvee_{\substack{I=I_{0}\amalg I_{1} \\ |I_{1}|=m}}
X_{0}^{\wedge I_{0}}\wedge \partial_{A}B^{\wedge I_{1}}\wedge Z}
  \ar[r]\ar[d]  &   
{\displaystyle 
\bigvee_{\substack{I=I_{0}\amalg I_{1} \\ |I_{1}|=m}}
X_{0}^{\wedge I_{0}}\wedge B^{\wedge I_{1}}\wedge Z}\ar[d] \\
\fil_{m-1}X_{1}\wedge Z  \ar[r]        & \fil_{m}X_{1}\wedge Z.
}
\end{equation}
Since $A\to B$ is a cofibration, the map in the top row is an
$h$-cofibration (Proposition~\ref{thm:266}) hence a closed inclusion.
It therefore suffices to show that $\Sigma$ acts freely away from the
base point on the upper right term (see Remark~\ref{rem:66}).
Induction on $|I|$ reduces this to $m=|I|$.  In this way the first
assertion of the proposition reduces to checking the special case
\[
X=p^{\vee}_{\ast}S^{-V}\wedge D^{k}_{+},
\]
with $p:\tilde I\to I$ a surjective map of $\Sigma\rtimes G$-sets and
$V$ an equivariant vector bundle over $\tilde I$ having a nowhere
vanishing invariant global section.   Since the factor
$(D^{k}_{+})^{\wedge I}$ can be absorbed into $Z$, we might as well
suppose
\[
X=p^{\vee}_{\ast}S^{-V}.
\]
The distributive law gives
\[
X^{\wedge I}=\bigvee_{\gamma\in\Gamma} S^{-V_{\gamma}}
\]
where $\Gamma$ is the $\Sigma\rtimes G$-set of sections $I\to\tilde I$, and 
\[
V_{\gamma} = \bigoplus_{i\in I}V_{\gamma(i)}.
\]
For an orthogonal $\Sigma\rtimes G$-representation $W$ we have, by
Lemma~\ref{thm:111},
\[
(X^{\wedge I}\wedge Z)_{W}=
\begin{cases}
\ast &\quad \dim W <\dim V_{\gamma} \\
\displaystyle\bigvee_{\gamma\in\Gamma}O(V_{\gamma}\oplus
U_{\gamma},W)_{+}\smashove{O(U_{\gamma})}Z_{U_{\gamma}}  &\quad\dim W\ge
\dim V_{\gamma}
\end{cases}
\]
in which $U=\{U_{\gamma}\}$ is any $\Sigma\rtimes G$-equivariant
vector bundle over $\Gamma$ with $\dim U_{\gamma}= \dim W-\dim
V_{\gamma}$.  We are interested in representations $W$ which are
pulled back from the projection map $\Sigma\rtimes G\to G$.  In the
first case there is nothing to prove.  In the second case the
complement of the base point is homeomorphic to
\[
\coprod_{\gamma\in\Gamma}O(V_{\gamma}\oplus U_{\gamma},W)\underset{O(U_{\gamma})}{\times} (Z_{U_{\gamma}}-\{\ast \})
\]
(see Remark~\ref{rem:66}).  The $\Sigma$-freeness then follows
from the fact that this space has an equivariant map to the disjoint
union of Stiefel-manifolds
\[
\coprod_{\gamma\in\Gamma}O(V_{\gamma},W) = \coprod_{\gamma\in\Gamma}O(V_{\gamma}\oplus U_{\gamma},W)/O(U_{\gamma}),
\]
which is $\Sigma$-free since each $V_{\gamma(i)}$ is non-zero, and
$\Sigma$ acts faithfully on $I$ but trivially on $W$.

With one additional observation, a similar argument reduces the
assertion about weak equivalences to the same case 
\begin{equation}
\label{eq:188}
X=p^{\vee}_{\ast}S^{-V}.
\end{equation}
To spell it out, abbreviate~\eqref{eq:123} as
\[
\xymatrix{
K  \ar[r]\ar[d] & L   \ar[d]\\
Y  \ar[r]       &Y'
}
\]
and form
\[
\xymatrix{ (\egsigma{G})_{+}\smashove{\Sigma} Y \ar[d]_{\sim} &\ar[l]
(\egsigma{G})_{+}\smashove{\Sigma} K
\ar[r]^{\flat}\ar[d]_{\sim} &
(\egsigma{G})_{+}\smashove{\Sigma} L \ar[d]^{\sim}\\
Y/\Sigma &\ar[l] K/\Sigma  \ar[r]^{\flat}  & L/\Sigma\mathrlap{\ .}
}
\]
By Remark~\ref{rem:73} the rightmost maps in both rows are
$h$-cofibrations, hence flat.  This means that if the vertical maps
are weak equivalences then the map of pushouts is a weak equivalence
(Remark~\ref{rem:52}).  With this in hand, one now reduces the second
claim to the cases $X=p^{\vee}_{\ast}S^{-V}\wedge S^{k-1}_{+}$ and
$X=p^{\vee}_{\ast}S^{-V}\wedge D^{k}_{+}$.  Absorbing the factors
$(S^{k-1}_{+})^{\wedge I}$ and $(D^{k}_{+})^{\wedge I}$ into $Z$
completes the reduction to~\eqref{eq:188}.  

With this $X$, the map on $W$-spaces induced by~\eqref{eq:131} is the
identity map of the terminal object if $\dim W<\dim V_{\gamma}$ and
otherwise the map of $\Sigma$-orbit spaces induced by
\[
(\egsigma{G})_{+}\wedge \bigvee_{\gamma\in\Gamma}O(V_{\gamma}\oplus
U_{\gamma},W)_{+}\smashove{O(U_{\gamma})}Z_{U_{\gamma}}\to 
\bigvee_{\gamma\in\Gamma}O(V_{\gamma}\oplus U_{\gamma},W)_{+}\smashove{O(U_{\gamma})}Z_{U_{\gamma}},
\]
in which $U=\{U_{\gamma}\}$ is any $\Sigma\rtimes G$-equivariant
vector bundle over $\Gamma$ with $\dim U_{\gamma}= \dim W-\dim
V_{\gamma}$.  The proposition then follows from the fact that
\[
\egsigma{G}\times \coprod_{\gamma\in\Gamma}O(V_{\gamma}\oplus U_{\gamma},W) \to 
\coprod_{\gamma\in\Gamma}O(V_{\gamma}\oplus U_{\gamma},W)
\]
is an equivariant homotopy equivalence for the compact Lie group
\[
\scr G = \big(\prod_{\gamma\in\Gamma}O(U_{\gamma})\rtimes\Sigma\big)\rtimes G.
\]
To see this, note that by~\cite{MR696520,MR0500993}, both sides are
$\scr G$-CW complexes so it suffices to check that the map is a weak
equivalence of $H$-fixed point spaces for all $H\subset \scr G$.  If
the image of $H$ in $\Sigma\rtimes G$ is not a subgroup of $\Sigma$
then $\egsigma{G}^{H}$ is contractible and the map of fixed points is
a homotopy equivalence.  If $H$ is a subgroup of $\prod O(U_{\gamma})$
then it acts trivially on $\egsigma{G}$, and once again
$\egsigma{G}^{H}$ is contractible.  Finally, suppose that there is an
element $h\in H$ whose image in $\Sigma\rtimes G$ is a non-identity of
$\Sigma$.  Since $W$ is pulled back from a $G$-representation, this
element acts trivially on $W$.  If $\gamma \in\Gamma$ is not fixed by
$h$ then no point of $O(V_{\gamma}\oplus U_{\gamma},W)$ can be fixed
by $h$.  If $\gamma\in\Gamma$ {\em is} fixed by $h$, then $h$ acts on
$V_{\gamma}$.  This action is non-trivial since $\Sigma$ acts
faithfully on $I$.  This means that $O(V_{\gamma}\oplus U_{\gamma},W)$
has no points fixed by $h$ since $h$ acts trivially on $W$.  Both
sides therefore have empty $H$-fixed points in this case.
\end{pf*}

\begin{pf*}{Proof of Proposition~\ref{thm:262}} 
The assertion that the upper arrow is a cofibration between cofibrant
objects and a weak equivalence if $X\to Y$ is, is contained in
Lemma~\ref{thm:201}.   Indeed consider the map of arrows 
\[
(X\to Y) \to (Y\to Y).
\]
If $X\to Y$ is a cofibration between cofibrant objects then both the
domain and range of the above map of arrows are cofibrant.  By
Lemma~\ref{thm:201} the map
\[
\big((\egsigma{G})_{+}\smashove{\Sigma}\partial_{X}Y^{\wedge I}
\to 
(\egsigma{G})_{+}\smashove{\Sigma}Y^{\wedge I}
\big)
\to 
\big((\egsigma{G})_{+}\smashove{\Sigma}Y^{\wedge I}
\to 
(\egsigma{G})_{+}\smashove{\Sigma}Y^{\wedge I}
\big)
\]
is a map of cofibrant objects, which is a weak equivalence if $X\to Y$
is.    This gives the assertion about the upper row.
The fact that the bottom row is an $h$-cofibration is
part of Remark~\ref{rem:73}. 

For the remaining assertions it will be helpful to reference the
expanded diagram
\[
\xymatrix{
(\egsigma{G})_{+}\smashove{\Sigma}\partial_{X}Y^{\wedge I}\wedge Z
\ar[r]\ar[d]  &
(\egsigma{G})_{+}\smashove{\Sigma}Y^{\wedge I} \wedge Z
\ar[r]\ar[d] &  (\egsigma{G})_{+}\smashove{\Sigma}(Y/X)^{\wedge
I}  \wedge Z
\ar[d]\\
\partial_{X}\sym^{I}_{\Sigma}Y \wedge Z \ar[r]        &
\sym^{I}_{\Sigma}Y \wedge Z
\ar[r] & \sym^{I}_{\Sigma}(Y/X) \wedge Z\mathrlap{\ ,}
}
\]
in which $Z$ is any $G$-spectrum.  By Proposition~\ref{thm:200} the
two right vertical maps are weak equivalences.  Since the left
horizontal maps are $h$-cofibrations, hence flat, this implies that
the left vertical map is a weak equivalences.  Taking $Z=S^{0}$ gives
the weak equivalence of the vertical arrows in the statement of
Proposition~\ref{thm:262}.  Letting $Z$ vary through a weak
equivalence and using the fact that cofibrant objects are flat gives
the flatness assertion.  By what we have already proved, when $X\to Y$
is a weak equivalence the vertical and top arrows in the left square
are weak equivalences, hence so is the bottom left map.  This
completes the proof.
\end{pf*}

\begin{pf*}{Proof of Proposition~\ref{thm:198}}
Suppose that $X\to Y$ is a weak equivalence of cofibrant objects, and
consider the diagram 
\[
\xymatrix{
(\egsigma{G})_{+}\smashove{\Sigma}X^{\wedge I}  \ar[r]\ar[d]  & (\egsigma{G})_{+}\smashove{\Sigma}Y^{\wedge I}
\ar[d] \\
\sym^{I}_{\Sigma}X  \ar[r]        & \sym^{I}_{\Sigma}Y\mathrlap{\ .}
}
\]
The vertical maps are weak equivalences by Proposition~\ref{thm:200}.
The top horizontal map is a weak equivalence by Lemma~\ref{thm:201}
(applied to, say, the map $(\ast \to X)\to (\ast \to Y)$).  The bottom
map is therefore a weak equivalence.
\end{pf*}

\subsubsection{Iterated indexed and symmetric powers}
\label{sec:iter-index-symm}

In our analysis of the norms of commutative rings in
\S\ref{sec:index-smash-prod-comm-rings} we will encounter iterated
indexed smash products and symmetric powers.  These work out just to
be other indexed smash or symmetric powers.  The point of this section
is to spell this out.

Suppose that $I$ and $J$ are $G$-sets and that $X$ is an equivariant
$I\times J$-diagram.   The factorization
\[
I\times J\to J\to \text{pt}
\]
gives an isomorphism 
\begin{equation}
\label{eq:189}
(X^{\wedge I})^{\wedge J} \approx X^{\wedge (I\times J)},
\end{equation}
in which $X^{\wedge I}$ is shorthand for $p^{\wedge}_{\ast}X$ with $p:I\times
J\to J$ the projection mapping.   Applying this to the arrow category
we get an isomorphism  of the corner map
\[
\partial_{X}Y^{\wedge (I\times J)} \to X^{\wedge(I\times J)}
\]
with the iterated corner map 
\[
\partial_{W}Z^{\wedge J} \to Z^{\wedge J}
\]
in which $W\to Z$ is the map 
\[
\partial_{X}Y^{\wedge I} \to Y^{\wedge I}.
\]

There is also a version with symmetric powers.  Suppose in addition
that $\Sigma\subset \Sigma_{I}$ is a $G$-stable subgroup.  Then the
action of $\Sigma^{J}$ on $I\times J$ by
\[
\phi\cdot (i,j) = (\phi(j)\cdot i,j)
\]
is $G$-stable, making $J\times I$ into a $\Sigma^{J}\rtimes G$-set,
and the projection map $I\times J\to J$ equivariant, with
$\Sigma^{J}\rtimes G$ acting on $J$ through $G$.  When $X$ is a
$\Sigma^{J}\rtimes G$-equivariant $J\times I$-diagram, the
isomorphism~\eqref{eq:189} is $\Sigma^{J}\rtimes G$-equivariant.
Passing to orbits gives an isomorphism of $G$-spectra 
\begin{equation}
\label{eq:191}
(\sym^{I}_{\Sigma}X)^{\wedge J} \approx
\sym^{I\times J}_{\Sigma^{J}}X.
\end{equation}
By working in the arrow category we get an isomorphism of the corner map 
\[
\partial_{X}\sym^{I\times J}_{\Sigma^{J}} Y \to \sym^{I\times J}_{\Sigma^{J}} Y
\]
with the iterated indexed corner map 
\begin{equation}
\label{eq:190}
\partial_{W}Z^{\wedge J} \to Z^{\wedge J}
\end{equation}
in which $W\to Z$ is the map 
\[
\partial_{X}\sym^{I}_{\Sigma}Y \to \sym^{I}_{\Sigma}Y.
\]

Our analysis of the homotopy properties of symmetric powers depended
on a convenient cofibrant approximation.   Let $\egsigma{G}$ be a
universal $G$-equivariant $\Sigma$ space.   The above discussion leads
to an isomorphism
\[
\big(\egsigma{G}_{+}\smashove{\Sigma} X^{\wedge
I}\big)^{\wedge J}
\approx 
\egsigma{G}^{J}_{+} \smashove{\Sigma^{J}} X^{\wedge(I\times J)},
\]
and an identification of the corner map
\[
\partial_{\tilde W}\tilde Z^{\wedge J} \to \tilde Z^{\wedge J}
\]
in which $\tilde W\to \tilde Z$ is the map
\[
\egsigma{G}_{+}\smashove{\Sigma}\big(\partial_{X}Y^{\wedge I} \to 
Y^{\wedge I}\big)
\]
with
\[
\big(\egsigma{G}^{J}\big)_{+} \smashove{\Sigma^{J}} \big(\partial_{X}Y^{\wedge(I\times J)}
\to Y^{\wedge(I\times J)}\big).
\]
To reduce this expression to one we have already considered we need to
know that $\egsigma{G}^{J}$ is a universal equivariant
$\Sigma^{J}$-space.

\begin{lem}
\label{thm:206} Let $J$ be a finite $G$-set.  If $\egsigma{G}$ is an
equivariant universal $\Sigma$-space then, under the product action,
$(\egsigma{G})^{J}$ is an equivariant universal $\Sigma^{J}$-space.
\end{lem}

\begin{pf}
The functor $T\mapsto T^{J}$ (from $\Sigma\rtimes G$-spaces to
$\Sigma^{J}\rtimes G$-spaces) has a left adjoint.  To describe it,
let $M$ be the set $\Sigma\rtimes G\times J$ and define a left action
of $\Sigma\rtimes G$ by the product of the translation action on
$\Sigma\rtimes G$ and the action of $G$ on $J$.  There is a commuting
right $\Sigma^{J}\rtimes G$-action
\[
\big(\Sigma\rtimes G\times J\big)\times \big(\Sigma^{J}\rtimes G\big)
\to \Sigma\rtimes G\times J.
\]
whose component in the second factor is just the projection, and in
the first factor is composed of the evaluation map
\[
J\times \Sigma^{J}\rtimes G\to \Sigma\rtimes G
\]
and the right action of $\Sigma\rtimes G$ on itself.  The functor
$T\mapsto T^{J}$ can be identified with
\[
\hom_{\Sigma\rtimes G}(M,T)
\]
and so its left adjoint is given by 
\[
S \mapsto M\underset{\Sigma^{J}\rtimes G}{\times}S.
\]
Breaking $M$ into right $\Sigma^{J}\rtimes G$-orbits gives the
decomposition
\[
M\underset{\Sigma^{J}\rtimes G}{\times}S =\coprod_{j\in J} S/\Sigma^{J-\{j\}}.
\]
In this latter expression, the action of $t\in\Sigma$ on $s\in
S/\Sigma^{J-\{j\}}$ can be computed as the orbit class of $\sigma s$,
where $\sigma\in\Sigma^{J}$ is any element with $\sigma(j)=t$.   For
example, the entire $\Sigma$-action can be computed by restricting to
the diagonal subgroup of $\Sigma^{J}$.

Observe that a $\Sigma^{J}\rtimes G$-set $S$ is $\Sigma^{J}$-free if
and only if $M\underset{\Sigma^{J}\rtimes G}{\times}S$ is
$\Sigma$-free.  Clearly if $S$ is $\Sigma^{J}$-free then for each
$j\in J$, $S/\Sigma^{J-\{j \}}$ is $\Sigma$-free.  On the other hand
if $\sigma\in\Sigma^{J}$ is a non-identity element fixing $s\in S$,
then there is a $j\in J$, with $\sigma(j)$ not the identity element.
For this $j$ we have $\sigma(j)\cdot \Sigma^{J-\{j \}}s=\Sigma^{J-\{j
\}}s$.

Now to the proof.  Let $S$ be a finite $\Sigma^{J}\rtimes G$-set.  We
need to show that the space of $\Sigma^{J}\rtimes G$-maps
\[
S\to \egsigma{G}^{J}
\]
is empty or contractible depending on whether or not $S$ has a point
fixed by a non identity element of $\Sigma^{J}$.  By adjunction, this
space can be identified with the space of $\Sigma\rtimes G$-maps from
\[
M\underset{\Sigma^{J}\rtimes G}{\times} S \to \egsigma{G}, 
\]
and so the result follows from the observation above.  
\end{pf}

We will be interested in the following case.  Suppose that
$X\to Y$ is a cofibration of cofibrant $\Sigma^{J}\rtimes
G$-equivariant $I\times J$-diagrams.
By Lemma~\ref{thm:206} and Proposition~\ref{thm:262}, in
the diagram 
\[
\xymatrix{
(\egsigma{G})^{J}_{+}\smashove{\Sigma^{J}}\partial_{X}Y^{\wedge(I\times
J)}  \ar[r]\ar[d]_{\sim}  &
(\egsigma{G})^{J}_{+}\smashove{\Sigma^{J}}Y^{\wedge (I\times J)}
\ar[d]^{\sim} \\
\partial_{X}\sym^{I\times J}_{\Sigma^{J}} Y  \ar[r]        & \sym^{I\times J}_{\Sigma^{J}} Y
}
\]
every object is flat, the top row is a cofibration of cofibrant
objects, the bottom row is an $h$-cofibration, and the vertical maps
are weak equivalences and remain so after smashing with any spectrum.
The same conclusion therefore holds for the corresponding diagram of
iterated indexed (symmetric) powers 
\[
\xymatrix{
\partial_{\tilde W}(\tilde Z^{\wedge(J)})  \ar[r]\ar[d]  &
 \tilde Z^{\wedge(J)} \ar[d] \\
\partial_{W}(Z^{\wedge(J)})  \ar[r]        & Z^{\wedge(J)}
}
\]
in which 
\[
\xymatrix{
\tilde W  \ar[r]\ar[d]  & \tilde Z
\ar[d] \\ 
W  \ar[r]        & Z
}
\]
is the diagram 
\[
\xymatrix{
\egsigma{G}_{+}\smashove{\Sigma} \partial_{X}Y^{\wedge I}
\ar[r]\ar[d]  & \egsigma{G}_{+}\smashove{\Sigma} Y^{\wedge I}
\ar[d] \\
\partial_{X}\sym^{I}_{\Sigma}Y  \ar[r]        &
\sym^{I}_{\Sigma}Y\mathrlap{\ .}
}
\]

Working fiberwise leads to an analogous result about the indexed smash
product along a map $q:J\to K$ of finite $G$-sets.  It plays an
important role in our analysis of the homotopy properties of the norms
of commutative rings.  Aside from the map $J\to K$ of finite $G$-sets,
the situation is the same as what we have been discussing in this
section.  We have fixed a finite $G$-set $I$, a $G$-stable subgroup
$\Sigma\subset \Sigma_{I}$, and a universal $G$-equivariant
$\Sigma$-space $\egsigma{G}$.

\begin{prop}
\label{thm:238} 
Let $X\to Y$ be a cofibration of cofibrant $\Sigma^{J}\rtimes
G$-equivariant $I\times J$-diagrams and write
\[
\xymatrix{
\tilde W  \ar[r]\ar[d]  & \tilde Z
\ar[d] \\ 
W  \ar[r]        & Z
}
\]
for the diagram 
\[
\xymatrix{
\egsigma{G}_{+}\smashove{\Sigma} \partial_{X}Y^{\wedge I}
\ar[r]\ar[d]  & \egsigma{G}_{+}\smashove{\Sigma} Y^{\wedge I}
\ar[d] \\
\partial_{X}\sym^{I}_{\Sigma}Y  \ar[r]        & \sym^{I}_{\Sigma}Y.
}
\]
In the $G$-equivariant $K$-diagram of corner maps 
\[
\xymatrix{
\partial_{\tilde W}(\tilde Z^{\wedge(J/K)})  \ar[r]\ar[d]  &
 \tilde Z^{\wedge(J/K)} \ar[d] \\
\partial_{W}(Z^{\wedge(J/K)})  \ar[r]        & Z^{\wedge(J/K)}
}
\]
every object is flat, the vertical maps are weak equivalences after
smashing with any object, the upper map is a cofibration of
cofibrant objects and the lower map is an $h$-cofibration.  The
horizontal maps are weak equivalences if $X\to Y$ is.  \qed
\end{prop}

\begin{rem}
\label{rem:75} The actual hypothesis on $X\to Y$ required for the
fiberwise argument is that for each $k\in K$, the map $X\to Y$ is a
cofibration of $\Sigma^{J_{k}}\rtimes G_{k}$-equivariant $I\times
J_{k}$-diagram, where $J_{k}\subset J$ is the inverse image of $k$,
and $G_{k}$ is the stabilizer of $k$.  For the sake of a cleaner
statement we have made the slightly stronger assumption that it is a
cofibration of cofibrant $\Sigma^{J}\rtimes G$-equivariant $I\times
J$-diagrams.  That this implies the ``fiberwise'' hypothesis is a
consequence of Proposition~\ref{thm:241}.
\end{rem}

\begin{rem}
\label{rem:74} As in Remark~\ref{rem:72}, Proposition~\ref{thm:238}
applies to the situation which $X\to Y$ is a cofibration of cofibrant
$G$-equivariant $J$-diagrams, regarded as a $\Sigma\rtimes
G$-equivariant $I\times J$ diagram by pulling back along the
projection mappings $\Sigma\rtimes G\to G$ and $I\times J\to J$.
\end{rem}

\subsection{Rings and modules}
\label{sec:rings-modules}

Aside from the alteration in model structure, the following is stated
as~\cite[Theorem III.8.1]{MR1922205}.  The proof depends on our
analysis of symmetric powers, which, as mentioned in
Remark~\ref{rem:53}, makes essential use of the norm.

\begin{prop}
\label{thm:114}
The forgetful functor 
\[
\ugeinftycat{G}\to \ugspectra{G}
\]
creates a topological model category structure on $\ugeinftycat{G}$ in
which the fibrations and weak equivalences in $\ugeinftycat{G}$ are the
maps that are fibrations and weak equivalences in $\ugspectra{G}$.
\end{prop}

\begin{pf}
Most of the proof is  formal.   One takes as generating 
cofibrations the maps $\sym A\to\sym B$ where $A\to B\in \gencof$
and generating acyclic cofibrations the maps $\sym A\to \sym B$ with
$A\to B\in \genacyclic$.     The only real point to check is that if
\begin{equation}
\label{eq:178}
\xymatrix{
\sym A  \ar[r]\ar[d]  &  \sym B \ar[d] \\
X  \ar[r]        & Y
}
\end{equation}
is a pushout diagram in which $A\to B$ is a generating acyclic
cofibration, then $X\to Y$ is a weak equivalence.  That is contained
in Lemma~\ref{thm:273} below.  The rest of the proof is left to the
reader.
\end{pf}

\begin{lem}
\label{thm:273}
Suppose that $A\to B$ is a map of $G$-spectra, and 
\[
\xymatrix{
\sym A  \ar[r]\ar[d]  &  \sym B \ar[d] \\
X  \ar[r]        & Y
}
\]
is a pushout diagram of equivariant commutative rings.   If $A\to B$
is an acyclic cofibration of cofibrant objects, then $X\to Y$ is a
weak equivalence.
\end{lem}

The proof of Lemma~\ref{thm:273} involves a filtration of $Y$ by
$X$-modules which we will use again in
\S\ref{sec:index-smash-prod-comm-rings}.  For a map $A\to B$ of
$G$-spectra define
\[
\fil_{m}\sym B = \bigvee_{n}\fil_{m}\sym^{n}B
\]
where the $\fil_{m}\sym^{n}B$ is obtained from the filtration
described in \S\ref{sec:index-mono-prod-pushout} by passing to
$\Sigma_{n}$-orbits, and fits into a pushout square
\[
\xymatrix{
\sym^{n-m}A \wedge \partial_{A}\sym^{m}B
\ar[r]\ar[d]  &   \sym^{n-m}A\wedge \sym^{m}B 
\ar[d] \\
\fil_{m-1}\sym^{n}(B)  \ar[r]     &  \fil_{m}\sym^{n} B\mathrlap{\ ,}
}
\]
with
\[
\partial_{A}\sym^{m}B = \left(\partial_{A}B^{\wedge m}  \right)/\Sigma_{m}.
\]
Wedging over $n$ one sees that the $\fil_{m}B$ are $\sym A$-submodules, and
that there is a pushout square of $A$-modules
\[
\xymatrix{
\sym A\wedge \partial_{A}\sym^{m}B  \ar[r]\ar[d]  & \sym A\wedge \sym^{m}B  \ar[d] \\
\fil_{m-1}\sym B  \ar[r]        &\fil_{m}\sym B\mathrlap{\ .}
}
\]
If a map $X\to Y$ of commutative rings fits into a pushout diagram
\[
\xymatrix{
\sym A  \ar[r]\ar[d]  & \sym B  \ar[d] \\
X  \ar[r]        & Y,
}
\]
we can define a filtration of $Y$ by $X$-modules by 
\[
\fil_{m}Y = X\smashove{\sym A}\fil_{m} \sym B.
\]
Evidently these $\fil_{m}Y$ are related by the pushout square of
$X$-modules
\begin{equation}
\label{eq:132}
\xymatrix{
X\wedge \partial_{A}\sym^{m}B  \ar[r]\ar[d]  & X\wedge \sym^{m}B  \ar[d] \\
\fil_{m-1}Y  \ar[r]        &\fil_{m}Y\mathrlap{\ .}
}
\end{equation}

\begin{pf*}{Proof of Lemma~\ref{thm:273}}
We use the filtration just described.  In the diagram~\eqref{eq:132},
if $A\to B$ is an acyclic cofibration between cofibrant objects, then
\[
\partial_{A}\sym^{m}B\to \sym^{m}B
\]
is a weak equivalence and an $h$-cofibration of flat spectra by
Proposition~\ref{thm:262}.  It follows that the bottom map is a weak
equivalence.
\end{pf*}

\begin{cor}
\label{thm:17}
For $H\subset G$, the adjoint functors 
\[
\ugeinftycat{H} \leftrightarrows \ugeinftycat{G}
\]
form a Quillen pair.
\end{cor}

\begin{pf}
The restriction functor obviously preserves the classes of fibrations
and weak equivalences.
\end{pf}

\begin{cor}
\label{thm:27}
The norm functor on commutative algebras
\[
\norm_{H}^{G}:\ugeinftycat{H}\to \ugeinftycat{G}
\]
is a left Quillen functor.  It preserves the classes of cofibrations
and acyclic cofibrations, hence weak equivalences between cofibrant
objects.
\end{cor}

\begin{pf}
This is immediate from Corollary~\ref{thm:17} and
Proposition~\ref{thm:184}.  The assertion about weak equivalences is
Ken Brown's Lemma (see, for example~\cite[Lemma
1.1.12]{hovey99:_model}).
\end{pf}

The category $\rmod{R}$ of left modules over an equivariant
associative algebra $R$ as defined in \S\ref{sec:equiv-einfty-ainfty}.
As pointed out there, when $R$ is commutative, a left $R$-module can
be regarded as a right $R$-module, and $\rmod{R}$ becomes a symmetric
monoidal category under the operation
\begin{equation}
\label{eq:39}
M\smashove{R} N.
\end{equation}

The following result is a consequence of
Propositions~\ref{thm:95},~\ref{thm:259}
and~\cite[Theorem~4.1]{MR1734325}.  Except for the slight change of
model structure, it is~\cite[Theorem III.7.6]{MR1922205}.

\begin{prop}
\label{thm:115} 
The forgetful functor 
\[
\rmod{R} \to \ugspectra{G}
\]
creates a model structure on the category $\rmod{R}$ in which the
fibrations and weak equivalences are the maps which become fibrations
and weak equivalences in $\ugspectra{G}$.  When $R$ is commutative, the
operation~\eqref{eq:39} satisfies the pushout product and monoid
axioms making $\rmod{R}$ into a symmetric monoidal model category. \qed
\end{prop}

Though not explicitly stated, the following formal result was surely
known to the authors of~\cite{MR1734325} (see the proof
of~\cite[Theorem~4.3]{MR1734325}.)

\begin{cor}
\label{thm:13}
Let $f:R\to S$ be a map of equivariant associative  algebras.   The
functors 
\[
S\smashove{R}(\slot): \rmod{R} \leftrightarrows \rmod{S}:U
\]
given by restriction and extension of scalars form a Quillen pair.  If
$S$ is cofibrant as a left $R$-module, then the restriction functor is
also a left Quillen functor.
\end{cor}

\begin{pf}
Proposition~\ref{thm:115} implies that the restriction functor
preserves fibrations and acyclic fibrations.  This gives the first
assertion.  The second follows from the fact that the restriction
functor preserves colimits, and the consequence of
Proposition~\ref{thm:115} that the generating (acyclic) cofibrations
for $\rmod{S}$ are formed as the smash product of $S$ with the
generating (acyclic) cofibrations for $\ugspectra{G}$.
\end{pf}

The following result is~\cite[Proposition~III.7.7]{MR1922205}.  Using
the fact that $h$-cofibrations are flat, the proof reduces to checking
the case $M=G_{+}\smashove{H}S^{-V}\wedge R$, which is
Proposition~\ref{thm-copied:250}.

\begin{prop}
\label{thm:48} Suppose that $R$ is an associative algebra, and $M$ is
a cofibrant right $R$-module.  The functor $M\smashove{R}(\slot)$
preserves weak equivalences.  \qed
\end{prop}

In other words, the functor $M\smashove{R}(\slot)$ is flat if $M$ is
cofibrant, and so need not be derived.   

\begin{cor}
\label{thm:102}
Suppose that $R$ is an associative algebra, $M$ a cofibrant right
$R$-module.   If $N\to N'$ a map of left $R$-modules whose underlying
map of spectra is an $h$-cofibration, then the sequence 
\[
M\smashove{R}N \to M\smashove{R}N' \to M\smashove{R}(N'/N)
\]
is weakly equivalent to a cofibration sequence.   
\end{cor}

Note that the assumption is {\em not} that $N\to N'$ is an
$h$-cofibration in the category of left $R$-modules.  In that case the
result would not require any hypothesis on $M$.

\begin{pf}
We must show that the map from the mapping cone of  
\begin{equation}
\label{eq:95}
M\smashove{R}N \to
M\smashove{R}N'
\end{equation}
to $M\smashove{R}(N'/N)$ is a weak equivalence.  The mapping cone
of~\eqref{eq:95} is isomorphic to 
\[
M\smashove{R}(N'\cup CN), 
\]
and the spectrum underlying the $R$-module mapping cone $N'\cup N$ is
the mapping cone formed in spectra.  Since $N\to N'$ is an
$h$-cofibration, the map $N'\cup CN\to N'/N$ is a weak equivalence
(Proposition~\ref{thm:116}).  The result now follows from
Proposition~\ref{thm:48}.
\end{pf}

Corollary~\ref{thm:102} can be used to show that many constructions
derived from the formation of monomial ideals have good homotopy
theoretic properties.   It is used in~\S\ref{sec:quotient-rings}
and in~\S\ref{sec:from-reduct-theor}.  In those cases
the map of spectra underlying $N\to N'$ is the inclusion of a wedge
summand, and so obviously an $h$-cofibration.

\subsection{Indexed smash products of commutative rings}
\label{sec:index-smash-prod-comm-rings}

\subsubsection{Description of the problem}
\label{sec:description-problem} 
Proposition~\ref{thm:199} asserts that
the indexed smash product functor
\[
(\slot)^{\wedge J}:\spectra^{\cat{B}_{J}G} \to \ugspectra{G}
\]
has a left derived functor
\[
(\slot)^{\lwedge J}:\ho\spectra^{\cat{B}_{J}G} \to \ho\ugspectra{G}
\]
which can be computed by applying the norm to a cofibrant
approximation.  We also know from Corollary~\ref{thm:17} (and the fact
that coproducts of weak equivalences are weak equivalences) that the
restriction functor and its left adjoint form a Quillen pair
\[
p_{!}:\comm\spectra^{\cat{B}_{J}G} \leftrightarrows
\comm\ugspectra{G}:p^{\ast}.
\]
Furthermore, the following diagram commutes in which the vertical
functors are the forgetful functors (Corollary~\ref{thm:184})
\[
\xymatrix{
\comm\spectra^{\cat{B}_{J}G}  \ar[r]^{p_{!}}\ar[d]  &  \comm\ugspectra{G} \ar[d] \\
\spectra^{\cat{B}_{J}G}   \ar[r]_{(\slot)^{\wedge J}}        & \ugspectra{G}\mathrlap{\ .}
}
\]
However, what we really want is the commutativity of the following
diagram
\[
\xymatrix{
\ho\comm\spectra^{\cat{B}_{J}G}  \ar[r]^{\lder p_{!}}\ar[d]  &  \ho\comm\ugspectra{G} \ar[d] \\
\ho\spectra^{\cat{B}_{J}G}   \ar[r]_{(\slot)^{\lwedge J}}        & \ho\ugspectra{G}
}
\]
in which the vertical maps are the forgetful functors (which are
homotopical, so don't need to be derived), and the horizontal arrows
are the left derived functors indicated.  The point of this section is
to establish this.

To clarify the issue, suppose that $R\in \comm\spectra^{\cat{B}_{J}G}$
is a cofibrant $J$-diagram of commutative rings.  Let $\tilde R\to R$
be a cofibrant approximation of the underlying $J$-diagram of spectra.
What needs to be checked is that the map
\begin{equation}
\label{eq:36}
(\tilde R)^{\wedge J}\to (R)^{\wedge J}
\end{equation}
is a weak equivalence.  The proof involves an elaboration of the
notion of flatness.   To motivate it we describe a bit of the argument.

The main point in the proof is to investigate the situation of a pushout
diagram of equivariant $J$-diagrams of commutative rings 
\[
\xymatrix{
\sym A  \ar[r]\ar[d]  & \sym B  \ar[d] \\
R_{1}  \ar[r]        & R_{2}
}
\]
in which the top row is constructed by applying the symmetric algebra
functor $\sym$ to a generating cofibration $A\to B$, and in which one
knows that the map~\eqref{eq:36} is a weak equivalence for $R=R_{1}$.
One would like to conclude that~\eqref{eq:36} is a weak equivalence
for $R=R_{2}$.   

To pass from $R_{1}$ to $R_{2}$ we use the $R_{1}$-module filtration
described after the statement of Lemma~\ref{thm:273}, whose stages fit
into a pushout square
\begin{equation}
\xymatrix{
R_{1}\wedge \partial_{A}\sym^{m}B  \ar[r]\ar[d]  & R_{1}\wedge \sym^{m}B  \ar[d] \\
\fil_{m-1}R_{2}  \ar[r]        &\fil_{m}R_{2}\mathrlap{\ ,}
}
\end{equation}
where 
\[
\partial_{A}\sym^{m}B = \left(\partial_{A}B^{\wedge m}  \right)/\Sigma_{m}.
\]
The filtration of \S\ref{sec:index-mono-prod-pushout} mediates between
$(\fil_{m-1}R_{2})^{\wedge J}$ and $(\fil_{m} R_{2})^{\wedge J}$ by
another sequence of pushout squares.  The upper right hand corner of a
typical stage is an indexed wedge of terms of the form
\begin{equation}
\label{eq:102}
(\fil_{m-1}R_{2})^{\wedge J_{0}} \wedge (R_{1}\wedge \sym^{m}B)^{\wedge J_{1}},
\end{equation}
indexed by the set-theoretic decompositions $J=J_{0}\amalg J_{1}$.

We need to know two things about this expression.  One is that the
left derived functor of its formation (in all variables) is computed
in terms of the expression itself, and the other is that formation of
each of the pushout squares we encounter is homotopical.    Motivated
by this we are led to consider a technical condition slightly stronger
than the requirement that~\eqref{eq:36} be a weak equivalence.  That
is the subject of the next section.

\subsubsection{Very flat diagrams}
\label{sec:very-flat-diagrams}

As in \S\ref{sec:iter-index-symm}, to make the diagrams more readable
we will use the notation
\[
X^{\wedge(K/L)} = q^{\wedge}_{\ast}X
\]
for the indexed smash product along a map $q:K\to L$ of finite  $G$-sets.

\begin{defin}
\label{def:60} An equivariant $J$-diagram $X$ {\em very flat} if it
has the following property: for every cofibrant approximation $\tilde
X\to X$, every diagram of finite $G$-sets
\[
J\xleftarrow{p} K \xrightarrow{q}{} L,
\]
and every weak equivalence of equivariant $L$-diagrams $\tilde Z\to Z$, the map
\begin{equation}
\label{eq:126}
(p^{\ast}\tilde X)^{\wedge(K/L)}\wedge \tilde Z \to 
(p^{\ast} X)^{\wedge(K/L)}\wedge Z 
\end{equation}
is a weak equivalence.
\end{defin}

Our main goal is to establish the following result.

\begin{prop}
\label{thm:16} If $R\in\spectra^{\cat{B}_{J}G}$ is cofibrant
commutative ring, then the equivariant $J$-diagram of spectra
underlying $R$ is very flat.
\end{prop}

The condition that $R$ be very flat certainly implies that~\eqref{eq:36} is a
weak equivalence.  Proposition~\ref{thm:16} therefore implies
\begin{cor}
\label{thm:256} The following diagram of left derived functors
commutes up to natural isomorphism 
\[
\xymatrix{
\ho\comm\spectra^{\cat{B}_{J}G}  \ar[r]^{\lder p_{!}}\ar[d]  &  \ho\comm\ugspectra{G} \ar[d] \\
\ho\spectra^{\cat{B}_{J}G}   \ar[r]_{(\slot)^{\lwedge J}}        & \ho\ugspectra{G}\mathrlap{\ .}
}
\]
\qed
\end{cor}

\begin{rem}
\label{rem:40} 
Since identity maps are weak equivalences, the
condition of being very flat implies that every arrow in the diagram
\[
\xymatrix{
(p^{\ast}\tilde X)^{\wedge(K/L)}\wedge \tilde Z   \ar[r]\ar[d]  & (p^{\ast}\tilde X)^{\wedge(K/L)}\wedge Z 
\ar[d] \\
(p^{\ast} X)^{\wedge(K/L)}\wedge \tilde Z   \ar[r]        & (p^{\ast} X)^{\wedge(K/L)}\wedge Z 
}
\]
is a weak equivalence.  In particular it implies that $(p^{\ast}
X)^{\wedge(K/L)}$ is flat.
\end{rem}

\begin{rem}
\label{rem:17} Since $\tilde X^{\wedge (K/L)}$ is cofibrant
(Proposition~\ref{thm:195}), and cofibrant objects are flat
(Proposition~\ref{thm-copied:250}), the top arrow in the above diagram
is always a weak equivalence.  It therefore suffices to check the very
flat condition when $\tilde Z\to Z$ is the identity map.
\end{rem}

\begin{rem}
\label{rem:16} If~\eqref{eq:126} is a weak equivalence for one
cofibrant approximation it is a weak equivalence for any cofibrant
approximation.  It therefore suffices to check the ``very flat''
condition for a single cofibrant approximation $\tilde X\to X.$
\end{rem}

\begin{lem}
\label{thm:121} Arbitrary wedges of very flat spectra are very flat.
Smash products of very flat spectra are very flat.  Filtered colimits
of very flat $G$-diagrams along $h$-cofibrations are very flat.
\end{lem}

\begin{pf}
The first assertion follows easily from the distributive law, and the
fact that the formation of indexed wedges is homotopical.  The second
follows from the fact that the formation of indexed smash products is
symmetric monoidal.  The third makes use of Proposition~\ref{thm:210}.
The details are left to the reader.
\end{pf}

\begin{eg}
\label{eg:1}
Here is one motivation for the definition of ``very flat.''
Suppose we are given a pushout square of equivariant $J$-diagrams
\[
\xymatrix{
A  \ar[r]\ar[d]  & B
\ar[d] \\
X  \ar[r]        & Y
}
\]
and we are interested in the filtration of $Y^{\wedge (K/L)}$
described in \S\ref{sec:index-mono-prod-pushout}, whose stages are
related by pushout squares 
\begin{equation}
\label{eq:100}
\entrymodifiers={+!!<0pt,\fontdimen22\textfont2>}
\xymatrix{
{\displaystyle 
\bigvee_{(\ell,K_{1})\in G_{n}}
X^{\wedge K_{0}}\wedge \partial_{A}B^{\wedge K_{1}}}
  \ar[r]\ar[d]  &   
{\displaystyle 
\bigvee_{(\ell,K_{1})\in G_{n}}
X^{\wedge K_{0}}\wedge B^{\wedge K_{1}}}\ar[d] \\
\fil_{n-1}Y^{\wedge K/L}  \ar[r]        & \fil_{n}Y^{\wedge K/L}\mathrlap{\ ,}
}
\end{equation}
where $G_{n}=G_{n}(K/L)$ is the $G$-set of pairs $(\ell, K_{1})$ with
$\ell\in L$ and $K_{1}\subset q^{-1}(\ell)$ a subset of cardinality
$n$, and the map $G_{n}\to L$ sends $(\ell, K_{1})$ to $\ell$.  For
$(\ell,K_{1})\in G_{n}$ we have written $K_{0}$ to denote the
complement of $K_{1}$ in $q^{-1}(\ell)$.

The condition that $B$ be very flat gives some control over the upper
right term.  To see this let $V_{n}=V_{n}(K/L)$ be the set of triples
$(\ell, K_{1}, k)$ with $(\ell, K_{1})\in G_{n}$ and $k \in K_{1}$.
We define maps
\[
J \xleftarrow{f} V_{n} \xrightarrow{g} G_{n}
\]
by 
\begin{align*}
f(\ell, K_{1},k) &=  q(k) \\
g(\ell, K_{1},k) &= (\ell,K_{1}).
\end{align*}
The spectra $X^{\wedge K_{0}}$ form an
equivariant $G_{n}$-diagram, which we denote $Z$.  The $B^{\wedge
K_{1}}$ are the constituents of $(f^{\ast} B)^{\wedge(V_{n}/G_{n})}$, and so the
indexed wedge occurring in the pushout square is
\[
\bigvee_{G_{n}} Z\wedge (f^{\ast} B)^{\wedge(V_{n}/G_{n})}.
\]
Since the formation of indexed wedges is homotopical, its homotopy
properties come down to understanding the homotopy properties of the
equivariant $G_{n}$-diagram $Z\wedge f^{\ast}
B^{\wedge(V_{n}/G_{n})}$, some of which are specified by the condition
that $B$ be very flat.  
\end{eg}

By replacing the category of equivariant $J$-diagrams with its arrow
category, we arrive at the notion of a {\em very flat object} of 
$\aspectra{\cat B_{J}G}$.   The formal properties of being very flat
persist in this context, and in particular 
the analogues of Remarks~\ref{rem:40},~\ref{rem:17},~\ref{rem:16}
and Lemma~\ref{thm:121} hold.  

To get a feel for the more particular aspects of very flat arrows,
suppose that $(A\to B)$ is an object of $\aspectra{\cat B_{J}G}$, and
$(\tilde A\to \tilde B)$ is a cofibrant approximation.  Consider a
weak equivalence of the form
\[
(\tilde X\to \ast)\to (X\to \ast).
\]
In this case the very flat condition becomes that
\[
(p^{\ast}(\tilde B/\tilde A)^{\wedge(K/L)}\wedge \tilde X \to
\ast) \to 
(p^{\ast} (B/A)^{\wedge(K/L)}\wedge  X \to
\ast)
\]
is a weak equivalence.  This is so if and only if $B/A$ is very flat.

Next consider a weak equivalence  of the form 
\[
(\ast \to \tilde X)\to (\ast\to X).
\]
The very flat condition in this case is that 
\begin{multline*}
(\partial_{p^{\ast}\tilde A}p^{\ast}\tilde B^{\wedge(K/L)}\wedge \tilde X \to
p^{\ast}\tilde B^{\wedge(K/L)}\wedge \tilde X ) \\
\to 
(\partial_{p^{\ast} A}p^{\ast}B^{\wedge(K/L)}\wedge  X \to
p^{\ast} B^{\wedge(K/L)}\wedge X )
\end{multline*}
is a weak equivalence.   This holds if and only if $B$ is very
flat and $(A\to B)$ satisfies the condition that 
\begin{equation}
\label{eq:192}
\partial_{p^{\ast}\tilde A}p^{\ast}\tilde B^{\wedge(K/L)}\wedge \tilde
X \to \partial_{p^{\ast} A}p^{\ast}B^{\wedge(K/L)}\wedge  X 
\end{equation}
is a weak equivalence.   If we happen to know that the indexed corner
maps
\begin{align*}
\partial_{p^{\ast}\tilde A}p^{\ast}\tilde B^{\wedge(K/L)} &\to \tilde
B^{\wedge(K/L)} \\
\intertext{and}
\partial_{p^{\ast} A}p^{\ast} B^{\wedge(K/L)} &\to 
B^{\wedge(K/L)} 
\end{align*}
are $h$-cofibrations then the leftmost horizontal maps in 
\[
\xymatrix{
\partial_{p^{\ast}\tilde A}p^{\ast}\tilde B^{\wedge(K/L)} \wedge \tilde X  \ar[r]\ar[d]  &   p^{\ast}\tilde B^{\wedge(K/L)} \wedge \tilde X \ar[r]\ar[d]   & p^{\ast}(\tilde B/\tilde A)^{\wedge(K/L)} \wedge \tilde X  \ar[d]\\
\partial_{p^{\ast}A}p^{\ast}B^{\wedge(K/L)}\wedge X   \ar[r]        &
p^{\ast}B^{\wedge(K/L)}\wedge X  \ar[r]         & p^{\ast}(B/A)^{\wedge(K/L)}\wedge X
}
\]
are $h$-cofibrations, hence flat.  Thus the middle and left vertical
arrows are weak equivalences if and only if the middle and right
vertical arrows are, or in other words if and only if both $B$ and
$B/A$ are very flat.  So in the presence of the condition above, a
necessary condition that $(A\to B)$ be a very flat arrow is that $B$
and $B/A$ are very flat.  This turns out to be sufficient.   We single
out the condition.

\begin{cond}
\label{cond:1}
For every $J\xleftarrow{p} K \xrightarrow{q} L$ the corner map 
\[
\partial_{p^{\ast}A}(p^{\ast}B)^{\wedge(K/L)} \to (p^{\ast}B)^{\wedge(K/L)}
\]
is an $h$-cofibration.  
\end{cond}

\begin{rem}
\label{rem:39} By Proposition~\ref{thm:267} and the monoid axiom for
$\aspectra{\cat B_{L}G}$, a cofibrant object $(A\to B)$ of
$\aspectra{\cat{B}_{J}G}$ is very flat and satisfies
Condition~\ref{cond:1}.
\end{rem}

\begin{lem}
\label{thm:244}
If $A_{1}\to A_{2}$ satisfies Condition~\ref{cond:1}, and both $A_{1}$ and $A_{2}/A_{1}$
are very flat, then  $A=(A_{1}\to A_{2})$ is very flat.
\end{lem}

\begin{pf}
Fix a diagram of finite $G$-sets 
\[
J\xleftarrow{p}K\xrightarrow{q}{}L
\]
let $\tilde A=(\tilde A_{1}\to\tilde A_{2})$ be a cofibrant
approximation to $A=(A_{1}\to A_{2})$, and
\begin{gather*}
\tilde X \to X \\
\tilde X = (\tilde X_{1}\to\tilde X_{2}) \\
X = (X_{1}\to X_{2})
\end{gather*}
a weak equivalence in $\aspectra{\cat B_{L}G}$.  By
Remark~\ref{rem:39}, $\tilde A$ also satisfies the conditions of the
lemma.   Let 
\[
X' \to X\to X'' 
\]
be the sequence
\[
(\ast\to X_{2}) \to (X_{1}\to X_{2})\to (X_{1}\to \ast).
\]
and $\tilde X'\to \tilde X\to\tilde X'$ the analogous sequence for
$\tilde X$.  The maps $X'\to X$ and $\tilde X'\to \tilde X$ are not
$h$-cofibrations, but they are so objectwise, and hence flat.

Consider the diagram
\begin{equation}
\label{eq:194}
\xymatrix{
p^{\ast}\tilde A^{\wedge (K/L)}\wedge \tilde X'  \ar[r]\ar[d]  & p^{\ast}\tilde A^{\wedge (K/L)}\wedge \tilde X  \ar[r]\ar[d]   &  p^{\ast}\tilde A^{\wedge (K/L)}\wedge \tilde X''  \ar[d]\\
p^{\ast}A^{\wedge (K/L)}\wedge X'   \ar[r]        & p^{\ast}A^{\wedge
(K/L)}\wedge X   \ar[r]         & p^{\ast}A^{\wedge (K/L)}\wedge X''\mathrlap{\ .}
}
\end{equation}
Our aim is to show that the middle vertical map is a weak equivalence.

The first step is to show that the left horizontal maps are flat.
This reduces us to checking that the left and right vertical maps are
weak equivalences.  For this, let's examine the bottom left horizontal
map in more detail.  It is given by
\begin{equation}
\label{eq:193}
(\partial_{p^{\ast}A_{1}}p^{\ast} A_{2}^{\wedge(K/L)}\wedge  X_{2} \to
p^{\ast} A_{2}^{\wedge(K/L)}\wedge X_{2} ) \to 
(C\to p^{\ast} A_{2}^{\wedge(K/L)}\wedge X_{2} ) 
\end{equation}
in which $C$ is defined by the pushout diagram
\begin{equation}
\label{eq:120}
\xymatrix{
\partial_{p^{\ast}A_{1}}p^{\ast} A_{2}^{\wedge(K/L)}\wedge  X_{1}   \ar[r]\ar[d]  &
p^{\ast} A_{2}^{\wedge(K/L)}\wedge  X_{1}\ar[d] \\
\partial_{p^{\ast}A_{1}}p^{\ast} A_{2}^{\wedge(K/L)}\wedge  X_{2}  \ar[r]
& C\mathrlap{\ .}
}
\end{equation}
When $A_{1}\to A_{2}$ satisfies Condition~\ref{cond:1} the top map
in~\eqref{eq:120} is an $h$-cofibration, hence so is the bottom map.
This means that~\eqref{eq:193} is an objectwise $h$-cofibration, and
so flat.  Since $\tilde A_{1}\to\tilde A_{2}$ also satisfies
Condition~\ref{cond:1} the upper left horizontal map in~\eqref{eq:194}
is also flat.  Thus we are reduced to checking that the maps
\begin{gather*}
p^{\ast}\tilde A^{\wedge(K/L)}\wedge \tilde X'\to p^{\ast}A^{(\wedge
K/L)}\wedge X' \\ 
p^{\ast}\tilde A^{\wedge(K/L)}\wedge \tilde X''\to p^{\ast}A^{(\wedge
K/L)}\wedge X''
\end{gather*}
are weak equivalences.  As described above, this fact for the second
map follows from the assumption that $A_{2}/A_{1}$ is very flat.   The
assertion in the case of the first map is that the middle and left vertical
arrows in
\[
\xymatrix{
\partial_{p^{\ast}\tilde A_{1}}p^{\ast}\tilde A_{2}^{\wedge(K/L)}\wedge \tilde
X_{2}  \ar[r]\ar[d]_{\therefore\sim}  & p^{\ast}\tilde A_{2}^{\wedge(K/L)}\wedge \tilde
X_{2}  \ar[r]\ar[d]^{\sim}   &   p^{\ast}(\tilde A_{2}/\tilde A_{1})^{\wedge(K/L)}\wedge \tilde
X_{2} \ar[d]^{\sim}\\
\partial_{p^{\ast}A_{1}}p^{\ast}A_{2}^{\wedge(K/L)}\wedge X_{2}
\ar[r]        &  p^{\ast}A_{2}^{\wedge(K/L)}\wedge X_{2}  \ar[r]
&  p^{\ast}(A_{2}/A_{1})^{\wedge(K/L)}\wedge X_{2}
}
\]
are weak equivalences.  Since $A_{2}$ and $A_{2}/A_{1}$ are very flat,
the middle and right vertical maps are weak equivalences.
Condition~\ref{cond:1} shows that the left horizontal maps are
$h$-cofibrations, hence flat.  It follows that the left vertical map
is a weak equivalence.
\end{pf}

We can now establish an important technical fact used in the proof of
Proposition~\ref{thm:16}.

\begin{lem}
\label{thm:242}
Suppose that $A\to B$ is a cofibrant object of
$\aspectra{\cat{B}_{J}G}$, $I$ is a $G$-set and $\Sigma\subset
\Sigma_{I}$ a $G$-stable subgroup.   Then 
\[
\sym^{I}_{\Sigma} (A\to B) = (\partial_{A}\sym^{I}_{\Sigma}B\to \sym^{I}_{\Sigma}B)
\]
is very flat. 
\end{lem}

\begin{pf}
Proposition~\ref{thm:238} implies that in this situation the map
$\sym^{I}_{\Sigma}(A\to B)$ satisfies Condition~\ref{cond:1}, and that
for {\em every} cofibrant $B$, $\sym^{I}_{\Sigma}B$ is very flat (so
both $\sym^{I}_{\Sigma}B$ and $\sym^{I}_{\Sigma}(B/A)$ are very flat).
The result then follows from Lemma~\ref{thm:244}.
\end{pf}

\begin{eg}
\label{eg:2}
Continuing with Example~\ref{eg:1}, the top map in~\eqref{eq:100}
arises naturally in the arrow category as 
\[
\bigvee_{G_{n}(K/L)} Z\wedge (p^{\ast} (A\to B)^{\wedge(K/L)}),
\]
where $Z$ is the identity arrow of the diagram $X^{\wedge J_{0}}$.
Since the formation of indexed wedges is homotopical, the information
in the homotopy type of this expression is contained in the
$G_{n}(K/L)$-diagram $Z\wedge (p^{\ast} (A\to B)^{\wedge(K/L)})$.  The
condition that $(A\to B)$ be very flat thus specifies good homotopical
properties of the top map in~\eqref{eq:100}.
\end{eg}

\begin{lem}
\label{thm:193}
Consider a pushout square 
\begin{equation}
\label{eq:125}
\xymatrix{
A  \ar[r]\ar[d]  & B  \ar[d] \\
X  \ar[r]        & Y
}
\end{equation}
in which $(A\to B)$ is a very flat object of $\aspectra{\cat B_{J}G}$
satisfying Condition~\ref{cond:1}.  If $X$ is very flat, then so is
$Y$.
\end{lem}

\begin{pf}
Using the fact that cofibrations are flat, we can arrange things so that
the cofibrant approximation $\tilde Y\to Y$ fits into a pushout square
\begin{equation}
\xymatrix{
\tilde A  \ar[r]\ar[d]  & \tilde B  \ar[d] \\
\tilde X  \ar[r]        & \tilde Y
}
\end{equation}
of cofibrant approximations to~\eqref{eq:125}, in which $\tilde
A\to\tilde B$ is a cofibration.  We give $\tilde
Y^{\wedge (K/L)}$ and $Y^{\wedge (K/L)}$ the filtration described in
\S\ref{sec:index-mono-prod-pushout}.  We will prove by induction on
$n$ that for any weak equivalence $\tilde Z\to Z$ of equivariant
$J$-diagrams, the map
\begin{equation}
\label{eq:129}
\fil_{n}\tilde Y^{\wedge K}\wedge \tilde Z \to \fil_{n}Y^{\wedge K}\wedge Z
\end{equation}
is a weak equivalence.  The case $n=0$ is the assertion that $X$ is
very flat, which is true by assumption.  For the inductive step,
consider the diagram
\[
\xymatrix@C=1em@!R{
\fil_{n-1}\tilde Y^{\wedge K}\wedge \tilde Z \ar[d] &{\displaystyle 
\bigvee_{G_{n}(K/L)}}
\tilde X^{\wedge K_{0}}\wedge \partial_{\tilde A}\tilde B^{\wedge K_{1}}\wedge \tilde Z
 \ar[l] \ar[r]\ar[d]  &   
{\displaystyle 
\bigvee_{G_{n}(L)}}
\tilde X^{\wedge K_{0}}\wedge \tilde B^{\wedge K_{1}}\wedge \tilde Z\ar[d] 
\\
\fil_{n-1}Y^{\wedge K}\wedge Z &{\displaystyle 
\bigvee_{G_{n}(K/L)}}s
X^{\wedge K_{0}}\wedge \partial_{A}B^{\wedge K_{1}}\wedge Z
 \ar[l] \ar[r]  &   
{\displaystyle 
\bigvee_{G_{n}(K/L)}}
X^{\wedge K_{0}}\wedge B^{\wedge K_{1}}\wedge Z \mathrlap{\ .}
}
\]
The map from the pushout of the top row to the pushout of the bottom
row is~\eqref{eq:129}.  The rightmost horizontal maps are
$h$-cofibrations by Condition~\ref{cond:1}.  The left vertical map is
a weak equivalence by induction, and the other two vertical maps are
weak equivalences since $(A\to B)$ is very flat (Example~\ref{eg:2}).
The map of pushouts is therefore a weak equivalence since
$h$-cofibrations are flat.
\end{pf}

\subsubsection{Proof of Proposition~\ref{thm:16}}
\label{sec:proof-prop-refthm:16}

Since the class of very flat $G$-diagrams is closed under the
formation of filtered colimits along $h$-cofibrations
(Lemma~\ref{thm:121}), it suffices to show that if $A\to B$ is a
generating cofibration in $\spectra^{\cat{B}_{J}G}$,
\[
\xymatrix{
\sym A  \ar[r]\ar[d]  & \sym B  \ar[d] \\
X  \ar[r]        & Y
}
\]
is a pushout square of commutative $J$-algebras, and $X$ is very flat,
then $Y$ is very flat.  Working fiberwise, the filtration described
after the statement of Lemma~\ref{thm:273} gives a filtration of $Y$
by $X$-modules, whose stages are related by the pushout squares
\begin{equation}
\label{eq:20}
\xymatrix{
X\wedge \partial_{A}\sym^{m}B  \ar[r]\ar[d]  & X\wedge \sym^{m}B  \ar[d] \\
\fil_{m-1}Y  \ar[r]        &\fil_{m}Y\mathrlap{\ .}
}
\end{equation}
We show by induction on $m$ that each $\fil_{m}Y$ is very flat.  Since
$\fil_{0}Y=X$, the induction starts.  The arrow
$(\partial_{A}\sym^{m}B\to \sym^{m}B)$ is very flat by
Lemma~\ref{thm:242}.  This means that the top row of~\eqref{eq:20} is
a very flat arrow, since smash products of very flat objects are very
flat (Lemma~\ref{thm:121}).  This places us in the situation of
Lemma~\ref{thm:193} which completes the inductive step.

\subsection{The slice tower, symmetric powers and the norm}
\label{sec:slice-tower-symm} 

The main goal of this section is to show that if $R$ is an equivariant
commutative ring in $\geslice{0}$, and $n\ge 0$ is an integer, then
the slice section $P^{n}R$ is also an equivariant commutative ring in
$\geslice{0}$.  The proof makes use of the technology used to show
that cofibrant commutative rings are very flat, and so has been
deferred to this appendix.  The reader may wish to look through the
first three subsections of \S\ref{sec:ths-slice-filtration} for the
basic definitions concerning the slice tower.

Our presentation of the slice tower was made in a context where the
emphasis was on homotopy theory, and the slice sections $P^{n}$, etc
were characterized by homotopy theoretic properties.  Here we will be
making use of some explicit constructions, and some care needs to be
taken to ensure that the derived functors we are ultimately interested
in can be computed on the objects that arise.  Using the fact that
indexed smash products of cofibrant objects are cofibrant, and that
indexed symmetric powers of cofibrant spectra are flat, one can check
that this is indeed the case.  We will take as the definition of
$P^{n}$ the colimit of the inductive construction described in
\S\ref{sec:slice-tower}, using the cofibrant approximations
$S^{-1}\wedge S^{1}\wedge G_{+}\smashove{H}S^{k\rho_{H}}$ and
$S^{-1}\wedge S^{1}\wedge G_{+}\smashove{H}S^{k\rho_{H}-1}$ for the
slice cells.  This particular choice of $P^{n}$ is homotopical, and
the natural map $X\to P^{n}X$ is a cofibration.  Our task will be to
show that something functorially weakly equivalent to $P^{n}$ takes
commutative rings in $\geslice{0}$ to commutative rings in
$\geslice{0}$.

We begin with the interaction of the slice filtration with the
formation of indexed smash products.  As in
\S\ref{sec:homot-prop-norm} we fix a finite $G$-set and work with the
homotopy theory of equivariant $J$-diagrams.  We define slice cells
and the slice filtration in the evident manner, so that the slice
filtration on equivariant $J$-diagrams corresponds to the product of
slice filtrations on $G_{t}$-spectra under the equivalence
\[
\spectra^{\cat{B}_{J}G} \approx \prod_{t}\ugspectra{G_{t}},
\]

The proposition below follows easily from Proposition~\ref{thm:67}.  

\begin{prop}
\label{thm:99} Suppose that $J$ is a non-empty $G$-set, $X$ is a
cofibrant equivariant $J$-diagram, and $n$ is an integer.  If each
$X_{j}$ is slice $(n-1)$ positive, then the indexed wedge
\[
\bigvee_{j\in J} X_{j}
\]
is slice $(n-1)$ positive.  \qed
\end{prop}

The next two results make use of the implication
\begin{equation}
\label{eq:93}
X\ge 0\quad\text{and}\quad Y\ge k\implies X\wedge Y\ge k
\end{equation}
proved in \S\ref{sec:mult-prop-slice} (Proposition~\ref{thm:6}).

\begin{prop}
\label{thm:98} Suppose that $J$ is a non-empty $G$-set, $X$ is a
cofibrant equivariant $J$-diagram, and $n\ge 0$ is an integer.  If
each $X_{j}$ is slice $(n-1)$ positive, then the indexed smash product
\[
\bigwedge_{j\in J} X_{j}
\]
is slice $(n-1)$ positive.  
\end{prop}

\begin{pf}
By induction on $|G|$ we may suppose that $i_{H}^{\ast}X^{\wedge J}$
is slice $(n-1)$-positive for any proper subgroup $H\subset G$.  This
implies that $T\wedge X^{\wedge J} \ge n$ if $T$ is any $G$-CW complex
built entirely from induced $G$-cells.  Since the formation of indexed
smash products commutes with filtered colimits, it suffices by
Proposition~\ref{thm:35} to consider a cofibration $A\to B$ of
equivariant $J$-diagrams in which $B/A$ is a wedge of slice cells of
dimension greater than $n$, and show that
\begin{equation}
\label{eq:63}
A^{\wedge J}\ge n \implies B^{\wedge J} \ge n.
\end{equation}
Using the filtration of \S\ref{sec:index-mono-prod-pushout} for the
identity pushout square
\[
\xymatrix{
A  \ar[r]\ar[d]  & B
\ar[d] \\
A  \ar[r]        & B\mathrlap{\ ,}
}
\]
gives a filtration of $B^{\wedge J}$ whose stages fit into cofibration
sequences
\begin{equation}
\label{eq:62}
\fil_{m-1}B^{\wedge J} \to \fil_{m}B^{\wedge J} \to \bigvee A^{\wedge
J_{0}}\wedge (B/A)^{\wedge J_{1}}
\end{equation}
in which the indexing $G$-set for the coproduct is the set of all set
theoretic decomposition $J=J_{0}\coprod J_{1}$ with $|J_{1}|=m$.  The
implication~\eqref{eq:93} and Proposition~\ref{thm:99} above reduce
the claim to showing that if $J_{1}\ne\emptyset$, then $(B/A)^{\wedge
J_{1}}$ (regarded as an equivariant spectrum for the stabilizer of
$J_{1}$) is slice $(n-1)$-positive.  In other words, it suffices to
prove the proposition when $X$ is a wedge of slice cells of dimension
greater than or equal to $n$.

Making use of the distributive law, and once again~\eqref{eq:93} and
Proposition~\ref{thm:99}, one reduces to the case in which $J=G/H$ is
a single orbit, and $X$ corresponds to $S^{k\rho_{H}}$ with $k|H|\ge n$
or $S^{k\rho_{H}-1}$, with $k|H| -1 \ge n$.  In the first case
\[
X^{\wedge J } \approx S^{k\rho_{G}} 
\]
has dimension $k|G|\ge k|H|\ge n$.  In the second case
\[
X^{\wedge J } \approx S^{(n-1)\rho_{G}+V}
\]
where $V=\rho_{G}-\ind_{H}^{G}1$.   Write $W=\ind_{H}^{G}-1$ so that
$S^{W}\wedge S^{V} \approx S^{\rho_{G}-1}$, and there is a cofibration
sequence 
\begin{equation}
\label{eq:45}
S(W)_{+}\wedge X^{\wedge J} \to X^{\wedge J} \to S^{(n-1)\rho_{G}+(\rho_{G}-1)}.
\end{equation}
The $G$-space $S(W)$ is homeomorphic to the boundary of the simplex
with vertices $G/H$, and has no $G$-fixed points.  The barycentric
subdivision gives $S(W)_{+}$ the structure of a $G$-CW complex built
entirely from induced $G$-cells.  It therefore follows from our induction
hypothesis that
\[
S(W)_{+}\wedge X^{\wedge J}
\]
is slice $(n-1)$-positive.  The rightmost term in~\eqref{eq:45} is a
slice cell of dimension
\[
k|G| -1 \ge k|H|-1 \ge n
\]
It follows that $X^{\wedge J}$ is slice $(n-1)$-positive.
\end{pf}

\begin{rem}
\label{rem:15} We will later need to know that in the situation of
Proposition~\ref{thm:98}, one has
\begin{equation}
\label{eq:94}
\Sigma^{-1} (\Sigma X)^{\wedge J} \ge n.
\end{equation}
To see this, rewrite the spectrum in~\eqref{eq:94} as
\[
\big(\Sigma^{-1} (S^{1})^{\wedge J})\wedge (X^{\wedge J}).
\]
The factor $\Sigma^{-1} (S^{1})^{\wedge J}$ is weakly equivalent to
the sphere $S^{V}$ with $V=\R^{J}-1$.  This gives
\[
\Sigma^{-1} (S^{1})^{\wedge J}\ge
0
\]
and the relation~\eqref{eq:94} then follows from
Proposition~\ref{thm:98} and~\eqref{eq:93}.
\end{rem}

We next turn to indexed symmetric powers.  As in
\S\ref{sec:symmetric-powers} we consider a finite $G$-set $I$, a
$G$-stable subgroup $\Sigma\subset \Sigma_{I}$, and the indexed
symmetric power
\[
\sym^{I}_{\Sigma}X = X^{\wedge I}/\Sigma.
\]

\begin{prop}
\label{thm:91} Let $n\ge 0$ be an integer, $I$ a non-empty $G$-set,
and $X$ a cofibrant equivariant $I$-diagram.  If $X$ is slice
$(n-1)$-positive then both the indexed symmetric power
$\sym^{I}_{\Sigma}X$ and $\Sigma^{-1}\sym^{I}_{\Sigma}(\Sigma X)$
are slice $(n-1)$-positive.
\end{prop}

\begin{pf}
Using the equivalences
\begin{align*}
(\egsigma{G})_{+}\smashove{\Sigma} X^{\wedge I} &
\approx\sym^{I}_{\Sigma}X  \\
\Sigma^{-1} (\egsigma{G})_{+}\smashove{\Sigma}
(\Sigma X)^{\wedge I} & \approx \Sigma^{-1}\sym^{I}_{\Sigma}(\Sigma X) 
\end{align*}
of Proposition~\ref{thm:200} and working through an equivariant cell
decomposition of $\egsigma{G}$  reduces the claim to showing that 
\begin{equation}
\label{eq:92}
S_{+}\smashove{\Sigma}X^{\wedge I}  \quad\text{and}\quad
\Sigma^{-1}S_{+}\smashove{\Sigma}(\Sigma X)^{\wedge I} 
\end{equation}
are slice $(n-1)$-positive when $S$ is a finite $\Sigma$-free
$\Sigma\rtimes G$-set.  But the first spectrum in~\eqref{eq:92} is an
indexed wedge of indexed smash products of $X$ (see the proof of
Lemma~\ref{thm:201}), hence slice $(n-1)$-positive by
Propositions~\ref{thm:98} and~\ref{thm:99}.  The second spectrum is an
indexed wedge of desuspensions of indexed smash products of $\Sigma
X$, hence slice $(n-1)$-positive by Remark~\ref{rem:15} and
Proposition~\ref{thm:99}.
\end{pf}

We can now investigate the slice sections of commutative rings.  Let
$\pnalg{n}:\ugeinftycat{G}\to\ugeinftycat{G}$ be the multiplicative
analogue of $P^{n}$, constructed as the colimit of a sequence of
functors
\[
W^{\text{alg}}_{0}R\to W^{\text{alg}}_{1}R\to\cdots.
\]
The $W^{\text{alg}}_{i}R$ are defined inductively starting with $W^{\text{alg}}_{0}R=R$, and
in which $W^{\text{alg}}_{k}R$ is defined by the pushout square
\[
\xymatrix{
{\displaystyle\sym\big(\bigvee_{I} \Sigma^{t}\slicecell\,\big)}
\ar[r]\ar[d]  & {\displaystyle\sym\big(\bigvee_{I} C\,\Sigma^{t}\slicecell\,\big)}
\ar[d] \\
W^{\text{alg}}_{k-1}R  \ar[r]        & W^{\text{alg}}_{k}R
}
\]
in which the indexing set $I$ is the set of maps
$\Sigma^{t}\slicecell\to W^{\text{alg}}_{k-1}R$ with $\slicecell>n$ a
cofibrant slice cell and $t\ge 0$.  The functor $\pnalg{n}$ is
homotopical and for any $R$, the map $R\to \pnalg{n}R$ is a
cofibration of equivariant commutative rings.  The arrow
$R\to\pnalg{n}R$ is characterized up to weak equivalence by the
following universal property: if $S$ is an equivariant commutative
ring whose underlying spectrum is slice $(n+1)$-null then the map
\[
\ho\ugeinftycat{G}(\pnalg{n}R,S) \to 
\ho\ugeinftycat{G}(R,S) 
\]
is an isomorphism.   

For clarity let's temporarily denote by $U$, the
forgetful functor
\[
U:\ugeinftycat{G}\to \ugspectra{G}.
\]
By the small object argument, the spectrum
$U\pnalg{n}R$  is slice $(n+1)$-null, so there is a natural
transformation 
\[
P^{n}UR\to U\pnalg{n}R
\]
of functors to $\ugspectra{G}$. 

\begin{prop}
\label{thm:79} If $R$ is a slice $(-1)$-positive cofibrant equivariant
commutative ring, then for all $n\in\Z$, the map
\[
P^{n}UR\to U\pnalg{n}R
\]
is a weak equivalence.
\end{prop}

\begin{pf}
When $n$ is negative, $P^{n}UR$ is contractible, and $\pnalg{n}R$ is a
commutative ring whose unit is nullhomotopic, hence also contractible.
We may therefore assume $n$ is non-negative.

It suffices to show that each of the maps
\[
UW^{\text{alg}}_{k-1}R \to UW^{\text{alg}}_{k}R
\]
is a $P^{n}$-equivalence.  We do this by working through the
filtration used in the proof of Lemma~\ref{thm:273}, whose successive
terms are related by the homotopy cocartesian square
\[
\xymatrix{
U W^{\text{alg}}_{k-1}R\wedge \partial_{A}\sym^{m}B  \ar[r]\ar[d]  & U W^{\text{alg}}_{k-1}R\wedge \sym^{m}B  \ar[d] \\
\fil_{m-1}W^{\text{alg}}_{k}R  \ar[r]        &\fil_{m}W^{\text{alg}}_{k}R\mathrlap{\ ,}
}
\]
in which $A\to B$ is the map
\begin{equation}
\label{eq:96}
\bigvee_{I} \Sigma^{t}\slicecell \to \bigvee_{I} C\,\Sigma^{t}\slicecell.
\end{equation}
By induction we may assume that the maps 
\[
U R\to U W_{k-1}^{\text{alg}}R \to \fil_{m-1}W^{\text{alg}}_{k}R
\]
are $P^{n}$ equivalences, and so among other things that the three
spectra are all in $\geslice{0}$.  The homotopy fiber of
$\fil_{m-1}W^{\text{alg}}_{k}R\to \fil_{m}W^{\text{alg}}R$ is
$U W_{k-1}^{\text{alg}}R\wedge \Sigma^{-1}\sym^{m}(B/A)$.  Now $B/A$ is
the suspension of the left term in~\eqref{eq:96} which is slice
$n$-positive.  It follows (Proposition~\ref{thm:91}) that
$\Sigma^{-1}\sym^{m}(B/A)$ is also slice $n$-positive hence so is
$U W_{k-1}^{\text{alg}}R\wedge \Sigma^{-1}\sym^{m}(B/A)$ since
$U W_{k-1}^{\text{alg}}R\ge 0$.  The fact that
$\fil_{m-1}W^{\text{alg}}_{k}R\to \fil_{m}W^{\text{alg}}R$ is a
$P^{n}$-equivalence is now a consequence of Lemma~\ref{thm:104}.
\end{pf}

\subsection{Geometric and  monoidal geometric fixed points}
\label{sec:geom-fixed-points}

The geometric fixed point functor was defined and its main properties
summarized in \S\ref{sec:isotropy-separation}.  In this section we
gives proofs of some of these properties, and describe the variation
constructed in Mandell-May~\cite[\S{}V.4]{MR1922205}.  We refer to the
Mandell-May construction as the {\em monoidal geometric fixed point
functor} and denote it $\phigm$, in order not to confuse it with the
usual geometric fixed point functor.

\subsubsection{Geometric fixed points}
\label{sec:geom-fixed-points-2}

The geometric fixed point functor was defined in~\S\ref{sec:isotropy-separation} by
\[
\phig(X) = \big((\tilde E\pfamily\wedge X)_{f}\big)^{G},
\]
in which the $G$-CW complex $\tilde E\pfamily$ is the one characterized up to
equivariant homotopy equivalence by the property
\[
\big(\tilde E\pfamily\big)^{H} \sim 
\begin{cases}
S^{0}& H=G \\
\ast & H\ne G.
\end{cases}
\]
The characterizing property of $\tilde E\pfamily$ implies that for any
{\em $G$-space} $Z$ and any $G$-CW complex $A$, the restriction map
\[
[A,\tilde E\pfamily \wedge Z]^{G} \to 
[A^{G},\tilde E\pfamily \wedge Z]^{G}
\]
is an isomorphism.   Since $G$-acts trivially on $A^{G}$, the right
hand side is isomorphic to 
\[
[A^{G},\big(\tilde E\pfamily \wedge Z\big)^{G}] = 
[A^{G},Z^{G}].
\]
Combining these gives the isomorphism 
\begin{equation}
\label{eq:7}
[A,\tilde E\pfamily \wedge Z]^{G} \approx
[A^{G},X^{G}].
\end{equation}
This isomorphism is the foundation for our investigation into $\phig$.

Let $\iota:\spectra\to\ugspectra{G}$
be the functor which regards a spectrum as a $G$-spectrum with
trivial action.  As described in \S\ref{sec:fixed-point-spectra}, the
fixed point functor $(\slot)^{G}$ is right adjoint to $\iota$
\[
\iota:\spectra \leftrightarrows \ugspectra{G}: (\slot)^{G}
\]
and together they form a Quillen morphism in the positive complete
model structures.

Since smashing with $\tilde E\pfamily$ is homotopical, and the fixed
point functor $(\slot)^{G}$ is a right Quillen functor, the functor
$\phig$ is homotopical.   

\begin{prop}
\label{thm:246} 
For a spectrum $X\in \spectra$, the map
\begin{equation}
\label{eq:128}
X\to \phig(\iota X)
\end{equation}
adjoint to
\[
\iota X\to \tilde E\pfamily \wedge \iota X \to(\tilde E\pfamily \wedge \iota X)_{f}.
\]
is a weak equivalence.
\end{prop}

\begin{pf}
We have
\begin{align*}
\pi_{k}\phig(i_{\ast}X) &\approx \ho\spectra(S^{k},\big(\tilde
E\pfamily\wedge X)^{G}_{f}\big) \\
&\approx \ho\spectra(S^{k},\big(\tilde
E\pfamily\wedge X)^{G}\big) \\
&\approx \ho\ugspectra{G}(S^{k},\tilde
E\pfamily\wedge X) \\
&\approx \varinjlim_{W> -k} \pi^{G}_{k+W}\tilde E\pfamily\wedge X_{W} \\
&\approx \varinjlim_{W> -k} \pi^{G}_{k+W^{G}}X_{W^{G}}
\end{align*}
with the last isomorphism coming from~\eqref{eq:7}.  Under the
composite isomorphism, the map on stable homotopy groups induced
by~\eqref{eq:128} is
\[
\varinjlim_{V>-k} \pi_{k+V}X_{V} \to 
\varinjlim_{W> -k} \pi^{G}_{k+W^{G}}X_{W^{G}},
\]
in which $V$ is ranging through the poset of finite dimensional
orthogonal vector spaces and $W$ through the poset of
$G$-representations.  This is clearly an isomorphism.
\end{pf}

Since $\tilde E\pfamily$ is $H$-equivariantly contractible when $H$ is
a proper subgroup of $G$, the smash product $\tilde E\pfamily\wedge X$
is contractible if $X$ is a cellular spectrum built entirely from
$G$-cells induced from a proper subgroup of $G$.  More generally
\begin{lem}
\label{thm:245} 
Let $A$ and $Y$ be $G$-spectra.  If $X$ is constructed from $A$ by
attaching $G$-cells induced from proper subgroups, then the inclusion
$A\to X$ induces a weak equivalence
\[
\tilde E\pfamily \wedge A\wedge Y \xrightarrow{\sim}{}
\tilde E\pfamily \wedge X \wedge Y
\]
hence a weak equivalence 
\[
\phig \big(A\wedge Y\big)\xrightarrow{\sim}{} \phig \big(X\wedge Y\big). 
\]
\qed
\end{lem}

\begin{cor}
\label{thm:205}
Let $V$ be a $G$-representation and $A$ a $G$-CW complex.  The maps 
\[
S^{-V^{G}}\wedge A^{G} \to S^{-V^{G}} \wedge A \leftarrow S^{-V} \wedge A
\]
induce weak equivalences  
\[
S^{-V^{G}}\wedge A^{G} \sim \phig \big(S^{-V^{G}}\wedge A^{G}
\big)\xrightarrow{\sim}{} \phig \big(S^{-V^{G}}\wedge A
\big)\xleftarrow{\sim} \phig \big(S^{-V}\wedge A\big),
\]
giving a zig-zag of weak equivalences
\[
\phig \big(S^{-V}\wedge A\big) \xleftrightarrow{\sim} S^{-V^{G}}\wedge A^{G}.
\]
\end{cor}

\begin{pf}
We work our way from the left.  The weak equivalence  $S^{-V^{G}}\wedge A^{G}\approx \phig
(A^{G}\wedge S^{-V^{G}})$ is Proposition~\ref{thm:246}.   The next map
is a weak equivalence by Lemma~\ref{thm:245} since $A$ is
constructed from $A^{G}$ by adding induced $G$-cells.   The last map can
be constructed by applying $\phig$ to the composition
\[
S^{-V}\wedge A
\to
S^{-V}\wedge S^{V-V^{G}}\wedge A
\to
S^{-V^{G}}\wedge A.
\]
The right arrow is a weak equivalence.  Since $S^{V-V^{G}}$ is a
$G$-CW complex with fixed point space $S^{0}$, it is constructed from
$S^{0}$ by adding induced $G$-cells.  The left map therefore induces an
equivalence of geometric fixed points by Lemma~\ref{thm:245}.
\end{pf}

\subsubsection{Motivation and definition of the monoidal geometric
fixed point functor}
\label{sec:motivation}

For an orthogonal representation $V$ of $G$ let $\fixg{V}\subset V$ be
the space of invariant vectors, and $V^{\perp}$ the orthogonal
complement of $\fixg{V}$.  Note that
\begin{equation}
\label{eq:145}
\igcat{G}(V,W)^{G} \approx \icat(\fixg{V},\fixg{W})\wedge O(V^{\perp},
W^{\perp})^{G}_{+},
\end{equation}
so that there is a canonical map 
\[
\igcat{G}(V,W)^{G} \to \icat{}(\fixg{V},\fixg{W}),
\]
given in terms of~\eqref{eq:145} by smashing the identity map with the
map $O(V^{\perp}, W^{\perp})^{G}\to \text{pt}$.

We wish to define a functor $\phigm$ with the property that
\begin{equation}
\label{eq:69}
\phigm(S^{-V}\wedge A) = S^{-\fixg{V}}\wedge A^{G}
\end{equation}
and which commutes with colimits as far as is possible.  A value needs
to be assigned to the effect of $\phigm$ on the map
\[
S^{-W}\wedge \igcat{G}(V,W)\to S^{-V}.
\]
The only obvious choice is to take 
\[
\phigm(S^{-W}\wedge \igcat{G}(V,W))\to \phigm(S^{-V})
\]
to be the composite
\begin{equation}
\label{eq:68}
S^{-\fixg{W}}\wedge \igcat{G}(V,W)^{G}  \to
S^{-\fixg{W}}\wedge \icat(\fixg{V},\fixg{W})  \to
S^{-\fixg{V}}.
\end{equation}

If $\phigm$ actually {\em were} to commute with colimits, it would be
determined by the specifications given by~\eqref{eq:69}
and~\eqref{eq:68}.  Indeed, using the tautological presentation to write
a general equivariant orthogonal spectrum $X$ as a reflexive
coequalizer
\[
\bigvee_{V,W} S^{-W}\wedge\igcat G(V,W)\wedge
X_{V}\rightrightarrows \bigvee_{V} S^{-V}\wedge X_{V} \to X,
\]
the value of $\phigm(X)$ would be given by the reflexive coequalizer diagram
\begin{equation}
\label{eq:70}
\bigvee_{V,W} S^{-\fixg{W}}\wedge\igcat G(V,W)^{G}_{+}\wedge
X_{V}^{G}\rightrightarrows \bigvee_{V} S^{-\fixg{V}}\wedge X_{V}^{G} \to
\phigm X.  
\end{equation}
We take this as the definition of $\phigm(X)$.

\begin{defin}
\label{def:3}
The {\em monoidal geometric fixed point functor} 
\[
\phigm:\ugspectra{G} \to \spectra
\]
is the functor defined by the coequalizer diagram~\eqref{eq:70}.
\end{defin}

\begin{rem}
\label{rem:13}
In case $X=S^{-V}\wedge A$, the tautological presentation is a split
coequalizer, and one recovers both~\eqref{eq:69} and~\eqref{eq:68}.
\end{rem}

A fundamental property of the usual geometric fixed point functor
$\phig$ is that for proper $H\subset G$, the spectrum
$\phig(G_{+}\smashove{H}X)$ is contractible.  The monoidal geometric
fixed point functor has this property on the nose.

\begin{prop}
\label{thm:73}
Suppose that $J$ is a $G$-set and $X$ an equivariant $J$-diagram.   If
$J$ has no $G$-fixed points then the map 
\[
\phigm\big(\bigvee_{j\in J}X_{j}) \to \ast
\]
is an isomorphism.  In particular, if $H\subset G$ is a proper
subgroup and $X$ an orthogonal $H$-spectrum, then the map
\[
\phigm(G_{+}\smashove{H}X) \to \ast
\]
is an isomorphism.
\end{prop}

\begin{pf}
Since indexed wedges are computed componentwise, the assumption that
$J$ has no fixed points implies that for all representations $W$
of $G$,
\[
\big(\bigvee_{j\in J}X_{j}\big)_{W}^{G} = \big(\bigvee_{j\in
J}(X_{j})_{W}\big)^{G} = \ast.
\]
The claim then follows from the definition of $\phigm$.
\end{pf}

Working through an equivariant cell decomposition gives

\begin{cor}
\label{thm:54} Let $A$ and $Y$ be $G$-spectra.  If $X$ is constructed
from $A$ by attaching $G$-cells induced from proper subgroups, then
the map $\phigm(A\wedge Y)\to \phigm(X\wedge Y)$ is an isomorphism.
\qed
\end{cor}

There is a natural map 
\begin{equation}
\label{eq:134}
X^{G}\to \phigm X
\end{equation}
from the fixed point spectrum of $X$ to the monoidal geometric fixed
point spectrum.   To construct it note that the fixed point spectrum
of $X$ is computed termwise, and so is given by the coequalizer
diagram 
\begin{equation}
\label{eq:133}
\bigvee_{V,W\in \icat} S^{-W}\wedge\icat(V,W)_{+}\wedge
X_{V}^{G}\rightrightarrows \bigvee_{V\in\icat } S^{-V}\wedge X_{V}^{G} \to
X^{G}.  
\end{equation}
The map~\eqref{eq:134} is given by the evident inclusion
of~\eqref{eq:133} into~\eqref{eq:70}.

The functor $\phigm$ cannot commute with all colimits.  However, since
colimits of orthogonal $G$-spectra are computed objectwise, the
definition implies that $\phigm$ commutes with whatever enriched
colimits are preserved by the fixed point functor on $G$-spaces.  This
means that there is a functorial isomorphism
\begin{equation}
\label{eq:76}
\phigm(X\wedge A) \approx \phigm(X)\wedge A^{G}
\end{equation}
for each pointed $G$-space $A$, and that $\phigm$ commutes with the
formation of wedges, directed colimits and cobase change along a
closed inclusion.  Because $h$-cofibrations and cofibrations are
objectwise closed inclusion (Lemma~\ref{thm:105} and
Remark~\ref{rem:36}), the functor $\phigm$ has good homotopy
theoretic properties.

\subsubsection{Homotopy properties of $\phigm$}
\label{sec:homot-prop-phig}

Several variations on the following appear in 
in~\cite[\S{}V.4]{MR1922205}.  

\begin{prop}
\label{thm:214}
The functor $\phigm$ sends cofibrations to cofibrations and acyclic
cofibrations to acyclic cofibrations.   It therefore sends weak
equivalences between cofibrant objects to weak equivalences.
\end{prop}

\begin{pf}
That $\phigm$ sends cofibrations to cofibrations follows from the fact
that it preserves cobase change along closed inclusions and sends
generating cofibrations to generating cofibrations.   A similar
argument applies to the acyclic cofibrations, once one checks that $\phigm$
sends both maps in the factorization~\eqref{eq:147} 
\[
S^{-V\oplus W}\wedge S^{W}\to \tilde S^{-V} \to S^{-V}
\]
to weak equivalences.  But the second map is a homotopy equivalence
and the composite map is sent to a weak equivalence by~\eqref{eq:69}.
The last assertion is a consequence of Ken Brown's Lemma.
\end{pf}

Proposition~\ref{thm:214} implies that the monoidal geometric fixed
point functor has a left derived functor which can be computed on any
cofibrant approximation.  A similar argument with a slightly different
model structure could be used to show that the left derived functor
can be computed on a cellular approximation.    We will show in
\S\ref{sec:relat-with-isotr} that the left derived functor
$\lder\phigm$ is the geometric fixed point functor $\phig$.

\subsubsection{Monoidal geometric fixed points and smash product}
\label{sec:geom-fixed-points-1}

The properties~\eqref{eq:69} and~\eqref{eq:68} give an identification
\[
\phigm(S^{-V}\wedge A\wedge S^{-W}\wedge B) \approx
\phigm(S^{-V}\wedge A)\wedge \phigm(S^{-W}\wedge B)
\]
making the diagram
\[
\xymatrix{
\phigm(S^{-V_{1}}\wedge \igcat{G}(W_{1},V_{1}))\wedge
\phigm(S^{-V_{2}}\wedge \igcat{G}(W_{2},V_{2}))   \ar[r] \ar[d]  
&  \phigm(S^{-W_{1}})\wedge \phigm(S^{-W_{2}})\ar[d] \\
\phigm(S^{-V_{1}}\wedge \igcat{G}(W_{1},V_{1})\wedge
S^{-V_{2}}\wedge \igcat{G}(W_{2},V_{2}))     \ar[r] &
\phigm(S^{-W_{1}}\wedge S^{-W_{2}})
}
\]
commute.  Applying $\phigm$ termwise to the smash product of the
tautological presentations of $X$ and $Y$, and using the above
identifications, gives a natural transformation
\begin{equation}\label{eq:67}
\phigm(X)\wedge \phigm(Y) \to \phigm(X\wedge Y),
\end{equation}
making $\phigm$ lax monoidal.  From the formula~\eqref{eq:69}  this map
is an isomorphism if $X=S^{-V}\wedge A$ and $Y=S^{-W}\wedge B$.   This
leads to

\begin{prop}[\cite{MR1922205}, Proposition V.4.7]
\label{thm:55} The functor $\phigm$ is weakly monoidal: the
map~\eqref{eq:67} is a weak equivalence (in fact an isomorphism) if
$X$ and $Y$ are cellular.
\end{prop}

\begin{pf}
The class of spectra $X$ and $Y$ for which~\eqref{eq:67} is an
isomorphism is stable under smashing with a $G$-space, the formation
of wedges, directed colimits, and cobase change along an objectwise
closed inclusion.  Since~\eqref{eq:67} is an isomorphism when
$X=G_{+}\smashove{H} S^{-V}\wedge A$ and $Y=G_{+}\smashove{H}
S^{-W}\wedge B$ this implies it is an isomorphism when $X$ and $Y$ are
cellular.  Since isomorphisms are weak equivalences, the result
follows.
\end{pf}

\begin{rem}
\label{rem:18} 
Blumberg and
Mandell~\cite[Appendix~A]{blumberg_cylcotomic} have shown that
Proposition~\ref{thm:55} remains true under the assumption that only
one of $X$ or $Y$ is cellular.  This implies that
Proposition~\ref{thm:174} below remains true if only one of $N$ or
$N'$ is cofibrant.
\end{rem}

\subsubsection{Relation with the geometric fixed point functor}
\label{sec:relat-with-isotr}

We now turn to identifying the left derived functor $\lder\phigm$ with
the geometric fixed point functor $\phig$.  The inclusion $S^{0}\to
\tilde E\pfamily$ and the fibrant replacement functor give maps
\[
X \to \tilde E\pfamily\wedge X \to (\tilde E\pfamily\wedge X)_{f}.
\]

\begin{prop}[\cite{MR1922205}, Proposition V.4.17]
\label{thm:14}
If $X$ is cofibrant, then the maps 
\[
\phig X = (\tilde E\pfamily\wedge X_{f})^{G} \to \phigm((\tilde
E\pfamily\wedge X)_{f}) \leftarrow \phigm(X)
\]
are weak equivalences.
\end{prop}

\begin{pf*}{Sketch of proof}
For the arrow on the left, note that both functors
are homotopical and, up to weak equivalence, preserve filtered
colimits along $h$-cofibrations.  Using the canonical homotopy
presentation, it suffices to check that the arrow on the left is a
weak equivalence when $X=S^{-V}\wedge A$, with $A$ a $G$-CW complex.
This follows from Corollary~\ref{thm:205}, the identity~\eqref{eq:69},
and a little diagram chasing to check compatibility.

The right arrow is the composition of
\[
\phigm(X)\to \phigm(\tilde E\pfamily\wedge X)
\]
which is an isomorphism by~\eqref{eq:76}, and 
\[
\phigm(\tilde E\pfamily\wedge X)\to \phigm((\tilde E\pfamily\wedge X)_{f}),
\]
which is an acyclic cofibration by Proposition~\ref{thm:214}.
\end{pf*}

\subsubsection{The relative monoidal geometric fixed point functor}
\label{sec:relat-mono-geom}

The functor $\phigm$ can be formulated relative to an equivariant
commutative or associative algebra $R$.  As described below, care must
be taken in using the theory in this way.

Because it is lax monoidal, the functor $\phigm$ gives a functor 
\[
\phigm:\rmod{R}\to \rmod{\phigm R}
\]
which is lax monoidal in case $R$ is commutative.

\begin{prop}
\label{thm:15}
The functor 
\[
\phigm:\rmod{R}\to \rmod{\phigm R}
\]
commutes with cobase change along a cofibration and preserves the
classes of cofibrations and acyclic cofibrations.  
\end{prop}

\begin{pf}
This follows easily from the fact that the maps of equivariant
orthogonal spectra underlying the generating cofibrations for
$\rmod{R}$ are $h$-cofibrations.
\end{pf}

\begin{prop}
\label{thm:174} 
When $R$ is commutative, the functor
\[
\phigm:\rmod{R}\to \rmod{\phigm R}
\]
is weakly monoidal, and in fact
\begin{equation}
\label{eq:104}
\phigm(N')\smashove{\phigm(R)}\phigm(N) \to 
\phigm(N'\smashove{R}N)
\end{equation}
is an isomorphism if $N'$ and $N$ are cofibrant.  \qed
\end{prop}

\begin{pf}
The proof is the same as that of Proposition~\ref{thm:55} once one
knows that the class of modules $N'$ and $N$ for which~\eqref{eq:104}
is an isomorphism is stable under cobase change along a generating
cofibration.  This, in turn, is a consequence of the fact that both
sides of~\eqref{eq:104} preserve $h$-cofibrations in each variable,
since $h$-cofibrations are closed inclusions.  The functor $\phigm$
does so since it commutes with the formation of mapping cylinders, and
$N'\smashove{R}(\slot)$ does since $\rmod{R}$ is a closed symmetric
monoidal category.
\end{pf}

As promising as it looks, it is not so easy to make use of
Proposition~\ref{thm:174}.  The trouble is that unless $X$ is
cofibrant, $\phigm(X)$ may not have the weak homotopy type of
$\phig(X)$.  So in order to use Proposition~\ref{thm:174} one needs a
condition guaranteeing that $N'\smashove{R}N$ is a cofibrant spectrum.
The criterion of Proposition~\ref{thm:190} below was suggested to us
by Mike Mandell.

\begin{prop}
\label{thm:190} Suppose $R$ is an associative algebra with the
property that $S^{-1}\wedge R$ is cofibrant.  If $N'$ is a cofibrant
right $R$-module, and $S^{-1}\wedge N$ is a cofibrant left $R$-module,
then
\[
N'\smashove{R}N
\]
is cofibrant.
\end{prop}

\begin{pf}
First note that the condition on $R$ guarantees that for every
representation $U$ with $\dim U^{G}>0$ and every cofibrant $G$-space
$T$, the spectrum
\begin{equation}
\label{eq:5}
S^{-U}\wedge R \wedge T
\end{equation}
is cofibrant.  Since the formation of $N'\smashove{R}N$ commutes with
cobase change in both variables, the result reduces to the case
$N'=S^{-V}\wedge R\wedge X$ and $N=S^{-W}\wedge R\wedge Y$ with $V$
having a non-zero fixed point space, and $X$ and $Y$ cofibrant
$G$-spaces.  But in that case
\[
N'\smashove{R}N \approx S^{-V\oplus W}\wedge R\wedge X\wedge Y
\]
which is of the form~\eqref{eq:5}, and hence cofibrant.
\end{pf}

\begin{cor}
\label{thm:57}
Suppose $R$ is an associative algebra with the property that
$S^{-1}\wedge R$ is cofibrant.  If $N'$ is a cofibrant right
$R$-module, then the equivariant orthogonal spectrum underlying $N'$ is cofibrant.
\end{cor}

\begin{pf}
Just take $N=R$ in Proposition~\ref{thm:190}.
\end{pf}

The following result plays an important role in determining $\phig
R(\infty)$ (\S\ref{sec:proof-theor-refthm:21}).

\begin{prop}
\label{thm:191} 
Suppose that $R$ is an equivariant associative algebra
whose underlying $G$-spectrum is cellular, and that $R\to S^{0}$ is an
equivariant associative algebra map.  If $N'$ is a cofibrant right
$R$-module, then $N'\smashove{R}S^{0}$ is a cofibrant spectrum, and
the map
\[
\phigm(N')\smashove{\phigm R} S^{0} \to \phigm (N'\smashove{R}S^{0})
\]
is an isomorphism.
\end{prop}

\begin{pf}
One easily reduces to the case $N'=S^{-V}\wedge X\wedge R$, in which
$V$ is a representation with $\fixg{V}\ne 0$, and $X$ is a cofibrant
$G$-space.  In this case $N'\smashove{R}S^{0}$ is isomorphic to
$S^{-V}\wedge X$ which is cofibrant.  The assertion about monoidal
geometric fixed points follows easily from Proposition~\ref{thm:55}.
\end{pf}

\subsection{Geometric fixed points and the norm}
\label{sec:geom-fixed-points-norm}

Our original version of the following result merely concluded that the
transformation in question was a weak equivalence on cofibrant
objects.  Andrew Blumberg and Mike Mandell pointed out that it is in
fact an isomorphism.  At their request we have included the stronger
statement.

\begin{prop}
\label{thm:101} Suppose $H\subset G$.  There is a natural
transformation
\[
\phihm(\slot) \to \phigm\circ \norm_{H}^{G}(\slot)
\]
which is an isomorphism, hence a weak equivalence on cellular objects.
\end{prop}

\begin{pf}
To construct the natural transformation, first note that there is a
natural isomorphism 
\[
A^{H}\approx (\norm_{H}^{G}A)^{G}
\]
for $H$-equivariant {\em spaces} $A$.  Next note that for an
orthogonal representation $V$ of $H$, Proposition~\ref{thm:12} and the
property~\eqref{eq:69} give isomorphisms
\[
\phigm\norm_{H}^{G}S^{-V}\approx \phigm
S^{-\ind_{H}^{G}V} \approx S^{-V^{H}}\approx \Phi^{H} S^{-V}.
\]
The monoidal properties of $\phigm$ and the norm then combine to give
an isomorphism
\begin{equation}
\label{eq:77}
\phih (S^{-V}\wedge A)\approx
\phig \norm_{H}^{G}(S^{-V}\wedge A)
\end{equation}
which one easily checks to be compatible with the maps
\[
S^{-V}\wedge\igcat{H}(W,V) \to S^{-W}.
\]
To construct the transformation, write a general $H$-spectrum $X$ in
terms of its tautological presentation
\[
\bigvee_{V,W}S^{-W}\wedge \igcat{H}(V,W)\wedge X_{V} \rightrightarrows
\bigvee_{V}S^{-V}\wedge X_{V} \to X,
\]
and apply~\eqref{eq:77} termwise to produce a diagram 
\[
\bigvee_{V,W}S^{-W^{H}}\wedge \igcat{H}(V,W)^{H}\wedge X_{V}^{H} \rightrightarrows
\bigvee_{V}S^{-V^{H}}\wedge X_{V}^{H} \to \phig\norm_{H}^{G} X.
\]
The coequalizer of the two arrows is, by definition,  $\phihm(X)$.
This gives the natural transformation.   

The isomorphism assertion for cellular $X$ reduces to the special
case~\eqref{eq:77}, once one show that $\phigm\circ
\norm_{H}^{G}(\slot)$ commute with the formation of wedges, cobase
change along cofibrations between cofibrant objects, and filtered
colimits along closed inclusions.  The last property is clear since
both of the functors being composed commutes with filtered colimits
along closed inclusions.  For the other two assertions it will be
easier to work in terms of equivariant $J$-diagrams for $J=G/H$.

Suppose that $T$ is an indexing set, and $X_{t}$, $t\in T$ a set of
equivariant $J$-diagrams.  We wish to show that the natural map
\begin{equation}
\label{eq:187}
\bigvee_{t\in T}\phigm X_{t}^{\wedge J} \to 
\phigm\big(\bigvee_{t\in T}X_{t}\big)^{\wedge J}
\end{equation}
is an isomorphism.  For this use the distributive law to rewrite the
argument of the right hand side as
\[
\bigvee_{\gamma\in \Gamma} X^{\wedge \gamma}
\]
where $\gamma$ is the $G$-set of functions $J\to T$ and 
\[
X^{\wedge \gamma} = \bigwedge_{j\in J} X_{\gamma(j)}.
\]
The map asserted to be an isomorphism on monoidal geometric fixed
points is the inclusion of the summand indexed by the constant
functions.  But since $G$ acts trivially on $T$, the other summands
form an indexed wedge over a $G$-set with no fixed points.  The claim
then follows from Proposition~\ref{thm:73}.

The cobase change property is similar.  Suppose we are given a pushout
square of equivariant $J$-diagrams
\[
\xymatrix{
A  \ar[r]\ar[d]  & B
\ar[d] \\
X  \ar[r]        & Y
}
\]
in which $A\to B$ is a cofibration and $A$ is cofibrant.  We consider
the filtration of $Y^{\wedge J}$ given in
\S\ref{sec:index-mono-prod-pushout} whose stages fit into a pushout
square
\[
\xymatrix{
{\displaystyle 
\bigvee_{\substack{J=J_{0}\amalg J_{1} \\ |J_{1}|=m}}
X^{\wedge J_{0}}\wedge \partial_A{B^{\wedge J_{1}}}}
  \ar[r]\ar[d]  &   
{\displaystyle 
\bigvee_{\substack{J=J_{0}\amalg J_{1} \\ |J_{1}|=m}}
X^{\wedge J_{0}}\wedge B^{\wedge J_{1}}}\ar[d] \\
\fil_{m-1}Y^{\wedge J}  \ar[r]        & \fil_{m}Y^{\wedge J}\mathrlap{\ .}
}
\]
By Proposition~\ref{thm:266}, the upper arrow is an $h$-cofibration,
so the resulting diagram of monoidal geometric fixed points is a
pushout.   But since $J$ is a transitive $G$-set, unless $m=|J|$ the
group $G$ has no fixed points on the $G$-set indexing the wedges.
Applying Proposition~\ref{thm:73} then shows that for $m<|J|$ the map 
\[
\phigm X^{\wedge J} \to \phigm \fil_{m}Y^{\wedge J}
\]
is an isomorphism, and that the pushout square when $m=|J|$ becomes
\[
\xymatrix{
\phigm \partial_{A}B^{\wedge J}  \ar[r]\ar[d]  & \phigm B^{\wedge J}
\ar[d] \\
\phigm X^{\wedge J}  \ar[r]        & \phigm Y^{\wedge J}\mathrlap{\ .}
}
\]
However the term $\partial_{A}B^{\wedge J}$ is the term
$\fil_{|J|-1}B^{\wedge J}$ in the case in which $X=A$ and $Y=B$, and
so $\phigm A^{\wedge J}\to \phigm \partial_{A}B^{\wedge J}$ is an
isomorphism.  This completes the proof.  
\end{pf}

Thinking in terms of left derived functors one can get a slightly
better result.  As long as $X$ has the property that the map
$(\lder\norm_{H}^{G})X\to \norm_{H}^{G}X$ is a weak equivalence, there
will be a weak equivalence between $\phih X$ and
$\phig\norm_{H}^{G}X$.  Since it plays an important role in our work,
we spell it out.  Start with $X\in\ugspectra{H}$ and let $X_{c}\to X$
be a cofibrant approximation.  Now consider the diagram
\begin{equation}
\label{eq:66}
\xymatrix{ \phih X_{c} \ar[d]_{\sim}
\ar@{<->}[r]^{\sim}_{\text{zig zag}} & 
\phihm X_{c} \ar[r]^-{\sim} & \phigm\norm_{H}^{G}X_{c} 
\ar@{<->}[r]^{\sim}_{\text{zig zag}} & 
\phig\norm_{H}^{G}X_{c} \ar[d] \\
\phih X   &&& \phig\norm_{H}^{G}X
}
\end{equation}
The left vertical arrows are weak equivalences since the geometric
fixed point functor preserves weak equivalences.  The weak
equivalences in the top row are given by Propositions~\ref{thm:14},
\ref{thm:195}, and~\ref{thm:101}.  Since $\phig$ is homotopical we
have 

\begin{prop}
\label{thm:97}
Suppose that $X\in\ugspectra{H}$ has the property that for some (hence
any) cofibrant approximation $X_{c}\to X$ the map 
\[
\norm_{H}^{G}X_{c}\to \norm_{H}^{G}X
\]
is a weak equivalence.  Then the functorial relationship between
$\phih X$ and $\phig\norm_{H}^{G}X$ given by~\eqref{eq:66} is a weak
equivalence.  \qed
\end{prop}

\begin{rem}
\label{rem:68}
Proposition~\ref{thm:97} can be proved without reference to $\phigm$
by using the canonical homotopy presentation.
\end{rem}

\begin{rem}
\label{rem:33} Proposition~\ref{thm:97} applies in particular when $X$
is {\em very flat} in the sense of \S\ref{sec:very-flat-diagrams}.  By
Proposition~\ref{thm:16} this means that if $R\in\ugspectra{H}$ is a
cofibrant commutative ring, then $\phih R$ and $\phig\norm_{H}^{G}R$
are related by a functorial zig-zag of weak equivalences.  The case of
interest to us is when $H=\zt$, $G=\ztn$ and $R=\mur$.  In this case
$\norm_{H}^{G}R = \mutn{G}$, and we get an equivalence
\[
\phig\mutn{G}\approx \gfp{\zt}\mur \approx MO.
\]
\end{rem}

\begin{rem}
\label{rem:10}
Proposition~\ref{thm:97} also applies to the suspension spectra of
cofibrant $H$-spaces.  Indeed, if $X$ is a cofibrant $H$-space then
$S^{-1}\wedge S^{1}\wedge X\to X$ is a cofibrant approximation.
Applying $\norm_{H}^{G}$ leads to the map 
\[
S^{-V}\wedge S^{V}\wedge\norm_{H}^{G}(X) \to
\norm_{H}^{G}(X)
\]
with $V=\ind_{H}^{G}\R$, which is a weak equivalence (in fact a
cofibrant approximation).  This case is used to show that
$\phig\circ\norm_{H}^{G}$ is a ring homomorphism on the $RO(G)$-graded
cohomology of $G$-spaces (Proposition~\ref{thm:125}).
\end{rem}

\subsection{Real bordism}
\label{sec:real-bordism}
In this section we give a construction of the real bordism spectrum
$\mur$ as a commutative algebra in $\ugspectra{\zt}$.  As will be
apparent to the reader, this construction owes a great deal to the
Stefan Schwede's construction of $MU$ in~\cite[Chapter
2]{schwede:_untit_book_projec_about_symmet_spect}. We are indebted to
Schwede for some very helpful correspondence concerning these matters.

Our goal is to construct a $\zt$-equivariant commutative ring $\mur$
admitting the canonical homotopy presentation
\begin{equation}
\label{eq:141}
\mur \approx \hocolim S^{-\C^{n}}\wedge MU(n),
\end{equation}
in which $MU(n)$ is the Thom complex of the universal bundle over
$BU(n)$.  The group $\zt$ is acting on everything by complex
conjugation, so we could also write this expression as 
\begin{equation}
\label{eq:143}
\mur \approx \hocolim S^{-n\rho_{2}}\wedge MU(n).
\end{equation}
The map
\[
S^{-\rho_{2}}\wedge MU(1) \to \mur
\]
defines a real orientation.  These things form the basis for
everything we proved about $\mur$.

The most natural construction of $\mur$ realizes
this structure in the category $\rspectra$ of {\em real spectra},
which is related to the category of $\zt$-equivariant orthogonal
spectra by a multiplicative Quillen equivalence
\[
i_{!}:\rspectra \leftrightarrows \ugspectra{\zt}:i^{\ast}.
\]
We will construct a commutative algebra $\MUR\in\comm\rspectra$, whose
underlying real spectrum has a canonical homotopy presentation of the form
\begin{equation}
\label{eq:197}
\MUR \xleftarrow{\sim} \hocolim S^{-n\C}\wedge MU(n)\xrightarrow{\sim}
\hocolim S^{-n\C}\wedge MU(n)_{\text{cf}}.
\end{equation}
Applying $i_{!}$ to~\eqref{eq:197} and making the identification
$i_{!}S^{-\C}=S^{-\rho_{2}}$ leads to the diagram
\begin{equation}
\label{eq:42}
i_{!}\MUR \leftarrow \hocolim S^{-n\rho_{2}}\wedge MU(n)\rightarrow 
\hocolim S^{-n\rho_{2}}\wedge MU(n)_{\text{cf}}.
\end{equation}
We define $\mur$ to be the spectrum $i_{!}\MUR'$, where $\MUR'\to
\MUR$ is a cofibrant commutative algebra approximation.  The functor
$i_{!}$ is strictly monoidal, so $\mur$ is a commutative ring in
$\ugspectra{\zt}$.  The map on the right in~\eqref{eq:42} is a weak
equivalence since $i_{!}$ is a left Quillen functor.  The problem is
to show that the one on the left is.

This involves two steps.  The first is to show that the forgetful
functor
\[
\comm\rspectra\to\rspectra
\]
creates a model category structure on $\comm\rspectra$.  This involves
analyzing the symmetric powers of cofibrant real spectra, which, as
pointed out in Remark~\ref{rem:53} depends in an essential way on
understanding the homotopy theoretic properties of indexed symmetric
powers.  The second is to show that the functor $i_{!}$ is homotopical
on a subcategory of $\rspectra$ containing the real spectra underlying
cofibrant real commutative rings.  As in our analysis of norms of
commutative rings, this involves a generalized notion of flatness.
There is no real way to short circuit the model structure on
$\comm\rspectra$.  Its role is to identify the cofibrant real
commutative algebras.  But the only real work in establishing the
model structure is showing that what one thinks is a cofibrant
approximation is actually a weak equivalence, and that is what is
needed to show that every real commutative algebra is weakly
equivalent to a cofibrant one.  

\subsubsection{Real and complex spectra}

In this section we describe the basics of {\em real and complex
spectra}.  The additive results are more or less all a special case of
the results of~\cite{MR1806878}, but the important multiplicative
properties require a separate analysis.

For finite dimensional complex Hermitian vector spaces $A$ and $B$ let
$U(A,B)$ be the Stiefel manifold of unitary embeddings
$A\hookrightarrow B$.  There is a natural Hermitian inner product on
the complexification $V_{\C}$ of a real orthogonal vector space $V$,
so there is a natural map
\[
O(V,W)\to U(V_{\C},W_{\C}).  
\]
The group $\zt$ acts on $U(V_{\C},W_{\C})$ by complex conjugation, and
the fixed point space is $O(V,W)$.

\begin{defin}
\label{def:22} 
The category $\icatc$ is the 
topological category whose objects are finite dimensional Hermitian
vector spaces, and whose morphism space $\icatc(A,B)$ is the
Thom complex
\[
\icatc(A,B) = \thom(U(A,B); B-A).
\]
The category $\icatr$ is the $\zt$-equivariant
topological category whose objects are finite dimensional orthogonal
real vector spaces $V$, and with
\[
\icatr(V,W) = \icatc(V_{\C},W_{\C}),
\]
on which $\zt$ acts by complex conjugation.
\end{defin}

\begin{defin}
\label{def:51} 
The category $\cspectra$ of {\em complex spectra} is the topological
category of (continuous) functors
\[
\icatc\to\spaces. 
\]
The category $\rspectra$ of {\em real spectra} is the
topological category of $\zt$-enriched functors
\[
\icatr\to\gspaces{\zt},
\]
and equivariant natural transformations.   
\end{defin}

We will write 
\[
V\mapsto X_{V_{\C}}
\]
for a typical real spectrum $X$, and let $S^{-V_{\C}}\in\spectrar$ be
the functor co-represented by $V\in\icatr$.  From the Yoneda lemma
there is a natural isomorphism
\[
\rspectra(S^{-V_{\C}},X) = X_{V_{\C}}.
\]
As with equivariant orthogonal spectra, every real
spectrum $X$ has a tautological presentation
\begin{equation}
\label{eq:78}
\bigvee_{V,W\in\icatr}S^{-W_{\C}}\wedge \icatr(V,W)\wedge X_{W_{\C}}
\rightrightarrows \bigvee_{V\in\icatr} S^{-V_{\C}}\wedge X_{V_{\C}} \to X.
\end{equation}
A similar apparatus exist for complex spectra.

\begin{rem}
\label{rem:69}
The category $\icatr$ is equivalent to its full subcategory with
objects $\R^{n}$, and similarly $\icatc$ is equivalent to its full
subcategory with objects $\C^{n}$.  Thus a real spectrum $X$ is
specified by the spaces $X_{V_{\C}}$ with $V=\R^{n}$ together with the
structure maps between them, and an object $Y\in\cspectra$ is
specified by its spaces $Y_{\C^{n}}$, together with the structure maps
between them.
\end{rem}

The group $\zt$ acts on $\cspectra$ through its action on $\icatc$.
We write this as $X\mapsto \bar X$, where
\[
(\bar X)_{V} = X_{\bar V}.
\]
A fixed point for this action is a complex spectrum $X$ equipped with
an isomorphism $X\to \bar X$ having the property that $X\to \bar X\to
\bar{\bar X}=X$ is the identity map.  Restricting to the spaces
$X_{\C^{n}}$ and using the standard basis to identify $\C^{n}$ with
$\bar\C^{n}$ one sees that a fixed point for this $\zt$-action
consists of a sequence $\zt$-spaces $X_{\C^{n}}$, together with an
associative family $\zt$-equivariant maps
\[
\icatc(\C^{n},\C^{m})\wedge_{U(\C^{n})} X_{\C^{n}}\to X_{\C^{m}},
\]
where $\zt$ is acting by conjugation.   But this is the same thing as
giving a real spectrum indexed on the spaces $\R^{n}$.   This shows
that the category of fixed points for the $\zt$-action on $\cspectra$
is $\rspectra$.

\subsubsection{Smash product and indexed smash products}
\label{sec:smash-prod-index}

The orthogonal sum makes $\icatc$ into a symmetric monoidal category and
$\icatr$ an $\ugspaces{\zt}$-enriched symmetric monoidal category.
Using this one can define the smash product $X\wedge Y$  giving both
$\rspectra$ and $\cspectra$ the structure of symmetric monoidal
categories.   The smash product in $\rspectra$ is specified by the formula
\[
S^{-V_{\C}}\wedge S^{-W_{\C}} = S^{-(V\oplus W)_{\C}}
\]
and the fact that it commutes with colimits in each variable.  A
similar characterization holds for $\cspectra$.

There are indexed monoidal products in this context.  Let $J$ be a
finite set with a $\zt$-action.  The actions of $\zt$ on $J$ and on
$\cspectra$ combine to give an action on the product category
$\cspectra^{J}$.  The category of $\rspectra^{J}$ of {\em real
$J$-diagrams} is the category of fixed points for this action.  The
category of real $J$-diagrams for $J=\{\text{pt} \}$ is equivalent to
$\rspectra$.  When $J=\zt$, the category of real $J$-diagrams is
equivalent to $\cspectra$.  For general $J=n_{1}+ n_{2}\zt$, one has
an equivalence
\[
\rspectra^{J} \approx \rspectra^{n_{1}}\times \cspectra^{n_{2}}.
\]
There are indexed wedges and indexed smash products from
$\rspectra^{J}$ to $\rspectra$.

\subsubsection{Homotopy theory of real and complex spectra}
\label{sec:homotopy-theory-real}

We now turn to the homotopy theory of real and complex spectra.  We
describe the case of $\rspectra$ and leave the analogous case of
$\cspectra$ to the reader.

Suppose that $X$ is a real spectrum.  For $H\subset \zt$ and $k\in\Z$
set 
\[
\pi_{k}^{H}(X) = \varinjlim_{V}\pi_{k+V_{\C}}^{H}X_{V_{\C}}.
\]
The colimit is taken over the poset of finite dimensional orthogonal
vector spaces over $\R$, ordered (in agreement with
Definition~\ref{def:28}) by dimension.  A {\em stable weak
equivalence} in $\rspectra$ is a map $X\to Y$ inducing an isomorphism
$\pi^{H}_{k}X\to \pi^{H}_{k}Y$ for all $H\subset \zt$ and $k\in\Z$.
For fixed $k$, the groups $\pi^{H}_{k}$ form a Mackey functor which we
denote $\pim_{k}$.

Equipped with the stable weak equivalences, the category $\rspectra$
becomes a homotopical category.  We refine it to a model category by
defining a map to be a {\em fibration} if for each {\em non-zero} $V$,
the map $X_{V_{\C}}\to Y_{V_{\C}}$ is a fibration in $\ugspaces{\zt}$.
The cofibrations are the maps having the left lifting property against
the acyclic fibrations.  This is the {\em positive stable model
structure} on $\rspectra$.

The positive stable model structure is cofibrantly generated.  The
generating cofibrations can be taken to be the maps of the form
\begin{align*}
S^{-V_{\C}}\wedge \big(S^{n-1}_{+} &\to D^{n}_{+}\big) \\
\intertext{and}
(\zt)_{+} \wedge S^{-V_{\C}}\wedge \big(S^{n-1}_{+}&\to
D^{n}_{+}\big)
\end{align*}
with $V>0$.    The generating acyclic
cofibrations are the analogous maps
\begin{align*}
S^{-V_{\C}}\wedge \big(I^{n-1}_{+}&\to I^{n}\big) \\
\intertext{and}
(\zt)_{+}\wedge S^{-V_{\C}}\wedge \big(I^{n-1}_{+}&\to I^{n}  \big)
\end{align*}
together with the corner maps formed by smashing 
\begin{equation}
\label{eq:195}
S^{-V_{\C}\oplus W_{C}}\wedge S^{W_{\C}}\to \tilde S^{-V_{C}} 
\end{equation}
with the maps $S^{n-1}_{+}\to D^{n}_{+}$ and $(\zt)_{+}\wedge
(S^{n-1}_{+}\to D^{n}_{+})$.  We assume $V>0$, while $W$ need not be.
The map~\eqref{eq:195} is extracted from the factorization
\[
S^{-V_{\C}\oplus W_{\C}}\wedge S^{W_{\C}}\to \tilde S^{-V_{\C}} \to S^{-V_{\C}}
\]
formed by applying the small object construction with the generating
cofibrations.  As in the case of the complete positive stable model
structure on $\ugspectra{G}$, the map $\tilde S^{-V_{\C}}\to
S^{-V_{\C}}$ is a homotopy equivalence.  The verification of the model
category axioms is straightforward.  See
\S\ref{sec:strong-posit-stable} or~\cite{MR1806878}.  

\subsubsection{Real spectra and $\zt$-spectra}
\label{sec:real-spectra-zt}
Let
\[
i:\icatr  \to \igcat{\zt} 
\]
be the functor sending $V$ to 
\[
V_{\rho_{2}} = V\otimes\rho_{2}.
\]
Then the restriction functor
\[
i^{\ast}:\ugspectra{\zt}\to \rspectra
\]
has both a left and right adjoint which we denote $i_{!}$ and
$i_{\ast}$ respectively.  The left adjoint sends $S^{-V_{\C}}$ to
$S^{-V_{\rho_{2}}}$, and is described in general by applying the
functor termwise to the tautological presentation.

Since the functor $i$ is symmetric monoidal, the left adjoint $i_{!}$ is
strongly symmetric monoidal.  

\begin{prop}
\label{thm:209}
The functors 
\[
i_{!}:\rspectra \leftrightarrows \ugspectra{\zt}:i^{\ast}
\]
form a Quillen equivalence. 
\end{prop}

\begin{rem}
\label{rem:70} 
A similar discussion leads to a Quillen equivalence
\[
\cspectra \leftrightarrows \spectra.
\]
\end{rem}

\begin{pf}
Since $i_{!}$ is a left adjoint and
\[
i_{!}(S^{-V_{\C}}\wedge A) = S^{-V_{\rho_{2}}}\wedge A
\]
it is immediate that $i_{!}$ sends the generating (acyclic)
cofibrations to (acyclic) cofibrations, and hence is a left Quillen
functor.  Using the fact that the sequence $\{\R^{n}\otimes\rho_{2}\}$
is exhausting, one can easily check that a map $X\to Y$ in
$\ugspectra{\zt}$ is a weak equivalence if and only if $i^{\ast}X\to
i^{\ast}Y$ is.  This means that to show that $i_{!}$ and $i^{\ast}$
form a Quillen equivalence it suffices to show that the unit map
\begin{equation}
\label{eq:196}
X\to i^{\ast}i_{!}X
\end{equation}
is a weak equivalence for every cofibrant $X\in \rspectra$.  Since
$i^{\ast}$ is also a left adjoint, it preserves colimits, and
therefore so does $i^{\ast}i_{!}$.  Since both functors also commute
with smashing with a $\zt$-space, we are reduced to checking that for
each $0\ne V\in\icatr$, the map
\begin{equation}
\label{eq:138}
S^{-V_{\C}}\to i^{\ast}S^{-V_{\rho_{2}}}
\end{equation}
is a weak equivalence.   

For $W\in\icatr$, the $W_{\C}$-space of
$S^{-V_{\C}}$ is
\[
\icatr(V,W) = \thom(U(V_{\C},W_{\C}); W_{C}-V_{\C}) 
\]
and the $W$-space of $i^{\ast} S^{-V_{\rho_{2}}}$ is 
\[
\igcat{\zt}(V_{\rho_{2}},W_{\rho_{2}}) =
\thom(O(V_{\rho_{2}},W_{\rho_{2}});
W_{\rho_{2}}-V_{\rho_{2}}).
\]
The unit of the adjunction is derived from the inclusion
$U(V_{\C},W_{\C})\to O(V_{\rho_{2}}, W_{\rho_{2}})$.   
We must therefore show that for each $k$,  the map
\begin{equation}
\label{eq:140}
\varinjlim_{W\in\icatr}\pim_{k+W_{\C}}\icatr(V_{\C},W_{\C}) \to 
\varinjlim_{W\in\icatr}\pim_{k+W_{\C}}\igcat{\zt}(V_{\rho_{2}},W_{\rho_{2}})
\end{equation}
is an isomorphism.   

We may suppose that $\dim W>\dim V$.  For a fixed $W$ choose an
orthogonal embedding $V\subset W$, write $W=V\oplus U$, and consider
the diagram 
\[
\xymatrix{
S^{U_{\C}}  \ar[r]\ar[d]_{\approx}  &  \icatr(V,W) \ar[d] \\
S^{U_{\rho_{2}}}  \ar[r]        & \igcat{\zt}(V_{\rho_{2}},W_{\rho_{2}})\mathrlap{\ .}
}
\]
The left vertical map is an equivariant isomorphism.  A
straightforward argument using the connectivity of Stiefel manifolds
shows that for $\dim W\gg 0$ the horizontal maps are isomorphisms in
both $\piu_{k+W_{\C}}$ and $\pi^{\zt}_{k+W_{\C}}$.  It follows that
the right vertical map is as well, and hence so is~\eqref{eq:140}.
\end{pf}

For later reference, we record one fact that emerged in the proof of
Proposition~\ref{thm:209}.

\begin{lem}
\label{thm:227} The functor $i^{\ast}$ reflects weak equivalences: a
map $X\to Y\in \ugspectra{\zt}$ is a weak equivalence if and only if
$i^{\ast}X\to i^{\ast}Y$ is.  \qed
\end{lem}

\subsubsection{Multiplicative aspects of real spectra}
\label{sec:mult-aspects-real}

The multiplicative homotopy theory of real spectra is similar to that
of $\ugspectra{G}$.  Though there does not seem to be a simple
way to directly deduce the results from the case of $\ugspectra{\zt}$,
the proofs are very similar.

\begin{prop}
\label{thm:207} 
If $J$ is a set with a $\zt$-action and $X\to Y$ is a
cofibration of cofibrant real $J$-diagrams, then both the indexed
corner map $\partial_{X}Y^{\wedge J}\to Y^{\wedge J}$ and the absolute
map $X^{\wedge J}\to Y^{\wedge J}$ are cofibrations between
cofibrant objects.  They are weak equivalences if $X\to Y$ is.
\end{prop}

\begin{pf}
This is an analogue of Propositions~\ref{thm:267} and~\ref{thm:232},
and is proved in the same way, using the arrow category and the
filtration of \S\ref{sec:index-mono-prod-pushout}.
\end{pf}

For the symmetric powers, we fix a $\zt$-set $J$ and a $\zt$-stable
subgroup $\Sigma\subset \Sigma_{I}$.  The following is an analog of
Proposition~\ref{thm:200} and, making use of
Proposition~\ref{thm:207}, is proved in the same manner.

\begin{prop}
\label{thm:248} 
If $X\in\rspectra$ is cofibrant and $Z$ is any real
spectrum equipped with an action of $\Sigma\rtimes \zt$ extending the
$G$-action, then the map
\[
(\egsigma{\zt})_{+}\smashove{\Sigma} (X^{\wedge J}\wedge Z)\to
(X^{\wedge J}\wedge Z)/\Sigma.
\]
is a weak equivalence. \qed
\end{prop}

\begin{prop}
\label{thm:58}
If $A\to B$ is a cofibration of cofibrant real spectra and $J$ is a
finite set with a $\zt$-action, then in the diagram
\[
\xymatrix{
\egsigma{\zt}_{+}\smashove{\Sigma}\partial_{A}B^{\wedge J}
\ar[r]\ar[d]  &
\egsigma{\zt}_{+}\smashove{\Sigma}B^{\wedge J} 
\ar[d] \\
\partial_{A}\sym^{J}B  \ar[r]        & \sym^{J}B
}
\]
the upper row is a cofibration between cofibrant objects, the vertical
maps are weak equivalences and remain so after smashing
with any object, and the bottom row is an $h$-cofibration of
flat spectra.  The horizontal maps are weak equivalences if $A\to B$
is.
\end{prop}

\begin{pf}
This is an analogue of Proposition~\ref{thm:262} and is proved in the
same way, making use of Proposition~\ref{thm:248}.
\end{pf}

\begin{prop}
\label{thm:217} The forgetful functor 
\[
\comm\rspectra \to \rspectra
\]
creates a model category structure on commutative algebras in
$\rspectra$, in which a map of commutative algebras is a fibration or
weak equivalence if and only if the underlying map of real spectra is.
\end{prop}

\begin{pf}
This is proved in the same manner as Proposition~\ref{thm:114}, making
use of Proposition~\ref{thm:58}.   
\end{pf}

\subsubsection{Generalized flatness}
\label{sec:generalized-flatness}

Our next task is to show that the left derived functor of $i_{!}$ can
be computed on a subcategory of real spectra containing those which
underlie real commutative rings.

\begin{defin}
\label{def:50} A real spectrum $X\in\rspectra$ is {\em
$i_{!}$-flat} if it satisfies the following property: for every
cofibrant approximation $\tilde X\to X$ and every weak equivalence
$\tilde Z\to Z\in\ugspectra{\zt}$
the map
\begin{equation}
\label{eq:146}
i_{!}\tilde X\wedge \tilde Z \to 
i_{!} X\wedge Z
\end{equation}
is a weak equivalence.  
\end{defin}

\begin{rem}
\label{rem:71} 
Since $i_{!}$ is a left Quillen functor and cofibrant objects of
$\ugspectra{\zt}$ are flat, cofibrant objects of $\rspectra$ are
$i_{!}$-flat.
\end{rem}

\begin{rem}
\label{rem:50} If~\eqref{eq:146} is a weak equivalence for one
cofibrant approximation it is a weak equivalence for any cofibrant
approximation.   
\end{rem}

Our main result is

\begin{prop}
\label{thm:213}
If $R\in\rspectra$ is a cofibrant commutative algebra then $R$ is $i_{!}$-flat.
\end{prop}

The proof of Proposition~\ref{thm:213} follows the argument for the
proof of Proposition~\ref{thm:16}.

\begin{lem}
\label{thm:215}
If $A\in\rspectra$ is cofibrant, and $n\ge 1$, then $\sym^{n}A$ is $i_{!}$-flat.
\end{lem}

\begin{pf}
By Proposition~\ref{thm:248}, the map
\[
(\eggroup{\zt}{\Sigma_{n}})_{+}\smashove{\Sigma_{n}} A^{\wedge n}\to
\sym^{n}A
\]
is a cofibrant approximation.   Since $i_{!}$ is a continuous left
adjoint, we may identify
\begin{equation}
\label{eq:136}
i_{!}\big((\eggroup{\zt}{\Sigma_{n}})_{+}\smashove{\Sigma_{n}} A^{\wedge
n}\big)\wedge \tilde Z\to i_{!}\big(\sym^{n}A\big)\wedge Z
\end{equation}
with 
\begin{equation}
\label{eq:135} (\eggroup{\zt}{\Sigma_{n}})_{+}\smashove{\Sigma_{n}}
(i_{!}A)^{\wedge n}\wedge \tilde Z\to \sym^{n}(i_{!}A)\wedge Z.
\end{equation}
Since $i_{!}$ is a left Quillen functor, $i_{!}(A)$ is cofibrant, and
Proposition~\ref{thm:200} implies that~\eqref{eq:135},
hence~\eqref{eq:136} is a weak equivalence.  
\end{pf}

We also require an analogue of Lemma~\ref{thm:193}, though the
statement and proof are much simpler in this case, since $i_{!}$ is
a left adjoint.

\begin{lem}
\label{thm:226}
If $S\to T$ is an $h$-cofibration in $\rspectra$, and two of $S$, $T$,
$T/S$ are $i_{!}$-flat, then so is the third.
\end{lem}

\begin{pf}
We may choose a map $\tilde S\to \tilde T$ of cofibrant approximations
which is a cofibration, hence an $h$-cofibration.  Our assumption is
that two of the vertical maps in
\[
\xymatrix{
i_{!}\tilde S\wedge \tilde Z  \ar[r]\ar[d]  & i_{!}\tilde T\wedge \tilde Z  \ar[r]\ar[d]   &  i_{!}(\tilde T/\tilde S)\wedge \tilde Z \ar[d]\\
i_{!}S\wedge Z  \ar[r]        & i_{!}T\wedge Z   \ar[r]         & i_{!}(T/S)\wedge Z 
}
\]
are weak equivalences.  This implies that the third is, since the two
left horizontal maps are $h$-cofibrations hence flat.
\end{pf}

\begin{lem}
\label{thm:224}
Consider a pushout square in $\rspectra$, 
\begin{equation}
\label{eq:18}
\xymatrix{
S  \ar[r]\ar[d]  & T  \ar[d] \\
X  \ar[r]        & Y
}
\end{equation}
in which $S\to T$ is an $h$-cofibration.  If $T$, $T/S$ and $X$ are
$i_{!}$-flat, then so is $Y$.
\end{lem}

\begin{pf}
Since $T$ and $T/S$ are $i_{!}$-flat, so is $S$ by
Lemma~\ref{thm:226}.  We may choose cofibrant approximations of
everything fitting into a pushout diagram
\[
\xymatrix{
\tilde S  \ar[r]\ar[d] &\tilde T   \ar[d] \\
\tilde X  \ar[r]       &\tilde Y
}	
\]
in which the top row is an $h$-cofibration.   Now consider 
\[
\xymatrix{
i_{!}\tilde X\wedge \tilde Z  \ar[d]  &\ar[l]  i_{!}\tilde S\wedge \tilde Z
\ar[r]\ar[d]   
  &i_{!}\tilde T\wedge \tilde Z   \ar[d]\\
i_{!} X\wedge Z &\ar[l] i_{!} S\wedge Z\ar[r] &i_{!} T\wedge
Z 
}
\]
The left horizontal maps are $h$-cofibrations, hence flat, and the
vertical maps are weak equivalences by assumption.   It follows that
the map of pushouts is a weak equivalence.   
\end{pf}

\begin{pf*}{Proof of Proposition~\ref{thm:213}}
It suffices to show that if  $A\to B$ is a generating cofibration in
$\rspectra$ then
\[
\xymatrix{
\sym A  \ar[r]\ar[d]  & \sym B  \ar[d] \\
X  \ar[r]        & Y
}
\]
is a pushout square of commutative algebras in $\rspectra$ and $X$ is
$i_{!}$-flat, then $Y$ is $i_{!}$-flat.  We induct over the filtration
described in \S\ref{sec:index-mono-prod-pushout}.  Since
$\fil_{0}Y=X$, the induction starts.  For the inductive step, consider
the pushout square
\begin{equation}
\label{eq:137}
\xymatrix{
X\wedge \partial_{A}\sym^{m}B  \ar[r]\ar[d]  & X\wedge \sym^{m}B  \ar[d] \\
\fil_{m-1}Y  \ar[r]        &\fil_{m}Y\mathrlap{\ ,}
}
\end{equation}
and assume that $\fil_{m-1}Y$ is $i_{!}$-flat.  Both $\sym^{m}B$ and
\[
\sym^{m}B/\partial_{A}\sym^{m}B=\sym^{m}(B/A)
\]
are $i_{!}$-flat by Lemma~\ref{thm:215}.  Since smash products of
$i_{!}$-flat spectra are $i_{!}$-flat, both $X\wedge\sym^{m}B$ and
$X\wedge \sym^{m}(B/A)$ are $i_{!}$-flat.  The top row
of~\eqref{eq:137} is an $h$-cofibration, so Lemma~\ref{thm:224}
implies that $\fil_{m}Y$ is $i_{!}$-flat.  This completes the
inductive step, and the proof.
\end{pf*}

Though we don't quite need the following result, having come this far
we record it for future reference.

\begin{prop}
\label{thm:219}
The functors $i_{!}$ and $i^{\ast}$ restrict to a Quillen equivalence 
\[
i_{!}:\comm\rspectra\leftrightarrows \comm\ugspectra{\zt}:i^{\ast}.
\]
\end{prop}

\begin{pf}
It is immediate from the definition of the model structures on 
$\comm\rspectra$ and $\comm \ugspectra{\zt}$, and the fact that 
\[
i_{!}:\rspectra\leftrightarrows \ugspectra{\zt}:i^{\ast}
\]
is a Quillen pair, that 
\[
i^{\ast}:\comm\ugspectra{\zt}\to \comm\rspectra
\]
preserves the classes of fibrations and acyclic fibrations.    It
remains to show that if $A\in \comm\rspectra$ is cofibrant, then the
composition
\[
A \to i^{\ast} i_{!}A \to i^{\ast}(i_{!}A_{f})
\]
is a weak equivalence, where $i_{!}A \to i_{!}A_{f}$ is a fibrant
replacement.  Since $i^{\ast}$ reflects weak equivalences
(Lemma~\ref{thm:227}) this is equivalent to showing that
\[
A\to i^{\ast}i_{!}A
\]
is a weak equivalence.   Let $A'\to A$ be a cofibrant
approximation in $\rspectra$, and consider the following diagram in
$\rspectra$
\begin{equation}
\label{eq:142}
\xymatrix{
A' \ar[r]^{\sim} \ar[d]_{\sim} & i^{\ast}i_{!}A'\ar[d]^{\sim} \\
A \ar[r] & i^{\ast}i_{!}A\mathrlap{\ .}
}
\end{equation}
By Proposition~\ref{thm:213} the map $i_{!}A'\to i_{!}A$ is a weak
equivalence.  The rightmost arrow in~\eqref{eq:142} is therefore a
weak equivalence.  The top arrow is a weak equivalence by
Proposition~\ref{thm:209}, and the left arrow is a weak equivalence by
assumption.  This implies that the bottom arrow is a weak equivalence.
\end{pf}

\subsubsection{The real bordism spectrum}
\label{sec:real-bord-spectr}

For $V\in\icatr$ let 
\[
MU(V_{\C})=\thom(BU(V_{\C}),V_{\C})
\]
be the Thom complex of the bundle
$EU(V_{\C})\underset{U(V_{\C})}\times V_{\C}$ over $BU(V_{\C})$,
equipped with the $\zt$-action of complex conjugation.  We will take
our model of $BU(V_{\C})$ to be the one given by Segal's
construction~\cite{Segal:CSandSS}, so that
\begin{equation}
\label{eq:109}
V\mapsto \thom(BU(V_{\C}),V_{\C})
\end{equation}
is a lax symmetric monoidal functor $\icatr \to \gspaces{\zt}$, and so
defines a commutative ring $\MUR\in\comm\rspectra$.  Let $\MUR'\to
\MUR$ be a cofibrant approximation to $\MUR$ in $\comm\rspectra$.

\begin{defin}
\label{def:53}
The {\em real bordism spectrum} is the spectrum $\mur = i_{!}\MUR'$.
\end{defin}

To get at the homotopy type of $\mur$, we examine the canonical
homotopy presentation of $\MUR$ using the exhausting sequence
$V_{n}=\R^{n}$.  This gives a weak equivalence
\begin{equation}
\label{eq:112}
\hocolim S^{-\C^{n}}\wedge MU(n) \xrightarrow{\sim}\MUR'
\end{equation}
in which $MU(n)=MU(\C^{n})$.  Applying $i_{!}$ and using
Proposition~\ref{thm:213} gives 
\[
\hocolim S^{-n\rho_{2}}\wedge MU(n) \xrightarrow{\sim}\mur.
\]
In this presentation the universal real orientation of $\mur$
(Example~\ref{eg:16}) is given by restricting to the term $n=1$
\[
S^{-\rho_{2}}\wedge MU(1) \to \mur.
\]

The next result summarizes some further consequences of the presentation~\eqref{eq:112}.

\begin{prop}
\begin{thmList}
\item The non-equivariant spectrum underlying $\mur$ is the usual
complex cobordism spectrum $MU$.  
\item The equivariant cohomology theory represented by $\mur$
coincides with the one studied
in~\cite{MR0222890,MR0420597,MR614829,MR1808224}.
\item There is an equivalence 
\[
\Phi^{\zt}\mur \approx MO.
\]
\item The Schubert cell decomposition of Grassmannians leads to a
cofibrant approximation of $\mur$ by a $\zt$-CW complex with one
$0$-cell ($S^{0}$) and the remaining cells of the form
$e^{m\rho_{2}}$, with $m>0$.
\end{thmList}
\end{prop}

\def\cprime{$'$}
\providecommand{\bysame}{\leavevmode\hbox to3em{\hrulefill}\thinspace}
\providecommand{\MR}{\relax\ifhmode\unskip\space\fi MR }
\providecommand{\MRhref}[2]{%
  \href{http://www.ams.org/mathscinet-getitem?mr=#1}{#2}
}
\providecommand{\href}[2]{#2}

\end{document}